\documentstyle[11pt]{article}
\newtheorem{theorem}{Theorem}[section]
\newtheorem{lemma}[theorem]{Lemma}
\newtheorem{corollary}[theorem]{Corollary}

\newtheorem{proposition}[theorem]{Proposition}

\makeindex \makeglossary
\begin{document}                
%
%

\long\def\ig#1{\relax}
\ig{Thanks to Roberto Minio for this def'n.  Compare the def'n of
\comment in AMSTeX.}

\newcount \coefa
\newcount \coefb
\newcount \coefc
\newcount\tempcounta
\newcount\tempcountb
\newcount\tempcountc
\newcount\tempcountd
\newcount\xext
\newcount\yext
\newcount\xoff
\newcount\yoff
\newcount\gap%
\newcount\arrowtypea
\newcount\arrowtypeb
\newcount\arrowtypec
\newcount\arrowtyped
\newcount\arrowtypee
\newcount\height
\newcount\width
\newcount\xpos
\newcount\ypos
\newcount\run
\newcount\rise
\newcount\arrowlength
\newcount\halflength
\newcount\arrowtype
\newdimen\tempdimen
\newdimen\xlen
\newdimen\ylen
\newsavebox{\tempboxa}%
\newsavebox{\tempboxb}%
\newsavebox{\tempboxc}%

\makeatletter
\setlength{\unitlength}{.01em}%
\def\settypes(#1,#2,#3){\arrowtypea#1 \arrowtypeb#2 \arrowtypec#3}
\def\settoheight#1#2{\setbox\@tempboxa\hbox{#2}#1\ht\@tempboxa\relax}%
\def\settodepth#1#2{\setbox\@tempboxa\hbox{#2}#1\dp\@tempboxa\relax}%
\def\settokens[#1`#2`#3`#4]{%
     \def\tokena{#1}\def\tokenb{#2}\def\tokenc{#3}\def\tokend{#4}}
\def\setsqparms[#1`#2`#3`#4;#5`#6]{%
\arrowtypea #1
\arrowtypeb #2
\arrowtypec #3
\arrowtyped #4
\width #5
\height #6
}
\def\setpos(#1,#2){\xpos=#1 \ypos#2}

\def\bfig{\begin{picture}(\xext,\yext)(\xoff,\yoff)}
\def\efig{\end{picture}}

\def\putbox(#1,#2)#3{\put(#1,#2){\makebox(0,0){$#3$}}}

\def\settriparms[#1`#2`#3;#4]{\settripairparms[#1`#2`#3`1`1;#4]}%

\def\settripairparms[#1`#2`#3`#4`#5;#6]{%
\arrowtypea #1
\arrowtypeb #2
\arrowtypec #3
\arrowtyped #4
\arrowtypee #5
\width #6
\height #6
}

\def\resetparms{\settripairparms[1`1`1`1`1;500]\width 500}

\resetparms

\def\mvector(#1,#2)#3{
\put(0,0){\vector(#1,#2){#3}}%
\put(0,0){\vector(#1,#2){30}}%
}
\def\evector(#1,#2)#3{{
\arrowlength #3
\put(0,0){\vector(#1,#2){\arrowlength}}%
\advance \arrowlength by-30
\put(0,0){\vector(#1,#2){\arrowlength}}%
}}

\def\horsize#1#2{%
\settowidth{\tempdimen}{$#2$}%
#1=\tempdimen
\divide #1 by\unitlength
}

\def\vertsize#1#2{%
\settoheight{\tempdimen}{$#2$}%
#1=\tempdimen
\settodepth{\tempdimen}{$#2$}%
\advance #1 by\tempdimen
\divide #1 by\unitlength
}

\def\vertadjust[#1`#2`#3]{%
\vertsize{\tempcounta}{#1}%
\vertsize{\tempcountb}{#2}%
\ifnum \tempcounta<\tempcountb \tempcounta=\tempcountb \fi
\divide\tempcounta by2
\vertsize{\tempcountb}{#3}%
\ifnum \tempcountb>0 \advance \tempcountb by20 \fi
\ifnum \tempcounta<\tempcountb \tempcounta=\tempcountb \fi
}

\def\horadjust[#1`#2`#3]{%
\horsize{\tempcounta}{#1}%
\horsize{\tempcountb}{#2}%
\ifnum \tempcounta<\tempcountb \tempcounta=\tempcountb \fi
\divide\tempcounta by20
\horsize{\tempcountb}{#3}%
\ifnum \tempcountb>0 \advance \tempcountb by60 \fi
\ifnum \tempcounta<\tempcountb \tempcounta=\tempcountb \fi
}

\ig{ In this procedure, #1 is the paramater that sticks out all the way,
#2 sticks out the least and #3 is a label sticking out half way.  #4 is
the amount of the offset.}

\def\sladjust[#1`#2`#3]#4{%
\tempcountc=#4
\horsize{\tempcounta}{#1}%
\divide \tempcounta by2
\horsize{\tempcountb}{#2}%
\divide \tempcountb by2
\advance \tempcountb by-\tempcountc
\ifnum \tempcounta<\tempcountb \tempcounta=\tempcountb\fi
\divide \tempcountc by2
\horsize{\tempcountb}{#3}%
\advance \tempcountb by-\tempcountc
\ifnum \tempcountb>0 \advance \tempcountb by80\fi
\ifnum \tempcounta<\tempcountb \tempcounta=\tempcountb\fi
\advance\tempcounta by20
}

\def\putvector(#1,#2)(#3,#4)#5#6{{%
\xpos=#1
\ypos=#2
\run=#3
\rise=#4
\arrowlength=#5
\arrowtype=#6
\ifnum \arrowtype<0
    \ifnum \run=0
        \advance \ypos by-\arrowlength
    \else
        \tempcounta \arrowlength
        \multiply \tempcounta by\rise
        \divide \tempcounta by\run
        \ifnum\run>0
            \advance \xpos by\arrowlength
            \advance \ypos by\tempcounta
        \else
            \advance \xpos by-\arrowlength
            \advance \ypos by-\tempcounta
        \fi
    \fi
    \multiply \arrowtype by-1
    \multiply \rise by-1
    \multiply \run by-1
\fi
\ifnum \arrowtype=1
    \put(\xpos,\ypos){\vector(\run,\rise){\arrowlength}}%
\else\ifnum \arrowtype=2
    \put(\xpos,\ypos){\mvector(\run,\rise)\arrowlength}%
\else\ifnum\arrowtype=3
    \put(\xpos,\ypos){\evector(\run,\rise){\arrowlength}}%
\fi\fi\fi
}}

\def\putsplitvector(#1,#2)#3#4{
\xpos #1
\ypos #2
\arrowtype #4
\halflength #3
\arrowlength #3
\gap 140
\advance \halflength by-\gap
\divide \halflength by2
\ifnum \arrowtype=1
    \put(\xpos,\ypos){\line(0,-1){\halflength}}%
    \advance\ypos by-\halflength
    \advance\ypos by-\gap
    \put(\xpos,\ypos){\vector(0,-1){\halflength}}%
\else\ifnum \arrowtype=2
    \put(\xpos,\ypos){\line(0,-1)\halflength}%
    \put(\xpos,\ypos){\vector(0,-1)3}%
    \advance\ypos by-\halflength
    \advance\ypos by-\gap
    \put(\xpos,\ypos){\vector(0,-1){\halflength}}%
\else\ifnum\arrowtype=3
    \put(\xpos,\ypos){\line(0,-1)\halflength}%
    \advance\ypos by-\halflength
    \advance\ypos by-\gap
    \put(\xpos,\ypos){\evector(0,-1){\halflength}}%
\else\ifnum \arrowtype=-1
    \advance \ypos by-\arrowlength
    \put(\xpos,\ypos){\line(0,1){\halflength}}%
    \advance\ypos by\halflength
    \advance\ypos by\gap
    \put(\xpos,\ypos){\vector(0,1){\halflength}}%
\else\ifnum \arrowtype=-2
    \advance \ypos by-\arrowlength
    \put(\xpos,\ypos){\line(0,1)\halflength}%
    \put(\xpos,\ypos){\vector(0,1)3}%
    \advance\ypos by\halflength
    \advance\ypos by\gap
    \put(\xpos,\ypos){\vector(0,1){\halflength}}%
\else\ifnum\arrowtype=-3
    \advance \ypos by-\arrowlength
    \put(\xpos,\ypos){\line(0,1)\halflength}%
    \advance\ypos by\halflength
    \advance\ypos by\gap
    \put(\xpos,\ypos){\evector(0,1){\halflength}}%
\fi\fi\fi\fi\fi\fi
}

\def\putmorphism(#1)(#2,#3)[#4`#5`#6]#7#8#9{{%
\run #2
\rise #3
\ifnum\rise=0
  \puthmorphism(#1)[#4`#5`#6]{#7}{#8}{#9}%
\else\ifnum\run=0
  \putvmorphism(#1)[#4`#5`#6]{#7}{#8}{#9}%
\else
\setpos(#1)%
\arrowlength #7
\arrowtype #8
\ifnum\run=0
\else\ifnum\rise=0
\else
\ifnum\run>0
    \coefa=1
\else
   \coefa=-1
\fi
\ifnum\arrowtype>0
   \coefb=0
   \coefc=-1
\else
   \coefb=\coefa
   \coefc=1
   \arrowtype=-\arrowtype
\fi
\width=2
\multiply \width by\run
\divide \width by\rise
\ifnum \width<0  \width=-\width\fi
\advance\width by60
\if l#9 \width=-\width\fi
\putbox(\xpos,\ypos){#4}
{\multiply \coefa by\arrowlength
\advance\xpos by\coefa
\multiply \coefa by\rise
\divide \coefa by\run
\advance \ypos by\coefa
\putbox(\xpos,\ypos){#5} }%
{\multiply \coefa by\arrowlength
\divide \coefa by2
\advance \xpos by\coefa
\advance \xpos by\width
\multiply \coefa by\rise
\divide \coefa by\run
\advance \ypos by\coefa
\if l#9%
   \put(\xpos,\ypos){\makebox(0,0)[r]{$#6$}}%
\else\if r#9%
   \put(\xpos,\ypos){\makebox(0,0)[l]{$#6$}}%
\fi\fi }%
{\multiply \rise by-\coefc
\multiply \run by-\coefc
\multiply \coefb by\arrowlength
\advance \xpos by\coefb
\multiply \coefb by\rise
\divide \coefb by\run
\advance \ypos by\coefb
\multiply \coefc by70
\advance \ypos by\coefc
\multiply \coefc by\run
\divide \coefc by\rise
\advance \xpos by\coefc
\multiply \coefa by140
\multiply \coefa by\run
\divide \coefa by\rise
\advance \arrowlength by\coefa
\ifnum \arrowtype=1
   \put(\xpos,\ypos){\vector(\run,\rise){\arrowlength}}%
\else\ifnum\arrowtype=2
   \put(\xpos,\ypos){\mvector(\run,\rise){\arrowlength}}%
\else\ifnum\arrowtype=3
   \put(\xpos,\ypos){\evector(\run,\rise){\arrowlength}}%
\fi\fi\fi}\fi\fi\fi\fi}}

\def\puthmorphism(#1,#2)[#3`#4`#5]#6#7#8{{%
\xpos #1
\ypos #2
\width #6
\arrowlength #6
\putbox(\xpos,\ypos){#3\vphantom{#4}}%
{\advance \xpos by\arrowlength
\putbox(\xpos,\ypos){\vphantom{#3}#4}}%
\horsize{\tempcounta}{#3}%
\horsize{\tempcountb}{#4}%
\divide \tempcounta by2
\divide \tempcountb by2
\advance \tempcounta by30
\advance \tempcountb by30
\advance \xpos by\tempcounta
\advance \arrowlength by-\tempcounta
\advance \arrowlength by-\tempcountb
\putvector(\xpos,\ypos)(1,0){\arrowlength}{#7}%
\divide \arrowlength by2
\advance \xpos by\arrowlength
\vertsize{\tempcounta}{#5}%
\divide\tempcounta by2
\advance \tempcounta by20
\if a#8 %
   \advance \ypos by\tempcounta
   \putbox(\xpos,\ypos){#5}%
\else
   \advance \ypos by-\tempcounta
   \putbox(\xpos,\ypos){#5}%
\fi}}

\def\putvmorphism(#1,#2)[#3`#4`#5]#6#7#8{{%
\xpos #1
\ypos #2
\arrowlength #6
\arrowtype #7
\settowidth{\xlen}{$#5$}%
\putbox(\xpos,\ypos){#3}%
{\advance \ypos by-\arrowlength
\putbox(\xpos,\ypos){#4}}%
{\advance\arrowlength by-140
\advance \ypos by-70
\ifdim\xlen>0pt
   \if m#8%
      \putsplitvector(\xpos,\ypos){\arrowlength}{\arrowtype}%
   \else
      \putvector(\xpos,\ypos)(0,-1){\arrowlength}{\arrowtype}%
   \fi
\else
   \putvector(\xpos,\ypos)(0,-1){\arrowlength}{\arrowtype}%
\fi}%
\ifdim\xlen>0pt
   \divide \arrowlength by2
   \advance\ypos by-\arrowlength
   \if l#8%
      \advance \xpos by-40
      \put(\xpos,\ypos){\makebox(0,0)[r]{$#5$}}%
   \else\if r#8%
      \advance \xpos by40
      \put(\xpos,\ypos){\makebox(0,0)[l]{$#5$}}%
   \else
      \putbox(\xpos,\ypos){#5}%
   \fi\fi
\fi
}}

\def\topadjust[#1`#2`#3]{%
\yoff=10
\vertadjust[#1`#2`{#3}]%
\advance \yext by\tempcounta
\advance \yext by 10
}
\def\botadjust[#1`#2`#3]{%
\vertadjust[#1`#2`{#3}]%
\advance \yext by\tempcounta
\advance \yoff by-\tempcounta
}
\def\leftadjust[#1`#2`#3]{%
\xoff=0
\horadjust[#1`#2`{#3}]%
\advance \xext by\tempcounta
\advance \xoff by-\tempcounta
}
\def\rightadjust[#1`#2`#3]{%
\horadjust[#1`#2`{#3}]%
\advance \xext by\tempcounta
}
\def\rightsladjust[#1`#2`#3]{%
\sladjust[#1`#2`{#3}]{\width}%
\advance \xext by\tempcounta
}
\def\leftsladjust[#1`#2`#3]{%
\xoff=0
\sladjust[#1`#2`{#3}]{\width}%
\advance \xext by\tempcounta
\advance \xoff by-\tempcounta
}
\def\adjust[#1`#2;#3`#4;#5`#6;#7`#8]{%
\topadjust[#1``{#2}]
\leftadjust[#3``{#4}]
\rightadjust[#5``{#6}]
\botadjust[#7``{#8}]}

\def\putsquarep<#1>(#2)[#3;#4`#5`#6`#7]{{%
\setsqparms[#1]%
\setpos(#2)%
\settokens[#3]%
\puthmorphism(\xpos,\ypos)[\tokenc`\tokend`{#7}]{\width}{\arrowtyped}b%
\advance\ypos by \height
\puthmorphism(\xpos,\ypos)[\tokena`\tokenb`{#4}]{\width}{\arrowtypea}a%
\putvmorphism(\xpos,\ypos)[``{#5}]{\height}{\arrowtypeb}l%
\advance\xpos by \width
\putvmorphism(\xpos,\ypos)[``{#6}]{\height}{\arrowtypec}r%
}}

\def\putsquare{\@ifnextchar <{\putsquarep}{\putsquarep%
   <\arrowtypea`\arrowtypeb`\arrowtypec`\arrowtyped;\width`\height>}}
\def\square{\@ifnextchar< {\squarep}{\squarep
   <\arrowtypea`\arrowtypeb`\arrowtypec`\arrowtyped;\width`\height>}}
\def\squarep<#1>[#2`#3`#4`#5;#6`#7`#8`#9]{{
\setsqparms[#1]
\xext=\width                                          
\yext=\height                                         
\topadjust[#2`#3`{#6}]
\botadjust[#4`#5`{#9}]
\leftadjust[#2`#4`{#7}]
\rightadjust[#3`#5`{#8}]
\begin{picture}(\xext,\yext)(\xoff,\yoff)
\putsquarep<\arrowtypea`\arrowtypeb`\arrowtypec`\arrowtyped;\width`\height>%
(0,0)[#2`#3`#4`#5;#6`#7`#8`{#9}]%
\end{picture}%
}}

\def\putptrianglep<#1>(#2,#3)[#4`#5`#6;#7`#8`#9]{{%
\settriparms[#1]%
\xpos=#2 \ypos=#3
\advance\ypos by \height
\puthmorphism(\xpos,\ypos)[#4`#5`{#7}]{\height}{\arrowtypea}a%
\putvmorphism(\xpos,\ypos)[`#6`{#8}]{\height}{\arrowtypeb}l%
\advance\xpos by\height
\putmorphism(\xpos,\ypos)(-1,-1)[``{#9}]{\height}{\arrowtypec}r%
}}

\def\putptriangle{\@ifnextchar <{\putptrianglep}{\putptrianglep
   <\arrowtypea`\arrowtypeb`\arrowtypec;\height>}}
\def\ptriangle{\@ifnextchar <{\ptrianglep}{\ptrianglep
   <\arrowtypea`\arrowtypeb`\arrowtypec;\height>}}

\def\ptrianglep<#1>[#2`#3`#4;#5`#6`#7]{{
\settriparms[#1]%
\width=\height                         
\xext=\width                           
\yext=\width                           
\topadjust[#2`#3`{#5}]
\botadjust[#3``]
\leftadjust[#2`#4`{#6}]
\rightsladjust[#3`#4`{#7}]
\begin{picture}(\xext,\yext)(\xoff,\yoff)
\putptrianglep<\arrowtypea`\arrowtypeb`\arrowtypec;\height>%
(0,0)[#2`#3`#4;#5`#6`{#7}]%
\end{picture}%
}}

\def\putqtrianglep<#1>(#2,#3)[#4`#5`#6;#7`#8`#9]{{%
\settriparms[#1]%
\xpos=#2 \ypos=#3
\advance\ypos by\height
\puthmorphism(\xpos,\ypos)[#4`#5`{#7}]{\height}{\arrowtypea}a%
\putmorphism(\xpos,\ypos)(1,-1)[``{#8}]{\height}{\arrowtypeb}l%
\advance\xpos by\height
\putvmorphism(\xpos,\ypos)[`#6`{#9}]{\height}{\arrowtypec}r%
}}

\def\putqtriangle{\@ifnextchar <{\putqtrianglep}{\putqtrianglep
   <\arrowtypea`\arrowtypeb`\arrowtypec;\height>}}
\def\qtriangle{\@ifnextchar <{\qtrianglep}{\qtrianglep
   <\arrowtypea`\arrowtypeb`\arrowtypec;\height>}}

\def\qtrianglep<#1>[#2`#3`#4;#5`#6`#7]{{
\settriparms[#1]
\width=\height                         
\xext=\width                           
\yext=\height                          
\topadjust[#2`#3`{#5}]
\botadjust[#4``]
\leftsladjust[#2`#4`{#6}]
\rightadjust[#3`#4`{#7}]
\begin{picture}(\xext,\yext)(\xoff,\yoff)
\putqtrianglep<\arrowtypea`\arrowtypeb`\arrowtypec;\height>%
(0,0)[#2`#3`#4;#5`#6`{#7}]%
\end{picture}%
}}

\def\putdtrianglep<#1>(#2,#3)[#4`#5`#6;#7`#8`#9]{{%
\settriparms[#1]%
\xpos=#2 \ypos=#3
\puthmorphism(\xpos,\ypos)[#5`#6`{#9}]{\height}{\arrowtypec}b%
\advance\xpos by \height \advance\ypos by\height
\putmorphism(\xpos,\ypos)(-1,-1)[``{#7}]{\height}{\arrowtypea}l%
\putvmorphism(\xpos,\ypos)[#4``{#8}]{\height}{\arrowtypeb}r%
}}

\def\putdtriangle{\@ifnextchar <{\putdtrianglep}{\putdtrianglep
   <\arrowtypea`\arrowtypeb`\arrowtypec;\height>}}
\def\dtriangle{\@ifnextchar <{\dtrianglep}{\dtrianglep
   <\arrowtypea`\arrowtypeb`\arrowtypec;\height>}}

\def\dtrianglep<#1>[#2`#3`#4;#5`#6`#7]{{
\settriparms[#1]
\width=\height                         
\xext=\width                           
\yext=\height                          
\topadjust[#2``]
\botadjust[#3`#4`{#7}]
\leftsladjust[#3`#2`{#5}]
\rightadjust[#2`#4`{#6}]
\begin{picture}(\xext,\yext)(\xoff,\yoff)
\putdtrianglep<\arrowtypea`\arrowtypeb`\arrowtypec;\height>%
(0,0)[#2`#3`#4;#5`#6`{#7}]%
\end{picture}%
}}

\def\putbtrianglep<#1>(#2,#3)[#4`#5`#6;#7`#8`#9]{{%
\settriparms[#1]%
\xpos=#2 \ypos=#3
\puthmorphism(\xpos,\ypos)[#5`#6`{#9}]{\height}{\arrowtypec}b%
\advance\ypos by\height
\putmorphism(\xpos,\ypos)(1,-1)[``{#8}]{\height}{\arrowtypeb}r%
\putvmorphism(\xpos,\ypos)[#4``{#7}]{\height}{\arrowtypea}l%
}}

\def\putbtriangle{\@ifnextchar <{\putbtrianglep}{\putbtrianglep
   <\arrowtypea`\arrowtypeb`\arrowtypec;\height>}}
\def\btriangle{\@ifnextchar <{\btrianglep}{\btrianglep
   <\arrowtypea`\arrowtypeb`\arrowtypec;\height>}}

\def\btrianglep<#1>[#2`#3`#4;#5`#6`#7]{{
\settriparms[#1]
\width=\height                         
\xext=\width                           
\yext=\height                          
\topadjust[#2``]
\botadjust[#3`#4`{#7}]
\leftadjust[#2`#3`{#5}]
\rightsladjust[#4`#2`{#6}]
\begin{picture}(\xext,\yext)(\xoff,\yoff)
\putbtrianglep<\arrowtypea`\arrowtypeb`\arrowtypec;\height>%
(0,0)[#2`#3`#4;#5`#6`{#7}]%
\end{picture}%
}}

\def\putAtrianglep<#1>(#2,#3)[#4`#5`#6;#7`#8`#9]{{%
\settriparms[#1]%
\xpos=#2 \ypos=#3
{\multiply \height by2
\puthmorphism(\xpos,\ypos)[#5`#6`{#9}]{\height}{\arrowtypec}b}%
\advance\xpos by\height \advance\ypos by\height
\putmorphism(\xpos,\ypos)(-1,-1)[#4``{#7}]{\height}{\arrowtypea}l%
\putmorphism(\xpos,\ypos)(1,-1)[``{#8}]{\height}{\arrowtypeb}r%
}}

\def\putAtriangle{\@ifnextchar <{\putAtrianglep}{\putAtrianglep
   <\arrowtypea`\arrowtypeb`\arrowtypec;\height>}}
\def\Atriangle{\@ifnextchar <{\Atrianglep}{\Atrianglep
   <\arrowtypea`\arrowtypeb`\arrowtypec;\height>}}

\def\Atrianglep<#1>[#2`#3`#4;#5`#6`#7]{{
\settriparms[#1]
\width=\height                         
\xext=\width                           
\yext=\height                          
\topadjust[#2``]
\botadjust[#3`#4`{#7}]
\multiply \xext by2 
\leftsladjust[#3`#2`{#5}]
\rightsladjust[#4`#2`{#6}]
\begin{picture}(\xext,\yext)(\xoff,\yoff)%
\putAtrianglep<\arrowtypea`\arrowtypeb`\arrowtypec;\height>%
(0,0)[#2`#3`#4;#5`#6`{#7}]%
\end{picture}%
}}

\def\putAtrianglepairp<#1>(#2)[#3;#4`#5`#6`#7`#8]{{
\settripairparms[#1]%
\setpos(#2)%
\settokens[#3]%
\puthmorphism(\xpos,\ypos)[\tokenb`\tokenc`{#7}]{\height}{\arrowtyped}b%
\advance\xpos by\height
\advance\ypos by\height
\putmorphism(\xpos,\ypos)(-1,-1)[\tokena``{#4}]{\height}{\arrowtypea}l%
\putvmorphism(\xpos,\ypos)[``{#5}]{\height}{\arrowtypeb}m%
\putmorphism(\xpos,\ypos)(1,-1)[``{#6}]{\height}{\arrowtypec}r%
}}

\def\putAtrianglepair{\@ifnextchar <{\putAtrianglepairp}{\putAtrianglepairp%
   <\arrowtypea`\arrowtypeb`\arrowtypec`\arrowtyped`\arrowtypee;\height>}}
\def\Atrianglepair{\@ifnextchar <{\Atrianglepairp}{\Atrianglepairp%
   <\arrowtypea`\arrowtypeb`\arrowtypec`\arrowtyped`\arrowtypee;\height>}}

\def\Atrianglepairp<#1>[#2;#3`#4`#5`#6`#7]{{%
\settripairparms[#1]%
\settokens[#2]%
\width=\height
\xext=\width
\yext=\height
\topadjust[\tokena``]%
\vertadjust[\tokenb`\tokenc`{#6}]
\tempcountd=\tempcounta                       
\vertadjust[\tokenc`\tokend`{#7}]
\ifnum\tempcounta<\tempcountd                 
\tempcounta=\tempcountd\fi                    
\advance \yext by\tempcounta                  
\advance \yoff by-\tempcounta                 %
\multiply \xext by2 
\leftsladjust[\tokenb`\tokena`{#3}]
\rightsladjust[\tokend`\tokena`{#5}]%
\begin{picture}(\xext,\yext)(\xoff,\yoff)%
\putAtrianglepairp
<\arrowtypea`\arrowtypeb`\arrowtypec`\arrowtyped`\arrowtypee;\height>%
(0,0)[#2;#3`#4`#5`#6`{#7}]%
\end{picture}%
}}

\def\putVtrianglep<#1>(#2,#3)[#4`#5`#6;#7`#8`#9]{{%
\settriparms[#1]%
\xpos=#2 \ypos=#3
\advance\ypos by\height
{\multiply\height by2
\puthmorphism(\xpos,\ypos)[#4`#5`{#7}]{\height}{\arrowtypea}a}%
\putmorphism(\xpos,\ypos)(1,-1)[`#6`{#8}]{\height}{\arrowtypeb}l%
\advance\xpos by\height
\advance\xpos by\height
\putmorphism(\xpos,\ypos)(-1,-1)[``{#9}]{\height}{\arrowtypec}r%
}}

\def\putVtriangle{\@ifnextchar <{\putVtrianglep}{\putVtrianglep
   <\arrowtypea`\arrowtypeb`\arrowtypec;\height>}}
\def\Vtriangle{\@ifnextchar <{\Vtrianglep}{\Vtrianglep
   <\arrowtypea`\arrowtypeb`\arrowtypec;\height>}}

\def\Vtrianglep<#1>[#2`#3`#4;#5`#6`#7]{{
\settriparms[#1]
\width=\height                         
\xext=\width                           
\yext=\height                          
\topadjust[#2`#3`{#5}]
\botadjust[#4``]
\multiply \xext by2 
\leftsladjust[#2`#3`{#6}]
\rightsladjust[#3`#4`{#7}]
\begin{picture}(\xext,\yext)(\xoff,\yoff)%
\putVtrianglep<\arrowtypea`\arrowtypeb`\arrowtypec;\height>%
(0,0)[#2`#3`#4;#5`#6`{#7}]%
\end{picture}%
}}

\def\putVtrianglepairp<#1>(#2)[#3;#4`#5`#6`#7`#8]{{
\settripairparms[#1]%
\setpos(#2)%
\settokens[#3]%
\advance\ypos by\height
\putmorphism(\xpos,\ypos)(1,-1)[`\tokend`{#6}]{\height}{\arrowtypec}l%
\puthmorphism(\xpos,\ypos)[\tokena`\tokenb`{#4}]{\height}{\arrowtypea}a%
\advance\xpos by\height
\putvmorphism(\xpos,\ypos)[``{#7}]{\height}{\arrowtyped}m%
\advance\xpos by\height
\putmorphism(\xpos,\ypos)(-1,-1)[``{#8}]{\height}{\arrowtypee}r%
}}

\def\putVtrianglepair{\@ifnextchar <{\putVtrianglepairp}{\putVtrianglepairp%
    <\arrowtypea`\arrowtypeb`\arrowtypec`\arrowtyped`\arrowtypee;\height>}}
\def\Vtrianglepair{\@ifnextchar <{\Vtrianglepairp}{\Vtrianglepairp%
    <\arrowtypea`\arrowtypeb`\arrowtypec`\arrowtyped`\arrowtypee;\height>}}

\def\Vtrianglepairp<#1>[#2;#3`#4`#5`#6`#7]{{%
\settripairparms[#1]%
\settokens[#2]
\xext=\height                  
\width=\height                 
\yext=\height                  
\vertadjust[\tokena`\tokenb`{#4}]
\tempcountd=\tempcounta        
\vertadjust[\tokenb`\tokenc`{#5}]
\ifnum\tempcounta<\tempcountd%
\tempcounta=\tempcountd\fi
\advance \yext by\tempcounta
\botadjust[\tokend``]%
\multiply \xext by2
\leftsladjust[\tokena`\tokend`{#6}]%
\rightsladjust[\tokenc`\tokend`{#7}]%
\begin{picture}(\xext,\yext)(\xoff,\yoff)%
\putVtrianglepairp
<\arrowtypea`\arrowtypeb`\arrowtypec`\arrowtyped`\arrowtypee;\height>%
(0,0)[#2;#3`#4`#5`#6`{#7}]%
\end{picture}%
}}

\def\putCtrianglep<#1>(#2,#3)[#4`#5`#6;#7`#8`#9]{{%
\settriparms[#1]%
\xpos=#2 \ypos=#3
\advance\ypos by\height
\putmorphism(\xpos,\ypos)(1,-1)[``{#9}]{\height}{\arrowtypec}l%
\advance\xpos by\height
\advance\ypos by\height
\putmorphism(\xpos,\ypos)(-1,-1)[#4`#5`{#7}]{\height}{\arrowtypea}l%
{\multiply\height by 2
\putvmorphism(\xpos,\ypos)[`#6`{#8}]{\height}{\arrowtypeb}r}%
}}

\def\putCtriangle{\@ifnextchar <{\putCtrianglep}{\putCtrianglep
    <\arrowtypea`\arrowtypeb`\arrowtypec;\height>}}
\def\Ctriangle{\@ifnextchar <{\Ctrianglep}{\Ctrianglep
    <\arrowtypea`\arrowtypeb`\arrowtypec;\height>}}

\def\Ctrianglep<#1>[#2`#3`#4;#5`#6`#7]{{
\settriparms[#1]
\width=\height                          
\xext=\width                            
\yext=\height                           
\multiply \yext by2 
\topadjust[#2``]
\botadjust[#4``]
\sladjust[#3`#2`{#5}]{\width}
\tempcountd=\tempcounta                 
\sladjust[#3`#4`{#7}]{\width}
\ifnum \tempcounta<\tempcountd          
\tempcounta=\tempcountd\fi              
\advance \xext by\tempcounta            
\advance \xoff by-\tempcounta           %
\rightadjust[#2`#4`{#6}]
\begin{picture}(\xext,\yext)(\xoff,\yoff)%
\putCtrianglep<\arrowtypea`\arrowtypeb`\arrowtypec;\height>%
(0,0)[#2`#3`#4;#5`#6`{#7}]%
\end{picture}%
}}

\def\putDtrianglep<#1>(#2,#3)[#4`#5`#6;#7`#8`#9]{{%
\settriparms[#1]%
\xpos=#2 \ypos=#3
\advance\xpos by\height \advance\ypos by\height
\putmorphism(\xpos,\ypos)(-1,-1)[``{#9}]{\height}{\arrowtypec}r%
\advance\xpos by-\height \advance\ypos by\height
\putmorphism(\xpos,\ypos)(1,-1)[`#5`{#8}]{\height}{\arrowtypeb}r%
{\multiply\height by 2
\putvmorphism(\xpos,\ypos)[#4`#6`{#7}]{\height}{\arrowtypea}l}%
}}

\def\putDtriangle{\@ifnextchar <{\putDtrianglep}{\putDtrianglep
    <\arrowtypea`\arrowtypeb`\arrowtypec;\height>}}
\def\Dtriangle{\@ifnextchar <{\Dtrianglep}{\Dtrianglep
   <\arrowtypea`\arrowtypeb`\arrowtypec;\height>}}

\def\Dtrianglep<#1>[#2`#3`#4;#5`#6`#7]{{
\settriparms[#1]
\width=\height                         
\xext=\height                          
\yext=\height                          
\multiply \yext by2 
\topadjust[#2``]
\botadjust[#4``]
\leftadjust[#2`#4`{#5}]
\sladjust[#3`#2`{#5}]{\height}
\tempcountd=\tempcountd                
\sladjust[#3`#4`{#7}]{\height}
\ifnum \tempcounta<\tempcountd         
\tempcounta=\tempcountd\fi             
\advance \xext by\tempcounta           %
\begin{picture}(\xext,\yext)(\xoff,\yoff)
\putDtrianglep<\arrowtypea`\arrowtypeb`\arrowtypec;\height>%
(0,0)[#2`#3`#4;#5`#6`{#7}]%
\end{picture}%
}}

\def\setrecparms[#1`#2]{\width=#1 \height=#2}%
%

\def\recursep<#1`#2>[#3;#4`#5`#6`#7`#8]{{%
\width=#1 \height=#2
\settokens[#3]
\settowidth{\tempdimen}{$\tokena$}
\ifdim\tempdimen=0pt
  \savebox{\tempboxa}{\hbox{$\tokenb$}}%
  \savebox{\tempboxb}{\hbox{$\tokend$}}%
  \savebox{\tempboxc}{\hbox{$#6$}}%
\else
  \savebox{\tempboxa}{\hbox{$\hbox{$\tokena$}\times\hbox{$\tokenb$}$}}%
  \savebox{\tempboxb}{\hbox{$\hbox{$\tokena$}\times\hbox{$\tokend$}$}}%
  \savebox{\tempboxc}{\hbox{$\hbox{$\tokena$}\times\hbox{$#6$}$}}%
\fi
\ypos=\height
\divide\ypos by 2
\xpos=\ypos
\advance\xpos by \width
\xext=\xpos \yext=\height
\topadjust[#3`\usebox{\tempboxa}`{#4}]%
\botadjust[#5`\usebox{\tempboxb}`{#8}]%
\sladjust[\tokenc`\tokenb`{#5}]{\ypos}%
\tempcountd=\tempcounta
\sladjust[\tokenc`\tokend`{#5}]{\ypos}%
\ifnum \tempcounta<\tempcountd
\tempcounta=\tempcountd\fi
\advance \xext by\tempcounta
\advance \xoff by-\tempcounta
\rightadjust[\usebox{\tempboxa}`\usebox{\tempboxb}`\usebox{\tempboxc}]%
\bfig
\putCtrianglep<-1`1`1;\ypos>(0,0)[`\tokenc`;#5`#6`{#7}]%
\puthmorphism(\ypos,0)[\tokend`\usebox{\tempboxb}`{#8}]{\width}{-1}b%
\puthmorphism(\ypos,\height)[\tokenb`\usebox{\tempboxa}`{#4}]{\width}{-1}a%
\advance\ypos by \width
\putvmorphism(\ypos,\height)[``\usebox{\tempboxc}]{\height}1r%
\efig
}}

\def\recurse{\@ifnextchar <{\recursep}{\recursep<\width`\height>}}

\def\puttwohmorphisms(#1,#2)[#3`#4;#5`#6]#7#8#9{{%
%
\puthmorphism(#1,#2)[#3`#4`]{#7}0a
\ypos=#2
\advance\ypos by 20
\puthmorphism(#1,\ypos)[\phantom{#3}`\phantom{#4}`#5]{#7}{#8}a
\advance\ypos by -40
\puthmorphism(#1,\ypos)[\phantom{#3}`\phantom{#4}`#6]{#7}{#9}b
}}

\def\puttwovmorphisms(#1,#2)[#3`#4;#5`#6]#7#8#9{{%
%
%
%
\putvmorphism(#1,#2)[#3`#4`]{#7}0a
\xpos=#1
\advance\xpos by -20
\putvmorphism(\xpos,#2)[\phantom{#3}`\phantom{#4}`#5]{#7}{#8}l
\advance\xpos by 40
\putvmorphism(\xpos,#2)[\phantom{#3}`\phantom{#4}`#6]{#7}{#9}r
}}

\def\puthcoequalizer(#1)[#2`#3`#4;#5`#6`#7]#8#9{{%
%
\setpos(#1)%
\puttwohmorphisms(\xpos,\ypos)[#2`#3;#5`#6]{#8}11%
\advance\xpos by #8
\puthmorphism(\xpos,\ypos)[\phantom{#3}`#4`#7]{#8}1{#9}
}}

\def\putvcoequalizer(#1)[#2`#3`#4;#5`#6`#7]#8#9{{%
%
%
%
%
\setpos(#1)%
\puttwovmorphisms(\xpos,\ypos)[#2`#3;#5`#6]{#8}11%
\advance\ypos by -#8
\putvmorphism(\xpos,\ypos)[\phantom{#3}`#4`#7]{#8}1{#9}
}}

\def\putthreehmorphisms(#1)[#2`#3;#4`#5`#6]#7(#8)#9{{%
\setpos(#1) \settypes(#8)
\if a#9 %
     \vertsize{\tempcounta}{#5}%
     \vertsize{\tempcountb}{#6}%
     \ifnum \tempcounta<\tempcountb \tempcounta=\tempcountb \fi
\else
     \vertsize{\tempcounta}{#4}%
     \vertsize{\tempcountb}{#5}%
     \ifnum \tempcounta<\tempcountb \tempcounta=\tempcountb \fi
\fi
\advance \tempcounta by 60
\puthmorphism(\xpos,\ypos)[#2`#3`#5]{#7}{\arrowtypeb}{#9}
\advance\ypos by \tempcounta
\puthmorphism(\xpos,\ypos)[\phantom{#2}`\phantom{#3}`#4]{#7}{\arrowtypea}{#9}
\advance\ypos by -\tempcounta \advance\ypos by -\tempcounta
\puthmorphism(\xpos,\ypos)[\phantom{#2}`\phantom{#3}`#6]{#7}{\arrowtypec}{#9}
}}

\def\putarc(#1,#2)[#3`#4`#5]#6#7#8{{%
\xpos #1
\ypos #2
\width #6
\arrowlength #6
\putbox(\xpos,\ypos){#3\vphantom{#4}}%
{\advance \xpos by\arrowlength
\putbox(\xpos,\ypos){\vphantom{#3}#4}}%
\horsize{\tempcounta}{#3}%
\horsize{\tempcountb}{#4}%
\divide \tempcounta by2
\divide \tempcountb by2
\advance \tempcounta by30
\advance \tempcountb by30
\advance \xpos by\tempcounta
\advance \arrowlength by-\tempcounta
\advance \arrowlength by-\tempcountb
\halflength=\arrowlength \divide\halflength by 2
\divide\arrowlength by 5
\put(\xpos,\ypos){\bezier{\arrowlength}(0,0)(50,50)(\halflength,50)}
\ifnum #7=-1 \put(\xpos,\ypos){\vector(-3,-2)0} \fi
\advance\xpos by \halflength
\put(\xpos,\ypos){\xpos=\halflength \advance\xpos by -50
   \bezier{\arrowlength}(0,50)(\xpos,50)(\halflength,0)}
\ifnum #7=1 {\advance \xpos by
   \halflength \put(\xpos,\ypos){\vector(3,-2)0}} \fi
\advance\ypos by 50
\vertsize{\tempcounta}{#5}%
\divide\tempcounta by2
\advance \tempcounta by20
\if a#8 %
   \advance \ypos by\tempcounta
   \putbox(\xpos,\ypos){#5}%
\else
   \advance \ypos by-\tempcounta
   \putbox(\xpos,\ypos){#5}%
\fi
}}

\makeatother


\sloppy

\newcommand{\nl}{\hspace{2cm}\\ }

\def\nec{\Box}
\def\pos{\Diamond}
\def\diam{{\tiny\Diamond}}

\def\lc{\lceil}
\def\rc{\rceil}
\def\lf{\lfloor}
\def\rf{\rfloor}
\def\lk{\langle}
\def\rk{\rangle}
\def\lse{[\!|}
\def\rse{|\!]}
\def\le{(\!|}
\def\re{|\!)}

\def\homo{{\approx\!\! >}}
\def\inn{\in\!\!\!\!\ra}
\def\dsum{\stackrel{\cdot}{\sqcup}}
\def\dsum{\sqcup\!\!\!\!\cdot\;}

\def\tl{\triangleleft}
\def\tr{\triangleright}

\def\lhb{\lhd \hspace {-1mm}\bullet}

\newcommand{\pa}{\parallel}
\newcommand{\lra}{\longrightarrow}
\newcommand{\hra}{\hookrightarrow}
\newcommand{\hla}{\hookleftarrow}
\newcommand{\ra}{\rightarrow}
\newcommand{\lla}{\longleftarrow}
\newcommand{\da}{\downarrow}
\newcommand{\ua}{\uparrow}
\newcommand{\dA}{\downarrow\!\!\!^\bullet}
\newcommand{\uA}{\uparrow\!\!\!_\bullet}
\newcommand{\Da}{\Downarrow}
\newcommand{\DA}{\Downarrow\!\!\!^\bullet}
\newcommand{\UA}{\Uparrow\!\!\!_\bullet}
\newcommand{\Ua}{\Uparrow}
\newcommand{\Lra}{\Longrightarrow}
\newcommand{\Ra}{\Rightarrow}
\newcommand{\Lla}{\Longleftarrow}
\newcommand{\La}{\Leftarrow}
\newcommand{\nperp}{\perp\!\!\!\!\!\setminus\;\;}
\newcommand{\pq}{\preceq}

\def\lrao(#1){\stackrel{#1}{{\lra}\hskip -2.8mm\circ\hskip 1.5mm}}
\def\rao(#1){\stackrel{#1}{{\ra}\hskip -2.4mm\circ\hskip 1.5mm}}
\def\raob(#1){\stackrel{#1}{{\ra}\hskip -2.8mm\circ\hskip 1.5mm}}
\def\rab(#1){\stackrel{#1}{{\ra}}}
\def\raobp{{\ra}\hskip -3.5mm\circ\hskip 1.3mm} 

\newcommand{\lms}{\longmapsto}
\newcommand{\ms}{\mapsto}
\newcommand{\subseteqnot}{\subseteq\hskip-4 mm_\not\hskip3 mm}

\newcommand{\bth}{\begin{theorem}}
\newcommand{\eth}{\end{theorem}}

\def\phi{\varphi}
\def\ve{\varepsilon}
\def\o{{\omega}}

\def\bA{{\bf A}}
\def\bM{{\bf M}}
\def\bC{{\bf C}}
\def\bI{{\bf I}}
\def\bL{{\bf L}}
\def\bT{{\bf T}}
\def\bS{{\bf S}}
\def\bD{{\bf D}}
\def\bB{{\bf B}}
\def\bW{{\bf W}}
\def\bP{{\bf P}}
\def\bX{{\bf X}}
\def\bY{{\bf Y}}
\def\ba{{\bf a}}
\def\bb{{\bf b}}
\def\bc{{\bf c}}
\def\bd{{\bf d}}
\def\bh{{\bf h}}
\def\bi{{\bf i}}
\def\bj{{\bf j}}
\def\bk{{\bf k}}
\def\bm{{\bf m}}
\def\bn{{\bf n}}
\def\bp{{\bf p}}
\def\bq{{\bf q}}
\def\be{{\bf e}}
\def\br{{\bf r}}
\def\bi{{\bf i}}
\def\bs{{\bf s}}
\def\bt{{\bf t}}
\def\b1{{\bf 1}}

\def\bio{{\mbox{\boldmath $\nu$}}}

\def\bkappa{{\mbox{\boldmath $\kappa$}}}
\def\bpi{{\mbox{\boldmath $\pi$}}}
\def\brho{{\mbox{\boldmath $\rho$}}}
\def\bvrho{{\mbox{\boldmath $\varrho$}}}
\def\bbeta{{\mbox{\boldmath $\beta$}}}
\def\btheta{{\mbox{\boldmath $\theta$}}}
\def\bmu{{\mbox{\boldmath $\mu$}}}
\def\bcomp{{\mbox{ $\diamond$}}}

\def\fk{\mbox{\tiny{\frame{$k$}}}}
\def\fz{\mbox{\tiny{\frame{$0$}}}}
\def\fj{\mbox{\tiny{\frame{$1$}}}}
\def\fd{\mbox{\tiny{\frame{$2$}}}}
\def\fkj{\mbox{\tiny{\frame{$k\!\! +\!\! 1$}}}}
\def\fkd{\mbox{\tiny{\frame{$k\!\! +\!\! 2$}}}}
\def\fnj{\mbox{\tiny{\frame{$n\!\! -\!\! 1$}}}}
\def\fn{\mbox{\scriptsize{\frame{$n$}}}}

\def\sfk{\mbox{\tiny{\frame{$_k$}}}}
\def\sfz{\mbox{\tiny{\frame{$_0$}}}}
\def\sfj{\mbox{\tiny{\frame{$_1$}}}}
\def\sfd{\mbox{\tiny{\frame{$_2$}}}}
\def\sfkj{\mbox{\tiny{\frame{$_{k\! +\! 1}$}}}}
\def\sfkd{\mbox{\tiny{\frame{$_{k\! +\! 2}$}}}}
\def\sfnj{\mbox{\tiny{\frame{$_{n\! -\! 1}$}}}}
\def\sfn{\mbox{\scriptsize{\frame{$_n$}}}}

\def\cBL{{\cal BL}}
\def\cB{{\cal B}}
\def\cA{{\cal A}}
\def\cC{{\cal C}}
\def\cD{{\cal D}}
\def\cE{{\cal E}}
\def\cF{{\cal F}}
\def\cG{{\cal G}}
\def\cI{{\cal I}}
\def\cJ{{\cal J}}
\def\cK{{\cal K}}
\def\cL{{\cal L}}
\def\cN{{\cal N}}
\def\cM{{\cal M}}
\def\cP{{\cal P}}
\def\cQ{{\cal Q}}
\def\cR{{\cal R}}
\def\cS{{\cal S}}
\def\cT{{\cal T}}
\def\cU{{\cal U}}
\def\cV{{\cal V}}
\def\cX{{\cal X}}
\def\cY{{\cal Y}}

\def\cat{{\bf Cat}}

\def\Cat{{\bf Cat}}
\def\oC{{{\omega}Cat}}
\def\kC{{kCat}}
\def\nC{{n{\bf Cat}}}
\def\njC{{(n+1){\bf Cat}}}
\def\oG{{{\omega}Gr}}
\def\mts{{MltSet}}

\def\ghg{{\bf gHg}}
\def\oghg{{\bf ogHg}}
\def\mto{{m/1}}

\def\Ope{{\bf Ope}}
\def\pOpe{{\bf pOpe}}
\def\pOpeCard{{\bf pOpeCard}}
\def\wpOpeCard{{\bf wpOpeCard}}
\def\wpfs{{\bf wpOpeCard}}

\def\Opei{{\bf Ope}_\iota}
\def\pOpei{{\bf pOpe}_\iota}
\def\Opeo{{\bf Ope}_\omega}
\def\pOpeo{{\bf pOpe}_\omega}
\def\nOpe{{\bf nOpe}}
\def\npOpe{{\bf npOpe}}

\def\pHg{{\bf pHg}}
\def\pOpeSet{{\bf pOpeSet}}

\def\wpnfs{{\bf wpOpeCard}_n}
\def\wpnjfs{{\bf wpOpeCard}_{n+1}}

\def\pOpeo{{\bf pOpe_\omega}}
\def\pOpeCardo{{\bf pOpeCard_\omega}}
\def\pnOpeo{{\bf npOpe_\omega}}

\def\Poly{{\bf Poly}} 
\def\pnPoly{{{\bf pPoly}_n}} 
\def\pnjPoly{{{\bf pPoly}_{n+1}}} 
\def\pnmjPoly{{{\bf pPoly}_{n-1}}} 
\def\pPoly{{\bf pPoly}} 
\def\pPoly{{\bf pPoly}} 
\def\CP{{\bf CatPoly}} 

\pagenumbering{arabic} \setcounter{page}{1}

\title{On positive opetopes, positive opetopic cardinals \\
and positive opetopic sets}

\author{Marek Zawadowski \\
Instytut Matematyki, Uniwersytet Warszawski\\
ul. S.Banacha 2,\\
00-913 Warszawa, Poland\\
zawado@mimuw.edu.pl} 

\maketitle
\begin{abstract}  We introduce the notion of a positive opetope and  positive opetopic cardinals as certain finite combinatorial structures. The positive opetopic cardinals to positive-to-one polygraphs are like simple graphs to free $\o$-categories over $\o$-graphs, c.f. \cite{MZ}. In particular, they allow us to give an explicit combinatorial description of positive-to-one polygraphs. Using this description we show, among other things, that positive-to-one polygraphs form a presheaf category with the exponent category being the category of positive opetopes. We also show that the category of $\o$-categories is monadic over the category of positive-to-one polygraphs with the `free functor' being the inclusion $\pPoly\ra\oC$. \end{abstract}

\tableofcontents

\section{Introduction}
In this paper we present a combinatorial description of the category of the positive-to-one polygraphs $\pPoly$. We show that this category is a presheaf category and we describe its exponent category in a combinatorial way as the category of positive opetopes $\pOpe$, see section \ref{sec-pfs}. However the proof of that requires some extended studies of the category of all positive opetopic cardinals. Intuitively, the (isomorphism classes of) positive opetopic cardinals correspond to the types of arbitrary cells in positive-to-one polygraphs. The notion of a positive opetopic cardinal is the main notion introduced in this paper. We describe in a combinatorial way the embedding functor $\be:\pPoly \ra \oC$ of the category of positive-to-one polygraphs into the category of $\o$-categories\index{category!oC@$\oC$} as the left Kan extension along a suitable functor $\bj$, and its right adjoint as the restriction along $\bj$. We end by adopting an argument due to V.Harnik \cite{H} to show that the right adjoint to $\be$ is monadic. This approach does not cover the problem of the cells with empty domains which is important for both Makkai's multitopic categories and Baez-Dolan's opetopic categories. However, it keeps something from the simplicity of the Joyal's $\theta$-categories, i.e., the category of positive opetopic cardinals with $o$-omega functors as morphisma $\pOpeCardo$ is not much more complicated than the category of simple $\o$-categories, the dual of the category of disks, c.f. \cite{Joyal}, \cite{MZ}, \cite{Berger}. In this sense this paper may be considered as a step towards a comparison of cellular and opetopic approaches.

This paper is an extended and improved version of \cite{Z}. The terminology, notation, and proofs are changed and adjusted in many cases.

\subsection*{Positive opetopic cardinals}
Positive opetopic cardinals represent all possible shapes of cells in positive-to-one polygraphs. A positive opetopic cardinal $S$ of dimension $2$ can be pictured as a figure
\begin{center} \xext=2400 \yext=750
\begin{picture}(\xext,\yext)(\xoff,\yoff)
 \settriparms[-1`1`1;300]
 \putAtriangle(300,350)[s_5`s_6`s_4;x_8`x_7`x_6]
 \put(560,420){\makebox(100,100){$\Downarrow a_2$}}
 \put(560,120){\makebox(100,100){$\Downarrow a_1$}}
 \putmorphism(0,70)(1,0)[s_7`\phantom{s_7}`x_4]{1200}{1}b
 \putmorphism(110,180)(1,1)[\phantom{s_7}`\phantom{s_3}`x_9]{70}{1}l
 \putmorphism(900,370)(1,-1)[\phantom{s_7}`\phantom{s_3}`x_5]{300}{1}r
 \putAtriangle(1750,50)[s_1`s_2`s_0;x_2`x_1`x_0]
 \put(2010,120){\makebox(100,100){$\Downarrow a_0$}}
 \put(0,600){\makebox(100,100){$S:$}}
\putmorphism(1250,50)(1,0)[s_3`\phantom{s_2}`x_3]{500}{1}a
\end{picture}
\end{center}
and a positive opetopic cardinal $T$ of dimension $3$ can be pictured as a figure
\begin{center} \xext=2400 \yext=680
\begin{picture}(\xext,\yext)(\xoff,\yoff)
 \label{figure ppof}
 \settriparms[-1`1`1;300]
 \putAtriangle(300,350)[t_2`t_3`t_1;y_5`y_4`y_3]
 \put(520,450){$\Da\! b_3$}
 \settriparms[-1`0`0;300]
 \putAtriangle(0,50)[\phantom{t_3}`t_4`;y_6``]
 \settriparms[0`1`0;300]
 \putAtriangle(600,50)[\phantom{t_1}``t_0;`y_1`]
 \putmorphism(0,50)(1,0)[\phantom{t_4}`\phantom{t_0}`y_0]{1200}{1}b
 \putmorphism(350,150)(3,1)[\phantom{t_4}`\phantom{t_0}`]{300}{1}a
 \put(600,170){$y_2$}
 \put(260,195){$\Da\! b_2$}
 \put(800,140){$\Da\! b_1$}

\put(1240,350){$\Longrightarrow$}
  \put(1250,375){\line(1,0){135}}
   \put(1280,410){$\beta$}

\settriparms[-1`1`0;300]
 \putAtriangle(1800,350)[t_2`t_3`t_1;y_5`y_4`]
 \put(2020,300){$\Da\! b_0$}
 \settriparms[-1`0`0;300]
 \putAtriangle(1500,50)[\phantom{t_3}`t_4`;y_6``]
 \settriparms[0`1`0;300]
 \putAtriangle(2100,50)[\phantom{t_1}``t_0;`y_1`]
 \putmorphism(1500,50)(1,0)[\phantom{t_4}`\phantom{t_0}`y_0]{1200}{1}b
\end{picture}
\end{center}
They have faces of various dimensions that fit together so that it makes sense to compose them in a unique way. By $S_n$ we denote the set of faces of dimension $n$ in $S$. Each face $a$ has a face $\gamma(a)$ as its codomain and a {\em non-empty set} of faces $\delta(a)$ as its domain. In $S$ above we have for $a_1$
\[ \gamma(a_1)=x_4\;\;\;{\rm and} \;\;\; \delta(a_1)=\{ x_5,x_6,x_9\}\]
and in $T$ we have for $\beta$
\[ \gamma(\beta)=b_0\;\;\;{\rm and} \;\;\; \delta(\beta)=\{ b_1,b_2,b_3\}\]
This is all the data we need. Moreover, these (necessarily finite) data satisfy four conditions (see Section \ref{sec-pfs} for details). Below we explain them in an intuitive way.

{\em Globularity.}\index{globularity} This is the main condition. It relates the sets that are obtained by double application of $\gamma$ and $\delta$. They are
\[ \gamma\gamma(a) =\gamma\delta(a)-\delta\delta(a),\hskip 15mm \delta\gamma(a) =\delta\delta(a)-\gamma\delta(a).\]
Let us look how it works for two faces $a_1$ and $\beta$. In case of the face $a_1$ we have
\[ \gamma\delta(a_1)=\{s_3,s_4,s_6 \},\hskip 15mm \delta\delta(a_1)=\{ s_4,s_6,s_7 \}\]
\[ \gamma\gamma(a_1)= s_3,\hskip 15mm \delta\gamma(a_1)=\{s_7\}\]
So we have indeed
\[ \delta\delta(a_1)-\gamma\delta(a_1) = \{ s_4,s_6,s_7 \}-\{s_3,s_4,s_6 \}= \{s_7\} = \delta\gamma(a_1) \]
\[  \gamma\delta(a_1)-\delta\delta(a_1)= \{s_3,s_4,s_6\}-\{ s_4,s_6,s_7 \}=\{s_3\} =\{ \gamma\gamma(a_1)\} \]
Similarly for the face $\beta$ we have
  \[ \gamma\gamma(\beta)= y_0,\hskip 15mm \delta\gamma(\beta)=\{y_1,y_4,y_5,y_6\}\]
  \[ \gamma\delta(\beta)=\{y_0,y_2,y_3 \},\hskip 15mm \delta\delta(\beta)=\{
  y_1,y_2,y_3,y_4,y_5,y_6 \}\]
and hence
  \[  \gamma\delta(\beta)-\delta\delta(\beta)= \{y_0,y_2,y_3 \}-\{
  y_1,y_2,y_3,y_4,y_5,y_6 \}=\{y_0\} =\{ \gamma\gamma(\beta)\} \]
  \[ \delta\delta(\beta)-\gamma\delta(\beta) =\{
  y_1,y_2,y_3,y_4,y_5,y_6 \}-\{y_0,y_2,y_3 \}=\{y_1,y_4,y_5,y_6\} = \delta\gamma(\beta) \]
As we see, in both cases $a_1$ and $\beta$ the first actual formula is a bit more baroque (due to the curly brackets around $\gamma\gamma(a_1)$, $\gamma\gamma(\beta)$) than in the globularity condition stated above. However, in the following we omit curly bracket on purpose to have a simpler notation, with the hope that it will contribute to simplicity without messing things up.

Using $\delta$'s and $\gamma$'s we can define two binary relations $<^+$ and $<^-$ on faces of the same dimension which are transitive closures of the relations $\lhd^+$ and $\lhd^-$, respectively, defined as follows: $a\lhd^+b$ holds iff there is a face $\alpha$ such that $a\in\delta(\alpha)$ and $\gamma(\alpha)=b$, and $a\lhd^- b$ holds iff $\gamma(a)\in\delta(b)$. We call $<^+$ the {\em upper order}\index{order!upper} and $<^-$ the {\em lower order}\index{order!lower}.  The following three conditions refer to these relations.

{\em Strictness.}\index{strictness} In each dimension, the relation $<^+$ is a strict order. The relation $<^+$ on $0$-dimensional faces is required to be a linear order.

{\em Disjointness.}\index{disjointness} This condition says that no two faces   can be comparable with respect to both orders $<^+$ and $<^-$.

{\em Pencil linearity.}\index{pencil!linearity} This final condition says that the sets of cells with common codomain ($\gamma${\em -pencil}) and the sets of cells that have the same distinguished cell in the domain ($\delta${\em -pencil}) are linearly   ordered by $<^+$.

The morphisms of positive opetopic cardinals are functions that preserve dimensions and operations $\gamma$ and $\delta$. The size of a positive opetopic cardinal $S$ is defined as an infinite sequence of natural numbers $size(S)=\{ size(S)_k\}_{k\in\o}=\{S_k-\delta(S_{k+1})\}_{k\in\o}$ (almost all equal $0$). We order the sequences lexicographically with higher dimensions being more important. The induction on the size of positive opetopic cardinals is a convenient way of reasoning about them. Dimension of a positive opetopic cardinal $S$ is the index of the largest non-zero number in the sequence $size(S)$. If for $k\leq dim(S)$ ($k<dim(S)$), $size(S)_k=1$, then $S$ is {\em principal} ({\em normal}).  The normal positive opetopic cardinals play role of the pasting diagrams in \cite{HMP} and the principal positive opetopic cardinals play role of the (positive) multitopes. Note that, contrary to \cite{MZ},  we do not consider either the empty-domain multitopes or the pasting diagrams. The precise connection between these two approaches will be described elsewhere. On positive opetopic cardinals we define operations of the domain, codomain, and special pushouts which play the role of composition.  With these operations (isomorphisms classes of) the positive opetopic cardinals form the terminal positive-to-one polygraph, and at the same time the monoidal globular category in the sense of Batanin.

\subsection*{Categories and functors}
We shall define the following categories
\begin{center} \xext=800 \yext=1150
\begin{picture}(\xext,\yext)(\xoff,\yoff)
 \setsqparms[1`1`1`1;800`500]
 \putsquare(0,50)[\pOpeCard`\pOpeCardo`\pPoly`\oC;\bj`(-)^*``\be]
 \setsqparms[0`1`1`0;800`500]
 \putsquare(0,550)[\pOpe`\cS`\phantom{\pOpeCard}`\phantom{\pOpeCardo};`\bi`\bk`]
\end{picture}
\end{center}
where $\pOpe$ is the category of principal positive opetopic cardinals, $\pOpeCard$ is the category of positive opetopic cardinals, $\cS$ is the category of simple categories c.f.  \cite{MZ}, $(-)^*$ is the embedding functor of positive opetopic cardinals  into positive-to-one polygraphs, $\be$ is the inclusion functor, $\pOpeCardo$ is the full image of the composition functor $(-)^*;\be$, with the non-full embedding $\bj$.

Having these functors we can form the following diagram
\begin{center} \xext=1500 \yext=1300
\begin{picture}(\xext,\yext)(\xoff,\yoff)
 \setsqparms[1`1`-1`0;1500`600]
 \putsquare(0,0)[sPb((\pOpeCard)^{op},Set)`sPb((\pOpeCardo)^{op},Set)` Set^{(\pOpe)^{op}}`sPb(\cS^{op},Set);Lan_\bj`\bi^*`Ran_\bk`]
 \put(1230,250){\makebox(100,100){$\bk^*$}}
 \putmorphism(150,600)(0,1)[``Ran_\bi]{600}{-1}r
\putmorphism(1350,600)(0,1)[``]{600}{1}a

\setsqparms[1`1`-1`0;1500`600]
 \putsquare(0,600)[\pPoly`\oC`\phantom{sPb((\pOpeCard)^{op},Set)}`
 \phantom{sPb((\pOpeCardo)^{op},Set)};\be`\widehat{(-)}`\widetilde{(-)}`]
 \put(1200,850){\makebox(100,100){$\widehat{(-)}$}}
 \putmorphism(150,1200)(0,1)[``\widetilde{(-)}]{600}{-1}r
\putmorphism(1350,1200)(0,1)[``]{600}{1}a
\putmorphism(410,520)(1,0)[``\bj^*]{640}{-1}b
\end{picture}
\end{center}
in which all the vertical arrows come in pairs and are adjoint equivalences of categories. The unexplained categories in the diagram above are: $sPb((\pOpeCard)^{op},Set)$ - the category of the special pullback preserving functors from $(\pOpeCard)^{op}$ to $Set$ and natural transformations, $sPb(\cS^{op},Set)$ - the category of the special pullback preserving  functors from $\cS^{op}$ to $Set$ and natural transformations, $sPb((\pOpeCardo)^{op},Set)$ - the category of the special pullbacks preserving  functors from $(\pOpeCardo)^{op}$ to $Set$ and natural transformations.

The functors $$\widehat{(-)}: \pPoly\ra sPb((\pOpeCard)^{op},Set)$$
and
$$\widehat{(-)}: \oC\ra sPb((\pOpeCardo)^{op},Set)$$
are defined similarly, due to the embedding
$$(-)^*:\pOpeCard\lra \pPoly$$
from the previous diagram. For a polygraph $Q$ and an $\o$-category $C$, $\widehat{Q}$ and $\widehat{C}$ are presheaves so that for a positive opetopic cardinal $S$ we have
\[ \widehat{Q} =\pPoly(S^*,Q)   \hskip 15mm \widehat{C}=\oC(S^*,C) \]
The adjoint functors $\widetilde{(-)}$ that produce $\o$-categories are slightly more complicated. They are defined in Sections \ref{sec_pPoly} and \ref{sec-nerve}. The other functors are standard. The functors $\bi^*$, $\bj^*$, $\bk^*$  are inverse image functors. $Ran_\bi$ and $Ran_\bk$ are the right Kan extensions along $\bi$, $\bk$, respectively and $Lan_\bj$ is the left Kan extension along $\bj$.

Since we have $\be;\widehat{(-)}=\widehat{(-)};Lan_\bj$, and $\widehat{(-)}$'s are equivalences of categories, the functor $Lan_\bj$ is like $\be$ but moved into a more manageable context. In fact, we have a very neat description of this functor.

\subsection*{The content}Since the paper is quite long, I describe below the content of each section to help the reading. Sections 2 and 3 introduce the notion of a positive hypergraph and positive opetopic cardinal. Section 4 is concerned with establishing what kind of inclusions hold between iterated applications of $\gamma$'s and $\delta$'s. Section 5 contains many technical statements concerning positive opetopic cardinals. All of them are there because they are needed afterwards. But it is not recommended to read the whole section at once. One of the main technical tool is so called Path Lemma \ref{fullpath}. Section 6 describes the embedding $(-)^*:\pOpeCard \ra \oC$, i.e., it's main goal is to define an $\o$-category $S^*$ for any positive opetopic cardinal $S$.  Section 7 describes rather technical but useful properties of normal positive opetopic cardinals. In section 8 we study a way we can decompose positive opetopic cardinals if they are at all decomposable. Any positive opetopic cardinal is either principal or decomposable.  This provides a way of proving the properties of positive opetopic cardinals by induction on the size. Using this in section 9 we show that the $\o$-category $S^*$ and in fact the whole functor $(-)^*$ end up in $\pPoly$. The section 10 describes the inner-outer factorization and its refinements, i.e., a further factorization of inner maps into inner-epi and inner monos. These factorizations will play an important role in describing the strongly cartesian monad (c.f. \cite{BMW}) $T_\o$ on opetopic sets for $\o$-categories and its decomposition into two simpler monads ($T_\o=T_\iota\circ T_c$), together with a distributive law combining them.

The next short section \ref{sec-terminal-polygraph} describes just what is in their title: the terminal positive-to-one polygraph as an $\o$-category in terms of positive opetopic cardinals.  Section \ref{sec-description-of-polygraphs} gives an explicit description of all the cells in a given positive-to-one polygraph with the help of positive opetopic cardinals. In other words, it describes in concrete terms the functor $\overline{(-)}^n: Set\da D_n\lra \pnPoly$.  Section \ref{sec_pPoly} establishes the equivalence of categories between $\pPoly$ and the category of presheaves over  $\pOpe$. In Section 14 the principal pullbacks are introduced and the adaptation of the Victor Harnik's argument to the opetopic context is presented. The original argument due to Harnik was expressed in a different setting and was supposed to show the monadicity of the category of all $\o$-categories and $\o$-functors $\oC$. However, this original proof contains a gap \cite{H}, and it is still open question, c.f. \cite{M} whether the  category $\oC$ is monadic over all polygraphs. Section 15 describes a full nerve functor
   \[ \widehat{(-)} : \oC \lra Set^{(\pOpeCardo)^{op}} \]
and identifies its essential image as the special pullbacks preserving functors. Section 16 describes the inclusion functor as the left Kan extension
  \[ Lan_\bj : sPb((\pOpeCard)^{op},Set)  \lra sPb((\pOpeCardo)^{op},Set) \]
with the formulas involving just coproduct (and no other colimits). This gives as a corollary the fact that $\be :   \pPoly\ra \oC$ preserves connected limits.  Then we show that the right adjoint to $Lan_\bj$
  \[ \bj^* : sPb((\pOpeCardo)^{op},Set) \lra   sPb((\pOpeCard)^{op},Set)\]
(and hence the right adjoint to $\be : \pPoly\ra \oC$) is monadic. In Appendix we recall the definition of the category of positive-to-one polygraphs.\
In section \ref{sec-more-on-monadic-adj} we describe in detail the strongly cartesian monad $T_\o$ on opetopic sets and its decomposition into two other strongly cartesian monads. We finish the introduction stating an open problem.

\vskip 2mm
{\bf Problem.} What are the full subcategories of the category of polygraphs $\cX\hookrightarrow \pPoly$ that are coreflective as subcategories of $\cX\hookrightarrow \pPoly\hookrightarrow\oC$ with coreflector $\oC\ra \cX$ being monadic?
\vskip 2mm
This paper shows that many-to-positive polygraphs form one such category.

\subsection*{Acknowledgments}

I am very grateful to Victor Harnik and Mihaly Makkai for the conversations concerning matters of this paper.

\section{Positive hypergraphs }

$\o$ is the set of natural numbers. A {\em positive hypergraph}\index{hypergraph!positive -} $S$ is a family $\{ S_k\}_{k\in \o}$ of finite sets of faces, a family of functions $\{ \gamma_k : S_{k+1}\ra S_k \}_{k \in\o }$, and a family of total relations $\{ \delta_k : S_{k+1}\ra S_k\}_{0\leq k < n}$. Moreover, $\delta_0 : S_1\ra S_0$ is a function and only finitely many among sets $\{ S_k\}_{k\in \o}$ are non-empty. As it is always clear from the context, we shall never use the indices of the functions $\gamma$ and $\delta$.

A {\em morphism of positive hypergraphs}\index{morphism!positive hypergraph -}\index{hypergraph!positive - morphism} $f:S\lra T$ is a family of functions $f_k : S_k \lra T_k$, for $k \in\o$, such that the diagrams
\begin{center}
\xext=1800 \yext=500 \adjust[`I;I`;I`;`I]
\begin{picture}(\xext,\yext)(\xoff,\yoff)
 \setsqparms[1`1`1`1;600`400]
 \putsquare(0,0)[S_{k+1}`T_{k+1}`S_k`T_k; f_{k+1}`\gamma`\gamma`f_k]
 \setsqparms[1`1`1`1;600`400]
 \putsquare(1200,0)[S_{k+1}`T_{k+1}`S_k`T_k; f_{k+1}`\delta`\delta`f_k]
\end{picture}
\end{center}
commute, for $k\in\o$. The commutation of the left hand square is the commutation of the diagram of sets of functions but in case of the right hand square we mean more than commutation of a diagram of relations, i.e., we demand that for any $a\in S_{\geq 1}$, $f_a:\delta(a)\lra \delta(f(a))$ be a bijection, where $f_a$ is the restriction of $f$ to $\delta(a)$. The category of positive\index{category!of positive hypergraphs!$\pHg$} hypergraphs is denoted by $\pHg$.

\vskip 3mm

{\em Some notions and notation.} Let $S$ be a positive hypergraph.
\begin{enumerate}
\item When convenient and does not lead to confusions, if $a\in S_k$, i.e., $a$ is  $k$-dimensional face in $S$, we sometime treat $\gamma(a)$ as an element of $S_{k-1}$ and sometimes as a subset $\{ \gamma(a) \}$ of $S_{k-1}$.
\item The dimension of $S$ is $max\{ k\in\o : S_k\neq\emptyset \}$, and it is denoted by $dim(S)$.
\item The sets of faces of different dimensions are assumed to be disjoint (i.e., $S_k\cap S_l=\emptyset$, for $k\neq l$). $S$ is also used to mean the set of all faces of $S$, i.e.,  $\bigcup_{k=0}^nS_k$; the notation $A\subseteq S$ means that $A$ is a set of some faces of $S$; $A_k=A\cap S_k$, for $k\in\o$.
\item For $a\in S_{\geq 1}$, the set $\partial(a)=\delta(a)\cup\gamma(a)$ is {\em the boundary of $A$}\index{face!boundary of -}, i.e., the set of codimension $1$ faces in $a$.
\item $S_{\geq k} = \bigcup_{i\geq k} S_i $, $\;\; S_{\leq k} = \bigcup_{i\leq k} S_i $. The set $S_{\leq k}$ is closed under $\delta$ and $\gamma$ so it is a sub-hypergraph of $S$, called $k$-truncation of $S$.

\item The image of  $A\subseteq S$ under $\delta$ and $\gamma$ will be denoted by
  \[\delta(A)=\bigcup_{a\in A}\delta(a),\hskip 5mm \gamma(A)=\{ \gamma(a) \,:\, a\in   A\},\]
respectively. In particular, $\delta\delta(a)=\bigcup_{x\in\delta(a)}\delta(x)$, $\gamma\delta(a)=\{ \gamma(x) : x\in\delta(a) \}$.
  \item $\iota(a)=\delta\delta(a)\cap\gamma\delta(a)$ is {\em the set of internal faces}\index{face!internal} of the face $a\in S_{\geq 2}$.
\item On each set $S_k$ we introduce two binary relations $<^{S_k,-}$ and $<^{S_k,+}$, called {\em lower} and {\em upper order}\index{order!upper}\index{order!lower}, respectively. We usually omit $k$ and even $S$ in the superscript.
\begin{enumerate}
  \item $<^{S_0,-}$ is the empty relation. For $k>0$, $<^{S_k,-}$ is the transitive closure of the relation $\lhd^{S_k,-}$ on $S_k$, such that $a \lhd^{S_k,-} b$ iff $\gamma (a)\in \delta(b)$.  We write $a\perp^- b$ iff either $a <^- b$ or $b <^- a$, and we write $a \leq^- b$ iff either $a=b$ or $a <^- b$.

  \item $<^{S_k,+}$ is the transitive closure of the relation $\lhd^{S_k,+}$ on $S_k$, such that $a \lhd^{S_k,+} b$ iff there is $\alpha\in S_{k+1}$, such that $a\in \delta (\alpha)$ and  $\gamma(\alpha)=b$. We write $a\perp^+ b$ iff either $a <^+ b$ or $b <^+ a$, and we write $a \leq^+ b$ iff either $a=b$ or $a <^+ b$.
   \item  $a\not\perp b$ if  both conditions $a\not\perp^+ b$ and $a\not\perp^- b$ hold.
\end{enumerate}
\item   Let $a,b\in S_k$. A {\em lower path} $a_0, \ldots , a_m$ {\em from}\index{path!lower} $a$ {\em to} $b$ in $S$ is a sequence of faces $a_0, \ldots, a_m\in S_k$ such that $a=a_0$, $b=a_m$ and for $\gamma(a_{i-1})\in\delta(a_i)$, $i=1,\ldots ,m$.
\item   Let $x,y\in S_k$. An  {\em upper path} $x,a_0, \ldots , a_m,y$ {\em from}\index{path!upper} $x$ {\em to} $y$ in $S$ is a sequence of faces $a_0, \ldots , a_m\in S_{k+1}$ such that $x\in\delta(a_0)$, $y=\gamma(a_m)$ and $\gamma(a_{i-1})\in\delta(a_i)$, for $i=1,\ldots ,m$.
 \item The iterations of $\gamma$ and $\delta$ will be denoted in two different ways. By $\gamma^k$ and $\delta^k$ we mean $k$ applications of $\gamma$  and $\delta$, respectively. By $\gamma^{(k)}$ and $\delta^{(k)}$ we mean the application as many times $\gamma$  and $\delta$, respectively, to get faces of dimension $k$.  For example, if $a\in S_5$, then $\delta^3(a)=\delta\delta\delta(a)\subseteq S_2$  and $\delta^{(3)}(a)=\delta\delta(a)\subseteq S_3$.
 \item For $l\leq k$, $a,b\in S_k$, we define $a<_lb$  iff   $\gamma^{(l)}(a)<^- \gamma^{(l)}(b)$.
 \item A face $a$ is {\em unary}\index{face!unary -} iff $\delta(a)$ is a singleton.
\end{enumerate}

We have
\begin{lemma}
If $S$ is a hypergraph and $k\in\o$, then $<^{S_{k+1},-}$ is a strict partial order iff $<^{S_k,+}$ is a strict partial order.
\end{lemma}

\section{Positive opetopic cardinals}\label{sec-pfs}

To simplify the notation, we treat both $\delta$ and $\gamma$ as functions acting on faces as well as on sets of faces, which means that sometimes we confuse elements with singletons. Clearly, both $\delta$ and $\gamma$, when considered as functions on sets of faces, are monotone.

A positive hypergraph $S$ is a {\em positive opetopic cardinal}\index{opetopic cardinal!positive -} if it is non-empty, i.e., $S_0\neq\emptyset$ and
\begin{enumerate}
 \item {\em Globularity:}\index{globularity}  for  $a\in S_{\geq 2}$:
  \[ \gamma\gamma(a) =\gamma\delta(a)-\delta\delta(a),\hskip 15mm \delta\gamma(a) =\delta\delta(a)-\gamma\delta(a);\]

  \item {\em Strictness}:\index{strictness} for $k\in\o$, the relation $<^{S_k,+}$ is a strict order; $<^{S_0,+}$ is linear;
  \item {\em Disjointness}:\index{disjointness} for $k>0$,
\[\perp^{S_k,-}\cap \perp^{S_k,+}=\emptyset\]
  \item {\em Pencil linearity}:
  \index{pencil!linearity}\index{linearity!pencil -} for any $k>0$
  and $x\in S_{k-1}$, the sets
  \[ \{ a\in S_k \; | \; x=\gamma(a) \} \;\;\;\;{\rm
  and}\;\;\;\; \{ a\in S_k \; | \; x\in \delta(a) \} \]
  are linearly ordered by $<^{S_k,+}$.
\end{enumerate}

{\em Remarks.}
\begin{enumerate}
\item The reason why we call the first condition `globularity' is that it will imply the usual globularity condition in the $\o$-categories generated by positive opetopic cardinals.
\item Note that if we were to assume that each positive opetopic cardinal has a single cell of dimension $-1$,  then linearity of $<^{S_0,+}$ would become a special case of pencil linearity.
\item  The fact that, for $x\in S_{k-1}$, the set $\{ a\in S_k \; | \; x=\gamma(a) \}$ is linearly ordered is sometimes referred to as {\em $\gamma$-linearity}\index{linearity!$\gamma$- -} of $<^{S_k,+}$, and the fact that the set $\{ a\in S_k \; | \; \; x\in \delta(a) \}$ is linearly ordered is sometimes referred to as {\em $\delta$-linearity}\index{linearity!$\delta$- -} of $<^{S_k,+}$.

\item The {\em size of positive opetopic cardinal}\index{positive opetopic cardinal!size of -} $S$ is the sequence of natural numbers $size(S)=\{ | S_n - \delta (S_{n+1})| \}_{n\in\o}$, with almost all being equal $0$. We have an order $<$ on such sequences, so that $\{ x_n \}_{n\in\o} < \{ y_n \}_{n\in\o}$ iff there is $k\in\o$ such that $x_k< y_k$ and for all $l>k$, $x_l = y_l$. This order is well founded and many facts about positive opetopic cardinals will be proven by induction on the size.

\item The {\em category of positive opetopic cardinal}\index{category!of positive opetopic cardinals!$\pOpeCard$} is the full subcategory of positive hypergraphs $\pHg$ whose objects are the positive opetopic cardinal and is denoted by $\pOpeCard$.

\item  Let $S$ be a positive opetopic cardinal.  $S$ is $k$-principal\index{positive opetopic cardinal!!principal -} iff $size(S)_l=1$, for $l\leq k$. $S$ is a {\em positive opetope}  iff $S$ is $dim(S)$-principal. $S$ is {\em normal}\index{positive opetopic cardinal!normal -} iff $S$ is $(dim(S)-1)$-principal. By $\pOpe$\index{category!pOpe@$\pOpe$} we denote the full subcategory of $\pOpeCard$ whose objects are positive opetopes.
\end{enumerate}

\section{Atlas for $\gamma$ and $\delta$}

We have an easy

\begin{lemma}\label{atlas2} Let $S$ be a positive opetopic cardinal, $a\in S_n$, $n>1$. Then
\begin{enumerate}
  \item the sets $\delta\gamma(a)$, $\iota(a)$, and $\gamma\gamma(a)$ are disjoint;
  \begin{center}
\xext=450 \yext=300
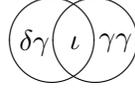
\begin{picture}(\xext,\yext)(\xoff,\yoff)
 \put(150,150){\circle{300}}
 \put(300,150){\circle{300}}
 \put(40,90){\makebox(100,100){$\delta\gamma$}}
 \put(180,90){\makebox(100,100){$\iota$}}
 \put(320,90){\makebox(100,100){$\gamma\gamma$}}
\end{picture}
\end{center}
  \item $\delta\delta(a)=\delta\gamma(a)\cup\iota(a)$;
  \item $\gamma\delta(a)=\gamma\gamma(a)\cup\iota(a)$.
\end{enumerate}
\end{lemma}

{\it Proof.}~ These are immediate consequences of globularity. $~~\Box$

Moreover

\begin{lemma}\label{atlas3} Let $S$ be a positive opetopic cardinal, $a\in S_n$, $n>2$. Then we have
\begin{enumerate}
  \item $\delta\gamma\gamma(a)\subseteq \delta\gamma\delta(a) \subseteq
  \delta\gamma\gamma(a)\cup \iota\gamma(a) = \delta\delta\gamma(a) =
  \delta\delta\delta(a)$;
  \item $\gamma\gamma\gamma(a)\subseteq \gamma\gamma\delta(a) \subseteq
  \gamma\gamma\gamma(a)\cup \iota\gamma(a) = \gamma\delta\gamma(a) =
  \gamma\delta\delta(a)$.
\end{enumerate}
\end{lemma}

{\it Proof.}~ From globularity we  have $ \gamma\gamma(\alpha)\subseteq \gamma\delta(\alpha)$. Thus by monotonicity of $\delta$ and $\gamma$ we get
\[ \gamma\gamma\gamma(\alpha)\subseteq \gamma\gamma\delta(\alpha) \hskip 5mm {\rm and}  \hskip 5mm \delta\gamma\gamma(\alpha)\subseteq \delta\gamma\delta(\alpha) \hskip 5mm {\rm and}  \hskip 5mm \gamma\gamma\delta(\alpha)\subseteq \gamma\delta\delta(\alpha).\]

Similarly, as we have from globularity: $\delta\gamma(\alpha)\subseteq \delta\delta(\alpha)$ it follows by monotonicity of $\delta$ and $\gamma$:
\[ \gamma\delta\gamma(\alpha)\subseteq \gamma\delta\delta(\alpha) \hskip 5mm {\rm and}  \hskip 5mm  \delta\delta\gamma(\alpha)\subseteq \delta\delta\delta(\alpha)
\hskip 5mm {\rm and}  \hskip 5mm \delta\gamma\delta(\alpha)\subseteq \delta\delta\delta(\alpha). \]

The equalities
\[ \delta\gamma\gamma(a)\cup \iota\gamma(a) = \delta\delta\gamma(a) \hskip 5mm {\rm and}  \hskip 5mm \gamma\gamma\gamma(a)\cup \iota\gamma(a) = \gamma\delta\gamma(a) \]
follow from Lemma \ref{atlas2}.

Thus it remains to show that:

\begin{enumerate}
  \item $\delta\delta\gamma(a) \supseteq \delta\delta\delta(a)$,
  \item $\gamma\delta\gamma(a) \supseteq \gamma\delta\delta(a)$.
\end{enumerate}
Both inclusions can be proven similarly. We shall show the first only.

Suppose to the contrary that there is $u\in \delta\delta\delta(a) - \delta\delta\gamma(a)$.  Let $x\in\delta(a)$ be $<^-$-minimal element in $\delta(a)$ such that there is $s\in\delta(x)$ with $u\in\delta(s)$. If $s\in\delta\gamma(a)$, then $u\in\delta\delta\gamma(a)$, contrary to the supposition. Thus $s\not\in\delta\gamma(a)$. Since $\delta\gamma(a)=\delta\delta(a)-\gamma\delta(a)$ it follows that $s\in\gamma\delta(a)$. Hence there is $x'\in\delta(a)$ with $\gamma(x')=s$. In particular, $x'<^-x$. Moreover
\[ u \in \delta(s)=\delta\gamma(x')\subseteq \delta\delta(x').\] Then there is $s'\in\delta(x')$ so that $u\in\delta(s')$. This contradicts the $<^-$-minimality of $x$. $~~\Box$

From Lemma \ref{atlas3} we get

\begin{corollary}\label{atlas} Let $S$ be a positive opetopic cardinal, $a\in S_n$, $n>2$, $k<n$. Then, with $\xi^{l}$ and $\xi'^{l}$ being two fixed strings of $\gamma$'s and $\delta$'s of length $l$, we have
\begin{enumerate}
  \item $\gamma^{k}(a) \subseteq \gamma\xi^{k-1}(a)$;
  \item $\delta\xi^{k-1}(a)\subseteq \delta^{k}(a)$;
  \item $\delta^{k}(a)\cap \gamma^{k}(a)=\emptyset$;
  \item $\xi^{k}(a)\subseteq \gamma^{k}(a)\cup \delta^{k}(a)$;
  \item $\delta^{2}\xi^{k-2}(a) = \delta^{2}\xi'^{k-2}(a)$, (e.g.
  $\delta^{k}(a) = \delta^{2}\gamma^{k-2}(a)$);
  \item $\gamma\delta\xi^{k-2}(a) = \gamma\delta\xi'^{k-2}(a)$,
  (e.g.  $\gamma\delta^{k-1}(a) = \gamma\delta\gamma^{k-2}(a)$);
  \item $\xi^{k-2}\delta\gamma(a) = \xi^{k-2}\delta^{2}(a)$, for $k>2$;
  \item  $\delta^{k}(a) = \delta\gamma^{k-1}(a)\cup \iota\gamma^{k-2}(a)$, for $k>1$.
\end{enumerate}
\end{corollary}

\section{Combinatorial properties of positive opetopic cardinals}
\subsection*{Local properties}
\begin{proposition}\label{cond2} Let $S$ be a positive opetopic cardinal, $k>0$ and $\alpha\in S_k$, $a_1,a_2\in\delta(\alpha)$,
$a_1\neq a_2$. Then  we have
\begin{enumerate}
  \item $a_1\not\perp^+a_2$;
  \item  $\delta(a_1)\cap\delta(a_2)=\emptyset$ and $\gamma(a_1)\neq\gamma(a_2)$.
\end{enumerate}
\end{proposition}

{\it Proof.}~Ad 1. Suppose contrary that there are $a_1,a_2\in\delta(\alpha)$ such that $a_1<^+a_2$. So we have an upper path
\[ a_1,\beta_1,\ldots,\beta_r,a_2 \]
and hence a lower path
\[  \beta_1,\ldots,\beta_r,\alpha. \]
In particular, $\beta_1<^-\alpha$. As $a_1\in\delta(\beta_1)\cap\delta(\alpha)$ by $\delta$-linearity, we have $\beta_1\perp^+\alpha$. But then $(\alpha,\beta_1)\in \perp^+\cap\perp^-\neq\emptyset$, i.e., $S$ does not satisfy the disjointness. This shows 1.

Ad 2.  This is an immediate consequence of 1. If $a_1,a_2\in\delta(\alpha)$ and either $\gamma(a_1)=\gamma(a_2)$ or $\delta(a_1)\cap\delta(a_2)\neq\emptyset$, then by pencil linearity we get that $a_1\perp^+a_2$, contradicting 1. $~~\Box$

After proving the above proposition we can introduce more notation. Let $S$ be a positive opetopic cardinal, $n\in\o$.

\begin{enumerate}
  \item For a face $\alpha\in S_{n+2}$, we shall denote by $\rho(\alpha)\in \delta(\alpha)$ the only face in $\delta(\alpha)$, such that $\gamma(\rho(\alpha))=\gamma\gamma(\alpha)$.
  \item $X\subseteq S_{n+1}$, $a,b\in S_n$ and $a,\alpha_1,\ldots ,\alpha_k, b$ be an upper path in $S$. We say that it is {\em a path in }$X$\index{path!in $X$} (or $X$-path) if $\{ \alpha_1,\ldots ,\alpha_k \}\subseteq X$.
\end{enumerate}

\begin{lemma}
\label{tech_lemma0} Let $S$ be a positive opetopic cardinal, $n\in\o$, $\alpha\in S_{n+2}$, $a,b\in S_{n+1}$,  $y\in \delta\delta(\alpha)$. Then
\begin{enumerate}
  \item there is a unique upper $\delta(\alpha)$-path from $y$ to   $\gamma\gamma(\alpha)$;
  \item there is a unique $x\in\delta\gamma(\alpha)$ and an upper $\delta(\alpha)$-path from $x$ to $y$ such that $\gamma(x)=\gamma(y)$;
 \item if $t\in\delta(y)$, there are a unique $x\in\delta\gamma(\alpha)$ such  $t\in\delta(x)$ and an upper $\delta(\alpha)$-path from $x$ to $y$;
  \item If $a<^+b$, then $\gamma(a)\leq^+\gamma(b)$.
\end{enumerate}
\end{lemma}

{\it Proof.}~ Ad 1. The uniqueness follows from Proposition \ref{cond2}.2. To show the existence, let us suppose contrary that there is no $\delta(\alpha)$-path from $y$ to $\gamma\gamma(\alpha)$. We shall construct an infinite upper $\delta(\alpha)$-path from $y$
\[ y,a_1,a_2, \ldots \]
As $y\in \delta\delta(\alpha)$, there is $a_1\in\delta(\alpha)$ such that $y\in\delta(a_1)$. So now suppose that we have already constructed $a_1,\ldots ,a_k$. By assumption $\gamma(a_k)\neq\gamma\gamma(\alpha)$ so, by globularity, $\gamma(a_k)\in\delta\delta(\alpha)$. Hence there is $a_{k+1}\in\delta(\alpha)$ such that $\gamma(a_k)\in\delta(a_{k+1})$. This ends the construction of the path.

As there are no infinite paths in positive opetopic cardinals, this is a contradiction and, in fact, there is a $\delta(\alpha)$-path from $y$ to $\gamma\gamma(\alpha)$.

Ad 2. Suppose not that there is no $x\in\delta\gamma(\alpha)$ as claimed. We shall construct an infinite descending lower $\delta(\alpha)$-path
\[ \ldots \lhd^-a_1\lhd^-a_0 \]
such that $\gamma(a_0)=y$, $\gamma\gamma(a_n)=\gamma(y)=t$, for $n\in\o$.

By assumption $y\not\in\delta\gamma(\alpha)=\delta\delta(\alpha)-\gamma\delta(\alpha)$. So $y\in\gamma\delta(\alpha)$. Hence there is $a_0\in\delta(\alpha)$ such that $\gamma(a_0)=y$. Now, suppose that the lower $\delta(\alpha)$-path
\[ a_k\lhd^-a_{k-1}\lhd^-\ldots \lhd^- a_0 \]
has been already constructed. By globularity, we can pick $z\in\delta(a_k)$ such that $\gamma(z)=t$. By assumption $z\not\in\delta\gamma(\alpha)=\delta\delta(\alpha)-\gamma\delta(\alpha)$. So $z\in\gamma\delta(\alpha)$. Hence there is $a_{k+1}\in\delta(\alpha)$ such that $\gamma(a_{k+1})=z\in\delta(a_k)$. Clearly, $\gamma\gamma(a_{k+1})=t$. This ends the construction of the path. But by strictness such a path has to be finite, so there is $x$ as needed. 

Ad 3. This case is similar. We put it for completeness.

Suppose not that there is no $x\in\delta\gamma(\alpha)$ as above. We shall construct an infinite descending lower $\delta(\alpha)$-path
\[ \ldots \lhd^-a_1\lhd^-a_0 \]
such that $\gamma(a_0)=y$, $t\in\delta\gamma(a_n)$, for $n\in\o$.

By assumption $y\not\in\delta\gamma(\alpha)=\delta\delta(\alpha)-\gamma\delta(\alpha)$. So $y\in\gamma\delta(\alpha)$. Hence there is $a_0\in\delta(\alpha)$ such that $\gamma(a_0)=y$. Now, suppose that the lower $\delta(\alpha)$-path
\[ a_k\lhd^-a_{k-1}\lhd^-\ldots \lhd^- a_0 \]
has been already constructed. By globularity, we can pick $z\in\delta(a_k)$, such that $t\in\delta(z)$. By assumption $z\not\in\delta\gamma(\alpha)=\delta\delta(\alpha)-\gamma\delta(\alpha)$. So $z\in\gamma\delta(\alpha)$. Hence there is $a_{k+1}\in\delta(\alpha)$ such that $\gamma(a_{k+1})=z\in\delta(a_k)$. Clearly, $t\in\delta\gamma(a_{k+1})$. This ends the construction of the path. But by strictness such a path has to be finite, so there is $x$ as needed.  The uniqueness again follows from Proposition \ref{cond2}.2.

Ad 4. The essential case is when $a\lhd^+b$. This follows from 1. Then use the induction on the length of the upper path from $a$ to $b$. $~~\Box$

\begin{lemma}\label{iota} Let $S$ be a positive opetopic cardinal, $n>1$, $\alpha\in S_{n+1}$, and $a,b\in S_n$ such that $a<^+b$.
Then
\begin{enumerate}
  \item $\iota\delta(\alpha)=\iota\gamma(\alpha)$ ;
  \item $\iota(a)\subseteq\iota(b)$;
  \item $\iota(a)\cup\gamma\gamma(a)\subseteq\iota(b)\cup\gamma\gamma(b)$;
  \item $\iota(a)\cup\delta\gamma(a)\subseteq\iota(b)\cup\delta\gamma(b)$;
  \item $\partial\partial(a)\subseteq\partial\partial(b)$.
\end{enumerate}
\end{lemma}

{\it Proof.}~ Ad 1. First we show $\iota\delta(\alpha)\subseteq\iota\gamma(\alpha).$ Fix $a\in\delta(\alpha)$ and $t\in\iota(a)$. Thus there are $x,y\in\delta(a)$ such that $\gamma(x)=t\in\delta(y)$. By Lemma \ref{tech_lemma0} 2,3 there are $x',y'\in\delta\gamma(\alpha)$ such that $x'\leq^+x$, $y'\leq^+y$ and $\gamma(x')=t\in\delta(y')$. Thus $t\in\iota\gamma(\alpha)$ and the first inclusion is proved.

\vskip 2mm
Now we prove the converse inclusion $\iota\delta(\alpha)\supseteq\iota\gamma(\alpha)$. Fix $t\in\iota\gamma(\alpha)$. In particular, there are $x,y\in\delta\gamma(\alpha)$, so that $\gamma(x)=t\in\delta(y)$. Suppose that $t\not\in\iota\delta(\alpha)$. We shall build an infinite $\delta(\alpha)$-path
\[  a_1\lhd ^- a_2 \ldots \]
such that $\gamma\gamma(a_i)=t$ for $i\in\o$.

Since $\delta\gamma(\alpha)\subseteq\delta\delta(\alpha)$, there is $a_1\in\delta(\alpha)$ such that $x\in\delta(a_1)$. Since $t\not\in\iota\delta(\alpha)$, it follows that $\gamma\gamma(a_1)=t$.  Suppose now that we have already
constructed the path
\[  a_1\lhd^-a_2 \lhd^-\ldots \lhd^- a_k\]
with the stated properties. We have $\gamma\gamma(a_k)=t\lhd^+\gamma(y)\leq^+\gamma\gamma\gamma(\alpha)$. So, by strictness, $\gamma(a_k)\neq\gamma\gamma(\alpha)$ and $\gamma(a_k)\in\delta\delta(\alpha)$. Then there is $a_{k+1}\in\delta(\alpha)$ such that $\gamma(a_k)\in\delta(a_{k+1})$. Again, as $t\not\in\iota\delta(\alpha)$, it follows that $\gamma\gamma(a_{k+1})=t$. This ends the construction of the path. Since, by strictness, such a path cannot exist, we get the other inclusion.

Ad 2.  Since the inclusion is transitive, it is enough to consider the case $a\lhd^+b$, i.e., there is an $\alpha\in S_{n+1}$ such that $a\in\delta(\alpha)$ and  $b=\gamma(\alpha)$. Then by 1. we have
\[ \iota(a)\subseteq
\iota\delta(\alpha)=\iota\gamma(\alpha)=\iota(b) \]

 Ad 3. Again it is enough to consider the case $a\lhd^+b$, i.e., that there is $\alpha\in S_{n+1}$ such that $a\in\delta(\alpha)$ and $\gamma(\alpha)=b$. By 2. we need to show that $\gamma\gamma(a)\in\iota(b)\cup\gamma\gamma(b)$. Using Lemma \ref{atlas3}.2 we have
\[ \gamma\gamma(a)\in\gamma\gamma\delta(\alpha)\subseteq \iota\gamma(\alpha)\cup\gamma\gamma\gamma(\alpha)=\iota(b)\cup\gamma\gamma(b) .\]

 Ad 4. This is similar to 3 and uses Lemma \ref{atlas3}.1.

 Ad 5. This follows from 3. and 4 and Lemma \ref{atlas2}.  $~~\Box$

\subsection*{Global properties}

\begin{lemma} \label{tech_lemma1} Let $S$ be a positive opetopic cardinal, $n\in\o$, $a,b\in S_n$, $a<^+b$. Then, there is an upper $S_{n+1}-\gamma(S_{n+2})$-path from $a$ to $b$.
\end{lemma}

{\it Proof.}~  Let $a,\alpha_1,\ldots ,\alpha_k, b$ be an upper path is $S$.  By Lemma \ref{tech_lemma0} we can replace each face $\alpha_i$ in this path which is not in $S-\gamma(S)$  by a sequence of faces which are $<^+$-smaller. Just take $\Gamma\in S_{n+2}$, such that $\gamma(\Gamma)=\alpha_i$  and take instead of $\alpha_i$ a path in $\delta(\Gamma)$ from $\gamma(\alpha_{i-1})$ (if $i=0$ then from $a$) to $\gamma(\alpha_i)$. Repeated application of this procedure will eventually yield the required path. $~~\Box$

\begin{lemma}\label{tech_lemma2} Let $S$ be a positive opetopic cardinal, $n>0$, $a\in S_n$, $\alpha\in S_{n+1}$, and either $\gamma(a)\in\iota(\alpha)$ or $\delta(a)\cap\iota(\alpha)\neq\emptyset$.  Then $a<^+\gamma(\alpha)$. Moreover, if $\alpha\in S-\gamma(S)$, then there is a unique $a'\in\delta(\alpha)$ such that $a\leq^+ a'$. \end{lemma}

{\it Proof.}~ If $a\in\delta(\alpha)$, there is nothing to prove. So we assume that $a\not\in\delta(\alpha)$.

We begin with  the second part of Lemma, i.e., we assume that $\alpha\in S_{n+1}-\gamma(S_{n+2})$.

Let $\gamma(a)\in\iota(\alpha)$. Thus there are $b,c\in\delta(\alpha)$ such that $\gamma(a)=\gamma(b)\in\delta(c)$. In particular, $a<^-c$. By $\gamma$-linearity either $b<^+a$ or  $a<^+b$. Suppose that $b<^+a$. Then we have an $(S-\gamma(S))$-upper path $b,\beta_0,\ldots,\beta_r, a$. As $b\in \alpha\cap\beta_0$ and $\alpha,\beta_0\in S-\gamma(S)$, we have $\alpha=\beta_0$. But then $c\in\delta(\alpha)=\delta(\beta_0)$ and hence $c<^+\gamma(\beta_0)\leq^+a$. So we get $a<^-c$ and $c<^+a$
contradicting the disjointness of $\perp^+$ and $\perp^-$. Thus we can put $a'=b$ and we have $a<^+a'$. The uniqueness of $a'$ follows from the fact that $\gamma(a)=\gamma(a')$.

The case $\delta(a)\cap\iota(\alpha)\neq\emptyset$ is similar and we put it for completeness. Thus there are $b,c\in\delta(\alpha)$ such that $\gamma(b)\in\delta(a)\cap\delta(c)$. In particular, $b<^-a$. By $\delta$-linearity either $c<^+a$ or  $a<^+c$. Suppose that $c<^+a$. Then we have an $(S-\gamma(S))$-upper path $c,\beta_0,\ldots,\beta_r, a$. As $c\in\alpha\cap\beta_0$ and $\alpha,\beta_0\in S-\gamma(S)$, we have $\alpha=\beta_0$. But then $b\in\delta(\alpha)=\delta(\beta_0)$ and hence $b<^+\gamma(\beta_0)\leq^+a$. So we get $b<^-a$ and $b<^+a$ contradicting the disjointness of $\perp^+$ and $\perp^-$. Thus we can put $a'=c$ and we have $a<^+a'$. The uniqueness of $a'$ follows from the fact that $\gamma(b)\in\delta(a')$ and $a'\in\delta(\alpha)$ and Proposition \ref{cond2}.

The first part of the Lemma follows from the above, Lemma \ref{tech_lemma0}.4 and the following Claim.

{\em Claim.} If $\alpha\in S_{n+1}$ and $x\in\iota(\alpha)$, then there is an $\alpha'\in S_{n+1}$ such that $\alpha'\leq^+\alpha$, $x\in\iota(\alpha')$ and $\alpha'\not\in \gamma(S_{n+2})$.

{\it Proof of the Claim.} Suppose that Claim is not true. To get a contradiction, we shall build an infinite descending $\gamma(S_{n+2})$-path
\[ \ldots \lhd^+\alpha_1\lhd^+\alpha_0=\alpha \]
such that $x\in\iota(\alpha_i)$, for $i\in\o$.

We put $\alpha_0=\alpha$. Suppose that we have already constructed $\alpha_0,\ldots, \alpha_k\in \gamma(S_{n+2})$. Hence there is $\beta\in S_{n+2}$ such that $\gamma(\beta)=\alpha_k$. Since, by Lemma \ref{iota}.1,  $\iota\delta(\beta)=\iota\gamma(\beta)=\iota(\alpha_k)$, there is $\alpha_{k+1}\in\delta(\beta)$ such that $x\in\iota(\alpha_{k+1})$. This ends the construction of the infinite path and the proof of the Claim and the Lemma. $~~\Box$

\begin{corollary}\label{gadeio} Let $S$ be a positive opetopic cardinal. If $a\in S-\delta(S)$, then $\gamma(a)\in S-\iota(S)$ and $\delta(a)\subseteq S-\iota(S)$.
\end{corollary}

{\it Proof.}~Let $a\in S_n$ and $\alpha\in S_{n+2}$.  If either $\gamma(a)\in\iota(\alpha)$ or $\delta(a)\cap\iota(\alpha)\neq\emptyset$, then by Lemma \ref{tech_lemma2} $a<^+\gamma(\alpha)$. Thus $a\in \delta(S)$. $~~\Box$
\vskip 2mm
A lower path $b_0,\ldots, b_m$ is a {\em maximal path}\index{path!maximal} if $\delta(b_0)\subseteq\delta(S)-\gamma(S)$ and $\gamma(b_m)\in\gamma(S)-\delta(S)$, i.e., if it can't be extended either way.

\begin{lemma}[Path Lemma]\label{fullpath} Let $k\geq 0$, $B=(a_0,\ldots, a_k )$ be a maximal lower $S_n$-path in a positive opetopic cardinal $S$,  $b\in S_n$, $0\leq s\leq k$, $a_s<^+b$. Then there are  $0\leq l\leq s\leq p \leq k$ such that
\begin{enumerate}
   \item $a_i<^+b$ for $i= l, \ldots , p$;
   \item $\gamma(a_p)=\gamma(b)$;
  \item either $l=0$ and $\delta(a_0)\subseteq \delta(b)$ or $l>0$ and $\gamma(a_{l-1})\in\delta(b)$;
  \item $\gamma(a_i)\in\iota(S)$, for $l\leq i <p$.
\end{enumerate}
\end{lemma}

{\it Proof.}~ Let $0\leq l\leq p\leq k$ be such that $a_i<^+b$, for $l\leq i\leq p$ and either $l=0$ or $a_{l-1}\not<^+b$ and either $p=k$ or $a_{p+1}\not<^+b$. We shall show that $l$ and $p$ have the properties stated in the Lemma. From the very definition the property $1$ holds.

We shall show 2. Take an upper $(S-\gamma(S))$-path from $a_p$ to $b$: $a_p,\beta_0,\ldots ,\beta_r,b$. If $\gamma(a_p)=\gamma\gamma(\beta_i)$, for $i=0,\ldots ,r$, then $\gamma(a_p)=\gamma\gamma(\beta_r)=\gamma(b)$ and we are done. So suppose contrary and let
\[ i_0=min\{ i : \gamma(a_p)\neq\gamma\gamma(\beta_i)\} \]
Then there are $a,c\in\delta(\beta_{i_0})$ such that $\gamma(a_p)=\gamma(a)\in\delta(c)$ (NB. $a=a_p$ if $i_0=0$ and $a=\gamma(\beta_{i_0-1})$, otherwise). In particular, $\gamma(a_p)\in\iota(\beta_{i_0})$. As $\gamma(a_p)\in\delta(S)$, we have $p<k$. Thus $\gamma(a_p)\in\delta(a_{p+1})\cap\iota(\beta_{i_0})$, and by Lemma \ref{tech_lemma2} $a_{p+1}<^+c<^+b$. But this contradicts the choice of $p$. So the property 2. holds.

Now we shall show 3. Take an upper $(S-\gamma(S))$-path from $a_l$ to $b$: $a_l,\beta_0,\ldots ,\beta_r,b$. We have two cases: $l=0$ and $l>0$.

If $l=0$, then there is no face $a\in S$ such that $\gamma(a)\in \delta(a_l)$. As $\delta(a_l)\subseteq\delta\delta(\beta_0)$, we must have $\delta(a_l)\subseteq\delta\gamma(\beta_i)$, for $i=0,\ldots , r$. Hence $\delta(a_l)\subseteq \delta\gamma(\beta_r)=\delta(b)$ and 3. holds in this case.

Now suppose that $l>0$.  If $\gamma(a_{l-1})\in\delta\gamma(\beta_i)$, for $i=0,\ldots , r$, then $\gamma(a_{l-1})\in \delta\gamma(\beta_r)=\delta(b)$ and 3. holds again. So suppose contrary , and let
\[ i_1 = min\{ i : \gamma(a_{l-1})\not\in\delta\gamma(\beta_i) \} \]
Then there are $a,c\in\delta(\beta_{i_1})$ such that $\gamma(a_{l-1})=\gamma(a)\in\delta(c)$ (NB: $c=a_l$ if $i_1=0$ and $c=\gamma(\beta_{i_1-1})$ otherwise). In particular, $\gamma(a_{l-1})\in\iota(\beta_{i_1})$, and by Lemma \ref{tech_lemma2} we have $a_{l-1}<^+a<^+b$ contrary to the choice of $l$. Thus 3. holds in this case as well.

Finally, we shall show 4. Let $l\leq j < p$ and $a_j, \beta_0,\ldots, \beta_r,b$ be an upper $(S-\delta(S))$-path from $a_j$ to $b$. As $a_j<^-a_p$ and $a_p<^+b$, we have $\gamma(a_j)\neq\gamma(b)$.  So we can put
\[ i_2 = min\{ i : \gamma(a_j)\neq\gamma\gamma(\beta_i) \} \]
But then $\gamma(a_j)\in\gamma\delta(\beta_{i_2})-\gamma\gamma(\beta_{i_2})=\iota(\beta_{i_2})$. Therefore $\gamma(a_j)\in\iota(S)$ and 4. holds. $~~\Box$

\begin{lemma}
\label{tech_lemma1'} Let $S$ be a positive opetopic cardinal, $n\in\o$, $x,y\in S_n$, $x<^+y$. If $x,y\not\in \iota(S_{n+2})$, then there is an upper $S_{n+1}-\delta(S_{n+2})$-path from $x$ to $y$.
\end{lemma}

{\it Proof.}~ Assume $x,y\in(S-\iota(S))$ and $x<^+y$. Let
\[ x,b_0,\ldots ,b_k,y \]
be an upper path from $x$ to $y$ with the longest possible initial segment $b_0,\ldots ,b_l$ in $S-\delta(S)$. As $x\not\in \iota(S_{n+2})$, such a non-empty path exists. We need to show that $k=l$.

Let $a$ be the $<^+$-largest element of the set $\{ b\in S : \gamma(b_l)\in\delta(b)\}$. Then $b_{l+1}\leq^+ a$.
and $a \not\in\delta(S)$. Since $y\not\in\iota(S)$, by Lemma \ref{fullpath}.4 there is $p$ such that $l+1\leq p\leq k$ such that $\gamma(b_p)=\gamma(a)$. Thus we have an upper path from $x$ to $y$, $x,b_0,\ldots,b_l,a,b_{p+1}, \ldots , b_k, y$ with a longer initial segment in $S-\delta(S)$.  But this is a contradiction with the choice of the path $ x,b_0,\ldots ,b_k,y$, and it means that in fact $l=k$, as required.  $~~\Box$

\subsection*{Order}

\begin{lemma}
\label{tech_lemma3} Let $S$ be a positive opetopic cardinal, $n\in\o$, $a,b\in S_n$. Then we have
\begin{enumerate}
  \item\label{tl1} If $a<^+b$, then for any   $x\in\delta(a)$ there is $y\in\delta(b)$ such that $y\leq^+x$;
  \item\label{tl2} If $a<^+b$ and $\gamma(a)=\gamma(b)$, then for any $y\in\delta(b)$ there is $x\in\delta(a)$ such that $y\leq^+x$;
  \item\label{tl3} If $\gamma(a)=\gamma(b)$, then either $a=b$ or $a\perp^+b$;
  \item\label{tl4} If $\gamma(a)<^+\gamma(b)$ then either $a<^+b$ or $a<^-b$;
  \item\label{tl5} If $a<^+b$ then $\gamma(a)\leq^+\gamma(b)$;
  \item\label{tl6} If $a<^-b$ then $\gamma(a)<^+\gamma(b)$;
  \item\label{tl7} If $\gamma(a)\perp^- \gamma(b)$ then  $a\not\perp^-b$ and $a\not\perp^+b$.
\end{enumerate}
\end{lemma}

{\it Proof.}~ Ad \ref{tl1}. Let $a<^+b$ and $x\in\delta(a)$. We have two cases: either $x\in\gamma(S)$ or $x\not\in\gamma(S)$.

In the first case there is $a'\in S-\gamma(S)$ such that $\gamma(a')=x$. Let $a_0,\ldots,a_k$ be a maximal path containing $a',a$, say $a_{s-1}=a'$ and $a_s=a$, where $0<s\leq k$. As $a_s<^+b$, by Lemma \ref{fullpath} there is $l\leq s$ and $y\in\delta(a_l)\cap\delta(b)$. Clearly, $y\leq^+x$.

In the second case consider an upper path from $a$ to $b$: $a,\beta_0,\ldots,\beta_r,b$. We have $x\in\delta(a)\subseteq\delta\delta(\beta_0)$. As $x\not\in\gamma(S)$ so $x\not\in\gamma\delta(\beta_0)$, and hence $x\in\delta\delta(\beta_0)-\gamma\delta(\beta_0)=\delta\gamma(\beta_0)$. Thus we can define
\[ r'=max\{ i : x\in\delta\gamma(\beta_i)  \} \]
If we had $r'<r$, then again we would have $x\in\delta\delta(\beta_{r'+1})-\gamma\delta(\beta_{r'+1})= \delta\gamma(\beta_{r'+1})$, contrary to the choice of $r'$. So $r'=r$ and $x\in\delta\gamma(\beta_r)=\delta(b)$. Thus we can put $y=x$.

 Ad \ref{tl2}. Fix  $a<^+b$ such that $\gamma(a)=\gamma(b)$ and $y\in\delta(b)$.  We need to find $x\in\delta(a)$ with $y\leq^+x$. Take a maximal $(S-\gamma(S))$-path $a_0,\ldots,a_k$ passing through $y$, i.e., there is $0\leq j\leq k$ such that $y\in\delta(a_j)$ and if $y\in \gamma(S)$, then moreover $j>0$ and $\gamma(a_{j-1})=y$. Since $a_j\not\in\gamma(S)$ by $\delta$-linearity $a_j<^+b$. Thus by Lemma \ref{fullpath} there is $j\leq p\leq k$ such that $\gamma(a_p)=\gamma(b)=\gamma(a)$. Since $a_p\not\in\gamma(S)$ by $\gamma$-linearity we have $a_p\leq^+a$. If $a_p=a$, then we can take as the face $x$ either $y$ if $p=0$ or $\gamma(a_{p-1})$ if $p>0$.  So assume now that $a_p<^+a$. Again by Lemma \ref{fullpath} there is $0\leq l\leq p$ such that either $l=0$ and $\delta(a_0)\subseteq\delta(a)$ or $l>0$ and $\gamma(a_{l-1})\in\delta(a)$. As $a_l$ is the first face in the path $a_0,\ldots,a_k$ such that $a_l<^+a$ and $a_j$ is the first face in the path $a_0,\ldots,a_k$ such that $a_j<^+b$ and moreover $a<^+b$, it follows that $j\leq l$.  Thus in this case we can take as the face $x$ either $y$ if $l=0$ or $\gamma(a_{l-1})$ if $l>0$.

 Ad \ref{tl3}. This is an immediate consequence of $\gamma$-linearity.

 Ad \ref{tl4}. Suppose $\gamma(a)<^+\gamma(b)$. So there is an upper path
\[ \gamma(a),c_1,\ldots , c_k, \gamma(b) \]
with $k>0$. We put $c_0=a$. We have $\gamma(c_k)=\gamma(b)$ so by $\gamma$-linearity $c_k\perp^+ b$ or $c_k=b$. In the later case $a<^-b$. In the former case, we have two possibilities: either $b<^+c_k$ or $c_k<^+b$.

If $b<^+c_k$, then by Lemma \ref{fullpath} for any maximal path that contains $b$ and the face $c_k$ we get that $c_{k-1}<^-b$. Thus we have $a<^-b$.

If $c_k<^+b$, then by Lemma \ref{fullpath} for any maximal path that extends $c_0,c_1,\ldots, c_k$ and face $b$ we get that either there is $0\leq i<k$ such that $\gamma(c_i)\in\delta(b)$ and then $a<^-b$ or else $a=c_0<^+b$.

 Ad \ref{tl5}. This is repeated from Lemma \ref{tech_lemma0}.

  Ad \ref{tl6}. Suppose $a<^-b$.  Then there is a lower path
 \[ a=a_0,a_1, \ldots , a_k =b \]
with $k>0$. Then we have an upper path
\[ \gamma(a)=\gamma(a_0),a_1, \ldots , a_k, \gamma(a_k) =\gamma(b). \]
Hence  $\gamma(a)<^+\gamma(b)$.

 Ad \ref{tl7}. Easily follows from \ref{tl5} and \ref{tl6}.
$~~\Box$

\begin{proposition}\label{comparison}
Let $S$ be a positive opetopic cardinal, $a,b\in S_n$, $a\neq b$. Let $\{ a_i\}_{0\leq i\leq n}$, $\{ b_i\}_{0\leq i\leq n}$ be two sequences of codomains of $a$ and $b$, respectively, so that
\[ a_i=\gamma^{(i)}(a)  \hskip 20mm b_i=\gamma^{(i)}(b)\]
(i.e., $dim(a_i)=i$), for $i=0,\ldots ,n$. Then there are two numbers $0\leq l\leq k \leq n$ such that either
\begin{enumerate}
  \item $a_i=b_i$, for $i<l$,
  \item $a_i<^+b_i$, for $l\leq i \leq k$,
  \item $a_i<^-b_i$, for $k+1 = i \leq n$,
  \item $a_i \not\perp b_i$, for $k+2\leq i \leq n$,
\end{enumerate}
or the roles of $a$ and $b$ are interchanged.
\end{proposition}

{\it Proof.}~ We can present the above conditions  more visually
as:
\[ a_0=b_0,\ldots , a_{l-1}=b_{l-1},\hskip 5mm a_{l}<^+b_{l}, \ldots  a_{k}<^+b_{k}, \]
\[  a_{k+1}<^-b_{k+1},\hskip 5mm  a_{k+2}\not\perp b_{k+2}, \ldots , a_{n}\not\perp b_{n}. \]
These conditions we will verify from the bottom up. Note that by strictness $<^{S_0,+}$ is a linear order. So either $a_0=b_0$ or $a_0\perp^+b_0$. In the later case $l=0$. As $a\neq b$, then there is $i\leq n$ such that $a_i\neq b_i$. Let $l$ be minimal such, i.e.,  $l=min\{ i : a_i\neq b_i \}$. By Lemma \ref{tech_lemma3} 3., $a_l\perp^+b_l$.  So assume that $a_l<^+b_l$. We put $k = max \{ i\leq n : a_i<^+b_i \}$. If $k=n$, we are done. If $k<n$, then by Lemma \ref{tech_lemma3} 4., we have $a_{k+1}<^-b_{k+1}$. Then if $k+1<n$, by Lemma \ref{tech_lemma3} 5. 6. 7., $a_i \not\perp b_i$ for $k+2\leq i \leq n$. This ends the proof. $~~\Box$

Having Proposition \ref{comparison} we can define a relation $<_l$ on $k$-faces of any positive opetopic cardinal S, $l<k$, as follows. For $a,b\in S_k$, $a<_l b$ iff $\gamma^{(l)}(a)<^- \gamma^{(l)}(b)$. We get immediately

\begin{corollary}\label{order} Let $S$ be a positive opetopic cardinal, $a,b\in S_n$, $a\neq b$. Then either $a\perp^+b$ or there is a unique $0\leq l\leq k$ such that  $a\perp^-_lb$, but not both.
\end{corollary}

The above Corollary allows us to define an order $<^S$ (also denoted $<$) on all cells of $S$ as follows. For $a,b\in S_n$,

\[ a<^S b \;\;\; {\rm iff}\;\;\; a<^+b \;\; {\rm or}\;\; \exists_l \;\; a<^-_l b. \]

\begin{corollary}
\label{linord} For any positive opetopic cardinal $S$, and $n\in\o$, the relation $<^S$ restricted to $S_n$ is a linear order.
\end{corollary}

{\it Proof.}~ We need to verify that $<^S$ is transitive.

Let $a,b,c\in S_n$. There are some cases to consider.

If $a<^+b<^+ c$, then clearly $a<^+c$.

If $a<^+b<^-_l c$, then, by Lemma \ref{tech_lemma0}.4., we have $\gamma^{(l)}(a)\leq^+ \gamma^{(l)}(b)<^- \gamma^{(l)}(c)$, and by transitivity of $<^-$ we have $\gamma^{(l)}(a)<^- \gamma^{(l)}(c)$. Hence $a<^-_l c$.

Now assume that $a<^-_lb<^+ c$ and consider a lower path in $S_l$ containing $\gamma^{(l)}(a)$ and $\gamma^{(l)}(b)$. By Lemma \ref{tech_lemma3}.\ref{tl5} $\gamma^{(l)}(b)<^+ \gamma^{(l)}(c)$, and hence by Lemma \ref{fullpath}, either $\gamma^{(l)}(a)<^+ \gamma^{(l)}(c)$ or $\gamma^{(l)}(a)<^- \gamma^{(l)}(c)$. In the later case, by transitivity of $<_l$ we have $a<_lc$, and we are done. In the former case, by Proposition \ref{comparison}, either $a<^+c$ and we are done, or there is $k>l$ such that $\gamma^{(k)}(a)<^- \gamma^{(k)}(c)$, i.e. $a<_kc$, as required.

The last case $a<^-_kb<^-_l c$ has three subcases.

If $k=l$, then clearly $a<_l c$.

If $k>l$, then $\gamma^{(l)}(a)\leq^+\gamma^{(l)}(b)<^- \gamma^{(l)}(c)$ and, by the previous argument, $\gamma^{(l)}(a)<^- \gamma^{(l)}(c)$, i.e., $a<^-_lc$.

Finally, assume that $k<l$. Then $\gamma^{(k)}(a)<^-\gamma^{(k)}(b)<^+ \gamma^{(k)}(c)$.  By Path Lemma, either $\gamma^{(k)}(a)<^-\gamma^{(k)}(c)$ or $\gamma^{(k)}(a)<^+\gamma^{(k)}(c)$. In the former case we are done. In the latter case, by Proposition \ref{comparison}, either $a<^+c$ or there is $k'$, such that $k<k'\leq n$ and $\gamma^{(k')}(a)<^+\gamma^{(k')}(c)$, as required.  $~~\Box$

From the proof of the above corollary we get

\begin{lemma}\label{le_linearity}
Let $S$ be a positive opetopic cardinal, $a\in S_n$. Then the set
\[ \{ b\in S_n\, :\, a\leq^+ b   \} \]
is linearly ordered by $\leq^+ $.
\end{lemma}

{\it Proof.}~ Suppose $a\leq^+ b,b'$. If we were to have $b<^-_l b'$ for some $l\leq n$ then, by Corollary \ref{linord} we would have $a<^-_l b'$ which is a contradiction. $~~\Box$

\begin{corollary}\label{mono} Any morphism of positive opetopic cardinals is one-to-one. Moreover, any automorphism of positive opetopic cardinals is an identity.
\end{corollary}

{\it Proof.}~ By Corollary \ref{linord}, the (strict, linear in each dimension) order $<^S$ is defined internally using relations $<^-$ and $<^+$ that are preserved by any morphism. Hence $<^S$ must be preserved by any morphism, as well. From this observation the Corollary follows. $~~\Box$

\begin{lemma}\label{LocLin2}
Let $S$ be a positive opetopic cardinal, $a,b\in S_n$. Then
\begin{enumerate}
  \item if $\iota(a)\cap\iota(b)\neq\emptyset$, then $a\perp^+b$;
  \item if $\emptyset\neq\iota(a)\subset\iota(b)\neq\iota(a)$, then $a<^+b$;
  \item if $a\perp^-b$, then  $\iota(a)\cap\iota(b)=\emptyset$.
\end{enumerate}
\end{lemma}

{\it Proof.}~ 2. is an easy consequence of 1. and Lemma \ref{iota}. 3. is an easy consequence of 1. and Disjointness. We shall show 1.

Assume that $u\in\iota(a)\cap\iota(b)$. Thus there are $x,y\in\delta(a)$ and $x',y'\in\delta(b)$ such that $\gamma(x)=\gamma(x')=u\in\delta(y)\cap\delta(y')$. If $x=x'$, then by pencil linearity $a\perp^+b$, as required.  So assume that $x\neq x'$. Again by pencil linearity $x\perp^+x'$, say $x'<^+x$. Thus there is an upper $T-\gamma(T)$-path $x',a_1,\ldots,a_k,x$. As, for $i=1,\ldots,k$, $\gamma\gamma(a_i)=u$ and $\gamma\gamma(b)\not\in\iota(b)\ni u$, we have that $\gamma(a_i)\neq\gamma(b)$ and $a_i\neq b$. Once again by pencil linearity $a_1\perp^+b$ and by Path Lemma $a_i<b$, for $i=1,\ldots,k$ with $\gamma(a_k)\neq\gamma(b)$. As $\gamma(a_k)=x\in\delta(a)$, again by Path Lemma $a<^+b$, as well. $~~\Box$

\begin{proposition}\label{successor} Let $S$ be a positive opetopic cardinal, $a,b\in S_k$, $\alpha\in S_{k+1}$, so that $\alpha$ is a $<^+$-minimal element in $S_{k+1}$, and $a\in\delta(\alpha)$, $b=\gamma(\alpha)$. Then $b$ is the $<^+$-successor of $a$.
\end{proposition}

{\it Proof.}~   Assume that $\alpha$ is a $<^+$-minimal element in $S_{k+1}$. Suppose that there is $c\in S_k$ such that $a<^+c<^+b$. Thus we have an upper path
\[ a,\beta_1,\ldots ,\beta_i,c,\beta_{i+1},\ldots ,\beta_l,b. \]
Hence $\beta_1<^-\beta_l$. Moreover, $a\in\delta(\beta_1)\cap\delta(\alpha)$ and $\gamma(\beta_l)=b=\gamma(\alpha)$. Thus both $\beta_1$ and $\beta_l$ are $<^+$-comparable with $\alpha$.  Since $\alpha$ is $<^+$-minimal, we have $\alpha<^+\beta_1, \beta_l$. By Lemma \ref{le_linearity}, $\beta_1\perp^+\beta_l$. But then we have $(\beta_1,\beta_l)\in \perp^+\cap\perp^-\neq\emptyset$, contradicting disjointness. $~~\Box$

\begin{proposition}\label{subhyp}
Let $T$ be a positive opetopic cardinal and $X\subseteq T$ a subhypergraph of $T$. Then $X$ is a positive opetopic cardinal iff the relation $<^{X,+}$ is the restriction of $<^{T,+}$ to $X$.
\end{proposition}

{\it Proof.}~ Assume that $X$ is a subhypergraph of a positive opetopic cardinal $T$. Then $X$ satisfies axioms of globularity, disjointness, and strictness of the relations $<^{X_k,+}$ for
$k>0$.

Clearly, if $<^{X_k,+}\; = \;<^{T_k,+}\cap (X_k)^2$, then the relation $<^{X_0,+}$ is linear, the relations $<^{X_k,+}$, for $k>0$, satisfy pencil linearity, i.e., $X$ is a positive opetopic cardinal.

Now we assume that the subhypergraph $X$ of positive opetopic cardinal $T$ is a positive opetopic cardinal.  We shall show that for $k\in\o$, $a,b\in X_k$, we have $a<^{X_k,+}b$ iff $a<^{T_k,+}b$. Since $X$ is a subhypergraph, $a<^{X_k,+}b$ implies $a<^{T_k,+}b$. Thus it is enough to show that if $a<^{T_k,+}b$ then $a\perp^{X_k,+}b$.  We shall prove this by induction on $k$. For $k=0$, it is obvious, since $<^{X_0,+}$ is linear. So assume that for faces $x,y\in X_l$, with $l<k$, we already know that $x<^{X_l,+}y$ iff $x<^{T_l,+}y$. Fix $a,b\in T_k$ such that $a<^{T_k,+}b$. Then by Lemma \ref{tech_lemma3}.\ref{tl2} $\gamma(a)\leq^{T_{k-1},+}\gamma(b)$ and hence by inductive hypothesis $\gamma(a)\leq^{X_{k-1},+}\gamma(b)$. Thus we have an upper $X-\gamma(X)$-path $a=a_r,\gamma(a), a_{r-1}\ldots, a_1,\gamma(b)$, with $r\geq 1$. Since $X$ is a positive opetopic cardinal and $\gamma(a_1)=\gamma(b)$, by pencil linearity we have $a_1\leq^+b$. By Path Lemma \ref{fullpath}, either $a<^{X,-}b$ or $a<^{X,+}b$. Since the first option is impossible, we have $a<^{X,+}b$, as required. $~~\Box$

\begin{lemma}\label{iota-mid}
Let $T$ be a positive opetopic cardinal, $a,b,\alpha\in T$. If $a\in \delta(\alpha)$ and $a<^+b<^+\gamma(\alpha)$, then $b\in\iota(T)$.
\end{lemma}

{\it Proof.}~Assume that $a,b,\alpha\in T$ are as in the assumption of the Lemma.  Thus we have an upper path $a,\alpha_0,\ldots,\alpha_r,b$. As $a\in \delta(\alpha)\cap\delta(\alpha_0)$, by pencil linearity we have $\alpha\perp^+\alpha_0$. If $\alpha<^+\alpha_0<^-\alpha_r$, then $\gamma(\alpha)\leq^+\gamma(\alpha_r)=b$ contradicting our assumption. Thus $\alpha_0<^+\alpha$. Then by Path Lemma \ref{fullpath}, since $b=\gamma(\alpha_r)<^+\gamma(\alpha)$, we have $\alpha_r<^+\alpha$ and $b\in\iota(T)$, as required. $~~\Box$

\subsection*{Some equations}

\begin{proposition}\label{equations_pfs}
Let $S$ be a positive opetopic cardinal $0<k\in\o$. Then
\begin{enumerate}
  \item $\iota(S_{k+1})=\iota(S_{k+1}-\delta(S_{k+2}))$ and $\iota(S_{k+1})=\iota(S_{k+1}-\gamma(S_{k+2}))$;
  \item $\delta(S_k)=\delta(S_k-\gamma({S_{k+1}}))$;
   \item $\gamma(S_k)=\gamma(S_k-\gamma({S_{k+1}}))$;
  \item $\delta(S_k)=\delta(S_k-\iota({S_{k+2}}))$;
  \item $\delta(S_k)=\delta(S_k-\delta({S_{k+1}}))\cup\iota(S_{k+1})$.
\end{enumerate}
\end{proposition}

{\it Proof.}~ In all the above equations the inclusion $\supseteq$ is obvious. So in each case we need to check the inclusion $\subseteq$ only.

Ad 1. Both equalities follow from Lemma \ref{iota}.

To prove the first equality, let $s\in\iota(S_{k+1})$, i.e., there is $a\in S_{k+1}$ such that $s\in\iota(a)$. By strictness, there is $b\in S_{k+1}$ such that $a\leq^+b$ and $b\not\in\delta(S_{k+1})$.  By Lemma \ref{iota}, we have
\[ s\in \iota(a)\subseteq \iota(b)\subseteq \iota(S_{k+1}-\delta(S_{k+2})) \]
 as required.

To prove the second equality, suppose contrary that there is $x\in\iota(S_{n+1})$ such that $x\not\in\iota(S_{n+1}-\gamma(S_n))$. Let $a\in S_{n+1}$ be a $<^+$-minimal face such that $x\in\iota(a)$. Since $x\not\in\iota(S_{n+1}-\gamma(S_n))$, there is $\alpha\in S_n$ such that $a=\gamma(\alpha)$. By Lemma \ref{iota} we have
\[ \iota(a)=\iota\gamma(\alpha)=\iota\delta(\alpha). \]
Therefore, there is $a'\in\delta(\alpha)$ such that $x\in\iota(a')$. Clearly $a'\lhd^+a$, and hence $a$ is not $<^+$-minimal, contrary to the supposition. This ends the proof of the first equality above.

Ad 2. Let $x\in\delta(S_k)$. Let $a\in S_k$ be the $<^+$-minimal element in $S_k$ such that $x\in\delta(a)$. We shall show that $a\in S_k-\gamma(S_{k+1})$. Suppose contrary that there is an $\alpha\in S_{k+1}$ such that $a=\gamma(\alpha)$. Then by globularity
 \[ x\in \delta(a)=\delta\gamma(\alpha)=\delta\delta(\alpha)-\gamma\delta(\alpha).\]
 So there is $b\in \delta(\alpha)$ such that $x\in\delta(b)$. As $b<^+a$, this contradicts the minimality of $a$.

Ad 3. This is similar to the previous one but simpler.

Ad 4. Since $\iota({S_{k+2}})\subseteq \gamma({S_{k+1}})$ 4. follows from 2.

Ad 5. Let $x\in\delta(S_k)$. Let $a\in S_k$ be the $<^+$-maximal element in $S_k$ such that $x\in\delta(a)$. If $a\not\in \delta(S_{k+1})$, then $x\in\delta(S_k-\delta({S_{k+1}}))$, as required. So assume that $a\in \delta(S_{k+1})$, i.e., there is $\alpha\in S_{k+1}$ such that $a\in\delta(\alpha)$.  Thus $x\in\delta\delta(\alpha)$. As $a<^+\gamma(\alpha)$, by choice of $a$ we have $x\not\in\delta\gamma(\alpha)\;(=\delta\delta(\alpha)-\gamma\delta(\alpha))$. So $x\in \gamma\delta(\alpha)$ and hence $x\in \iota(\alpha)\subseteq\iota(S_{k+1})$, as required. $~~\Box$

\section{The $\o$-categories generated by the positive opetopic cardinals}\label{sec-o-cat-generated}

The main objective of this section is to construct an embedding
\[ (-)^* : \pOpeCard \lra \oC \]
of the category of opetopic cardinals into the category of $\o$-categories. This embedding is not full. In Section \ref{sec-Sstar_is_pPoly} we shall show that the image of $(-)^*$ factorises through $\pPoly\ra \oC$ as a full functor.

Let $T$ be a positive opetopic cardinal. By $T^*_n$ we denote the set of all positive opetopic cardinals contained in $T$ of dimension at most $n$. If one wants to make these sets disjoint, one can think that an element of $T^*_n$ is  pair $\lk n, X \rk$, where $X$ is a positive opetopic cardinals contained in $T$.
We define below an $\o$-category, denoted $T^*$, whose set of $n$-cells is $T^*_n$.

We introduce operations
\[ \bd^{(k)}, \bc^{(k)}:T^*_n \lra T^*_k\]
of the $k$-th {\em domain} and the $k$-th {\em codomain} (of an $m$-dimensional cell), where $0\leq k\leq n$. For $S$ in $(T^*)_{m}$, the faces $k$-th {\em domain} $\bd^{(k)} S$ are:
\begin{enumerate}
   \item $ (\bd^{(k)} S)_l = \emptyset$, for $l> k$,
  \item $ (\bd^{(k)} S)_k= S_k - \gamma(S_{k+1})$,
  \item $ (\bd^{(k)} S)_l= S_l$, for $0\leq l<k$.
\end{enumerate}
and faces $k$-th {\em codomain} $\bc^{(k)} S$ are:
\begin{enumerate}
   \item $ (\bc^{(k)} S)_l= \emptyset$, for $l> k$,
  \item $(\bc^{(k)} S)_k= S_k - \delta (S_{k+1})$,
  \item $(\bc^{(k)} S)_{k-1} = S_{k-1} - \iota(S_{k+1})$, if $k>0$,
  \item $ (\bc^{(k)} S)_l = S_l$, for $0\leq l<k-1$.
\end{enumerate}

Note that the definitions of $\bd^{(k)}(S)$ and $\bc^{(k)}(S)$, for  $S\in T^*_n$ do not depend on the ambient opetopic cardinal $T$, nor even $dim(S)$.  Therefore we can write $\bd^{(k)}(S)$ and $\bc^{(k)}(S)$ without specifying either $T$ or $k$. We always assume that these cells are $k$-cells in $T^*$. If $n\in\o$ and we assume that $S\in T^*_{n+1}$, we write $\bd S$ for $\bd^{(n)}(S)$, and $\bc S$ for $\bc^{(n)}(S)$. We have

\begin{lemma}\label{o-cat-dom-codom} Let $S$ and $T$ be positive opetopic cardinals, $k\in \o$. Then
\begin{enumerate}

  \item  if $k<dim(S)$, then both $\bd^{(k)}(S)$, $\bc^{(k)}(S)$ are positive opetopic cardinals of dimension $k$; if $k\geq dim(S)$, then $\bd^{(k)}(S)=S=\bc^{(k)}(S)$;
  \item $\bd \bd^{(k+1)}(S)= \bd^{(k)}(S)$ , $\bc \bc^{(k+1)}(S)= \bc^{(k)}(S)$;
  \item if $S\in T^*_k$ and $k\geq 2$, then $\bd \bd S= \bd \bc S$, $\bc \bd S= \bc \bc S$;
   \item for any $\alpha\in S_k$, the least sub-hypergraph of $S$ containing the face $\alpha$ is again a positive opetopic cardinal of dimension $k$; it is denoted by $[\alpha]$. Moreover, if $k>0$,   then
  \[ \bc [\alpha] = [\gamma(\alpha)], \hskip 10mm \bd [\alpha] = [\delta(\alpha)] \]
  where $[\delta(\alpha)]$ is the least sub-hypergraph of $S$ containing the set of face $\delta(\alpha)$.
\end{enumerate}
\end{lemma}

{\it Proof.}~ Ad 1. It is obvious that $\bd^{(k)}S$ is a sub-hypergraph $S$. By Corollary \ref{gadeio} $\bc^{(k)}S$ is a sub-hypergraph $S$ as well. Any sub-hypergraph $T$ of a positive opetopic cardinal $S$ satisfies the conditions of globularity, strictness (possibly without $<^{T_0,+}$ being linear), and disjointness.

By Lemma \ref{tech_lemma1}, for $a,b\in\bd^{(k)} S_{l}$ we have $a<^{S_l,+}b$ iff $a<^{\bd^{(k)} S_l,+}b$. Moreover, by Lemma \ref{tech_lemma1'}, for $a,b\in{\bc^{(k)} S}_{l}$ we have $a<^{S_l,+}b$ iff $a<^{{\bc^{(k)} S}_l,+}b$. Hence by Lemma \ref{subhyp} both $\bd^{(k)} S$ and $\bc^{(k)} S$ are positive opetopic cardinals.

Ad 2. Fix a positive opetopic cardinal $S$ and $k\in\o$ such that $dim(S)> k$. Then the faces of $\bc^{(k+1)}(S)$, $\bc\bc^{(k+1)}(S)$, and $\bc^{(k)}(S)$ are as in the table
 \[ \begin{array}{|c|c|c|c|}
    dim & \bc^{(k+1)}(S)         & \bc\bc^{(k+1)}(S) & \bc^{(k)}(S) \\ \hline
    k+1 & S_{k+1}-\delta(S_{k+2}) & \emptyset & \emptyset \\ \hline
    k   & S_k-\iota(S_{k+2})    &  (S_k-\iota(S_{k+2}))-\delta(S_{k+1}-\delta(S_{k+2})) & S_k-\delta(S_{k+1})\\ \hline
    k-1 & S_{k-1}   & S_{k-1}-\iota(S_{k+1}-\delta(S_{k+2})) & S_{k-1}-\iota(S_{k+1}) \\ \hline
    l & S_l & S_l & S_l
  \end{array}
\] where $l<k-1$.
Moreover, the faces of $\bd^{(k+1)}(S)$, $\bd\bd^{(k+1)}(S)$, and $\bd^{(k)}(S)$ are as in the table
 \[ \begin{array}{|c|c|c|c|}
    dim & \bd^{(k+1)}(S)         & \bd\bd^{(k+1)}(S) & \bd^{(k)}(S) \\ \hline
    k+1 & S_{k+1}-\gamma(S_{k+2}) & \emptyset & \emptyset \\ \hline
    k   & S_k    &  S_k-\gamma(S_{k+1}-\gamma(S_{k+2})) & S_k-\gamma(S_{k+1})\\ \hline
    l & S_l & S_l & S_l
  \end{array}
\]
where $l<k$. Thus the equalities in question all follow from Lemma \ref{equations_pfs}.

Ad 3. Let $dim(S)=n>1$. Note that both $(\bd \bd S)_{n-2}$ and $(\bd \bc S)_{n-2}$ are the sets of all $<^+$-minimal elements in $S_{n-2}$, i.e., they are equal and the equation $\bd \bd S=\bd \bc S$ holds.

To see that $\bc \bd S= \bc \bc S$ holds, note first that both $(\bc \bd S)_{n-2}$ and  $(\bc \bc S)_{n-2}$ are the sets of all $<^+$-maximal elements in $S_{n-2}$. Moreover
\[ (\bc \bd S)_{n-3}=S_{n-3}-\iota(S_{n-1}-\gamma(S_n)), \]
\[ (\bc \bc S)_{n-3}=S_{n-3}-\iota(S_{n-1}-\delta(S_n)).\]
Now the equality $\bc \bd S= \bc \bc S$ follows from the following equalities
\[ \iota(S_{n-1}-\gamma(S_n))= \iota(S_{n-1}) =\iota(S_{n-1}-\delta(S_n)). \]
that follow from Lemma \ref{equations_pfs}.1.

Ad 4. Fix $\alpha\in S_k$.  We need to show that $[\alpha]$ is a positive opetopic cardinal. The globularity, strictness (except for linearity of $<^{[\alpha]_0,+}$), and disjointness are clear.

The linearity of $<^{[\alpha]_0,+}$. If $k\leq 2$, it is obvious. Put $a=\gamma^{(k+2)}(\alpha)$. Using Corollary \ref{atlas}, we have
\[ [\alpha]_k=\delta^{(k)}(\alpha)\cup\gamma^{(k)}(\alpha)= \]
\[ =\delta\delta(\gamma^{(k+2)}(\alpha))\cup\gamma\gamma(\gamma^{(k+2)}(\alpha))= \delta\delta(a)\cup\gamma\gamma(a). \]
Thus it is enough to check the linearity of $<^{[\alpha]_0,+}$ for $\alpha$ of dimension $k=2$. But in this case, as we mentioned, the  linearity of $<^{[\alpha]_0,+}$ is obvious.

The $\gamma$-linearity of $[\alpha]$. The proof proceeds by induction on $k=dim(\alpha)$. For $k\leq 2$, the $\gamma$-linearity is obvious. So assume that $k>2$ and that for $l<k$ and $a\in S_l$ $\gamma$-linearity holds in $[a]$.

First we shall show that $\bc([\alpha])=[\gamma(\alpha)]$. We have
\[ \bc([\alpha])_{k-1} =(\gamma(\alpha)\cup\delta(\alpha))-\delta(\alpha)=\gamma(\alpha)=[\gamma(\alpha)]_{k-1}, \]
\[ \bc([\alpha])_{k-2} = (\gamma\gamma(\alpha)\cup\delta\delta(\alpha))-\iota(\alpha)= \gamma\gamma(\alpha)\cup\delta\gamma(\alpha)=[\gamma(\alpha)]_{k-2}, \]
and for $l<k-2$
\[ \bc([\alpha])_l = \gamma^{(l)}(\alpha)\cup\delta^{(l)}(\alpha)= \gamma^{(l)}(\alpha)\cup\delta^{(l)}\gamma(\alpha)=[\gamma(\alpha)]_{l}. \]
Note that the definition of $\bc(H)$ makes sense for any positive hypergraph $H$ and in the above argument we haven't used the fact (which we don't know yet) that $[\alpha]$ is a positive opetopic cardinal.

Thus, for $l<k-2$, $[\alpha]_l=[\gamma(\alpha)]_l$. By induction hypothesis, $[\gamma(\alpha)]$ is a positive opetopic cardinal, and hence $[\alpha]_l$ is $\gamma$-linear for $l<k-2$.  Clearly $[\alpha]_l$ is $\gamma$-linear for $l=k-1,k$. Thus it remains to show the $\gamma$-linearity of $(k-2)$-cells in $[\alpha]$.

Fix $t\in [\alpha]_{k-3}$, and let
\[ \Gamma_t=\{ x\in [\alpha]_{k-2} : \gamma(x)=t \}. \]
We need to show that $\Gamma_t$ is linearly ordered by $<^+$. We can assume that  $t\in\gamma([\alpha]_{n-2}) =\gamma\delta\delta(\alpha)= \gamma\delta\gamma(\alpha)$  (otherwise $\Gamma_t=\emptyset$ is clearly linearly ordered by $<^+$). By Proposition \ref{cond2} there is a unique $x_t\in\delta\gamma(\alpha)$ such that $\gamma(x_t)=t$. From Lemma \ref{tech_lemma0}.2 we get easily the following Claim.

 {\em Claim 1.} For every $x\in\Gamma_t$ there is a unique upper $\delta(\alpha)$-path from $x_t$ to $x$.

Now fix $x,x'\in\Gamma_t$.  By the Claim 1, we have the unique upper $\delta(\alpha)$-path
\[ x_t,a_0,\ldots , a_l,x,\;\;\;\; x_t,a'_0,\ldots , a'_{l'},x'. \]
Suppose $l\leq l'$. By Proposition \ref{cond2}, for $i\leq l$, $a_i=a'_i$.  Hence either $l=l'$ and $x=x'$ or $l<l'$ and
\[ x, a_{l+1}, \ldots , a_{l'},x' \]
is a $\delta(\alpha)$-upper path.  Hence either $x=x'$ or $x\perp^+x'$ and $[\alpha]_{k-2}$ satisfy the $\gamma$-linearity, as required.

The proof is of the $\delta$-linearity of $[\alpha]$ is very similar to the one above.  For the same reasons the only non-trivial thing to check is the condition for $(k-2)$-faces. We pick $t\in\delta\delta(\alpha)$ and consider the set
\[ \Delta_t=\{ x\in [\alpha]_{k-2} : t\in\delta(x) \}. \]
Then we have a unique $y_t\in\delta\gamma(\alpha)$ such that $t\in\delta(y_t)$. From Lemma \ref{tech_lemma0}.3 we get the following Claim.

{\em Claim 2.} For every $y\in\Delta_t$ there is a unique upper $\delta(\alpha)$-path from $y_t$ to $y$.

The $\delta$-linearity of the $(k-2)$-faces in  $[\alpha]$ can be proven from Claim 2 similarly as the $\gamma$-linearity was proven from Claim 1.

It remains to verify the equalities
\[ \bc [\alpha] = [\gamma(\alpha)], \hskip 10mm \bd [\alpha] = [\delta(\alpha)]. \]
The first one we already checked on the way. To see that the second equality also holds we calculate
\[ \bd[\alpha]_{k-1}=(\gamma(\alpha)\cup\delta(\alpha))-\gamma(\alpha) = \delta(\alpha)=[\delta(\alpha)]_{k-1}, \]
\[ \bd[\alpha]_{k-2}=(\gamma\gamma(\alpha)\cup\delta\delta(\alpha)) = \gamma\delta(\alpha)\cup\delta\delta(\alpha)=[\delta(\alpha)]_{k-2}, \]
and for $l<k-2$
\[ \bd[\alpha]_{l}=\gamma^{(l)}(\alpha)\cup\delta^{(l)}(\alpha)) = \gamma^{(l)}\delta(\alpha)\cup\delta^{(l)}(\alpha)=[\delta(\alpha)]_{l}. \]
So the second equality holds as well.$~~\Box$

\vskip 2mm

{\em Remarks.}
\begin{enumerate}
  \item Inspired by the above Lemma \ref{o-cat-dom-codom}.4 we call $S$ a {\em weak  positive opetopic cardinal}\index{opetopic cardinal!weak positive -} if $S$ satisfies globularity, strictness, and disjointness as a positive hypergraph and moreover that for any face $\alpha$ in $S$ the sub-hypergraph $[\alpha]$ is an opetope.   (i.e., pencil linearity is required to hold only `locally'). The {\em category of weak  positive opetopic cardinal}\index{category!of positive opetopic cardinals!weak $\wpOpeCard$} is the full subcategory of positive hypergraphs $\pHg$ whose objects are the weak positive opetopic cardinals and is denoted by $\wpOpeCard$. For each $k\in\o$, the $k$-truncation of a weak positive opetopic cardinal $S$ is again a weak positive opetopic cardinal $S_{\leq k}$. In particular, any $k$-truncation of a positive opetopic cardinal $S$ is a weak positive opetopic cardinal $S_{\leq k}$, but it does not necessarily satisfy linearity condition.
  \item From Lemma \ref{o-cat-dom-codom}.1 we know for any positive opetopic cardinal $S$ the hypergraphs $\bc^{(k)}$ and $\bd^{(k)}$ are positive opetopic cardinals contained in $S$. We shall denote these embeddings by
      \begin{center} \xext=1200 \yext=120
\begin{picture}(\xext,\yext)(\xoff,\yoff)
  \putmorphism(0,50)(1,0)[\bc^{(k)}S`S`\bc^{(k)}_S]{600}{1}a
  \putmorphism(600,50)(1,0)[\phantom{S}`\bd^{(k)}S.`\bd^{(k)}_S]{600}{-1}a
 \end{picture}
\end{center}

\end{enumerate}

\begin{lemma}\label{o-cat-comp} Let $S$ and $T$ be positive opetopic cardinals such that $\bc^{(k)} S\subseteq \bd^{(k)} T$. Then the pushout
  $S+_k T$ in $\pOpeCard$ of $S$ and $T$ over ${\bc^{(k)} S}$  exists.   Moreover, if $\bc^{(k)} S=S\cap T$, then the diagram of inclusions in $\pOpeCard$
\begin{center} \xext=650 \yext=600
 \begin{picture}(\xext,\yext)(\xoff,\yoff)
 \setsqparms[-1`-1`-1`-1;650`500]
 \putsquare(0,50)[T\cup S`S`T`\bc^{(k)}S;``\bc^{(k)}_S`\bd^{(k)}_T]
\end{picture}
\end{center}
  is the pushout.
\end{lemma}

{\it Proof.}~  Assume that $\bc^{(k)} S=S\cap T\subseteq \bd^{(k)} T$. Let $T\cup S$ be the obvious sum of $S$ and $T$ as positive hypergraphs. The fact that $T\cup S$ is a pushout in $\pHg$ is obvious. Thus the only thing we need to verify is that $T\cup S$ is a positive opetopic cardinal.

First we write in detail the condition $\bc^{(k)} S=S\cap T\subseteq \bd^{(k)} T$:

\begin{enumerate}
   \item $S_l\cap T_l = \emptyset$, for $l> k$,
  \item $S_k-\delta(S_{k+1})\subseteq T_k-\gamma(T_{k+1})$,
  \item $S_{k-1}-\iota(S_{k+1})\subseteq T_{k-1}$,
  \item $S_l\subseteq T_l$, for $l<k-1$.
\end{enumerate}

Now we describe the orders $<^+$ in $T\cup S$:
\[ <^{(T\cup S)_l,+}  = \left\{ \begin{array}{ll}
        <^{S_l,+} +  <^{T_l,+} & \mbox{for  $l>k$,}  \\
        <^{S_l,+} +_{(S_k-\delta(S_{k+1}))}  <^{T_l,+} & \mbox{for  $l=k$,}  \\
        <^{S_l,+}+_{(S_{k-1}-\iota(S_{k+1}))} <^{T_l,+} & \mbox{for  $l=k-1$,}  \\
        <^{T_l,+}   & \mbox{for  $l<=k-1$.}
                                    \end{array}
                \right. \]
We shall comment on these formulas. For $l>k$, the formulas say that the order $<^+$ in $(T\cup S)_l$ is the disjoint sum of the orders in $S_l$ and $T_l$. This is obvious.

For $l<k-1$, the order $<^+$ in $(S\cap T)_l$ is just the order $<^{T_l,+}$.  The only case that requires an explanation is $l=k-2$. So suppose that $a,b\in T_{k-2}$ and $a<^{(T\cup S)_{k-2},+}b$. So we have an upper path
\[ a,\alpha_1,\ldots , \alpha_m, b \]
such that $\alpha_i\in (T\cup S)_{k-1}=\iota(S_{k+1})\cup T_{k-1}$. By Lemma \ref{tech_lemma1}, we can assume that if $\alpha_i\in S_{k-1}$, then $\alpha_i\not\in\gamma(S_{k})$. But then $\alpha_i\not\in\iota(S_{k+1})$. So in fact $\alpha_i\in T_{k-1}$, as required.

The most involved are the formulas for $<^{(S\cap T)_l,+}$, for $l=k$ and $l=k-1$. In both cases the comparison in $T\cup S$ involves orders both from $S$ and $T$. In the former case we have that, for $a,b\in (T\cup S)_k$, we have
\[ a<^{(T\cup S)_k,+}b  \mbox{ iff  }  \hskip 70mm \]
\[    \left\{ \begin{array}{ll}
          \mbox{either }  a,b\in T_k & \mbox{ and } a<^{T_k,+}b,    \\
          \mbox{or }  a,b\in S_k& \mbox{ and } a<^{S_k,+}b,  \\
          \mbox{or }  a\in \delta(S_{k+1}),\; b\in T_k& \mbox{ and }
          \exists_{a'\in S_k-\delta(S_{k+1})} a<^{S_k,+}a' \mbox{ and } a'\leq^{T_k,+}b.
                                    \end{array}
                \right. \]
The orders $<^{S_k,+}$ and $\leq^{T_k,+}$ are glued together along the set $S_k-\delta(S_{k+1})$ which is the set of $<^{S_k,+}$-maximal elements in $S_k$ and at the same time it is contained in the set of $<^{T_k,+}$-minimal elements $T_k-\gamma(T_{k+1})$. This is obvious when we realize that $\delta(S_{k+1})\cap\gamma(T_{k+1})=\emptyset$.

In the latter case, for $x,y\in (T\cup S)_{k-1}$, we have
\[ x<^{(T\cup S)_{k-1},+}y  \mbox{ iff  }  \hskip 70mm \]
\[    \left\{ \begin{array}{ll}
          \mbox{either }  x,y\in S_{k-1} & \mbox{ and } x<^{S_{k-1},+}y,    \\
          \mbox{or }  x,y\in T_{k-1}& \mbox{ and } x<^{T_{k-1},+}y,  \\
          \mbox{or }  x\in \iota(S_{k+1}),\; y\in T_k& \mbox{ and }
          \exists_{x'\in S_{k-1}-\iota(S_{k+1})} x<^{S_k,+}x' \mbox{ and }
          x'\leq^{T_k,+}y, \\
          \mbox{or }  x\in T_k,\; y\in \iota(S_{k+1}) & \mbox{ and }
          \exists_{x'\in S_{k-1}-\iota(S_{k+1})} x<^{T_k,+}x' \mbox{ and }
          x'\leq^{S_k,+}y.
                                    \end{array}
                \right. \]
The order $<^{S_{k-1},+}$ is `plugged into' the order $\leq^{T_{k-1},+}$, along the set $S_k-\iota(S_{k+1})$.

To show that these formulas hold true, we argue by cases. Assume that $x,y\in (T\cup S)_{k-1}$ and that $x<^{(T\cup S)_{k-1},+}y$, i.e., there is an upper path
\[ x, a_1, \ldots , a_m , y \]
with $a_i\in (T\cup S)_k,$ for $i=1,\ldots ,m$.

First suppose that $x,y\in S_{k-1}$ and that the set $\{a_i\}_{i}\not\subseteq S_k$. Let $a_{i_0},a_{i_0+1},\ldots ,a_{i_1}$ be a maximal subsequence of  consecutive elements of the path $a_1, \ldots ,a_m$ such that $\{ a_i\}_{i_0\leq i\leq i_1 }\subseteq T_k$. Thus it is an upper path in $T_k$ from $\bar{x}$ to $\bar{y}=\gamma(a_{i_1})$, where
\[   \bar{x} = \left\{ \begin{array}{ll}
          x & \mbox{if } i_0=1,    \\
          \gamma(a_{i_0-1})& \mbox{otherwise.}
                                    \end{array}
                \right. \]
Note that from maximality of the path $a_{i_0},\ldots , a_{i_1}$ follows that both $\bar{x},\bar{y}\in S_{k-1}-\iota(S_{k+1})$. As we have $\bar{x}<^{T_{k-1},+}\bar{y}$ from Corollary \ref{order}, we have $\bar{x}\not\perp^{T_l,-}\bar{y}$, for all $l< k-1$. Clearly $\perp^{S_l,-}\subseteq \perp^{T_l,-}$. Thus $\bar{x}\not\perp^{S_l,-}\bar{y}$, for all $l< k-1$, as well. But then again by Corollary \ref{order} we have that $\bar{x}\perp^{S_{k-1},+}\bar{y}$. If we were to have $\bar{y}<^{S_{k-1},+}\bar{x}$, then, as $\bar{x},\bar{y}\in S_{k-1}-\iota(S_{k+1})$, we would have $\bar{y}<^{T_{k-1},+}\bar{x}$.  But this would contradict the strictness of $<^{T_{k-1},+}$.  So we must have $\bar{x}<^{S_{k-1},+}\bar{y}$. In this way we can replace the upper path $a_1, \ldots , a_m$ in $(T\cup S)_k$ from $x$ to $y$ by an upper path from $x$ to $y$ in $S_k$.

Next, suppose that $x,y\in T_{k-1}$ and that the set $\{a_i\}_{i}\not\subseteq T_k$. Let $a_{i_0},a_{i_0+1},\ldots , a_{i_1}$ be a maximal subsequence of consecutive elements of the path $a_1, \ldots ,a_m$ such that $\{ a_i\}_{i_0\leq i\leq i_1 }\subseteq S_k$. Thus it is an upper path in $S_k$ from $\bar{x}$ to $\bar{y}=\gamma(a_{i_1})$, where
\[   \bar{x} = \left\{ \begin{array}{ll}
          x & \mbox{if } i_0=1,    \\
          \gamma(a_{i_0-1})& \mbox{otherwise.}
                                    \end{array}
                \right. \]
Note that from maximality of the sequence $a_{i_0},\ldots , a_{i_1}$ follows that both $\bar{x},\bar{y}\in S_{k-1}-\iota(S_{k+1})\subseteq T_{k-1}$. Thus by Lemma
\ref{tech_lemma1'} there is an upper path from $\bar{x}$ to $\bar{y}$ in $S_{k-1}-\delta(S_{k})\subseteq T_{k-1}$. In this way we can replace the upper path $a_1, \ldots , a_m$  in $(T\cup S)_k$ from $x$ to $y$ by an upper path from $x$ to $y$ in $T_k$.

Thus we have justified the first two cases of the above formula. The following two cases are easy consequences of these two. This end the description of the orders in $T\cup S$.

From these descriptions follows immediately that $<^{(T\cup S),+}$ is strict, for all $l$. It remains to show the pencil linearity. Both $\gamma$- and $\delta$-linearity of $l$-cells, for $l<k-1$ or $l>k$, are obvious.

To see the $\gamma$-linearity of $k$-cells assume, $a\in S_k$ and $b\in T_k$, such that $\gamma(a)=\gamma(b)$. Let $\bar{a}\in S_k$ be the $<^{S_k,+}$-maximal $k$-cells such that $\gamma(a)=\gamma(\bar{a})$. Then $\bar{a}\in\bc^{(k)}(S)_k \subseteq \bd^{(k)}(T)_k$. So $\bar{a}\in T_k$ is $<^{T_k,+}$-minimal $k$-cells such that $\gamma(\bar{a})=\gamma(b)$. Thus
\[ a\leq^{S_k,+}\bar{a}  \leq^{T_k,+}b.\]
Thus the $\gamma$-linearity of $k$-cells holds. The proof of $\delta$-linearity of $k$-cells is similar.

Finally, we need to establish the $\gamma$- and $\delta$-linearity of $(k-1)$-cells in $T\cup S$.

In order to prove the $\gamma$-linearity, let $x\in\iota(S_{k+1})$ and $y\in T_{k-1}$ such that $\gamma(x)=\gamma(y)$. We need to show that $x\perp^{(T\cup S)_{k-1},+}y$.

Let $\alpha_0\in S_{k+1}$ such that $x\in\iota(\alpha_0)$, $a\in\delta(\alpha_0)$ such that $x=\gamma(a)$ and let $\alpha_0,\ldots ,\alpha_l$ be a lower path in $S_{k+1}$ such that $\gamma(\alpha_l)\in T_k$. Since $x\in\iota(\alpha_0)$, then $x\in\gamma\delta(\alpha_0)$ and, by Lemma \ref{atlas3}
 \[ \gamma(x)\in\gamma\gamma\delta(\alpha_0)\subseteq \iota\gamma(\alpha_0)\cup \gamma\gamma\gamma(\alpha_0). \]
As $\gamma(\alpha_0)\leq^+\gamma(\alpha_l)$, by Lemma \ref{iota}.3, we have $\gamma(x)\in\iota\gamma(\alpha_l)\cup\gamma\gamma\gamma(\alpha_l)$. Thus we have two cases:
\begin{enumerate}
  \item $\gamma(x)\in\iota\gamma(\alpha_l)$,
  \item $\gamma(x)=\gamma\gamma\gamma(\alpha_l)$.
\end{enumerate}

Case 1: $\gamma(x)\in\iota\gamma(\alpha_l)$.

By Lemma \ref{tech_lemma0}.2, there is a unique $z\in\delta\gamma(\alpha_l)$ such that $\gamma(z)=\gamma(x)$ and $z<^+x$. As $\gamma(\alpha_l)\in T_k$, so $z\in T_{k-1}$. If
$y<^{T_{k-1},+}z$, then indeed $y<^{(T\cup S)_{k-1},+}z$, as required. By $\gamma$-linearity in $T_{k-1}$, it is enough to show that it is impossible to have $z<^{T_{k-1},+}y$.

Suppose contrary that there is an upper path $z,b_0,\ldots , b_r, y$ in $T$. Since $\gamma(\alpha_l)$ is $<^+$-minimal in $T$ (as $\alpha_l\in S$) and $z\in\delta\gamma(\alpha_l)\cap\delta(b_0)$, by $\delta$-linearity in $T_k$ we have $\gamma(\alpha_l)<^+b_0$. By Lemma \ref{iota}.2, we have
\[ \gamma(x)\in\iota\gamma(\alpha_l)\subseteq\iota(b_0)\subseteq \iota(b_r) \]
But $\gamma(b_r)=y$ so $\gamma\gamma(b_r)=\gamma(y)=\gamma(x)$. In particular, $\gamma(x)\not\in\iota(b_r)$ and we get a contradiction.

Case 2: $\gamma(x)=\gamma\gamma\gamma(\alpha_l)$.

By Lemma \ref{tech_lemma0}.2 there is $z\in\delta\gamma(\alpha_l)$ such that $\gamma(x)=\gamma(z)(=\gamma\gamma\gamma(\alpha_l))$, so that we have
\[ z<^{S_{k-1},+}x<^{S_{k-1},+}\gamma\gamma(\alpha_l). \]
As $\gamma(\alpha_l)\in T_{k}$ and it is $<^+$-minimal in $T_k$, by Proposition \ref{successor}, there is no face $y'\in T_{k-1}$ so that
\[ z<^{T_{k-1},+}y'<^{T_{k-1},+}\gamma\gamma(\alpha_l). \]
So if $y\in T_{k-1}$ and $\gamma(y)=\gamma(x)$, then either
\[ y\leq^{T_{k-1},+}z<^{S_{k-1},+}x\;\;\; \mbox{ or }\;\;\;  x<^{S_{k-1},+}\gamma\gamma(\alpha_l)\leq^{S_{k-1},+}y. \]
In either case $x\perp^{{(T\cup S)_{k-1},+}}y$, as required. This ends the proof of $\gamma$-linearity of $(k-1)$-faces in $(T\cup S)$.

Finally, we prove the $\delta$-linearity of $(k-1)$-faces in $T\cup S$. Let $x\in\iota(S_{k+1})$ and $y\in T_{k-1}$, $t\in T_{k-2}$ such that $t\in \delta(x)\cap\delta(y)$. We need to show that $x\perp^{(T\cup S)_{k-1},+}y$.

Let $\alpha_0\in S_{k+1}$ such that $x\in\iota(\alpha_0)$, $a\in\delta(\alpha_0)$ such that $x=\gamma(a)$, and let $\alpha_0, \ldots , \alpha_l$ be a lower path in $S_{k+1}$ such that $\gamma(\alpha_l)\in T_k$. As $x\in\iota(\alpha_0)$, using Lemma \ref{atlas3} we have
\[t\in\delta(x)\subseteq\delta\gamma\delta(\alpha_0)\subseteq \delta\gamma\gamma(\alpha_0)\cup\iota\gamma(\alpha_0).\]
As $\gamma(\alpha_0)<^+\gamma(\alpha_l)$, by Lemma \ref{iota}.4, we have two cases:
\begin{enumerate}
  \item $t\in\iota\gamma(\alpha_l)$,
  \item $t\in\delta\gamma\gamma(\alpha_l)$.
\end{enumerate}

Case 1: $t\in\iota\gamma(\alpha_l)$.

By Lemma \ref{tech_lemma0}.3, there is a unique $z\in\delta\gamma(\alpha_l)$ such that $t\in\delta(z)$ and $z<^+x$.  As $\gamma(\alpha_l)\in T_k$, so $z\in T_{k-1}$. If $y<^{T_{k-1},+}z$, then indeed $y<^{(T\cup S)_{k-1},+}z$, as required. By $\delta$-linearity in $T_k$, it is enough to show that it is impossible to have $z<^{T_{k-1},+}y$.

Suppose contrary that there is an upper path in $T$
\[ z,b_0,\ldots ,b_r,y .\]
Since $\gamma(\alpha_l)$ is $<^+$-minimal in $T_k$ and $z\in\delta\gamma(\alpha_l)\cap\delta(b_0)$, by $\delta$-linearity of $k$-faces in $T$ we have $\gamma(\alpha_l)<^+b_0$. By Lemma \ref{iota}.2, we have
\[ t\in \iota\gamma(\alpha_l)\subseteq\iota(b_0)\subseteq \ldots \subseteq \iota(b_r).\]
But $\gamma(b_r)=y$, so $t\in\delta(y)\subseteq\delta\gamma(b_r)$. In particular, $t\not\in\iota(b_r)$ and we get a contradiction.

Case 2: $t\in\delta\gamma\gamma(\alpha_l)$.

By Lemma \ref{tech_lemma0}.3 there is $z\in\delta\gamma\gamma(\alpha_l)$ such that $t\in\delta(z)$ and we have
\[ z <^{S_{k-1},+}x<^{S_{k-1},+}\gamma\gamma(\alpha_l). \]
As $\gamma(\alpha_l)\in T_k$, and it is $<^+$-minimal face in $T_k$, by Lemma \ref{successor}, there is no face $y'\in T_{k-1}$ such that
\[ z<^{T_{k-1},+}y'<^{T_{k-1},+}\gamma\gamma(\alpha_l). \]
So if $y\in T_{k-1}$ and $t\in\delta(y)$, then either
\[ y\leq^{T_{k-1},+}z<^{S_{k-1},+}x \;\;\; \mbox{ or }\;\;\;  x<^{S_{k-1},+}\gamma\gamma(\alpha_l)\leq^{S_{k-1},+}y. \]
In either case $x\perp^{{(T\cup S)_{k-1},+}}y$, as required. This ends the proof of $\delta$-linearity of $(k-1)$-faces in $(T\cup S)$ and the whole proof that $T\cup S$ is a positive opetopic cardinal. $~~\Box$

Let $S$ and $T$ be positive opetopic cardinals such that $\bc^{(k)} S= \bd^{(k)} T$. Then the pushout just described
\begin{center}
\xext=650 \yext=600
\begin{picture}(\xext,\yext)(\xoff,\yoff)
 \setsqparms[-1`-1`-1`-1;650`500]
 \putsquare(0,100)[T+_kS`S`T`\bc^{(k)} S;``\bc^{(k)}_S`\bd^{(k)}_T]
\end{picture}
\end{center}
is called {\em special pushout}\index{pushout!special -}\index{special!pushout} in $\pOpeCard$ (or {\em special pullback}\index{pullbackt!special -}\index{special!pullback} in $\pOpeCard^{op}$) and is denoted as $T\oplus_kS$.

Now we shall describe an $\o$-category $T^*$ generated by the positive opetopic cardinal $T$. The set of $m$-cell of $T^*$ is $T^*_m$, i.e., the set of all the positive opetopic cardinals contained in $T$ of dimension at most $m$, for $m\in\o$. The {\em $k$-th domain} and {\em $k$-th codomain} operations in $T^*$ are the operations
\[ \bd^{(k)}, \bc^{(k)} :T^*_{m} \lra T^*_k\]
defined above, with $m\geq k$. The {\em identity} operations
\[ \bi^{(m)}:T^*_k \lra T^*_{m}\]
are inclusions, and the composition map
 \[ \bm_{m,k,m} : T^*_{m}\times_{T^*_k} S^*_{m}\lra T^*_{m},\]
where $k<m$, is the sum, i.e., if $X$, $Y$ are positive opetopic cardinals contained in $T$ of dimension at most $m$ such that $\bc^{(k)} X=\bd^{(k)}Y$, then
 \[ \bm_{m,k,m}(X,Y) = X \oplus_k Y = X \cup Y. \]

\begin{corollary}\label{Star_functor}
Let $T$ be a weak positive opetopic cardinal. Then $T^*$ is an $\o$-category.  In fact, we have a functor
\[ (-)^* : \wpOpeCard \lra \oC \]
\end{corollary}

{\it Proof.}~ The fact that the operations on $T^*$ are well defined and satisfy the laws of $\o$-category follows from Lemmas \ref{o-cat-dom-codom} and \ref{o-cat-comp}.

If $f:S\ra T$ is a morphism of positive opetopic cardinals, $X\in S^*$ , then the image $f(X) \in T^*$ is isomorphic to $X$. Then again using Lemmas  \ref{o-cat-dom-codom} and \ref{o-cat-comp}, the association $X\mapsto f(X)$ is easily seen to be an $\o$-functor.  $~~\Box$

\section{Normal positive opetopic cardinals}

Let $S$ be a normal positive opetopic cardinal of dimension $k$, i.e.,  $S$ is $(k-1)$-principal. By $\bp_l^S$ we denote the unique element of the set $S_l-\delta(S_{l+1})$, for $l<k$. Moreover, as we shall show below, $\bp^S_{k-1}\in\gamma(S_k)$ and hence the set $\{ x \in S_k : \gamma(x)=\bp^S_{k-1}\}$ is not empty. We denote by $\bp^S_{k}$ the $<^+$-maximal element of this set. We shall omit the superscript $S$ if it does not lead to a confusion.

\begin{lemma}\label{norm_positive0}
Let $S$ be a $(k-1)$-principal opetope of dimension at least $k$, $k>0$. Then
\begin{enumerate}
 \item $S_l=\delta^{(l)}(S_k)\cup\gamma^{(l)}(S_k)=\delta^{(l)}(S_k)\cup\{\bp_l \}$, for $l<k$.
 \item $\delta(S_{l+1})=\delta^{(l)}(S_k)$, for $l<k$.
 \item $\bp_{k}$ is $<^-$-largest element in $S_k-\delta(S_{k+1})$.
 \item $\gamma(\bp_l)=\bp_{l-1}$, for $0<l\leq k$.
  \item $\delta(\bp_l)=\delta(S_l)-\gamma(S_l)$, for $0<l<k$.
 \item $S_l=\delta^{(l)}(\bp_{k-1})\cup \gamma^{(l)}(\bp_{k-1})$, for $l<k-2$.
\end{enumerate}
\end{lemma}

{\it Proof.}~Ad 1. If $H$ is a hypergraph of dimension greater than $l$ and  $\gamma(H_{l+1})\subseteq \delta(H_{l+1})$, then there is an infinite lower path in $H_{l+1}$, i.e., $<^{H_l,+}$ is not strict. Thus, if $S$ is a positive opetopic cardinal of dimension greater than $l$, we have $\delta(S_{l+1})\subseteqnot S_l$. A positive opetopic cardinal is normal iff this difference
\[ S_l-\delta(S_{l+1}) \]
is a singleton, for $l<k$. Thus, by the above, we must have
\begin{equation}\label{eq-sum}
  S_l=\delta(S_{l+1})\cup\gamma(S_{l+1})
\end{equation}
The first equation of the statement 1. we shall show by the downward induction on $l$. Suppose that we have $S_{l+1}=\delta^{(l)}(S_k)\cup \gamma^{(l)}(S_k)$ (for $l=k-2$, it is true by the above). Then
\[ S_l= \delta^{(l)}(S_k)\cup \gamma^{(l)}(S_k) = \]
\[ =\delta(\delta^{(l+1)}(S_k)\cup \gamma^{(l+1)}(S_k)) \cup  \gamma(\delta^{(l+1)}(S_k)\cup \gamma^{(l+1)}(S_k))= \]
\[ =\delta\delta^{(l+1)}(S_k)\cup \delta\gamma^{(l+1)}(S_k) \cup \gamma\delta^{(l+1)}(S_k)\cup\gamma\gamma^{(l+1)}(S_k)= \]
\[ =\delta^{(l)}(S_k)\cup \delta\gamma^{(l)}(S_k) \cup \gamma\delta^{(l)}(S_k)\cup\gamma^{(l)}(S_k)= \]
\[ =\delta^{(l)}(S_k)\cup\gamma^{(l)}(S_k) \]
where the last equation follows from Corollary \ref{atlas}.

The second equation of 1. is obvious, for $l=k-1$. So assume that $l<k-1$. We have
\[ \{ \bp_l \} = S_l -\delta(S_{l+1}) = \]
\[ = S_l -\delta(\delta^{(l+1)}(S_{k})\cup\gamma^{(l+1)}(S_{k})) =\]
\[ = S_l -(\delta^{(l)}(S_{k})\cup\delta\gamma^{(l+1)}(S_{k})) =\]
\[ S_l -\delta^{(l)}(S_{k}). \]
Thus \[ S_l= \delta^{(l)}(S_{k})\cup\{ \bp_l \}\] as required.

Ad 2. Let $l<k$. Then using 1. we have
\[ \delta^{(l)} (S_k) \subseteq \delta(S_{l+1})
\subseteqnot \delta^{(l)}(S_k)\cup\{ \bp_l \} \]
 Hence
 \[ \delta^{(l)}(S_k)=\delta(S_{l+1}). \]

Ad 3. First we shall show that $\bp_k\in S_k-\delta(S_{k+1})$. Suppose contrary that there is $\alpha\in S_{k+1}$ such that $\bp_k\in\delta(\alpha)$. Then $\gamma(\bp_k)\in\gamma\delta(\alpha)=\gamma\gamma(\alpha)\cup\iota(\alpha)$. If $\gamma(\bp_k) =\gamma\gamma(\alpha)$, then $\bp_k<^+\gamma(\alpha)$, i.e., $\gamma(\alpha)$ is $<^+$-smaller element than $\bp_k$ such that $\gamma(\gamma(\alpha))=\bp_{k-1}$. This contradicts the choice of $\bp_k$. If $\gamma(\bp_k)  =\iota(\alpha)$, then there is $a\in\delta(\alpha)$ such that $\gamma(\bp_k)\in\delta(a)$. But this means that $\bp_{k-1}=\gamma(\bp_k)\in \delta(S_k)$ contradicting the choice of $\bp_{k-1}\in  S_{k-1}-\delta(S_k)$. This shows that $\bp_k\in S_{k}-\delta(S_{k+1})$.

We need to prove that any maximal lower $(S_k-\delta(S_{k+1}))$-path ends at $\bp_k$. By strictness, it is enough to show that if $x\in S_k-\delta(S_{k+1})$ and $x\neq\bp_k$, then there is $x'\in S_k-\delta(S_{k+1})$ such that $\gamma(x)\in\delta(x')$. So fix $x\in S_k-\delta(S_{k+1})$. If we were to have $\gamma(x)\in\iota(\beta)$, for some $\beta\in S_{k+1}$, then by Lemma \ref{tech_lemma2} we would have $x<^+\gamma(\beta)$. In particular, $x\in \delta(S_{k+1})$, contrary to the assumption. Therefore $\gamma(x)\in S_{k-1}-\iota(S_{k+1})$. As $x,\bp_k\in S_{k}-\delta(S_{k+1})$, by $\gamma$-linearity we have $\gamma(x)\neq\gamma(\bp_k)=\bp_{k-1}$. Hence by 1. the set
\[ \Delta_{\gamma(x)} =\{ y\in S_k : \gamma(x)\in \delta(y) \} \]
is not empty. Let $x'$ be the $<^+$-maximal element of this set. It remains to show that $x'\not\in\delta(S_{k+1})$. Suppose contrary that there is $\alpha\in S_{k+1}$ such that $x'\in \delta(\alpha)$. As $\gamma(x)\not\in\iota(S_{k+1})$ and $\gamma(x)\in\delta(x')$, so $\gamma(x)\not\in\iota(\alpha)$ and $\gamma(x)\neq\gamma\gamma(\alpha)$. Thus $\gamma(x)\in\delta\gamma(\alpha)$. But this means that $x'<^+\gamma(\alpha)$ and $\gamma(\alpha)\in\Delta_{\gamma(x)}$. This contradicts the choice of $x'$. This ends the proof of 3.

Ad 4. $\gamma(\bp_k)=\bp_{k-1}$ by definition. Fix $0<l<k$. As $S_l=\delta(S_{l+1})\cup\{ \bp_l\}$, $\bp_l$ is $<^+$-greatest element in $S_l$.  Assume that $\gamma(\bp_l)\neq\bp_{l-1}$.  Thus $\gamma(\bp_l)<^+\bp_{l-1}$. Let $x\in S_l$.  Then $x\leq \bp_l$ and, by Lemma \ref{tech_lemma3}, $\gamma(x)\leq^+\gamma(\bp_l)<^+\bp_{l-1}$. Thus $\bp_{l-1}\not\in\gamma(S_l)$. So $\gamma(S_l)\subseteq\delta(S_l)$. But this is impossible in a positive opetopic cardinal, as we noticed in the proof of 1. This ends 4.

Ad 5. Fix $l<k$. First we shall show that
\begin{equation}\label{e1}
 \delta(\bp_l)\cap\gamma(S_l)=\emptyset
\end{equation}
Let $z\in\gamma(S_l)$, i.e., there is $a\in S_l$ such that $\gamma(a)=z$. By 1. $a\leq^+\bp_l$. By Lemma \ref{fullpath}, there are $x\in\delta(\bp_l)$ and $y\in\delta(a)$ such that $x\leq^+y$. Hence $x<^+\gamma(a)=z$.  By Proposition \ref{cond2}, since $x\in\delta(\bp_l)$, it follows that $z\not\in\delta(\bp_l)$. This shows (\ref{e1}).

By Lemma \ref{equations_pfs}, we have
\begin{equation}\label{e2}
\delta(S_l)=\delta(S_l-\delta({S_{l+1}}))\cup\iota(S_{l+1})
\end{equation}
Since $\delta(\bp_l)=\delta(S_l-\delta({S_{l+1}}))$ and $\iota(S_{l+1})\subseteq\gamma(S_l)$, we have by (\ref{e1})
\begin{equation}\label{e3}
\delta(S_l-\delta({S_{l+1}}))\cap\iota(S_{l+1})=\emptyset
\end{equation}
Next we shall show that
\begin{equation}\label{e4}
  \iota(S_{l+1})=\gamma(S_l)\cap\delta(S_l)
\end{equation}
The inclusion $\subseteq$ is obvious. Let $x\in\gamma(S_l)\cap\delta(S_l)$. Hence there are $a,b\in S_l$ such that $\gamma(a)=x\in\delta(b)$.  We can assume that $a$ is $<^+$-maximal with this property. As $a<^-b$, neither $a$ nor $b$ is equal to the $<^+$-greatest element $\bp_l\in S_l$. Therefore there is $\alpha\in S_{l+1}$ such that $a\in\delta(\alpha)$. If we were to have $x=\gamma(a)=\gamma\gamma(\alpha)$, then $\gamma(\alpha)$ would be a $<^+$-greater element than $a$ with $\gamma(\gamma(\alpha))=x$. So $\gamma(a)\neq\gamma\gamma(\alpha)$. Clearly, $x\in\gamma\delta(\alpha)$. By globularity, $x\in\delta\delta(\alpha)$, as well. Thus $x\in\iota(\alpha)$, and (\ref{e4}) is shown.

Using (\ref{e1}), (\ref{e2}), (\ref{e3}), and (\ref{e4}) we have
\[ \delta(\bp_l)=\delta(S_l-\delta(S_{l+1})) = \]
\[ =\delta(S_l)-\iota(S_{l+1}) = \]
\[ =\delta(S_l)-(\gamma(S_l)\cap\delta(S_l)) = \]
\[ =\delta(S_l)-\gamma(S_l) \]
as required.

Ad 6. By 1. and 2. it is enough to show
 \[ \delta^{(l)}(S_{k-1}) = \delta^{(l)}(\bp_{k-1}),\]
  for $l<k-2$. The inclusion $\supseteq$ is obvious.

Pick $x\in S_{k-1}$. We have an upper path $x,a_1,\ldots ,a_r,\bp_{k-1}$. By Corollary \ref{atlas}, as $\gamma(a_i)\in\delta(a_{i+1})$, we have
\[ \delta^{(l)}(a_i)= \delta^{(l)}\gamma(a_i)\subseteq \delta^{(l)}(\delta(a_{i+1})) = \delta^{(l)}(a_{i+1}) \]
for $i=0,\ldots , r-1$. Then, by transitivity of $\subseteq$ and  again Corollary \ref{atlas}, we get
\[ \delta^{(l)}(x)\subseteq \delta^{(l)}(a_1) \subseteq \delta^{(l)}(a_r)\subseteq \delta^{(l)}(\gamma(a_r))=\delta^{(l)}(\bp_{k-1}). \]
This ends the proof of the inclusion $\subseteq$ and 6. $~\Box$

\begin{lemma}\label{norm_positive1}
Let $S$ be a positive opetopic cardinal of dimension at least $k$.
Then
\begin{enumerate}
 \item $S$ is $(k-1)$-principal iff $\;\bd^{(k)}(S)$ is normal iff  $\;\bc^{(k-1)}(S)$ is principal,
 \item if $S$ is normal, so is $\bd(S)$,
 \item if $S$ is normal, $\bc(S)$ is principal.
\end{enumerate}
\end{lemma}

{\it Proof.}~ The whole Lemma is an easy consequence of Lemma \ref{equations_pfs}. We shall show 1., leaving 2. and 3. for the reader.

First note that all three conditions in 1. imply that $|S_l-\delta(S_{l+1}|=1$, for $l<k-2$. In addition to this these conditions say:
\begin{enumerate}
  \item $S$ is $(k-1)$-principal iff $|S_l-\delta(S_{l+1})|=1$, for $l= k-2, k-1$.
  \item $\bd^{(k)}(S)$ is normal iff
\begin{enumerate}
  \item $|S_{k-1}-\delta(S_k-\gamma(S_{k+1}))|=1$, and
  \item $|S_{k-2}-\delta(S_{k-1})|=1$.
\end{enumerate}
  \item $\bc^{(k-1)}(S)$ is principal iff
\begin{enumerate}
  \item $|(S_{k-1}-\iota(S_{k+1}))-\delta(S_k-\delta({S_{k+1}}))|=1$, and
  \item $|S_{k-2} -\delta(S_{k-1}-\iota({S_{k+1}}))|=1$.
\end{enumerate}
\end{enumerate}
So the equivalence of these conditions follows directly from Lemma \ref{equations_pfs}. $~~\Box$

Let $N$ be a normal positive opetopic cardinal of dimension $n$. We define a $(n+1)$-hypergraph $N^\bullet$ that contains two additional faces: $\bp_{n+1}^{N^\bullet}$ of dimension $n+1$, and $\bp_{n}^{N^\bullet}$ of dimension $n$. We shall drop superscripts if it does not lead to confusions.  We also put
 \[ \delta(\bp_{n+1})=N_n, \;\;\;\; \gamma(\bp_{n+1})=\bp_n, \]
  \[  \delta(\bp_{n})=\delta(N_n)-\gamma(N_n),     \;\;\;\;\gamma(\bp_n)=\bp_{n-1}(=\gamma(N_n)-\delta(N_n)). \]
As $N$ is normal, the $\gamma(N_n)-\delta(N_n)$ has one element so $\gamma(\bp_n)$ is well defined. This determines $N^\bullet$ uniquely. $N^\bullet$ is called {\em a simple extension of}\index{simple extension} $N$.

 {\em Example.} For a normal positive opetopic cardinal $N$ like this
\begin{center} \xext=1000 \yext=450
\begin{picture}(\xext,\yext)(\xoff,\yoff)
\putmorphism(220,400)(1,0)[x_2`x_1`f_1]{600}{1}a
\putmorphism(50,120)(1,2)[``f_2]{90}{1}l
\putmorphism(780,440)(1,-2)[``f_0]{230}{1}r
\putmorphism(0,0)(1,0)[x_3`x_0`]{1000}{0}b
\end{picture}
\end{center}
the hypergraph $N^\bullet$ looks like this
\begin{center} \xext=1000 \yext=500
\begin{picture}(\xext,\yext)(\xoff,\yoff)
\putmorphism(220,450)(1,0)[x_2`x_1`f_1]{600}{1}a
\putmorphism(50,170)(1,2)[``f_2]{90}{1}l
\putmorphism(780,490)(1,-2)[``f_0]{230}{1}r
\putmorphism(0,50)(1,0)[x_3`x_0`\bp^{N^\bullet}_1]{1000}{1}b
\put(400,200){\makebox(200,100){$\Downarrow\bp^{N^\bullet}_2$}}
\end{picture}
\end{center}

 We have

\begin{proposition}\label{bullet}
Let $N$ be a normal positive opetopic cardinal of dimension $n$. Then
\begin{enumerate}
 \item $N^\bullet$ is a positive opetope of dimension $n+1$.
 \item We have $\bd (N^\bullet)  \cong N $, $\bc (N^\bullet)  \cong (\bd N)^\bullet$.
 \item If $N$ is a principal, then $N\cong(\bd N)^\bullet$.
 \item If $T$ is a positive positive opetopic cardinal contained in $N^\bullet$, then either $T=N^\bullet$ or $T=\bc(N^\bullet)$ or $T\subseteq N$.
\end{enumerate}
\end{proposition}

{\it Proof.}~ Ad 1. We shall check globularity of the new added cells.  The other conditions are simple.

For $\bp_{n+1}$, we have:
\[ \gamma\gamma(\bp_{n+1})=\gamma(\bp_n)= \]
\[ = \gamma(N_n)-\delta(N_n) = \gamma\delta(\bp_{n+1})-\delta\delta(\bp_{n+1}) \]
 and
\[\delta\gamma(\bp_{n+1})= \delta(\bp_n)= \]
\[ = \delta(N_n)-\gamma(N_n) = \delta\delta(\bp_{n+1})-\gamma\delta(\bp_{n+1}). \]
 So globularity holds for $\bp_{n+1}$.

For $\bp_{n}$, using Lemmas \ref{norm_positive0}, \ref{equations_pfs} and normality of $N$, we have:
\[ \gamma\gamma(\bp_n)= \gamma(\bp_{n-1})= \bp_{n-2} = \]
\[ = \gamma(N_{n-1})-\delta(N_{n-1}) =  \]
\[ = \gamma(N_{n-1}-\gamma(N_n))-\delta(N_{n-1}-\gamma(N_n)) =  \]
\[ = \gamma(\delta(N_n)-\gamma(N_n))-\delta(\delta(N_n)-\gamma(N_n)) =  \]
\[ = \gamma\delta(\bp_{n})-\delta\delta(\bp_{n})  \]
and similarly
 \[ \delta\gamma(\bp_n)= \delta(\bp_{n-1})= \]
 \[ =\delta(N_n)-\gamma(N_n) = \]
 \[ = \delta(\delta(N_n)-\gamma(N_n))-\gamma(\delta(N_n)-\gamma(N_n)) = \]
 \[ = \delta\delta(\bp_n)-\gamma\delta(\bp_n) \]
So globularity  holds for $\bp_{n}$, as well.

Ad 2. The first isomorphism is obvious.

The faces of $(N^\bullet)$, $\bc (N^\bullet)$, $\bd N$, and $(\bd N)^\bullet$ are as in the tables
 \[ \begin{array}{|c|c|c|}
    dim & (N^\bullet)       & \bc (N^\bullet)\\ \hline
    n+1 & \{ \bp_{n+1}^{N^\bullet}\} & \emptyset   \\ \hline
    n   & N_n\cup \{ \bp_n^{N^\bullet}\} &  \{ \bp_n^{N^\bullet}\}  \\ \hline
    n-1 & N_{n-1}   &  N_{n-1}-(\gamma(N_n)\cap\delta(N_n))  \\ \hline
    n-2 & N_{n-2} & N_{n-2}
  \end{array}
\]
and
\[ \begin{array}{|c|c|c|}
    dim & \bd N & (\bd N)^\bullet\\ \hline
    n+1 & \emptyset & \emptyset \\ \hline
    n   & \emptyset &\{ \bp_n^{(\bd N)^\bullet}\} \\ \hline
    n-1 & N_{n-1}-\gamma(N_n) & (N_{n-1}-\gamma(N_n))\cup\{ \bp_{n-1}^{(\bd N)^\bullet}\} \\ \hline
    n-2 & N_{n-2} & N_{n-2}
  \end{array}
\]
We define the isomorphism $f:\bc (N^\bullet)\lra (\bd N)^\bullet$ as follows
 \[ f_n(\bp_{n+1}^{N^\bullet} )= \bp_{n+1}^{(\bd N)^\bullet}, \]
\[  f_{n-1}(x) = \left\{ \begin{array}{ll}
         \bp_{n-1}^{(\bd N)^\bullet}  & \mbox{if $x=\gamma(\bp_n^{N^\bullet})$,}  \\
         x & \mbox{otherwise.}
                                    \end{array}
                \right. \]

and $f_l=1_{N_l}$ for $l<n-1$. Clearly, all $f_i$'s are bijective. The preservation of the domains and codomains is left for the reader.

3. is left as an exercise.

Ad 4. If $\bp_{n+1}\in T_{n+1}$, then $T=N^\bullet$. If $\bp_n\not\in T_n$, then $T\subseteq N$.

Suppose that $\bp_{n+1}\not\in T_{n+1}$ but $\bp_n\in T_n$. Since $N^\bullet=[\bp_{n+1}]$, by Lemma \ref{o-cat-dom-codom} it is enough to show that $T=[\bp_n]$. Clearly $[\bp_n]\subseteq T$. As $[\bp_n]_l=N_l$, for $l<n-1$, we have $[\bp_n]_l=T_l$, for $l<n-1$, as well.

Fix $x\in N_n$.  As $x\in\delta(\bp_{n+1})$ and $\gamma(\bp_{n+1})=\bp_n$, we have $x<^{N^\bullet,+}\bp_n$. So by Corollary \ref{order} $x\not\perp^{N^\bullet,-}_l\bp_n$, for any
$l\leq n$.  Thus we cannot have $x\perp^{T,-}_l\bp_n$, for any $l\leq n$, as well. As $T$ is a positive opetopic cardinal, again by Corollary \ref{order}, $x\not\in T$. Since $x$ was an arbitrary element of $N_n$, we have $T_n=\{ \bp_n \}=[\bp_n]_n$.

It remains to show that $T_{n-1}=[\bp_n]_{n-1}$. Suppose that $x\in N_{n-1}-(\delta(\bp_n)\cup\gamma(\bp_n))$. Then $x<^{N,+}\gamma(\bp_n)$ and hence $x\not\perp^{N^\bullet,-}_l\gamma(\bp_n)$, for $l\leq n$. So $x$ and $\gamma(\bp_n)$ cannot be $<^{T,-}_l$ comparable, for $l\leq n$. Since, as we have shown, $N_n\cap T_n=\emptyset$, it follows that $x$ and $\gamma(\bp_n)$ cannot be $<^{T,+}$ comparable. So by Lemma \ref{order}, $x\not\in T_{n-1}$, i.e., $T_{n-1}=\delta(\bp_n)\cup\gamma(\bp_n)=[\bp_n]_{n-1}$. $~~\Box$

\section{Decomposition of positive opetopic cardinals}\label{sec-decomposition}

Let $T$ be a positive opetopic cardinal, $X\subseteq T$ a subhypergraph of $T$, $k\in\o$, $a\in (T_k -\iota(T_{k+2}))$. We define two subhypergraphs of $T$, $X^{\da a}$ and $X^{\ua a}$, as follows:

\[ X^{\da a}_l = \left\{ \begin{array}{ll}
           \{ \alpha\in X_l : \gamma^{(k)}(\alpha)\leq^+a \}  &  \mbox{for $l> k$,}  \\
            \{ b\in X_k : b\leq^+a \;{\rm or}\;\; b\not\in\gamma(X_{k+1}) \} & \mbox{for $l=k$} \\
            X_l &\mbox{for $l<k$.}
                                    \end{array}
                \right. \]

\[ X^{\ua a}_l = \left\{ \begin{array}{ll}
            \{ \alpha\in X_l : \gamma^{(k)}(\alpha)\not\leq^+a \} &  \mbox{for $l> k$,}  \\
            \{ b\in X_k : b\not<^+a \mbox{ or } b\not\in\delta (X_{k+1})\} &\mbox{for $l=k$} \\
            X_{k-1} - \iota(X^{\da a}_{k+1}) &\mbox{for $l=k-1$} \\
            X_l &\mbox{for $l<k-1$.}
                                    \end{array}
                \right. \]
Intuitively, if $X$ is a positive opetopic cardinal contained in $T$, $X^{\da a}$ is the least positive opetopic cardinal contained in $X$ that contains faces `smaller or equal' $a$ and can be $k$-pre-composed with the `rest' to get $X$. $X^{\ua a}$ is this `rest' or in other words it is the largest positive opetopic cardinal contained in $X$ that can be $k$-post-composed with $X^{\da a}$ to get $X$ (or the largest positive opetopic cardinal contained in $X$ that does not contain  faces `smaller' than $a$).
Note that $a$ does not need to be a face in $X$, in general.

{\em Examples.} If $X$ is a hypergraph $a\in T$, then $X^{\da a}$ is a hypergraph, as well.  However, this is not the case with $X^{\ua a}$, if $a\in\iota(T)$, as we can see below:
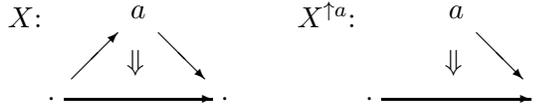
\begin{center}
\xext=1800 \yext=350
\begin{picture}(\xext,\yext)(\xoff,\yoff)
\settriparms[-1`1`1;300]
 \putAtriangle(150,0)[a`\cdot`\cdot;``]
  \put(410, 100){$\Downarrow$}
  \put(0, 250){$X$:}

  \settriparms[0`1`1;300]
 \putAtriangle(1250,0)[a`\cdot`\cdot;``]
  \put(1510, 100){$\Downarrow$}
  \put(1000,250){$X^{\uparrow a}$:}
\end{picture}
\end{center}
Here $X=T$. The faces in the domain of the 2-dimensional face are not in $X^{\ua a}$, i.e., $X^{\ua a}$ is not closed under $\delta$.

To see some real decompositions, let fix a positive opetopic cardinal $T$ as follows:
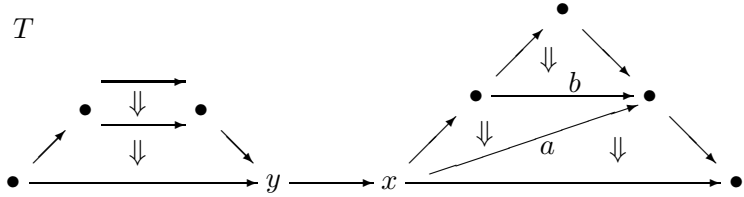
\begin{center} \xext=2400 \yext=640
\begin{picture}(\xext,\yext)(\xoff,\yoff)
 \put(0,500){$T$}
\settriparms[-1`0`0;250]
 \putAtriangle(0,0)[\bullet`\bullet`;``]
 \settriparms[0`1`0;250]
 \putAtriangle(400,0)[\bullet``y;``]
 \putmorphism(0,0)(1,0)[\phantom{\bullet}`\phantom{\bullet}`]{900}{1}b
  \put(400,80){$\Da$}
 \putmorphism(250,200)(1,0)[\phantom{\bullet}`\phantom{\bullet}`]{400}{1}b
 \putmorphism(250,350)(1,0)[\phantom{\bullet}`\phantom{\bullet}`]{400}{1}a
 \put(400,250){$\Da$}

\putmorphism(900,0)(1,0)[\phantom{\bullet}`\phantom{\bullet}`]{400}{1}b

 \put(1920,315){$b$} \put(1820,100){$a$}
 \settriparms[-1`1`1;300]
 \putAtriangle(1600,300)[\bullet`\bullet `\bullet;``]
 \put(1820,400){$\Da$}

 \settriparms[-1`0`0;300]
 \putAtriangle(1300,0)[\phantom{\bullet}`x`;``]
 \settriparms[0`1`0;300]
 \putAtriangle(1900,0)[\phantom{\bullet}``\bullet ;``]
 \putmorphism(1300,0)(1,0)[\phantom{\bullet}`\phantom{\bullet}`]{1200}{1}b
 \putmorphism(1650,100)(3,1)[\phantom{\bullet}`\phantom{\bullet}`]{300}{1}a
 \put(1600,145){$\Da$}
 \put(2060,90){$\Da$}
\end{picture}
\end{center}
Clearly $x,y,a,b\in T-\iota(T)$. Then
\begin{center} \xext=2400 \yext=640
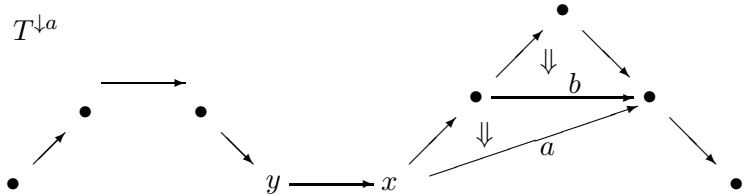
\begin{picture}(\xext,\yext)(\xoff,\yoff)
 \put(0,500){$T^{\da a}$}
\settriparms[-1`0`0;250]
 \putAtriangle(0,0)[\bullet`\bullet`;``]
 \settriparms[0`1`0;250]
 \putAtriangle(400,0)[\bullet``y;``]
 \putmorphism(250,350)(1,0)[\phantom{\bullet}`\phantom{\bullet}`]{400}{1}a

\putmorphism(900,0)(1,0)[\phantom{\bullet}`\phantom{\bullet}`]{400}{1}b

 \put(1920,315){$b$} \put(1820,100){$a$}
 \settriparms[-1`1`1;300]
 \putAtriangle(1600,300)[\bullet`\bullet `\bullet;``]
 \put(1820,400){$\Da$}

 \settriparms[-1`0`0;300]
 \putAtriangle(1300,0)[\phantom{\bullet}`x`;``]
 \settriparms[0`1`0;300]
 \putAtriangle(1900,0)[\phantom{\bullet}``\bullet ;``]
 \putmorphism(1650,100)(3,1)[\phantom{\bullet}`\phantom{\bullet}`]{300}{1}a
 \put(1600,145){$\Da$}
\end{picture}
\end{center}
and
\begin{center} \xext=2400 \yext=500
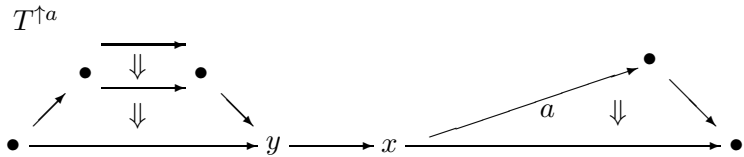
\begin{picture}(\xext,\yext)(\xoff,\yoff)
 \put(0,400){$T^{\ua a}$}
\settriparms[-1`0`0;250]
 \putAtriangle(0,0)[\bullet`\bullet`;``]
 \settriparms[0`1`0;250]
 \putAtriangle(400,0)[\bullet``y;``]
 \putmorphism(0,0)(1,0)[\phantom{\bullet}`\phantom{\bullet}`]{900}{1}b
  \put(400,80){$\Da$}
 \putmorphism(250,200)(1,0)[\phantom{\bullet}`\phantom{\bullet}`]{400}{1}b
 \putmorphism(250,350)(1,0)[\phantom{\bullet}`\phantom{\bullet}`]{400}{1}a
 \put(400,250){$\Da$}

\putmorphism(900,0)(1,0)[\phantom{\bullet}`\phantom{\bullet}`]{400}{1}b

 \put(1820,100){$a$}
 \settriparms[0`0`0;300]
 \putAtriangle(1600,300)[` `\bullet;``]
 \settriparms[0`0`0;300]
 \putAtriangle(1300,0)[\phantom{\bullet}`x`;``]
 \settriparms[0`1`0;300]
 \putAtriangle(1900,0)[\phantom{\bullet}``\bullet ;``]
 \putmorphism(1300,0)(1,0)[\phantom{\bullet}`\phantom{\bullet}`]{1200}{1}b
 \putmorphism(1650,100)(3,1)[\phantom{\bullet}`\phantom{\bullet}`]{300}{1}a
 \put(2060,90){$\Da$}
\end{picture}
\end{center}
Moreover, with
\begin{center} \xext=2400 \yext=500
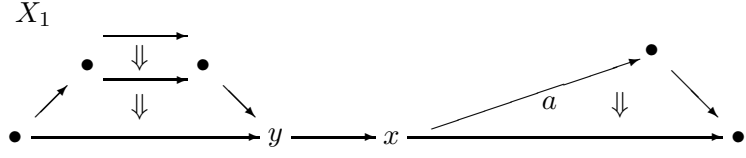
\begin{picture}(\xext,\yext)(\xoff,\yoff)
 \put(0,400){$X_1$}
\settriparms[-1`0`0;250]
 \putAtriangle(0,0)[\bullet`\bullet`;``]
 \settriparms[0`1`0;250]
 \putAtriangle(400,0)[\bullet``y;``]
 \putmorphism(0,0)(1,0)[\phantom{\bullet}`\phantom{\bullet}`]{900}{1}b
  \put(400,80){$\Da$}
 \putmorphism(250,200)(1,0)[\phantom{\bullet}`\phantom{\bullet}`]{400}{1}b
 \putmorphism(250,350)(1,0)[\phantom{\bullet}`\phantom{\bullet}`]{400}{1}a
 \put(400,250){$\Da$}

\putmorphism(900,0)(1,0)[\phantom{\bullet}`\phantom{\bullet}`]{400}{1}b

 \put(1820,100){$a$}
  \settriparms[0`0`0;300]
 \putAtriangle(1600,300)[` `\bullet;``]
 \settriparms[0`0`0;300]
 \putAtriangle(1300,0)[\phantom{\bullet}`x`;``]
 \settriparms[0`1`0;300]
 \putAtriangle(1900,0)[\phantom{\bullet}``\bullet ;``]
 \putmorphism(1300,0)(1,0)[\phantom{\bullet}`\phantom{\bullet}`]{1200}{1}b
 \putmorphism(1650,100)(3,1)[\phantom{\bullet}`\phantom{\bullet}`]{300}{1}a
 \put(2060,90){$\Da$}
\end{picture}
\end{center}
we have $X_1^{\ua b}=X_1$ and
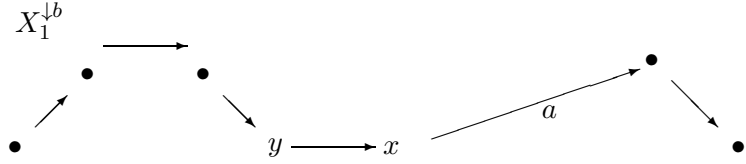
\begin{center} \xext=2400 \yext=500
\begin{picture}(\xext,\yext)(\xoff,\yoff)
 \put(0,400){$X_1^{\da b}$}
\settriparms[-1`0`0;250]
 \putAtriangle(0,0)[\bullet`\bullet`;``]
 \settriparms[0`1`0;250]
 \putAtriangle(400,0)[\bullet``y;``]
 \putmorphism(250,350)(1,0)[\phantom{\bullet}`\phantom{\bullet}`]{400}{1}a

\putmorphism(900,0)(1,0)[\phantom{\bullet}`\phantom{\bullet}`]{400}{1}b

 \put(1820,100){$a$}
  \settriparms[0`0`0;300]
 \putAtriangle(1600,300)[` `\bullet;``]
 \settriparms[0`0`0;300]
 \putAtriangle(1300,0)[\phantom{\bullet}`x`;``]
 \settriparms[0`1`0;300]
 \putAtriangle(1900,0)[\phantom{\bullet}``\bullet ;``]
 \putmorphism(1650,100)(3,1)[\phantom{\bullet}`\phantom{\bullet}`]{300}{1}a
\end{picture}
\end{center}
i.e., $X_1^{\da b}=\bd^{(1)}(X_1)$. For
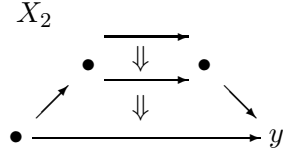
\begin{center} \xext=900 \yext=500
\begin{picture}(\xext,\yext)(\xoff,\yoff)
 \put(0,400){$X_2$}
\settriparms[-1`0`0;250]
 \putAtriangle(0,0)[\bullet`\bullet`;``]
 \settriparms[0`1`0;250]
 \putAtriangle(400,0)[\bullet``y;``]
 \putmorphism(0,0)(1,0)[\phantom{\bullet}`\phantom{\bullet}`]{900}{1}b
  \put(400,80){$\Da$}
 \putmorphism(250,200)(1,0)[\phantom{\bullet}`\phantom{\bullet}`]{400}{1}b
 \putmorphism(250,350)(1,0)[\phantom{\bullet}`\phantom{\bullet}`]{400}{1}a
 \put(400,250){$\Da$}
\end{picture}
\end{center}
we have $X_2^{\da x}=X_2$ and $X_2^{\ua x}=\{ y \}$.

We have

\begin{lemma}\label{decomp0}
 Let $T$ be a positive opetopic cardinal, $X\subseteq T$ a subhypergraph of $T$, $a\in (T -\iota(T))$, $a\in X_k$. Then
\begin{enumerate}
  \item $X^{\da a}$ and $X^{\ua a}$ are positive opetopic cardinals;
  \item $\bc^{(k)}(X^{\da a})=\bd^{(k)}(X^{\ua a})=X^{\da a}\cap X^{\ua   a}$;
  \item $\bd^{(k)}(X^{\da a})=\bd^{(k)}(X)$, $\bc^{(k)}(X^{\ua   a})=\bc^{(k)}(X)$;
  \item $X= X^{\ua a}\oplus_k X^{\da a}=X^{\ua a}\cup X^{\da a}$.
\end{enumerate}
\end{lemma}

{\it Proof.}~Ad 1.  The verification that both $X^{\da a}$ and $X^{\ua a}$ are closed under $\gamma$ and $\delta$ is routine.

For any $k$, if $x,y\in X^{\da a}_k$, then $x<^{+,X}y$ iff $x<^{+,X^{\da a}}y$. Similarly, for any $k$, if $x,y\in X^{\ua a}_k$, then $x<^{+,X}y$ iff $x<^{+,X^{\ua a}}y$. Thus by Lemma \ref{subhyp} both $X^{\da a}$ and $X^{\ua a}$ are positive opetopic cardinals.

Ad 2. Let us spell in detail both sides of the equation.

$\bc^{(k)}( X^{\da a})$ is:
\begin{enumerate}
  \item $\bc^{(k)}(X^{\da a})_l = \emptyset$, for $l>k$;
  \item $\bc^{(k)}(X^{\da a})_k = \\ (\{ b\in X_k : b\leq^+a \}\cup (X_k-\gamma(X_{k+1})))-\delta ( \{ \alpha\in X_{k+1} : \gamma(\alpha)\leq^+a \})$;
  \item $\bc^{(k)}(X^{\da a})_{k-1} = X_{k-1}-\iota(X^{\da a}_{k+1})$;
  \item $\bc^{(k)}(X^{\da a})_l = X_l$, for $l<k-1$.
\end{enumerate}
and $\bd^{(k)}( X^{\ua a})$ is:
\begin{enumerate}
  \item $\bd^{(k)}(X^{\ua a})_l =  \emptyset  $, for $l>k$;
   \item $\bd^{(k)}(X^{\ua a})_k = \\ \{ b\in X_k : b\not<^+a \mbox{ or } b\not\in\delta (X_{k+1})\} -   \gamma (X_{k+1}-\{ \alpha\in X_{k+1} : \gamma(\alpha)\leq^+a \}) $;
  \item $\bd^{(k)}(X^{\ua a})_{k-1} = X_{k-1}-\iota(X^{\da a}_{k+1})$;
  \item $\bd^{(k)}(X^{\ua a})_l = X_l$, for $l<k-1$.
\end{enumerate}

Thus to show that $\bc^{(k)}(X^{\da a})=\bd^{(k)}(X^{\ua a})$, we need to verify that $\bc^{(k)}(X^{\da a})_k=\bd^{(k)}(X^{\ua a})_k$. As both sets are contained in $X_k$, we can compare their complements. We have
\[ X_k-\bc^{(k)}(X^{\da a})_k = \{ b\in\delta(X_{k+1}):b<^+a \}\cup \gamma(X_{k+1}-\{ \alpha\in X_{k+1}:\gamma(\alpha)\not\leq^+a \}) \]
 and
\[ X_k-\bd^{(k)}(X^{\ua a})_k = \{ b\in\gamma(X_{k+1}):b\not\leq^+a \}\cup \delta(\{ \alpha\in X_{k+1}\gamma(\alpha)\leq^+a \}). \]
 But it easy to see that
\[  \{ b\in\delta(X_{k+1}):b<^+a \}= \delta(\{ \alpha\in X_{k+1}\gamma(\alpha)\leq^+a \}) \]
 and
\[ \gamma(X_{k+1}-\{ \alpha\in X_{k+1}:\gamma(\alpha)\not\leq^+a \})=\{ b\in\gamma(X_{k+1}):b\not\leq^+a \}. \]
The second equality
uses the fact that $a\not\in \iota(T)$. Thus  $\bc^{(k)}(X^{\da a})_k=\bd^{(k)}(X^{\ua a})_k$, as required.

Ad 3. To see that $\bc^{(k)}(X^{\ua a})=\bc^{(k)}(X)$, it is enough to note that $\iota(X_{k+1})=\iota(X^{\da a }_{k+1})\cup\iota(X^{\ua a }_{k+1})$. The equation $\bd^{(k)}(X^{\da a})=\bd^{(k)}(X)$ is even simpler.

Ad 4. Obvious. $~~\Box$

\begin{corollary} Let $T$ be a positive opetopic cardinal, $k\in\o$, $a\in (T_k-\iota(T_{k+2}))$. Then the square
\begin{center}
\xext=650 \yext=600
\begin{picture}(\xext,\yext)(\xoff,\yoff)
 \setsqparms[-1`-1`-1`-1;650`500]
 \putsquare(0,100)[T`T^{\da a}`T^{\ua a}`\bc^{(k)}(T^{\da a});``\bc^{(k)}_{T^{\da a}}`\bd^{(k)}_{T^{\ua a}}]
\end{picture}
\end{center}
is a special pushout in $\pOpeCard$.
\end{corollary}

 {\it Proof.}~ Follows immediately from Lemmas \ref{o-cat-comp} and \ref{decomp0}. $~~\Box$

\vskip 3mm

We need more notions and notations. Let $X$, $T$ be positive opetopic cardinals $X\subseteq T$, $a\in (T-\iota(T))$.  The decomposition $X=X^{\da a}\cup X^{\ua a}$ is said to be {\em proper}\index{decomposition!proper} iff $size(X^{\da a}),size(X^{\ua a})<size(X)$.  If the decomposition $X=X^{\da a}\cup X^{\ua a}$ is proper then $a$ is said to be a {\em saddle face}\index{face!saddle -} of $X$. $Sd(X)$  is the set of saddle faces of $X$, $Sd(X)_k=Sd(X)\cap X_k$.

\begin{lemma} \label{Sd}Let $X$, $S$, $T$ be positive opetopic cardinals, $X\subseteq T$, $l\in\o$. Then
\begin{enumerate}
  \item if $a\in (T_l-\iota(T))$, then  $a\in Sd(X)$ iff there are   $\alpha,\beta\in X_{l+1}$ such that $\gamma(\alpha)\leq^+a$ and $\gamma(\beta)\not\leq^+a$;
  \item if $\bc^{(k)}(S)=\bd^{(k)}(T)$, then
\[  size(S\oplus_kT)_l = \left\{ \begin{array}{ll}
          size(S)_l+ size(T)_l  & \mbox{ if $l>k$,}  \\
       size(T)_l & \mbox{ if $l\leq k$;}
                                    \end{array}
                \right. \]
  \item $size(S)_k\geq 1$ iff $k\leq dim(S)$;
  \item if $a\in Sd(S)_k$, then $size(S)_{k+1}\geq 2$;
  \item $S$ is principal iff $Sd(S)$ is empty.
\end{enumerate}
\end{lemma}

{\it Proof.}~ We shall show 5. The rest is easy.

If there is $a\in Sd(S)_k $, then by 2., 3. and Lemma \ref{decomp0} we have that $size(S)_{k+1}=size(S^{\da a})_{k+1}+size(S^{\ua a})_{k+1}\geq 1+1>1$. So in that case $S$ is not principal.

For the converse, assume that $S$ is not principal. Fix $k\in\o$ such that $size(S)_{k+1}>1$. Thus there are $a,b\in S_{k+1}$, such that $a\neq b$. Suppose $\gamma(a)\in\iota(\alpha)$, for some $\alpha\in S_{k+2}$. Then by Lemma \ref{tech_lemma2}, $a<^+\gamma(\alpha)$ contrary to the assumption on $a$. Hence $a\in S- \iota(S)$ and for similar reasons $b\in S- \iota(S)$. We have $a\not\perp^+b$ and, by pencil linearity, $\gamma(a)\neq\gamma(b)$. Then either $\gamma(a)\not<^+\gamma(b)$ and then $\gamma(b)\in Sd(S)_k$ or $\gamma(b)\not<^+\gamma(a)$ and then $\gamma(a)\in Sd(S)_k$. In either case $Sd(S)$ is not empty, as required. $~\Box$

\begin{lemma}\label{decomp1new}
Let $T$, $X$ be positive opetopic cardinals, $X\subseteq T$,  and $a,x\in X-\iota(X)$, $k=dim(x)<dim(a)=m$.
\begin{enumerate}
  \item We have the following equations of positive opetopic cardinals:
  \[ X^{\da x \da a}= X^{\da a \da x}, \;\;\; X^{\da x \ua a}= X^{\ua a \da x}, \;\;\; X^{\ua x \da a}= X^{\da a \ua x}, \;\;\; X^{\ua x \ua a}= X^{\ua a \ua x}, \]
   i.e., `the decompositions o different dimensions commute'.
   \item If $x\in Sd(X)$, then $x\in Sd(X^{\da a})\cap Sd(X^{\ua a})$.
   \item Moreover, we have the following equations concerning domains and codomains
   \[ \bc^{(k)}(X^{\da x \da a})= \bc^{(k)}(X^{\da x \ua a})= \bd^{(k)}(X^{\ua x \da a})= \bd^{(k)}(X^{\ua x \ua a})\]
   \[ \bc^{(m)}(X^{\da x \da a})= \bd^{(m)}(X^{\da x \ua a}),\;\;\;\; \bc^{(m)}(X^{\ua x \da a})= \bd^{(m)}(X^{\ua x \ua a}).\]
   \item Finally, we have the following equations concerning compositions
   \[  X^{\da x \ua a}\oplus_mX^{\da x \da a}=X^{\da x},\;\;\;\;  X^{\ua x \ua a}\oplus_mX^{\ua x \da a}=X^{\ua x}\]
    \[  X^{\ua x \da a}\oplus_k X^{\da x \da a}=X^{\da a},\;\;\;\; X^{\ua x \ua a}\oplus_kX^{\da x \ua a}=X^{\ua a}.\]
\end{enumerate}
\end{lemma}

{\it Proof.}~ Simple check. $~\Box$

\begin{lemma}\label{decomp2new}
Let $T$, $X$ be positive opetopic cardinals, $X\subseteq T$,  and $a,b\in X-\iota(X)$, $dim(a)=dim(b)=m$.
\begin{enumerate}
  \item We have the following equations of positive opetopic cardinals:
\[ X^{\da a \da b}= X^{\da b \da a}, \;\;\;
   X^{\ua a \ua b}= X^{\ua b \ua a},  \]
   i.e., `the decompositions in the same dimension and the same directions commute'.
  \item Assume $a<^+b$. Then we have the following further equations of positive opetopic cardinals:
  \[  X^{\da a}=  X^{\da a \da b},\;\;\; X^{\ua b}=  X^{\ua a \ua b},\;\;\;  X^{\da b \ua a}=  X^{\ua a \da b}. \]
  Moreover, if $a,b\in Sd(X)$, then $a\in Sd(X^{\da b})$ and $b\in Sd(X^{\ua a})$.
  \item Assume $a<^-_lb$, for some $l<m$. Then $X^{\ua b \da a}$, $X^{\ua a \da b}$,
  are positive opetopic cardinals, and
   \[ X^{\ua a \da b}\oplus_m X^{\da a}=X^{\ua b \da a}\oplus_m X^{\da b}\]
  Moreover, if $a,b\in Sd(X)$, then either there is $k$ such that  $l-1\leq k<m$ and $\gamma^{(k)}(a)\in Sd(X)$   or $a\in Sd(X^{\ua b})$ and $b\in Sd(X^{\ua a})$.
\end{enumerate}
\end{lemma}
{\it Proof.}~ Simple check. $~\Box$

\begin{lemma}\label{decomp4}
Let $T$, $X$ be positive opetopic cardinals, $X\subseteq T$, $dim(X)=n$, $l< n-1$, $a\in Sd(X)_l$. Then
\begin{enumerate}
  \item $a\in Sd(\bc X)\cap Sd(\bd X)$;
  \item $\bd(X^{\da a})=(\bd X)^{\da a}$;
  \item $\bd(X^{\ua a})=(\bd X)^{\ua a}$;
  \item $\bc(X^{\da a})=(\bc X)^{\da a}$;
  \item $\bc(X^{\ua a})=(\bc X)^{\ua a}$.
\end{enumerate}
\end{lemma}

{\it Proof.}~ The proof is again by a long and simple check. We shall check part of 5. We should consider separately cases: $l=n-2$, $l=n-3$, and $l<n-3$, but we shall check the case $l=n-3$ only. The other cases can be also shown by similar, but easier, check.

$(\bc X)^{\ua a}$ is:
\begin{enumerate}
  \item $(\bc X)^{\ua a}_l = \emptyset$,  for $l\geq n$;
  \item $(\bc X)^{\ua a}_{n-1} = \{ x\in X_{n-1} :
  \gamma^{(n-3)}(x)\not\leq^+a,\;\; x\not\in\delta(X_n) \} $;
  \item $(\bc X)^{\ua a}_{n-2} = \{ x\in X_{n-2} :
  \gamma(x)\not\leq^+a,\;\;  x\not\in\iota(X_n) \} $;
  \item $(\bc X)^{\ua a}_{n-3} = \{ x\in X_{n-3} : x\not<^+a \mbox{ or } x\not\in\delta(X_{n-2}-\iota(X_n)) \} $;
  \item $(\bc X)^{\ua a}_{n-4} = X_{n-4}- \iota(\{ x\in X_{n-2} :
  x\not\in\iota(X_n),\;  \gamma(x)\leq^+a \}) $;
  \item $X^{\da a}_l = X_l$, for $l<n-4$.
\end{enumerate}
and $ \bc(X^{\ua a})$ is:
\begin{enumerate}
  \item $ \bc(X^{\ua a})_l = \emptyset$,  for $l\geq n$;
  \item $ \bc(X^{\ua a})_{n-1} = \{ x\!\in\! X_{n-1} : \gamma^{(n-3)}(x)\not\leq^+a \}-\delta(\{ z\!\in\! X_{n} : \gamma^{(n-3)}(z)\not\leq^+a  \})$;
  \item $ \bc(X^{\ua a})_{n-2} = \{ x\!\in\! X_{n-2} : \gamma(x)\not\leq^+a \}-\iota(\{ z\!\in\! X_{n} : \gamma^{(n-3)}(z)\not\leq^+a  \})$;
  \item $ \bc(X^{\ua a})_{n-3} = \{ x\in X_{n-3} : x\not<^+a \mbox{ or } x\not\in\delta(X_{n-2})\} $;
  \item $ \bc(X^{\ua a})_{n-4} = X_{n-4}- \iota(X^{\da a}_{n-2}) $;
  \item $ \bc(X^{\ua a})_l = X_l$, for $l<n-4$.
\end{enumerate}
We need to verify the equality $(\bc X)^{\ua a}_l= \bc(X^{\ua a})_l$, for $l=n-1,\ldots, n-4$.

In dimension $n-1$, it is enough to show that if $x\in X_{n-1}$ and $z\in X_n$ so that $\gamma^{(n-3)}(x)\not\leq^+a$ and $x\in\delta(z)$, then $\gamma^{n-3}(z)\not\leq^+a$.

So assume that $x\in X_{n-1}$,  $\gamma^{(n-3)}(x)\not\leq^+a$, $z\in X_n$ such that $x\in\delta(z)$. Hence $x\lhd^+\gamma(z)$. By Lemma \ref{tech_lemma3}.5, $\gamma^{(n-3)}(x)\leq^+\gamma^{(n-3)}(z)$. Therefore $\gamma^{(n-3)}(z)\not\leq^+a$ (otherwise we would have $\gamma^{(n-3)}(x)\not\leq^+a$), as required.

In dimension $n-2$, it is enough to show that if $x\in X_{n-2}$ and $z\in X_n$ so that $x\not\leq^+a$ and $x\in\iota(z)$, then $\gamma^{n-3}(z)\not\leq^+a$.

So assume that $x\in X_{n-2}$, $z\in X_n$ so that $x\not\leq^+a$ and $x\in\iota(z)$. Hence $x\leq^+ \gamma\gamma(z)$. By Lemma \ref{tech_lemma3}.5,  $\gamma(x)\leq^+\gamma^{(n-3)}(z)$. Therefore $\gamma^{(n-3)}(z)\not\leq^+a$, as required.

The equality in dimension $n-3$ follows immediately from Lemma \ref{equations_pfs}.4.

To show that in dimension $n-4$, the above equation also holds, we shall show that
\[ \iota(X^{\da a}_{n-2})\subseteq\iota(\{ x\in X_{n-2} :  x\not\in\iota(X_n),\;  \gamma(x)\leq^+a \}) \]
Note that, by Lemma \ref{iota}.1, if $t\in X_{n-4}$ and $x\in X_{n-2}$, $y\in X_{n-1}$, $t\in\iota(x)$ and $\gamma(x)\leq^+a$ and $x=\gamma(y)$, then there is $x'\in\delta(y)$ (i.e., $x'\lhd^+x$ and hence $\gamma(x')\leq^+a$) such that $t\in\iota(x')$.

Thus, as $<^+$ is well founded, from the above observation follows that, for any $t\in X_{n-4}$ and $x\in X_{n-2}$ such that $t\in\iota(x)$ and $\gamma(x)\leq^+a$, there is $x''\leq^+x$ such that $t\in\iota(x'')$ and $x''\not\in\gamma(X)$. Then we clearly have that $x''\not\in\iota(X)$ and $\gamma(x'')\leq^+a$, as required. $~\Box$

\vskip 3mm
The following Lemma describes how one can express decompositions a special pushout in term of decompositions of its components.

\begin{lemma}\label{decomp3new}
Let $T, T_1,T_2$ be  positive opetopic cardinals, $dim(T_1),dim(T_2)>k$ such that $\bc^{(k)}(T_1)=\bd^{(k)}(T_2)$ and $T=T_2\oplus_kT_1$. Then $\bc^{(k)}(T_1)_k\cap\gamma(T_1)\neq\emptyset$. For any $a\in\bc^{(k)}(T_1)_k\cap\gamma(T_1)$, we have $a\in Sd(T)_k$ and either $T_1=T^{\da a}$ and $T_2=T^{\ua a}$ or $a\in Sd(T_1)_k$, $T^{\da a}=T_1^{\da a}$ and  $T^{\ua a}=T_2\oplus_kT_1^{\ua a}$.
\end{lemma}

{\it Proof.}~By assumption, $(T_1)_{k+1}\neq\emptyset$ and $(T_2)_{k+1}\neq\emptyset$. So $\bc^{(k)}(T_1)\cap\gamma(T_1)\neq\emptyset$. Fix $a\in\bc^{(k)}(T_1)\cap\gamma(T_1)\neq\emptyset$. Then $T^{\da a}_{k+1}\neq\emptyset$. As $T^{\da a}_{k+1}\cap (T_2)_{k+1} =\emptyset$, we must have $a\in Sd(T)_k$.

Assume that $T_1\neq T^{\da a}$. Then $T^{\da a}\subseteqnot T_1$. Hence $(T_1) - (T^{\da a})\neq\emptyset$.  But this means that $a\in Sd(T_1)_k$. The verification that the equalities $T^{\da a}=T_1^{\da a}$ and  $T^{\ua a}=T_2\oplus_kT_1^{\ua a}$ hold in this case is left as an exercise. $~\Box$

\section{Positive opetopic cardinals as positive-to-one polygraphs}\label{sec-Sstar_is_pPoly}

For the definition of positive-to-one polygraphs and related notation see Appendix. In this section we show that the image of the embedding defined in Section \ref{sec-o-cat-generated}
$$(-)*:\pOpeCard \lra \oC$$
is in fact contained in the category of polygraphs.

\begin{proposition}\label{Sstar} Let $S$ be a weak positive opetopic cardinal. Then $S^*$ is a positive-to-one polygraph whose $k$-indeterminates correspond to faces in $S_k$.
\end{proposition}
{\it Proof.}~ The proof is by induction on dimension $n$ of the weak positive opetopic cardinal $S$. For $n=0, 1$, the Proposition is obvious.

So assume that for any weak positive opetopic cardinal $T$ of dimension $n$, $T^*$ is a positive-to-one polygraph of dimension $n$, generated by faces in $T$. Suppose that $S$ is a weak positive opetopic cardinal of dimension $n+1$. We shall show that $S^*$ is a polygraph generated by faces in $S$.  Since $S_{\leq n}$ is a weak positive opetopic cardinal, by inductive assumption, $S_{\leq n}^*$ is the polygraph  generated by faces in $S_{\leq n}$. So we need to verify that, for any $\o$-functor $f:S_{\leq n}^*\lra C$ to any $\o$-category $C$ and any function $|f|:S_{n+1}\lra C_{n+1}$ such that for $a\in S_{n+1}$
\[ d_C(|f|(a)) = f(\bd([a])), \;\;\;\;\; c_C(|f|(a)) = f(\bc([a])), \]
there is a unique $\o$-functor $F:S^*\lra C$ such that
\[ F_{n+1}([a])=|f|(a), \;\;\;\;\; F_{\leq n}=f \]
as in the diagram
\begin{center}
\xext=650 \yext=1000
\begin{picture}(\xext,\yext)(\xoff,\yoff)
 \setsqparms[1`0`0`1;650`500]
 \putsquare(0,50)[S_{n+1}`S^*_{n+1}`S_{\leq n}`S_{\leq n}^*;```{[-]}]
  \putmorphism(-100,550)(0,-1)[``\delta ]{500}{1}l
  \putmorphism(0,550)(0,-1)[``\gamma ]{500}{1}r
  \putmorphism(600,550)(0,-1)[``\bd]{500}{1}l
  \putmorphism(700,550)(0,-1)[``\bc]{500}{1}r

\putmorphism(860,650)(2,1)[``]{150}{1}r
\put(860,520){\makebox(100,100){$F_{n+1}$}}

\putmorphism(360,650)(4,1)[``]{500}{1}l
\putmorphism(800,180)(2,3)[``f]{350}{1}r
\put(1180,780){\makebox(100,100){$C$}}
\put(265,430){\makebox(100,100){$[-]$}}
\put(565,730){\makebox(100,100){$|f|$}}
\end{picture}
\end{center}

We define $F_{n+1}$ as follows. For $X\in S^*_{n+1}$

\[ F_{n+1}(X) = \left\{ \begin{array}{ll}
                 id_{f(X)} & \mbox{ if $dim(X)\leq n$,} \\
                 |f|(a)  & \mbox{ if $dim(X)=n+1$, $X$ is principal and
                 $X=[a]$,} \\
                F_{n+1}(X^{\ua a})\circ_l  F_{n+1}(X^{\da a}) &
                            \mbox{ if $dim(X)=n+1$, $a\in Sd(X)_l$.}
                 \end{array}  \right. \]
$\circ_l$ refers to the composition in the $\o$-category $C$. Clearly $F_k=f_k$, for $k\leq n$. The above morphism, if well defined, clearly preserves identities. We need to verify, for $X\in S^*_{n+1}$ and $dim(X)=n+1$, three conditions:
\begin{itemize}
  \item[\bf I] $F$ is well defined, i.e.,  $F_{n+1}(X)= F_{n+1}(X^{\ua a})\circ_lF_{n+1}(X^{\da a})$ does not depend on the choice of the saddle face $a\in Sd(X)$;
  \item[\bf II] $F$ preserves the domains and codomains, i.e., $F_n(d(X))=d(F_{n+1}(X))$ and $F_n(c(X))=c(F_{n+1}(X))$;
  \item[\bf III] $F$ preserves compositions, i.e., $F_{n+1}(X)=F_{n+1}(X_2)\circ_kF_{n+1}(X_1)$   whenever $X=X_2\oplus_kX_1$ and $dim(X_1),dim(X_2)>k$.
\end{itemize}

We have an embedding $[-]:S_{\leq n}\lra S_{\leq n}^*$.

So let assume that for positive opetopic cardinals of $S$ of size less than $size(X)$ the above assumption holds. If  $size(X)_{n+1}=0$ or $X$ is principal, all three conditions are obvious. So assume that $X$ is not principal and $dim(X)=n+1$. To save on notation we write $F$ for $F_{n+1}$.

Ad {\bf I}. First we will consider two saddle faces $a,x\in Sd(X)$ of different dimension $k=dim(x)<dim(a)=m$. Using Lemma
\ref{decomp1new} we have
\begin{eqnarray*}
F(X^{\ua a})\circ_m  F(X^{\da a}) = & ind.\; hyp.\; I    \\
 = ( F(X^{\ua a\ua x})\circ_kF(X^{\ua a \da x}))\circ_m( F(X^{\da a\ua x})\circ_kF(X^{\da a\da x}))  = & MEL \\
 = (F(X^{\ua a\ua x})\circ_m F(X^{\da a\ua x}))\circ_k ( F(X^{\ua a \da x})\circ_mF(X^{\da a\da x})) = &  \\
 = ( F(X^{\ua x\ua a})\circ_mF(X^{\ua x\da a}))\circ_k  (F(X^{\da x \ua a})\circ_m F(X^{\da x\da a})) = & ind.\; hyp.\; III  \\
 =  F(X^{\ua x}) \circ_m F(X^{\da x}) &
\end{eqnarray*}
Now we will consider two saddle faces $a,b\in Sd(X)$ of the same dimension $dim(a)=dim(b)=m$. We shall use Lemma \ref{decomp2new}. Assume that $a<^-_l b$, for some $l<m$. If $\gamma^{(k)}(a)\in Sd(X)$, for some $k<m$, then this case reduces to the previous one for two pairs $a,\gamma^{(k)}(a)\in Sd(X)$ and $b,\gamma^{(k)}(a)\in Sd(X)$. Otherwise $a\in Sd(X^{\ua b })$ and $a\in Sd(X^{\ua b })$ and we have
\begin{eqnarray*}
  F(X^{\ua a})\circ_kF(X^{\da a}) = & ind.\; hyp\; I\\
 = (F(X^{\ua a\ua b}) \circ_k  F(X^{\ua a\da b}))\circ_k F(X^{\da a}) = &  \\
 = (F(X^{\ua a\ua b}) \circ_k  (F(X^{\ua a\da b})\circ_k F(X^{\da a})) = & ind\; hyp\; III\\
  =  F(X^{\ua b\ua a})\circ_k F( X^{\ua a\da b} \oplus_kX^{\da a}) = &  \\
  =  F(X^{\ua b\ua a})\circ_k F( X^{\ua b\da a} \oplus_kX^{\da b}) = &  ind\; hyp\;  III\\
  =  F(X^{\ua b\ua a})\circ_k (F(X^{\ua b\da a}) \circ_k F(X^{\da b})) = &  \\
  =  (F(X^{\ua b\ua a})\circ_k F(X^{\ua b\da a})) \circ_k F(X^{\da b}) = &  \\
 =  F(X^{\ua b}) \circ_kF(X^{\da b}) &
\end{eqnarray*}
Finally, we consider the case $a<^+b$. We have
\begin{eqnarray*}
  F(X^{\ua a})\circ_kF(X^{\da a}) = & ind.\; hyp\; I  \\
 =  ( F(X^{\ua a\ua b})\circ_kF(X^{\ua a\da b})) \circ_kF(X^{\da a}) = & \\
 =   (F(X^{\ua b})\circ_k F(X^{\da b\ua a}))\circ_k F(X^{\da b\da a}) = & \\
 = F(X^{\ua b})\circ_k (F(X^{\da b\ua a})\circ_k F(X^{\da b\da a})) = & ind\; hyp\; I \\
  =  F(X^{\ua b})\circ_kF(X^{\da b}) &
\end{eqnarray*}
This shows that $F(X)$ is well defined.

Ad {\bf II}. We shall show that the domains are preserved. The proof that  the codomains are preserved is similar.

The fact that if $Sd(X)=\emptyset$, then $F$ preserves domains and codomains follows immediately from the assumption on $f$ and $|f|$. So assume that $Sd(X)\neq\emptyset$ and  let $a\in Sd(X)$, $dim(a)=k$. We use Lemma \ref{decomp4}. We have to consider two cases  $k<n$, and $k=n$.

If $k<n$, then
\begin{eqnarray*}
F_n(d(X))= F_n(d( X^{\ua a} \oplus_kX^{\da a}))= &    \\
 = F_n( d((X^{\ua a})\oplus_k d(X^{\da a}))= &  \\
 = F_n( d(X)^{\ua a} \oplus_kd(X)^{\da a})= & ind\; hyp\; III \\
 =  F_n(d(X)^{\ua a}) \circ_k F_n(d(X)^{\da a})= & \\
 = F_n(d(X^{\ua a})) \circ_k F_n(d(X^{\da a}))= & ind\; hyp\; II \\
 =  d(F_{n+1}(X^{\ua a})) \circ_k d(F_{n+1}(X^{\da a}))= &  \\
 =  d(F_{n+1}(X^{\ua a}) \circ_kF_{n+1}(X^{\da a}))= & ind\; hyp\; I \\
 = d(F_{n+1}(X)) &
\end{eqnarray*}

If $k=n$, then
\begin{eqnarray*}
F_n(d(X))= F_n(d( X^{\ua a} \oplus_nX^{\da a}))= &   \\
 = F_n(d(X^{\da a}))= & ind\; hyp\; II \\
 = d(F_{n+1}(X^{\da a}))= & \\
 =d(F_{n+1}(X^{\ua a}))\circ_n F_{n+1}(X^{\da a}))= & ind\; hyp\; I \\
 = d(F_{n+1}(X)) &
\end{eqnarray*}

Ad {\bf III}. Suppose that $X=X_2\oplus_kX_1$ and $dim(X)\leq n+1$.  We shall show that $F$ preserves this composition.  If $dim(X_1)=k$, then $X=X_2$, $X_1=\bd^{(k)}(X_2)$.  We have
\[ F_{n+1}(X)=F_{n+1}(X_2)= \]
\[ = F_{n+1}(X_2) \circ_k 1^{(n+1)}_{F_k(\bd^{(k)}(X_2))}= \]
\[ =F_{n+1}(X_2) \circ_k1^{(n+1)}_{F_k(X_1)}=\]
\[ =F_{n+1}(X_2) \circ_k F_{n+1}(X_1)\]
The case $dim(X_2)=k$ is similar. So now assume that $dim(X_1),dim(X_2)>k$. We shall use Lemma \ref{decomp3new}. Fix $a\in\bc^{(k)}(X_1)_k\cap\gamma(X_1)$. So $a\in Sd(X)_k$. If $X_1=X^{\da a}$ and $X_2=X^{\ua a}$, then we have
\[ F(X)=F(X^{\ua a}) \circ_kF(X^{\da a})= F(X_2) \circ_kF(X_1). \]
If $a\in Sd(X_1)_k$, then

\begin{eqnarray*}
F(X)=F(X^{\ua a}) \circ_kF(X^{\da a})= & ind\; hyp\; II\\
 =(F(X_2) \circ_k F(X_1^{\ua a})) \circ_kF(X^{\da a})= &  \\
 =F(X_2) \circ_k( F(X_1^{\ua a})\circ_kF(X_1^{\da a}))= & ind\; hyp\; II \\
 F(X_2)\circ_kF(X_1) &
\end{eqnarray*}
So in any case the composition is preserved. This ends the proof of the Lemma. $~~\Box$

\vskip 2mm

Let $\wpnfs$ be the full subcategory of $\wpfs$ whose objects have dimension at most $n\geq 0$. For $n\in\o$, we have a functor
\[ (-)^{\sharp,n} : \wpnfs\lra Set\da D_{n-1} \]
such that, for $S$ in $\wpnfs$
 \[ S^{\sharp,n}=(S_n,S^*_{<n}, [\delta],[\gamma]) \]
and, for $f:S\ra T$ in $\wpnfs$, we have
\[ f^{\sharp,n}=(f_n,(f_{<n})^*). \]

\begin{corollary}
For every $n\in\o$, the functor $(-)^{\sharp,n}$ is well defined, full, faithful, and it preserves existing pushouts. Moreover, for $S$ in $\wpnfs$, we have $S^*=\overline{S^{\sharp,n}}^n$.
\end{corollary}
 {\it Proof.}~ The functor $\overline{(-)}^n: Set\da D_n \lra \pnPoly$ which is an equivalence of categories is defined in the Appendix.

Fullness and faithfulness of $(-)^{\sharp,n}$ is left for the reader.  We shall show simultaneously that for every $n\in\o$, both functors
\[  (-)^{\sharp,n} : \wpnfs\lra Set\da D_n,\hskip 10mm (-)^{*,n} : \wpnfs\lra \pnPoly \]
preserve existing pushouts. For $n=0$, there is nothing to prove. For $n=1$, this is obvious. So assume that $n\geq 1$ and that $(-)^{\sharp,n}$ preserves existing pushouts. Let
\begin{center} \xext=500 \yext=400
\begin{picture}(\xext,\yext)(\xoff,\yoff)
 \setsqparms[1`-1`-1`1;500`400]
 \putsquare(0,0)[S`S+_RT`R`T;```]
\end{picture}
\end{center}
be a pushout in $\wpnjfs$. Clearly its $n$-truncation is a pushout in $\wpnfs$. Hence by inductive hypothesis it is preserved by $(-)^{*,n}$. In dimension $n+1$, the functor $(-)^{\sharp,n+1}$ is an inclusion.  Hence, in dimension $n+1$,  this square is a pushout (of monos) in $Set$. So the whole square
\begin{center} \xext=800 \yext=450
\begin{picture}(\xext,\yext)(\xoff,\yoff)
 \setsqparms[1`-1`-1`1;800`400]
 \putsquare(0,0)[S^{\sharp,n+1}`(S+_RT)^{\sharp,n+1}`R^{\sharp,n+1}`T^{\sharp,n+1};```]
\end{picture}
\end{center}
is a pushout in $Set\da D_{n+1}$, i.e., $(-)^{\sharp,n+1}$ preserves pushouts. As $(-)^{*,n+1}$ is a composition of $(-)^{\sharp,n+1}$ with an equivalence of categories, it preserves the pushouts, as well. $~~\Box$

\begin{corollary}\label{ff-pres_sp} The functor
\[ (-)^* :\wpOpeCard \lra \pPoly \]
is full and faithful and preserves special pushouts. In particular, it is conservative.
\end{corollary}
{\it Proof.}~This follows from the previous Corollary and the fact that the functor $\overline{(-)}^n: Set\da D_n \lra \pnPoly$ is an equivalence of categories.  $~~\Box$
\vskip 3mm

\begin{corollary}\label{fc-pres_sp oC} The functor
\[ (-)^* :\wpOpeCard \lra \oC \]
is faithful, conservative and preserves special pushouts.
\end{corollary}
{\it Proof.}~ The faithfulness and preservation of special pushout follows from the previous Corollary and the fact that the functor $F_n: \pnPoly \lra \nC$ (see Appendix) is faithful and a left adjoint. Conservativity follows from previous Corollary and the fact that any isomorphic $\o$-functor between polygraphs preserves and reflects indeterminates, i.e., is a map of polygraphs.  $~~\Box$
\vskip 3mm

Let $P$ be a positive-to-one polygraph, $a$ a $k$-cell in $P$.  A {\em description of the cell} $a$\index{description of a cell}\index{cell!description of a -}  is a pair
  \[ <T_{a}, \tau_{a}: T_{a}^*\lra P> \]
where $T_{a}$ is a positive opetopic cardinal and $\tau_{a}$ is a polygraph map such that \[ \tau_{a}(T_{a})=a. \]

In the remainder of this section we shall define some specific positive opetopic cardinals that will be used later. First we define $\alpha^n$\label{def_alpha_n}, for $n\in\o$.  We put
\[ \alpha^n_l \;\;= \;\; \left\{ \begin{array}{ll}
                \emptyset    & \mbox{ if  $l> n$} \\
                \{ 2n \}     & \mbox{ if  $l= n$} \\
                \{ 2l+1,\, 2l  \}     & \mbox{ if  $0\leq l< n$}
                                    \end{array}
                            \right.\label{n.alpha0} \]

\[ d,\, c: \alpha^n_l \lra \alpha^n_{l-1} \]

\[ d(x)=\{ 2l-1 \} \;\;\;\;\;\;\;\;\; c(x)=2l-2 \]
for $x\in\alpha^n_l$, and $1\leq l\leq n$.

For example, $\alpha^4$\label{n.alpha1} can be pictured as follows:

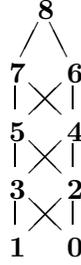
\begin{center}
\begin{picture}(1700,900)
\put(700, 0){${\bf 1}$} \put(900, 0){${\bf 0}$}
\put(720,110){\line(0,1){80}} \put(770,90){\line(1,1){110}}
\put(770,200){\line(1,-1){110}} \put(920,110){\line(0,1){80}}
\put(700,200){${\bf 3}$} \put(900, 200){${\bf 2}$}

\put(720,310){\line(0,1){80}} \put(770,290){\line(1,1){110}}
\put(770,400){\line(1,-1){110}} \put(920,310){\line(0,1){80}}

\put(700,400){${\bf 5}$} \put(900, 400){${\bf 4}$}

\put(720,510){\line(0,1){80}} \put(770,490){\line(1,1){110}}
\put(770,600){\line(1,-1){110}} \put(920,510){\line(0,1){80}}
\put(700,600){${\bf 7}$} \put(900,600){${\bf 6}$}

\put(740,690){\line(1,2){60}}
 \put(840,810){\line(1,-2){60}}
 \put(800,820){${\bf 8}$}
\end{picture}
\end{center}
i.e., $8$ is the unique cell of dimension $4$ in $\alpha^4$ that has $7$ as its domain and $6$ as its codomain, $7$ and $6$ have $5$ as its domain and $4$ as its codomain, and so on.  Note that, for any $k\leq n$, we have
\[ \bd^{(k)} \alpha^n = \alpha^k = \bc^{(k)} \alpha^n. \]

Let $n_1<n_0,n_2$ and $n_3<n_2,n_4$. We define the positive opetopic cardinals $\alpha^{n_0,n_1,n_2}$ and $\alpha^{n_0,n_1,n_2,n_3,n_4}$ as the following colimits in
$\pOpeCard$:
\begin{center} \xext=2800 \yext=650
\begin{picture}(\xext,\yext)(\xoff,\yoff)
 \setsqparms[1`-1`-1`1;700`500]
 \putsquare(0,100)[\alpha^{n_0},`\alpha^{n_0,n_1,n_2}`\alpha^{n_1}`\alpha^{n_2};
 \kappa_1`\bc^{(n_1)}_{\alpha^{n_0}}`\kappa_2`\bd^{(n_1)}_{\alpha^{n_2}}]

 \putsquare(1300,100)[\alpha^{n_0},`\alpha^{n_0,n_1,n_2,n_3,n_4}`\alpha^{n_1}`\alpha^{n_2};
 \kappa_1`\bc^{(n_1)}_{\alpha^{n_0}}`\kappa_2`\bd^{(n_1)}_{\alpha^{n_2}}]
 \setsqparms[-1`-1`-1`-1;700`500]
 \putsquare(2000,100)[\phantom{\alpha^{n_0,n_1,n_2,n_3,n_4}}`\alpha^{n_4}`
 \phantom{\alpha^{n_2}}`\alpha^{n_3};
 \kappa_3``\bd^{(n_3)}_{\alpha^{n_4}}`\bc^{(n_3)}_{\alpha^{n_2}}]
\end{picture}
\end{center}

We have

\begin{proposition}
The above colimits are preserved by the functor
\[ (-)^*:\pOpeCard \lra \pPoly. \]
Moreover, for any $\o$-category $C$, we have bijective correspondences
\[ \oC ((\alpha^n)^*,C)= C_n \]

\[ \oC ((\alpha^{n_0,n_1,n_2})^*,C) = \{(x,y)\in C_{n_0}\times C_{n_2} : c^{(n_1)}(x)=d^{(n_1)}(y) \} \]

\[ \oC ((\alpha^{n_0,n_1,n_2,n_3,n_4})^*,C) = \hskip 80mm \]
\[ = \{(x,y,z)\in C_{n_0}\times C_{n_2}\times C_{n_4} : c^{(n_1)}(x)=d^{(n_1)}(y) \;\; {\rm and }\;\; c^{(n_3)}(y)=d^{(n_3)}(z) \} \]
which are natural in $C$.
\end{proposition}

{\it Proof.}~ As both positive opetopic cardinals $\alpha^{n_0,n_1,n_2}$ and $\alpha^{n_0,n_1,n_2,n_3,n_4}$ are obtained via special pushout (in the second case applied twice), these  colimits are preserved by $(-)^*$. $~~\Box$
\vskip 3mm

Let $T$ be a positive opetopic cardinal. We have a functor
\[ \Sigma^T : \pOpe \da T \lra \pOpeCard \]
such that
\[ \Sigma^T (f:B\ra T) = B \]
and a cocone
\[ \sigma^T : \Sigma^T\lra T\]
such that
\[ \sigma^T_{(f:B\ra T)} = f : \Sigma^T (f:B\ra T) = B\lra T. \]

We have

\begin{lemma}\label{speclim} The cocone $\sigma^T: \Sigma^T \stackrel{\cdot}{\lra} T$ is a colimiting cocone in $\pOpeCard$. Such colimiting cones are called {\em special colimits}\index{colimit!special}\index{special!colimit}. Any functor from $\pOpeCard^{op}$ preserves special limits preserves special pullbacks as well.
\end{lemma}

{\it Proof.}~ To see that $\sigma^T: \Sigma^T \stackrel{\cdot}{\lra} T$ is a colimiting cone we proceed by induction on the size of $T$. If $T$ is a positive opetope, then  the category $\pOpe \da T$ has terminal object $id_T$ which is sent by $\Sigma^T$ to $T$. Thus in this case $T$ is the colimit of $\Sigma^T$. If $T$ is not a positive opetope, then, by Lemma \ref{Sd}.5, it can be presented as a special pushout $T=T_2\oplus_kT_1$
\begin{center}
\xext=650 \yext=500
\begin{picture}(\xext,\yext)(\xoff,\yoff)
 \setsqparms[-1`-1`-1`-1;650`400]
 \putsquare(0,50)[T_2\oplus_kT_1`T_1`T_2`\bd^{(k)} T_2;\kappa_1`\kappa_2`\bc^{(k)}_{T_1}`\bd^{(k)}_{T_2}]
\end{picture}
\end{center}
with both $T_1$ and $T_2$ of dimension larger than $k$ and size smaller than the size of $T$, for some $k\in \o$. By inductive assumption, the limits of $\Sigma^T_1$, $\Sigma^T_2$, and $\Sigma^{\bd^{(k)}T_2}$ are $T_1$, $T_2$ and $\Sigma^{\bd^{(k)}T_2}$, respectively. Each object $\pOpe \da T$ factorises (as a morphism) via either $\kappa_1$ or $\kappa_2$. If it factorises by both, it factorises by  $\bd^{(k)} T_2$. From this description it is easy to see that indeed in this case $T$ is also the colimit of the functor $\Sigma^T$. Moreover, if the limit of $\Sigma^T$ is preserved, then this special pushout is also preserved.  $~~\Box$

\vskip 3mm

{\em Remarks and notation.} The full image of the functor $(-)^*:\pOpeCard\lra \oC$ will be denoted by $\pOpeCardo$\index{category!ctypo@$\pOpeCardo$}. The objects of $\pOpeCardo$ are $\o$-categories isomorphic to those of form $S^*$ for $S$ being positive opetopic cardinal and the morphism in $\pOpeCardo$ are all $\o$-functors. In fact, when convenient, we shall think about positive opetopic cardinals $S$ as if they where $\o$-categories and talk about $\o$-functors between them. As the above embedding $(-)^*$ is conservative, see Corollary \ref{fc-pres_sp oC}, this will not lead to any confusions.

\section{The inner-outer factorization in $\pOpeCardo$} \label{sec_inner-outer}

Let $f:S^*\lra T^*$ be a morphism in $\pOpeCardo$. We say that $f$ is {\em outer}\index{map!outer -}\footnote{The names `inner' and `outer' are introduced in analogy with the morphism with the same name and role in the category of disks in \cite{Joyal}.} if there is a map of positive opetopic cardinals $g:S\lra T$ such that $g^*=f$. We say that $f$ is {\em inner}\index{map!inner -} iff $f_{dim(S)}(S)=T$. From Corollary \ref{ff-pres_sp} we have

\begin{lemma}\label{outer map}
An $\o$-functor $f:S^*\lra T^*$ is outer iff it is a polygraph map.  $~~\Box$
\end{lemma}

\begin{proposition}\label{inner map}
Let $f:S^*\lra T^*$ be an inner map, $dim(S)=dim(T)>0$.  The maps $\bd f : \bd S \lra \bd T$ and $\bc f : \bc S \lra \bc   T$, being the restrictions of $f$, are well defined, inner and the squares
\begin{center} \xext=1400 \yext=550
\begin{picture}(\xext,\yext)(\xoff,\yoff)
 \setsqparms[1`-1`-1`1;700`500]
 \putsquare(0,0)[(\bd T)^*`T^*`(\bd S)^*`S^*;
 \bd_T^*`\bd f`f`\bd_S^*]
 \setsqparms[-1`0`-1`-1;700`500]
 \putsquare(700,0)[\phantom{T^*}`(\bc T)^*`\phantom{S^*}`(\bc S)^*;
 \bc_T^*``\bc f`\bc_S^*]
\end{picture}
\end{center}
  commute.
\end{proposition}
{\it Proof.}~So suppose that $f:S^*\ra T^*$ is an inner map. So $f(S)=T$. Since $f$ is an $\o$-functor, we have
\[ f(\bd S)=\bd f( S)=\bd T\;\;\;\;{\rm and} \;\;\;\; f(\bc S)=\bc f(S)= \bc T. \]
This shows the proposition.  $~~\Box$

We have

\begin{proposition}\label{factorization}
The  inner and outer morphisms form a  factorization system in $\pOpeCardo$. So any $\o$-functor $f:S^*\lra T^*$ can be factored essentially uniquely by inner map   $\stackrel{\bullet}{f}$ followed by outer map $\stackrel{\circ}{f}$:
  \begin{center}
\xext=1000 \yext=500 \adjust[`I;I`;I`;`I]
\begin{picture}(\xext,\yext)(\xoff,\yoff)
\settriparms[1`1`-1;500]
\putVtriangle(20,0)[S^*`T^*`f(S)^*;f`\stackrel{\bullet}{f}`\stackrel{\circ}{f}]
\end{picture}
\end{center}
\end{proposition}
{\it Proof.}~ This is almost tautological. $~~\Box$

The inner maps between positive opetopic cardinals can be further factorized into inner epi and inner mono.

Let $P$ and $Q$ be positive opetopic cardinals, $f:P^*\ra Q^*$ an $\o$-functor between $\o$-categories generated by them. The {\em kernel of} $f$ is the set of faces of $P$ sent by $f$ to identities on cells of $Q^*$ of lower dimension. $\ker(f)$ denotes the kernel of $f$. We say that a set  $I$ of faces of $P_{>0}$ is an {\em ideal in} $P$ iff for any $b\in P_{>0}$
\begin{description}
      \item[$1_i.$]  if    $\gamma(b)\in I$, then $b\in I$; 
   \item[$2_i.$] if   $\delta(b)\subseteq I$, then $b\in I$; 
  \item[$3_i.$] if $b\in I$, then $|\delta(b) \setminus I| =  |\gamma(b)$. 
\end{description}

We have an easy

\begin{lemma} \label{ker->ideal}
Let $f:P^*\ra Q^*$ be an $\o$-functor between $\o$-categories generated by positive opetopic cardinals. Then $\ker(f)$ is an ideal.  $~~\Box$
\end{lemma}

\vskip2mm
We shall prove the converse of the above lemma, i.e., that any ideal is a kernel of an $\o$-functor, in fact an inner epi.

Recall that a face $u\in P_{>0}$ is {\em unary} if $\delta(u)$ contains one element. Let  $U(P)$ be the set of unary faces in $P$, $I \subseteq P$ an ideal in $P$, and $I\neq \emptyset$. The face $u\in P$ is called {\em safe for} $P$ iff $u \in U(P)-\gamma(P)-\delta(U(P)$, i.e., $u$ is a unary face in $P$ that is not a codomain of any other face in $P$ and it is not in the domain of a unary face in $P$.

The following Lemma says that we can always divide any opetopic cardinal by its safe face.

\begin{lemma} \label{safe faces}
Let $P$ be an opetopic cardinal, and $u$ a safe face for $P$. Then we can divide $P$ by $u$, i.e., we have a quotient $\o$-functor $q_u : P^* \lra P_{/u}^*$ whose kernel is $\{ u\}$. $q_u$ is an inner epi. \end{lemma}

{\em Proof.} The opetope $P_{/u}$ is obtained by gluing together $\delta(u)$ and $\gamma(u)$ (to a face $\{\delta(u), \gamma(u)\}$) and dropping $u$. The map $q_u$ is defined as follows (we describe it on faces of $P$ only)

\[ q_u(a)=
\left\{
  \begin{array}{ll}
    \{\delta(u), \gamma(u)\}, & \hbox{if $a = \delta(u),\; \gamma(u), \; u$;} \\
    a, & \hbox{otherwise.}
  \end{array}
\right.
\]
i.e., $q_u(\delta(u))=q_u(\gamma(u))= q_u(u) = \{\delta(u), \gamma(u)\}$, i.e., the equivalence class containing $\delta(u)$ and $\gamma(u)$.
$q_u(a)=a$, for other faces $a$ in $P$.

Then $\gamma$ and $\delta$ on $P_{/u}$ are so defined to make the quotient map $q_u: P*\lra P/u^*$ preserve both of them.

The only cell sent by $q_u$ to a(n identity on a) cell of a lower dimension is $u$. Thus  $\ker(q_u)=\{ u\}$. $~~\Box$
\vskip2mm

The following two lemmas show that in any non-empty ideal $I$ in $P$ there is always a safe face for $P$.

\begin{lemma}\label{safe faces exist 1} Let $I$ be a non-empty ideal in $P$. There is always a unary face $u\in I -\gamma(P)$. \end{lemma}

Proof. Suppose not, and let c be a cell of minimal dimension and $<^+$- minimal in $I$. If $c$ is not unary then this contradicts condition 3.
as c is of minimal dimension in I and hence $\delta(c)\cap I =\emptyset$.
If $c\in \gamma(P)$, then this contradicts the choice of $c$ as, if $\gamma(b)=c$, then $\delta(b)$ contains only unary faces, and as $c$ is in $I$, by 3., $\delta(b)\subseteq I$.
Thus $c \in (I \cap U(P))-\gamma(P)$. $~~\Box$

\begin{lemma}\label{safe faces exist 2}
Let $I$ be a non-empty ideal in $P$. There is always a face $u\in I$ safe for $P$.
\end{lemma}

{\em Proof.}  Take the unary face $u$ of maximal dimension in $I$. If $u \in \delta(v)$ with $v \in U(P)$, then, since $I$ satisfies 2., $v\in I\cap U(P)$. This contradicts the choice of $u$, as $dim(u)<dim(v)$. $~~\Box$

\vskip2mm
\begin{theorem} \label{ideal->ker}
 If  $I\subseteq P_{\geq 1}$ is an {\em ideal}, then there is a $\iota$-epi map $q_I:P\ra P_{/I}$ such that it has $I$ as its kernel. Moreover, $q_I$ is a universal map with this property, i.e., whenever there is an $\o$-functor $h: P^*\ra Q^*$ such that $I\subseteq \ker(f)$, then there is a unique map $f': P_{/I}*\ra Q$ such that $f=f'\circ q_I$.
\end{theorem}

{\em Proof.}
We can divide the opetope $P$ by unary cell $u\in I-\gamma(I)$ of maximal dimension, getting the map
\[ q_u: P\ra P_{/u}. \]
Then $q_u(I-\{u\})$, the image of $I-\{u\}$ in $P_{/u}$, is an ideal in $P_{/u}$. Thus we can iterate the construction until the resulting ideal will be empty.

To see that $q_I$ has the stated universal property, it is enough to notice that if $u\in \ker(f)$ is a safe face for $P$, then we have a factorization
\begin{center}
\xext=600 \yext=400 \adjust[`I;I`;I`;`I]
\begin{picture}(\xext,\yext)(\xoff,\yoff)
\settriparms[1`1`-1;400]
  \putVtriangle(0,50)[P`Q`P_{/u};f`q_u`f']
\end{picture}
\end{center}
Using the description of $q_u$, this is clear. $~~\Box$

\vskip2mm
Note that for an $\o$-functor $f:P^*\ra Q^*$ it might seem that to say `that it is mono (or epi)' is ambiguous, since this can be applied to either just faces of $P$ or all the cells of $P^*$. However, in both cases the notions of epi and of mono coincide. Thus, in fact, they are not ambiguous no matter how these notions are interpreted. We have

\begin{lemma}\label{mono <-> ker=0} Let $f:P^*\ra Q^*$ an $\o$-functor in $pOpeCard_\o$. Then $\ker(f)=\emptyset$ iff $f$ is mono.  \end{lemma}

Proof. Suppose that $f$ is not a mono. Let $a,b\in P_m$, $m\in\o$, be two different cells of minimal dimension such that $f(a)=f(b)$. If $m=0$, then $a\perp^+b$ by linearity of $<^+$ on $P_0$, and if  $m>0$, then since
\[  f(\gamma(a)) = \gamma(f(a)) = \gamma(f(b)) = f(\gamma(b)),\]
and by minimality of $m$, we have that  $\gamma(a)) =\gamma(b)$. Then by pencil linearity we have $a\perp^+b$, as well. Suppose $a<^+b$. Thus there is an upper path $a,\alpha_1,\ldots, \alpha_k,b$ in $P$. As $f(a)=f(b)$, we have $f(\alpha_i)=f(a)$, and hence, $\alpha\in \ker(f)$ for $i=,1,\ldots ,k$, i.e., $\ker(f)\neq\emptyset$.  $~~\Box$

\begin{theorem} \label{inner epi-mono}
The inner epis and inner monos form a factorization system on the category $\pOpeCard_{inn}$ of opetopic cardinals with inner maps.
\end{theorem}

{\em Proof.} Let $f: P^*\ra Q^*$ be an inner map. Then, by Lemma \ref{ker->ideal}, $\ker(f)$ is an ideal. By the universal property of $q_I$ stated in Lemma \ref{ideal->ker} we have a factorization
\begin{center}
\xext=600 \yext=400 \adjust[`I;I`;I`;`I]
\begin{picture}(\xext,\yext)(\xoff,\yoff)
\settriparms[1`1`-1;400]
  \putVtriangle(0,50)[P`Q`P_{/u};f`q_u`f']
\end{picture}
\end{center}
with $q_I$ inner epi. Since $\ker(f)=I=\ker(q_I)$, it follows that $\ker(f')=\emptyset$. $~~\Box$

\section{The terminal positive-to-one polygraph}\label{sec-terminal-polygraph}

In this section we shall describe the terminal positive-to-one polygraph $\cT$ as an $\o$-category.

The set of $n$-cell $\cT_n$ consists of (isomorphisms classes of) positive opetopic cardinals of dimension less than or equal to $n$. For $n>0$, the operations of domain and codomain $d^\cT,c^\cT:\cT_n\ra \cT_{n-1}$ are given, for $S\in\cT_n$, by

\[  d(S) = \left\{ \begin{array}{ll}
        S   & \mbox{if $dim(S)<n$,}  \\
        \bd S & \mbox{if $dim(S)=n$,}
                                    \end{array}
                \right. \]
and
\[  c(S) = \left\{ \begin{array}{ll}
        S   & \mbox{if $dim(S)<n$,}  \\
        \bc S & \mbox{if $dim(S)=n$.}
                                    \end{array}
                \right. \]
and, for $S,S'\in\cT_n$ such that $c^{(k)}(S)=d^{(k)}(S')$, the composition in $\cT$ is just the special pushout $$S'\circ_kS=S'\oplus_k,S$$
i.e.,
\begin{center} \xext=650 \yext=550
\begin{picture}(\xext,\yext)(\xoff,\yoff)
 \setsqparms[-1`-1`-1`-1;650`500]
 \putsquare(0,50)[S'\oplus_kS`S`S'`\bc^{(k)} S;``\bc^{(k)}`\bd^{(k)}]
\end{picture}
\end{center}
The identity $id_\cT : T_{n-1}\ra \cT_n$ is the inclusion map. The $n$-indeterminates in $\cT$ are positive opetopic cardinals of dimension $n$.

\begin{proposition} $\cT$ just described is the terminal positive-to-one polygraph.
\end{proposition}

{\it Proof.}~ The fact that $\cT$ is an $\o$-category is easy. The fact that $\cT$ is free with free $n$-indeterminates being positive opetopic cardinals of dimension $n$ can be shown much like the freeness of $S^*$ before. The fact that $\cT$ is terminal follows from the following observation.

{\em Observation.} For every pair of parallel positive opetopic cardinals of dimension $n$, $N$ and $B$ (i.e., $\bd N = \bd B$ and $\bc N = \bc B$) such that $B$ is principal,
it follows that $N$ is normal and there is a unique (up to an iso) principal  positive opetopic cardinal  $N^\bullet$ of dimension $n+1$ such that $\bd N^\bullet =N$ and $\bc N^\bullet =B$. $~~\Box$

\begin{lemma}\label{shape_pres} Let $S$ be a positive opetopic cardinal and $!:S^*\lra \cT$ the unique map from $S^*$ to $\cT$. Then, for $T\in S^*_k$, we have
\[ !_k(T)=T.\]
\end{lemma}

{\it Proof.}~ The proof is by induction on $k\in\o$ and the size of $T$ in $S^*_k$. For $k=0,1$, the Lemma is obvious. Let $k>1$ and assume that lemma holds for $i<k$.

If $dim(T)=l<k$, then, using the inductive hypothesis and the fact that $!$ is an $\o$-functor, we have
\[ !_k(T)=!_k(1^{(k)}_T) = 1^{(k)}_{!_l(T)} = 1^{(k)}_{T} = T \]

Suppose that $dim(T)=k$ and $T$ is principal. As $!$ is a polygraph map, $!_k(T)$ is an indeterminate, i.e., it is principal, as well. Using  again the inductive hypothesis and the fact that $!$ is an $\o$-functor, we obtain
\[ d(!_k(T))=!_{k-1}(\bd T)= \bd T \]
\[ c(!_k(T))=!_{k-1}(\bc T)= \bc T \]
As $T$ is the only (up to a unique iso) positive opetopic cardinal with the domain $\bd T$ and the codomain $\bc T$, it follows that $!_k(T)=T$, as required.

Finally, suppose that $dim(T)=k$, $T$ is not principal, and fo, the positive opetopic cardinals of size smaller than the size of $T$, the Lemma holds. Thus there are $l\in\o$ and $a\in Sd(T)_l$ so that
\[ !_k(T)= !_k(T^{\da a}\oplus_lT^{\ua a}) = !_k(T^{\da a})\oplus_l!_k(T^{\ua a})= T^{\da a}\oplus_lT^{\ua a} = T, \]
as required  $~~\Box$

\section{A description of the positive-to-one polygraphs}\label{sec-description-of-polygraphs}

In this section we shall describe all the cells in positive-to-one polygraphs using positive opetopic cardinals, in other words we shall describe in concrete terms the functor: \[ \overline{(-)}^n: Set\da D_n \lra \pnPoly \]

More precisely, the positive-to-one polygraphs of dimension 1 (and all polygraphs, as well) are free polygraphs over graphs and are well understood. So suppose that $n>1$, and we are given an object of $Set\da D_n$, i.e., a quadruple $(|P|_n,P,d,c)$ such that
\begin{enumerate}
  \item a positive-to-one $(n-1)$-polygraph $P$;
  \item a set $|P|_n$ with two functions $c: |P|_n\lra |P|_{n-1}$
  and $d: |P|_n\lra P_{n-1}$ such that, for $x\in |P|_n$, $cc(x)=cd(x)$ and
  $dc(x)=dd(x)$; we assume that $d(x)$ is not an identity, for any
  $x\in |P|_n$.
\end{enumerate}
If the maps $d$ and $c$ in the object $(|P|_n,P,d,c)$ are understood from the context, we can abbreviate notation to $(|P|_n,P)$.

Recall that for a positive opetopic cardinal $S$, with $dim(S)\leq n$, we denote by $S^{\sharp,n}$ the object  $(S_n,(S_{<n})^*,[\delta],[\gamma])$ in $Set\da D_n$. In fact, we have an obvious functor
\[ (-)^{\sharp,n}: \pOpeCard \lra  Set\da D_n\]
such that
 \[ S\mapsto (S_n,(S_{<n})^*,[\delta],[\gamma]) \]
Any positive-to-one polygraph $P$ can be restricted to its part in $Set\da D_n$.  So we have a forgetful functor
\[ U_n :\pPoly_n \lra Set\da D_n \]
such that
\[ P \mapsto (|P|_n,P_{<n},d,c), \]
Thus the above functor $U_n$ is the same functor as the one with the same name defined in the Appendix, but restricted to the category $\pPoly_n$. It is an essential inverse of
$\overline{(-)}^n$.

We shall describe the positive-to-one $n$-polygraph $\overline{P}(=\overline{P}^n)$ whose $(n-1)$-truncation is $P_{<n}$ and whose $n$-indeterminates are $|P|_n$ with the domains and codomains given by maps $c$ and $d$.

{\bf n-cells} of $\overline{P}$. An $n$-cell in $\overline{P}_n$ is a(n equivalence class of) pair(s)  $(S,f)$ where
\begin{enumerate}
  \item $S$ is a positive opetopic cardinal, $dim(S)\leq n$;
  \item $f: (S_n,(S_{<n})^*,[\delta],[\gamma])\lra (|P|_n,P,d,c)$ is a morphism in
  $Set\da D_n$, i.e.,
   \begin{center}
\xext=650 \yext=600
\begin{picture}(\xext,\yext)(\xoff,\yoff)
 \setsqparms[1`0`0`1;650`500]
 \putsquare(0,50)[S_n`|P|_n`S^*_{n-1}`P_{n-1};|f|_n```f_{n-1}]
  \putmorphism(-50,550)(0,-1)[``{[}\delta {]}]{500}{1}l
  \putmorphism(50,550)(0,-1)[``{[}\gamma {]}]{500}{1}r
  \putmorphism(600,550)(0,-1)[``d]{500}{1}l
  \putmorphism(700,550)(0,-1)[``c]{500}{1}r
\end{picture}
\end{center}
commutes.
\end{enumerate}

We identify two pairs $(S,f)$, $(S',f')$ if there is an isomorphism $h:S\lra S'$ such that the triangles of sets and of $(n-1)$-polygraphs
 \begin{center}
\xext=2200 \yext=500
\begin{picture}(\xext,\yext)(\xoff,\yoff)
 \settriparms[1`1`1;400]
 \putVtriangle(0,0)[S_n`S'_n`|P|_n;h_n`f_n`f'_n]
 \putVtriangle(1400,0)[(S_{< n})^*`(S'_{<n})^*`P;(h_{<n})^*`f_{<n}`f'_{<n}]
\end{picture}
\end{center}
commute. Clearly, such an $h$, if exists, is unique. Even if formally cells in $P_n$ are equivalence classes of triples, we will work on triples themselves as if they were cells understanding that equality between such cells is an isomorphism in the sense defined above.

{\bf Domains and codomains}. The domain and codomain functions
\[ d^{(k)},c^{(k)}: \overline{P}_n \lra \overline{P}_{k} \]
are defined for an $n$-cell $(S,f)$ as follows:
\[ d^{(k)}(S,f)=(\bd^{(k)}S,\bd^{(k)}f) \]
 where, for $x\in (\bd^{(k)}S)_k$
\[ (\bd^{(k)}f)_k(x)= f_k([x])(x) \]
(i.e., we take the  positive opetopic cardinals $[x]$ contained in $S$, then the value of $f$ on it, and then we evaluate the map in $Set\da D_n$ on $x$, the
only element of  $[x]_k$),
\[ (\bd^{(k)}f)_l= f_l \]
for $l<k$;

\[ c^{(k)}(S,f)=(\bc^{(k)}S,\bc^{(k)}f) \]
 where, for $x\in (\bc^{(k)}S)_k$
\[ (\bc^{(k)}f)_k(x)= f_k([x])(x) \]
and
\[ (\bd^{(k)}f)_l= f_l \]
for $l<k$, i.e., we calculate the $k$-th domain and $k$-th codomain of an $n$-cell $(S, f)$ by taking $\bd^{(k)}$ and $\bc^{(k)}$ of the domain $S$ of the cell $f$, respectively, and by restricting the maps $f$ accordingly.

{\bf Identities}. The identity function
\[  \bi : \overline{P}_{n-1}\lra \overline{P}_n \]
is defined, for an $(n-1)$-cell $((S,f)$ in $P_{n-1}$, as follows:
\[ \bi(S, f) = \left\{ \begin{array}{ll}
                            (S,f) & \mbox{ if $dim(S)<n-1$,} \\
                            (S,\overline{f})  & \mbox{ if dim(S)=n-1}
                           \end{array}
                   \right. \]
Note that $\overline{f}$ is the map $\pnmjPoly$ which is the value of the functor $\overline{(-)}$ on a map $f$ from $Set\da D_{n+1}$. So it is in fact defined as `the same $(n-1)$-cell' but considered as an $n$-cell.

{\bf Compositions}. Suppose that $(S^i,f^i)$ are $n$-cells for $i=0,1$ such that
\[ c^{(k)}(S^0, f^0)= d^{(k)}(S^1, f^1). \]
Then their composition is define, via pushout in $Set\da D_n$, as
\[ (S^1, f^1)\circ_k(S^0, f^0)= (S^1\oplus_kS^0,[f^1,f^0]) \]
i.e.,
\begin{center}
\xext=1500 \yext=700
\begin{picture}(\xext,\yext)(\xoff,\yoff)
 \setsqparms[1`0`0`1;1450`500]
 \putsquare(0,100)[S^0_n\dsum S^1_n`|P|_n`
 ((S^0\oplus_kS^1)_{\leq n-1})^*_{n-1}`P_{n-1};
 [f^0_n,f^1_n]```{[f^0_{n-1},f^1_{n-1}]}]
  \putmorphism(-50,600)(0,-1)[``{[}\delta {]}]{500}{1}l
  \putmorphism(50,600)(0,-1)[``{[}\gamma {]}]{500}{1}r
  \putmorphism(1400,600)(0,-1)[``d]{500}{1}l
  \putmorphism(1500,600)(0,-1)[``c]{500}{1}r
\end{picture}
\end{center}
This ends the description of the polygraph $\overline{P}$.

Now let $h:P\ra Q$ be a morphism in $Set\da D_n$, i.e., a function $h_n:|P|_n\lra |Q|_n$ and a $(n-1)$-polygraph morphism $h_{<n}: P_{<n}\lra Q_{<n}$ such that the square
\begin{center} \xext=1500 \yext=700
\begin{picture}(\xext,\yext)(\xoff,\yoff)
 \setsqparms[1`0`0`1;950`500]
 \putsquare(0,100)[|P|_n`|Q|_n`P_{n-1}`Q_{n-1};
 h_n```h_{n-1}]
  \putmorphism(-50,600)(0,-1)[``d]{500}{1}l
  \putmorphism(50,600)(0,-1)[``c]{500}{1}r
  \putmorphism(900,600)(0,-1)[``d]{500}{1}l
  \putmorphism(1000,600)(0,-1)[``c]{500}{1}r
\end{picture}
\end{center}
commutes serially. We define
\[ \bar{h}: \bar{P}\lra \bar{Q} \]
by putting $\bar{h}_k=h_k$, for $k<n$, and, for $(S,f)\in \bar{P}_n$, we put
\[ \bar{h}(S,f)= (S,h\circ f).\]

{\em Notation.} Let $x=(S,f)$ be a cell in $\bar{P}_n$ as above, and $a\in Sd(S)$. Then by  $x^{\da a}=(S^{\da a},f^{\da a})$ and $x^{\ua a}=(S^{\ua a},f^{\ua a})$  we denote the cells in $\overline{P}_n$ that are the obvious restriction  $x$.  Clearly, we have $c^{(k)}(x^{\da a})=d^{(k)}((x^{\ua a})$ and $x=x^{\ua a}\circ_kx^{\da a}$, where $k=dim(a)$.
\vskip 2mm

The following Proposition we collect several statements concerning the above construction.  This includes that the above construction is correct.  We have put all these statement together as we need to prove them together, that is, by simultaneous induction.

\begin{proposition}\label{type} Let $n\in\o$. We have
\begin{enumerate}
  \item Let $P$ be an object of $Set\da D_n$. We define the function \[ \eta_P: |P|_n\lra \overline{P}_n\] as follows.   Let $x\in |P|_n$. As $c(x)$ is an indeterminate, $d(x)$ is a normal cell of dimension $n-1$.  Thus there is a unique description of the cell $d(x)$
  \[ <T_{d(x)}, \tau_{d(x)}: T_{d(x)}^*\lra P_{<n}> \]
   with $T_{d(x)}$ being normal positive opetopic cardinal. Then we have a unique $n$-cell in $\overline{P}$:
   \[ \bar{x} = <T_{d(x)}^{\bullet},\;\; |\overline{\tau}_x|_n :
   \{ T_{d(x)}^{\bullet} \} \ra |P|_n,\;\;
    (\overline{\tau}_x)_{<n} : (T_{d(x)}^{\bullet})_{<n}^*\ra P_{<n}> \]
   (note: $|T_{d(x)}^{\bullet}|_n=\{ T_{d(x)}^{\bullet} \}$) such that
   \[ |\overline{\tau}_x|_n(T_{d(x)}^{\bullet})=x\]
   and
   \[  (\overline{\tau}_x)_{n-1}(S) = \left\{ \begin{array}{ll}
                            c(x) & \mbox{ if $S=\bc(T_{d(x)}^{\bullet})$} \\
                            (\tau_{dx})_{n-1}(S)  & \mbox{ if $S\subseteq T_{dx}$}
                           \end{array}
                   \right. \]
  and $(\overline{\tau}_x)_{<(n-1)}=(\tau_{dx})_{<(n-1)}$. We put $\eta_P(x)=\bar{x}$.

  Then $\overline{P}$ is a positive-to-one polygraph with $\eta_P$ the inclusion of $n$-indeterminates. Then any positive-to-one $n$-polygraph $Q$ is equivalent to a polygraph $\overline{P}$, for some $P$ in $Set\da D_n$.

  \item  Let $P$ be an object of $Set\da D_n$, $!:\overline{P}\lra \cT$ the unique morphism into the terminal object $\cT$ and $f:S^{\sharp,n}\ra P$ a cell in $\overline{P}_n$. Then
  \[ !_n(f:S^{\sharp,n}\ra P  )=S. \]

  \item Let $h:P\ra Q$ be an object of $Set\da D_n$. Then $\bar{h}: \overline{P}\lra \overline{Q}$ is a polygraph morphism,

  \item Let $S$ be a positive opetopic cardinal of dimension at most $n$. Moreover, for a morphism $f:S^{\sharp,n}\lra P$ in $Set\da D_n$, we have that
  \[ \overline{f}_k(T)=f\circ (i_T)^{\sharp,n} \]
  where $k\leq n$, $T\in S^*_k$ and $i_T:T \lra S$ is the inclusion.

 \item Let $S$ be a positive opetopic cardinal of dimension $n$, $P$ positive-to-one polygraph, $g,h:S^*\lra P$ polygraph maps. Then
 \[ g=h \hskip 10mm {\rm iff} \hskip 10mm g_n(S)=h_n(S). \]

 \item Let $S$ be a positive opetopic cardinal of dimension at most $n$, $P$ be an object in $Set\da D_n$. Then we have a bijective correspondence
  $$
\begin{array}{c}
f:S^{\sharp,n}\lra P \;\;\in Set\da D_n
\\ \hline
\overline{f}:S^*\lra \overline{P} \;\;\in \pnPoly
\end{array}
$$
such that $\overline{f}_n(S)=f$, and, for $g:S^*\lra \overline{P}$, we have $g=\overline{g_n(S)}$.
  \item We have a bijection
  \[ \kappa^P_n:\coprod_{S} \pPoly (S^*,\overline{P}) \lra \overline{P}_n \]
  \[ g:S^*\ra \overline{P}\;\;\; \mapsto\;\;\; g_n(S)\]
  where coproduct is taken over all (up to iso) positive opetopic cardinals $S$ of dimension at most $n$. In other words, any cell in $\bar{P}$ has a unique description.

\end{enumerate}
\end{proposition}

{\it Proof.} ~Ad 1. We have to verify that $\overline{P}$ satisfies the laws of $\o$-categories and that it is free in the appropriate sense.

The laws for $\o$-categories are left for the reader, as they easily follow from the fact that $S^*$ is a positive to one-polygraph for any positive opetopic cardinal $S$. We shall show that $\overline{P}$ is free in the appropriate sense.

Let $C$ be an $\o$-category, $g_{<n}:P_{<n}\ra C_{<n}$ and $(n-1)$-functor and $g_n:|P|_n\ra C_n$ a function so that the diagram
\begin{center} \xext=1500 \yext=700
\begin{picture}(\xext,\yext)(\xoff,\yoff)
 \setsqparms[1`0`0`1;950`500]
 \putsquare(0,100)[|P|_n`C_n`P_{n-1}`C_{n-1};
 g_n```g_{n-1}]
  \putmorphism(-50,600)(0,-1)[``d]{500}{1}l
  \putmorphism(50,600)(0,-1)[``c]{500}{1}r
  \putmorphism(900,600)(0,-1)[``d]{500}{1}l
  \putmorphism(1000,600)(0,-1)[``c]{500}{1}r
\end{picture}
\end{center}
commutes serially. We shall define an $n$-functor $\overline{g}:\overline{P}\ra C$ extending $g_{<n}$ and $g_n$. For $x=(S,f)\in \overline{P}_n$, we put
\[  \overline{g}_n(x)= \left\{ \begin{array}{ll}
               1_{g_{n-1}\circ f_{n-1}(S)} & \mbox{if $dim(S)<n$,} \\
               g_{n}\circ f_{n}(m_{S}) & \mbox{if $dim(S)=n$, $S$ is principal, $S_n=\{m_S\}$ } \\
                \overline{g}_n(x^{\ua a})\circ_k\overline{g}_n(x^{\da a}) & \mbox{if $dim(S)=n$, $a\in Sd(S)_k$}
                           \end{array}
                   \right. \]

We need to check that $\overline{g}$ is well defined, it unique that extends $g$, and it preserves domains, codomains, compositions and identities.

All these calculations are similar, and they are very much like those in the proof of Proposition \ref{Sstar} and use facts from Section \ref{sec-decomposition}. We shall check, assuming that we already know that  $\overline{g}$ is well defined and preserves identities, that compositions are preserved. The proof is by induction on the size of the composition and uses Lemma \ref{decomp3new}. So let $T$, $T_1$, $T_2$ be positive opetopic cardinals such that $T=T_2\oplus_kT_1$. Since $\overline{g}$ preserves identities, we can restrict to the case $dim(T_1),dim(T_2)>k$.

Fix $a\in\bc^{(k)}(T_1)_k\cap\gamma(T_1)$. So $a\in Sd(T)_k$. If $T_1=T^{\da a}$ and $T_2=T^{\ua a}$, then we have
\[ \overline{g}(T)=\overline{g}(T^{\ua a})\circ_k\overline{g}(T^{\da a})= \overline{g}(T_2)\circ_k\overline{g}(T_1). \]
If $a\in Sd(T_1)_k$, then
\begin{eqnarray*}
\overline{g}(T)=\overline{g}(T^{\ua a})\circ_k\overline{g}(T^{\da a})= \\
 =(\overline{g}(T_2)\circ_k\overline{g}(T_1^{\ua a})) \circ_k\overline{g}(T^{\da a})=   \\
 =\overline{g}(T_2) \circ_k(\overline{g}(T_1^{\ua a}) \circ_k\overline{g}(T_1^{\da a}))=  \\
 \overline{g}(T_2)\circ_k\overline{g}(T_1)
\end{eqnarray*}
The remaining verifications are similar.

 Ad 2. Let $!:\overline{P}\lra \cT$ be the unique polygraph map into the terminal object, $S$ a positive opetopic cardinal such that $dim(S)=l\leq n$, $f:S^{\sharp,n}\lra P$ a cell in $\overline{P}_n$.

If $l<n$, then by induction we have $!_n(f)=S$. If $l=n$ and $S$ is principal, then we have, by induction
 \[ !_n(d(f):(\bd S)^{\sharp,n}\ra P  )=\bd S, \hskip 10mm
 !_n(c(f):(\bc S)^{\sharp,n}\ra P  )=\bc S. \]
As $f$ is an indeterminates in $\overline{P}$, $!_n(f)$ is a positive opetope. But the only (up to an iso) positive opetope $B$ such that
\[ \bd B =\bd S, \hskip 10mm \bd B =\bd S \]
is $S$ itself. Thus, in this case, $!_n(f)=S$.

Now assume that $l=n$, and $S$ is not principal, and that for positive opetopic cardinals $T$ of smaller size than $S$ the
statement holds. Let $a\in Sd(S)_k$. We have
\[ !_n(f)= f^{\ua a})\circ_k !_n(f^{\da a} =   !_n (f^{\ua a}) \oplus_k!_n(f^{\da a}) =  S^{\ua a} \oplus_k S^{\da a} = S \]
 where $f^{\da a}=f\circ (\kappa^{\da a})^{\sharp,n}$ and $f^{\ua a}=f\circ (\kappa^{\ua a})^{\sharp,n}$ and $\kappa^{\da a}$ and $\kappa^{\ua
a}$ are the maps, as in the following pushout
\begin{center}
\xext=650 \yext=400
\begin{picture}(\xext,\yext)(\xoff,\yoff)
 \setsqparms[-1`-1`-1`-1;550`400]
 \putsquare(0,0)[S`S^{\da a}`S^{\ua a}`\bc^{(k)} S;\kappa^{\da a}`\kappa^{\ua a}``]
\end{picture}
\end{center}

Ad 3. The main thing is to show that $\overline{h}$ preserves compositions. This follows from the fact that the functor
\[ (-)^{\sharp,n} :\pOpeSet_n \lra Set\da D_n \]
preserves special pullbacks.

Ad 4. This is an immediate consequence of 3.

Ad 5. Let $S$ be a positive opetopic cardinal $S$ of dimension at most $n$. To prove 5., we are going to use the description of the $n$-cells in positive-to-one polygraphs given in 1. Moreover, note that by 3. and Lemma \ref{shape_pres} we have that for $T\in S^*_k$, the value of $g$ at $T$ is a map in $Set\da D_k$ such that $g_k(T):T^{\sharp,k}\lra U_k(P)$, i.e., the domain of $g_k(T)$ is necessarily $T^{\sharp,k}$.

The implication $\Ra$ is obvious. So assume that $g,h: S^* \lra P$ are different polygraph maps. Then there is $k\leq n$ and $x\in S_k$ such that $g_k([x])\neq h_k([x])$.  We shall show, by induction on size of $T$,  that for any $T\in S^*_l$ such that $x\in T$ , we have
\begin{equation}\label{eq gh}  g_k(T)\neq h_k(T)\end{equation}

$T=[x]$ has the least size among those positive opetopic cardinals that contain $x$. Clearly, (\ref{eq gh}) holds in this case by assumption.

Suppose that (\ref{eq gh}) holds for all $U\in S^*_{l'}$ whenever $l'<l$ and $x\in U$. Suppose that $T=[y]$, for some  $y\in S_l$, and $x\in [y]$. Then either $x\in \bd [y]$ or $x\in \bc [y]$. In the former case we have, by inductive hypothesis, that $g_k(\bd T)\neq h_k(\bd T)$.Thus

\[   d (g_k(T))= g_k(\bd T)\neq h_k(\bd T) = d (g_k(T)) \]
But then (\ref{eq gh}) holds as well. The later case ($x\in \bc [y]$) is similar.

Now suppose that $T$ is not principal $x\in T$ and that, for $U$ of a smaller size with $x\in U$, the condition (\ref{eq gh}) holds. Let $a\in Sd(T)_r$. Then either $x\in T^{\da a}$ or $x\in T^{\ua a}$. Both cases are similar, so we will consider the first one only. Thus, as $T^{\da a}$ has a smaller size than $T$, by inductive hypothesis we have
\begin{equation}\label{eq gha} g_k(T^{\da a})\neq h_k(T^{\da a})\end{equation}
As the compositions in $P$ are calculated via pushouts, we have that
\[ g_l(T^{\ua a})\circ_r g_l(T^{\da a}) = [g_l(T^{\ua a}),g_l(T^{\da a})] \]
 where $[g_l(T^{\da a}), g_l(T^{\ua a})]$ is the unique morphism  from the pushout as in the following diagram:
\begin{center} \xext=1500 \yext=1300
\begin{picture}(\xext,\yext)(\xoff,\yoff)
 \setsqparms[1`-1`-1`1;800`500]
 \putsquare(0,100)[(T^{\ua a})^{\sharp,l}`T^{\sharp,l}`
 (\bc T^{\ua a})^{\sharp,l}`(T^{\da a})^{\sharp,l}; ```]
 \putmorphism(950,700)(1,1)[``{[ g_l(T^{\ua a}),g_l(T^{\da a})]}]{300}{1}l
 \put(1350,1100){\makebox(100,100){$U_n(P)$}}
 \putmorphism(1050,200)(1,3)[``g_l(T^{\ua a})]{250}{1}r
 \put(5,700){\line(1,3){60}}
 \putmorphism(350,950)(4,1)[``g_l(T^{\da a})]{450}{1}l
\end{picture}
\end{center}
 Similarly
 \[  h_l(T^{\ua a})\circ_r h_l(T^{\da a}) = [ h_l(T^{\ua a}),h_l(T^{\da a})] \]
 As morphisms from the pushout are equal if and only if their both  components are equal we have
 \[ g_l(T)= g_l(T^{\ua a}\oplus_r T^{\da a})= g_l(T^{\ua a}) \circ_r g_l(T^{\da a}) = \]
 \[  [g_l(T^{\ua a}),g_l(T^{\da a})] \neq [h_l(T^{\ua a}),h_l(T^{\da a})]=\]
 \[ = h_l(T^{\ua a}) \circ_r  h_l(T^{\da a})= h_l(T^{\ua a} \oplus_r T^{\da a}) =h_l(T) \]

Thus (\ref{eq gh}) holds for all $T\in S^*$ such that $x\in T$. As $x\in S$, we get that
\[ g_n(S)\neq h_n(S), \]
as required.

 Ad 6. Fix a positive opetopic cardinal $S$ of dimension at most $n$.

Let $f:S^{\sharp,n}\lra P$ be a cell in $\overline{P}_n$. By 4, we have
\[ \overline{f}_n(S)=f\circ (i_S)^{\sharp,n}=f\circ (1_S)^{\sharp,n}=f\circ (1_S^{\sharp,n}) = f. \]

Let $g:S^*\lra \overline{P}$ be a polygraph map. To show that $g=\overline{g_n(S)}$, by 5, it is enough to show that
\[ (\overline{g_n(S)})_n(S) = g_n(S).\]
Using 4 again, we have,
\[ (\overline{g_n(S)})_n(S) = g_n(S)\circ (i_S)^{\sharp,n} = \]
\[ =g_n(S)\circ i_{S^{\sharp,n}}= g_n(S)\circ 1_{S^{\sharp,n}} =
g_n(S).\]
 Thus, by 5,  $(\overline{g_n(S)})=g$.

 Ad 7. It follows immediately from 6. $~~\Box$

From the Proposition \ref{type}.7 we know that each cell in a positive-to-one polygraph has (up to an isomorphism) a unique description. The following Proposition is a bit more specific.

\begin{proposition}\label{type1}
Let $P$ be a positive-to-one polygraph, $n\in\o$, and $a\in P_n$. Let $T_a$ be $!^P_n(a)$ (where $!^P:P\lra \cT$ is the unique morphism into the terminal polygraph). Then there is a unique polygraph map $\tau_a :T_a^* \lra P$ such that $(\tau_a)_n(T_a)=a$, i.e., each cell has an essentially unique description. Moreover, we have:
\begin{enumerate}
  \item for any $a\in P$, we have
  \[ \tau_{da}= d(\tau_a) = \tau_{da}= \tau_a\circ (\bd_{T_a})^*,\hskip 10mm   \tau_{c(a)}= c(\tau_a) = \tau_{c(a)}= \tau_a\circ (\bc_{T_a})^*, \]
  \[ \tau_{1_a} = \tau_a \]
  \item for any $a,b\in P$ such that $c^{(k)}(a)=d^{(k)}(b)$, we have
  \[ \tau_{a;_kb} = [\tau_a,\tau_b]: T_a^*+_{\bc^{(k)}T_a^*}T_b^* \lra P,\]
  \item for any positive opetopic cardinal $S$, for any polygraph map $f:S^*\lra P$,
  \[ \overline{\tau_{f_n(S)}}=f. \]
  \item for any positive opetopic cardinal $S$, any $\o$-functor $f:S^*\lra P$   can be essentially uniquely factorized as
  \begin{center}
\xext=600 \yext=300 \adjust[`I;I`;I`;`I]
\begin{picture}(\xext,\yext)(\xoff,\yoff)
 \settriparms[1`1`-1;300]
 \putVtriangle(0,0)[S^*`P`T_{f(S)}^*;f`f^{in}`\tau_{f(S)}]
\end{picture}
\end{center}
  where $f^{in}$ is an inner map and $(\tau_{f(S)},T_{f(S)})$ is the   description of the cell $f(S)$.
  \end{enumerate}
\end{proposition}

{\it Proof.}~ Using the above description of the positive-to-one polygraph $P$ we have that $a:(T_a)^{\sharp,n}\lra U_n(P)$. We put $\tau_a=\overline{a}$. By Proposition \ref{type} point 6, we have that $(\tau_a)_n(T_a)=\overline{a}_n(T_a)=a$, as required.

The uniqueness of $(T_a,\tau_a)$ follows from Proposition \ref{type} point 5.

The remaining part is left for the reader. $~~\Box$

\section{Positive-to-one polygraphs form a presheaf category}\label{sec_pPoly}

In this section we want to prove that the category $\pPoly$ is equivalent to the presheaf category $Set^{(\pOpe)^{op}}$. In fact, we will show that both categories are equivalent to the category $sPb((\pOpeCard)^{op},Set)$ of special pullbacks preserving functors from $(\pOpeCard)^{op}$ to $Set$.

First note that the inclusion functor $\bi : \pOpe \lra \pOpeCard$ induces the adjunction
\begin{center} \xext=1500 \yext=250
\begin{picture}(\xext,\yext)(\xoff,\yoff)
\putmorphism(0,150)(1,0)[\phantom{Set^{{(\pOpe)}^{op}}}`
\phantom{Set^{{(\pOpeCard)}^{op}}}`Ran_\bi]{1500}{1}a
\putmorphism(0,50)(1,0)[\phantom{Set^{{(\pOpe)}^{op}}}`
\phantom{Set^{{(\pOpeCard)}^{op}}}`\bi^*]{1500}{-1}b
\putmorphism(0,100)(1,0)[Set^{{(\pOpe)}^{op}}`Set^{{(\pOpeCard)}^{op}}`]{1500}{0}b
\end{picture}
\end{center}
where $\bi^*$ is the functor of composing with $\bi$ and $Ran_\bi$ is the right Kan extension along $\bi$. Recall that for $F$ in $Set^{{(\pOpe)}^{op}}$, $S$ in $\pOpeCard$, it is defined as the following limit
\[ (Ran_\bi F)(S) = Lim F\circ \Sigma^{S,op} \]
where $\Sigma^{S,op}$ is the dual of the functor $\Sigma^S$ defined before Lemma \ref{speclim}. Note that, as $(\pOpe \da S)^{op} = S\da (\pOpe)^{op}$, we have
\[ \Sigma^{S,op} : S\da (\pOpe)^{op}\lra (\pOpe)^{op}. \]

As $\bi$ is full and faithful, the right Kan extension $Ran_\bi (F)$ is an extension. Therefore the counit of this adjunction
\[ \varepsilon_F : (Ran_\bi\, F)\circ \bi \lra F \]
is an isomorphism. The functor $Ran_\bi F$ is so defined that it preserves special limits. Hence, by Lemma \ref{speclim}, it preserves special pullbacks.
As any positive opetopic cardinal can be constructed from positive opetopes via special pushouts, for $G$ in $Set^{{(\pOpeCard)}^{op}}$, the unit of adjunction
\[  \eta_G : G \lra Ran_\bi (G\circ\bi) \]
is an isomorphism iff $G$ preserves special pullbacks.  Thus we have

\begin{proposition}\label{presheaf_eq1}
The above adjunction restricts to the following equivalence of categories
\begin{center} \xext=1500 \yext=250
\begin{picture}(\xext,\yext)(\xoff,\yoff)
\putmorphism(0,150)(1,0)[\phantom{Set^{{(\pOpe)}^{op}}}`
\phantom{sPb((\pOpeCard)^{op},Set)}`Ran_\bi]{1500}{1}a
\putmorphism(0,50)(1,0)[\phantom{Set^{{(\pOpe)}^{op}}}`
\phantom{sPb((\pOpeCard)^{op},Set)}`\bi^*]{1500}{-1}b
\putmorphism(0,100)(1,0)[Set^{{(\pOpe)}^{op}}`sPb((\pOpeCard)^{op},Set)`]{1500}{0}b
\end{picture}
\end{center}
$~~\Box$
\end{proposition}

Now we will set up the adjunction
\begin{center} \xext=2500 \yext=300
\begin{picture}(\xext,\yext)(\xoff,\yoff)
\putmorphism(0,200)(1,0)[\phantom{sPb((\pOpeCard)^{op},Set)}`
\phantom{\pPoly}`\widetilde{(-)}]{2500}{1}a
\putmorphism(0,100)(1,0)[\phantom{sPb((\pOpeCard)^{op},Set)}`
\phantom{\pPoly}`\widehat{(-)}=\pPoly((\simeq)^*, - )]{2500}{-1}b
\putmorphism(0,150)(1,0)[sPb((\pOpeCard)^{op},Set)`\pPoly`]{2500}{0}b
\end{picture}
\end{center}
which will be later proved to be an equivalence of categories. The functor $\widehat{(-)}$ sends a positive-to-one polygraph $P$ to a functor
\[ \widehat{P} = \pPoly((-)^*,P) : (\pOpeCard)^{op} \lra Set \]
$\widehat{(-)}$ is defined on morphism in the obvious way, by composition. We have

\begin{lemma} Let $P$ be a positive-to-one polygraph. Then $\widehat{P}$ defined above is a special pullbacks preserving functor.
\end{lemma}

{\it Proof.}~This is an immediate consequence of the fact that the functor $(-)^*$ preserves special pushouts. $~~\Box$

Now suppose we have a special pullbacks preserving functor $F : (\pOpeCard)^{op} \lra Set$. We shall define a positive-to-one polygraph $\widetilde{F}$.

As $n$-cells of $\widetilde{F}$ we put
\[ \widetilde{F}_n = \coprod_S F(S) \]
where the coproduct is taken over all\footnote{ In fact, we think about such a coproduct $\coprod_S F(S)$ as if it were to be taken over sufficiently large (so that each isomorphism type of positive opetopic cardinals is represented) set of positive opetopic cardinals $S$ of dimension at most $n$. Then, if positive opetopic cardinals
$S$ and $S'$ are isomorphic via (necessarily unique) isomorphism $h$, then the cells $x\in F(S)$ and $x'\in F(S')$ are considered equal iff $F(h)(x)=x'$.} (up to iso) positive opetopic cardinals $S$ of dimension at most $n$.

If $k\leq n$, the identity map
\[ 1^{(n)} : \widetilde{F}_k \lra \widetilde{F}_n \]
is the obvious embedding induced by identity maps on the components of the coproducts.

Now we shall describe the domains and codomains in $\widetilde{F}$. Let $S$ be a positive opetopic cardinal of dimension at most $n$, $a\in F(S)\hra \widetilde{F}_n$. We have in $\pOpeCard$ the $k$-th domain and the $k$-th codomain morphisms:
\begin{center}
\xext=800 \yext=400 \adjust[`I;I`;I`;`I]
\begin{picture}(\xext,\yext)(\xoff,\yoff)
 \settriparms[-1`-1`0;400]
 \putAtriangle(0,0)[S`\bd^{(k)} S`\bc^{(k)} S;\bd^{(k)}_S`\bc^{(k)}_S`]
\end{picture}
\end{center}
We put
\[ \bd^{(k)}(a) = F(\bd^{(k)}_S)(a)\in F(\bd^{(k)}S) \hra\widetilde{F}_k ,\]
\[ \bc^{(k)}(a) = F(\bc^{(k)}_S)(a)\in F(\bc^{(k)}S) \hra\widetilde{F}_k. \]

Finally, we define the compositions in $\widetilde{F}$. Let $n_1,n_2\in \o$, $n=max(n_1,n_2)$, $k<min(n_1,n_2)$, and
\[ a\in F(S) \hra\widetilde{F}_{n_1} \hskip 10mm b\in F(T) \hra\widetilde{F}_{n_2},\]
 such that
 \[ \bc^{(k)}(a)=F(\bc^{(k)}_S)(a)=F(\bd^{(k)}_T)(b) = \bd^{(k)}(b). \]
We shall define the cell $b\circ_ka\in\widetilde{F}_n$. We take a special pushout in $\pOpeCard$:
\begin{center} \xext=1500 \yext=680
\begin{picture}(\xext,\yext)(\xoff,\yoff)
 \setsqparms[-1`-1`-1`-1;800`500]
 \putsquare(0,100)[T\oplus_kS`S`T`{\bc^{(k)}S}; \kappa_1`\kappa_2`\bc^{(k)}_S`\bd^{(k)}_T]
\end{picture}
\end{center}

As $F$ preserves special pullbacks (in $(\pOpeCard)^{op}$), it follows that the square
\begin{center} \xext=1500 \yext=700
\begin{picture}(\xext,\yext)(\xoff,\yoff)
 \setsqparms[1`1`1`1;800`500]
 \putsquare(0,100)[F(T\oplus_kS)`F(S)`F(T)`F({\bc^{(k)}S});
 F(\kappa_1)`F(\kappa_2)`F(\bc^{(k)}_S)`F(\bd^{(k)}_T)]
\end{picture}
\end{center}
is a pullback in $Set$. Thus there is a unique element
\[ x\in F(S\oplus_kT)\hra \widetilde{F}_n \]
such that
\[ F(\kappa_1)(x)=a, \hskip 10mm   F(\kappa_2)(x)=b. \]
 We put
\[ b\circ_ka = x. \]
This ends the definition of $\widetilde{F}$.

 For a morphism $\alpha :F\lra G$ in $sPb((\pOpeCard)^{op},Set)$, we put
 \[ \widetilde{\alpha} = \{ \widetilde{\alpha}_n :
 \widetilde{F}_n\lra \widetilde{G}_n \}_{n\in\o} \]
 such that
 \[ \widetilde{\alpha}_n = \coprod_S \alpha_S : \widetilde{F}_n\lra \widetilde{G}_n\]
where the coproduct is taken over all (up to iso) positive opetopic cardinals $S$ of dimension at most $n$. This ends the definition of the functor $\widetilde{(-)}$.

We have

\begin{proposition}
The functor
\[ \widetilde{(-)} : sPb((\pOpeCard)^{op},Set)\lra \pPoly \]
is well defined.
\end{proposition}

{\it Proof.}~The verification that $\widetilde{(-)}$ is a functor into $\oC$ is left for the reader.  We shall verify that, for any special pullbacks preserving functor $F:{\pOpeCard}^{op}\lra Set$, $\widetilde{F}$ is a positive-to-one polygraph whose $n$-indeterminates are
 \[ |\widetilde{F}|_n = \coprod_{B\in \pOpe, dim(B)=n} F(B)\; \hra \;
 \coprod_{S\in \pOpeCard, dim(S)\leq n} F(S) = \widetilde{F}_n.\]

Let $P$ be the truncation of $\widetilde{F}$ in  $Set\da D_n$, i.e., $P=U_n(\widetilde{F})$. We shall show that $\widetilde{F}_n$ is in a bijective correspondence with $\overline{P}_n$ described in the previous section. We define a function
\[ \varphi : \overline{P}_n \lra \widetilde{F}_n\]
so that, for a cell $f:S^{\sharp,n}\lra P$ in $\overline{P}_n$, we put
\[  \varphi(f) = \left\{ \begin{array}{ll}
        1_{f_{n-1}(S)}   & \mbox{if $dim(S)<n$,}  \\
        f_{n}(m_S)   & \mbox{if $dim(S)=n$, $S$ principal, $S_n=\{ m_S \}$}  \\
       \varphi(f^{\ua a}) \circ_k \varphi(f^{\da a}) & \mbox{if $dim(S)=n$, $a\in Sd(S)_k$.}
                                    \end{array}
                \right. \]
and the morphisms in $\varphi(f^{\da a})$ and $\varphi(f^{\ua a})$ in $Set\da D_n$ are obtained by compositions  so that the diagram
\begin{center}
\xext=1000 \yext=700
\begin{picture}(\xext,\yext)(\xoff,\yoff)
\settriparms[0`1`-1;300]
 \putDtriangle(0,0)[(S^{\da a})^{\sharp,n}`S^{\sharp,n}`(S^{\ua a})^{\sharp,n};``]
 \putmorphism(400,300)(1,0)[`P`f]{600}{1}a
 \putmorphism(0,700)(3,-1)[``f^{\da a}]{1100}{1}r
 \putmorphism(400,50)(3,1)[``f^{\ua a}]{280}{1}r
\end{picture}
\end{center}
commutes. We need to verify, by induction on $n$, that $ \varphi$ is well defined, bijective and that it preserves compositions, domains, and codomains.

We shall only verify (partially) that $\phi$ is well defined, i.e., the definition $\phi$ for any positive opetopic cardinal $S$ of dimension $n$ that is not a positive opetope does not depend on the choice of the saddle point of $S$. Let $a,x\in Sd(S)$ so that $k=dim(x)<dim(a)=m$. Using  Lemma \ref{decomp1new} and the fact that $(-)^{\sharp,n}$ preserves special pushouts, we have
\begin{eqnarray*}
 \varphi(f^{\ua a})\circ_m\varphi(f^{\da a}) =    \\
 = (\varphi(f^{\ua a\ua x})\circ_k \varphi(f^{\ua a \da x}))\circ_m( \varphi(f^{\da a\ua x})\circ_k \varphi(f^{\da a\da x})) =  \\
 =  (\varphi(f^{\ua a\ua x})\circ_m\varphi(f^{\da a\ua x}))\circ_k (\varphi(f^{\ua a \da x})\circ_m \varphi(f^{\da a\da x})) =   \\
 =  (\varphi(f^{\ua x\ua a})\circ_m \varphi(f^{\ua x\da a}))\circ_k (\varphi(f^{\da x \ua a})\circ_m \varphi(f^{\da x\da a})) =   \\
 = \varphi(f^{\ua x})\circ_m \varphi(f^{\da x}),
\end{eqnarray*}
as required in this case. The reader can compare these calculations with those, in the same case, of Proposition \ref{Sstar} ($F$ is replaced by $\varphi$ and $T$ is replaced by $f$). So there is no point to repeat the other calculations. $~~\Box$

 For $P$ in $\pPoly$, we define a polygraph map
\[ \eta_P : P\lra \widetilde{\widehat{P}} \]
so that, for $x\in P_n$, we put
\[ \eta_{P,n}(x)=\tau_x : T_x^*\ra P. \]

For $F$ in $sPb((\pOpeCard)^{op},Set)$, we define a natural transformation
\[ \varepsilon_F : \widehat{\widetilde{F}} \lra F, \]
such that, for a positive opetopic cardinal $S$ of dimension $n$,
\[ (\varepsilon_F)_S : \widehat{\widetilde{F}}(S) \lra F(S) \]
and $g : S^*\ra \widetilde{F} \in \widehat{\widetilde{F}}(S)$, we put
\[ (\varepsilon_F)_S(g) = g_n(S).\]

\begin{proposition}\label{presheaf_eq2}
The functors
\begin{center} \xext=2500 \yext=300
\begin{picture}(\xext,\yext)(\xoff,\yoff)
\putmorphism(0,200)(1,0)[\phantom{sPb((\pOpeCard)^{op},Set)}`
\phantom{\pPoly}`\widetilde{(-)}]{2500}{1}a
\putmorphism(0,100)(1,0)[\phantom{sPb((\pOpeCard)^{op},Set)}`
\phantom{\pPoly}`\widehat{(-)}=\pPoly((\simeq)^*, - )]{2500}{-1}b
\putmorphism(0,150)(1,0)[sPb((\pOpeCard)^{op},Set)`\pPoly`]{2500}{0}b
\end{picture}
\end{center}
together with the natural transformations $\eta$ and $\varepsilon$ defined above form an adjunction ($\widehat{(-)}\dashv\widetilde{(-)}$). It establishes the equivalence of categories $sPb((\pOpeCard)^{op},Set)$ and $\pPoly$.
\end{proposition}

{\it Proof.}~ The fact that both $\eta$ and $\varepsilon$ are bijective on each component follows immediately from Proposition \ref{type}.6. So we shall verify the triangular equalities only.

Let $P$ be a polygraph, and $F$ be a functor in $sPb((\pOpeCard)^{op},Set)$. We need to show that the triangles
\begin{center}
\xext=2000 \yext=400 \adjust[`I;I`;I`;`I]
\begin{picture}(\xext,\yext)(\xoff,\yoff)
 \settriparms[-1`1`1;400]
 \putAtriangle(0,0)[\widehat{\widetilde{\widehat{P}}}`\widehat{P}`\widehat{P};
 \widehat{\eta_P}`\varepsilon_{\widehat{P}}`1_{\widehat{P}}]
 \putAtriangle(1200,0)[\widetilde{\widehat{\widetilde{F}}}`\widetilde{F}`\widetilde{F};
 \eta_{\widetilde{F}}`\widetilde{\varepsilon_{F}}`1_{\widetilde{F}}]
\end{picture}
\end{center}
commute. So let $f:S^*\ra P \in \widehat{P}(S)$, $dim(S)=n$. Then we have

\[ \varepsilon_{\widehat{P}}\circ \widehat{\eta_P} (f) = \varepsilon_{\widehat{P}}(\eta_P\circ f)= (\eta_P\circ f)_n(S) = \]

\[ =   (\eta_P)_n(f_n(S)) = \tau_{f_n(S)} = f \]

So let $x \in F(S)\lra \widetilde{F}_n$. Then we have

\[ \widetilde{\varepsilon_F}\circ \eta_{\widetilde{F}}(x)= \widetilde{\varepsilon_F}(\tau_x) = (\tau_x)_n(1_{T_x}) =x  \] So
both triangles commutes, as required. $~~\Box$

From Propositions \ref{presheaf_eq1} and \ref{presheaf_eq2} we get immediately

\begin{corollary}\label{xxx}
The functor
\[ \widehat{(-)} :  \pPoly \lra Set^{(\pOpe)^{op}} \]
such that, for a positive-to-one polygraph $X$,
\[ \widehat{X} = \pPoly((-)^*,X): (\pOpe)^{op}\lra Set  \]
is an equivalence of categories.
\end{corollary}

\section{The principal pushouts}\label{sec-principal-pushouts}

Recall the positive opetopic cardinal $\alpha^n$ from section \ref{sec-Sstar_is_pPoly}. A {\em total composition map}\index{map!total composition -}\index{composition!total - map}
is an inner $\o$-functor whose domain is of form $(\alpha^n)^*$, for some $n\in\o$. If $S$ is a positive opetopic cardinal of dimension $n$, then the total composition of $S$ (in fact $S^*$) is denoted by
\[ \bmu^{S^*} : \alpha^{n,*}\lra S^*. \]
It is uniquely determined by the condition $\bmu^{S^*}_n(\alpha^n)=S$. We have the following

\begin{proposition}\label{principal pushout}
Let $N$ be a normal positive opetopic cardinal. With the notation as above, the square
\begin{center} \xext=700 \yext=650
\begin{picture}(\xext,\yext)(\xoff,\yoff)
 \setsqparms[1`-1`-1`1;700`500]
 \putsquare(0,50)[N^*`N^{\bullet,*}`\alpha^{n,*}`\alpha^{n+1,*};
 \bd_{N^\bullet}^*`\bmu^{N^*}`\bmu^{N^{\bullet,*}}`\bd_{\alpha^{n+1}}^*]
\end{picture}
\end{center}
is a pushout in $\pOpeCardo$.
\end{proposition}
{\it Proof.}~ This is an easy consequence of Proposition \ref{bullet}, particularly point 4.  $~~\Box$

\vskip 3mm
Pushouts described in the above Proposition are called {\em principal pushouts}\index{pushout!principal}\index{principal
pushout}.

From the above proposition we immediately get

\begin{corollary}\label{principal pushout1}
If $n>0$ and $P$ is a positive opetope of dimension $n$, then the square
\begin{center} \xext=700 \yext=650
\begin{picture}(\xext,\yext)(\xoff,\yoff)
 \setsqparms[1`-1`-1`1;700`500]
 \putsquare(0,50)[\bd P^*`P^*`\alpha^{n-1,*}`\alpha^{n,*};
 \bd_P^*`\bmu^{\bd P}`\bmu^P`\bd_{\alpha^{n,*}}]
\end{picture}
\end{center}
is a (principal) pushout in $\pOpeCardo$. $~~\Box$
\end{corollary}

\begin{theorem}[V.Harnik]\label{principal pushout2}\footnote{The original statement of V.Harnik is saying that the nerve functor from $\o$-categories to (all) polygraphs is monadic. However, in the present context the argument given by V.Harnik, c.f. \cite{H}, is directly proving the present statement, i.e., that the principal pullbacks are preserved whenever the special ones are. This statement is used to show that the category of $\o$-categories is equivalent to the category of special pullback preserving functors from $(\pOpeCardo)^{op}$ to $Set$, c.f. Corollary \ref{sp-nevre2}. From that statement, the monadicity of the nerve functor is an easy corollary, c.f. Theorem \ref{monadicity}. In the remainder of this section the Harnik's argument, adopted to the present context, is presented.} Let $F:(\pOpeCardo)^{op}\lra Set$ be a special pullback preserving functor. Then $F$ preserves the principal pullbacks, as well.
\end{theorem}

Theorem \ref{principal pushout2} is a special case of Lemma \ref{principal pushout5}, for $k=n-1$.

The proof of the above theorem will be divided into three Lemmas. Theorem \ref{principal pushout2} is a special case of Lemma \ref{principal pushout5}, for $k=n-1$.

\vskip 3mm

Before we even formulate the next three Lemmas, we need to introduce some constructions on positive opetopic cardinals and define some $\o$-functors between positive opetopic cardinals. Introducing these constructs and notation for them, we shall make some comments how they are going to be interpreted by special pullback preserving morphisms from $(\pOpeCardo)^{op}$ to $Set$.

\vskip 3mm

{\em Notation for presheaves.} To simplify the notation concerning presheaf $F:\pOpeCardo^{op}\ra Set$, for a morphism $g:P\ra  Q\in \pOpeCardo$, and an element $a\in F(Q)$, we will write $a\cdot g$ for $F(g)(a)$, i.e., we treat $F$ as a family of sets with a right action of the category $\pOpeCardo$. This convention is explained, for example, in \cite{SGL} page 121. Such notation suppresses the name of the presheaf $F$ but it will be always clear from the context. In this section we deal only with one presheaf named $F$ that preserves special pullbacks. We also have $(a\cdot g)\cdot f =a\cdot (g\circ f)$ and $a\cdot id_Q=a$ whenever these formulas are well defined.

\vskip 2mm

Fix $0<k\leq n$, and a $P$ positive opetope of dimension $n$.  We say that $P$ is $k${\em -globular}\index{globular!k-} iff $\bd^{(l)}P$ is a positive opetope, for $k\leq l\leq n$, i.e., $\delta^{(l)}(\bp_n)$ is a singleton, for $k\leq l\leq n$, where $P_n=\{ \bp_n\}$. The $k${\em -globularization}\index{globularization!k-} $\fk P$ of $P$ is the $k$-globular positive opetope of dimension $n$ defined as follows. 
We put
  \[ \fk P_l  = \left\{ \begin{array}{ll}
        \{ \bp_n\} & \mbox{for  $l=n$,}  \\
        \{ \bq_l,\bp_l\} & \mbox{for  $k\leq l<n$,}  \\
        \delta^{(k-1)}(\bp_n)\cup \{\bp_{k-1} \} & \mbox{for  $l=k-1$,}  \\
        P_l & \mbox{otherwise.}
                                    \end{array}
                \right. \]

For $x\in\fk P$,
\[ \gamma^{\fk P}(x)  = \left\{ \begin{array}{ll}
        \bp_{l-1} & \mbox{if $x=\bq_l$, for some  $k\leq l<n$,}  \\
        \gamma^P(x) & \mbox{otherwise.}
                                    \end{array}
                \right. \]
and
\[ \delta^{\fk P}(x)  = \left\{ \begin{array}{ll}
        \bq_l & \mbox{if $x\in \fk P_{l+1}$, for some  $k< l<n$,}  \\
        \delta^P(\bp_k) & \mbox{if $x\in \fk P_k$,}  \\
        \delta^P(x) & \mbox{otherwise.}
                                    \end{array}
                \right. \]
Note that $\fn P$ is $P$ itself and $\fn P$ is $\alpha^n$. Thus the elements of the shape $\fk P^*$ are $k$-globularized versions of the elements of the shape $P^*$. As the following positive opetopic cardinals
\[ \bc^{(k)} P\;\cong \;  \bc^{(k)} \fk P\;\cong\; \bc\bc^{(k+1)} \fk P\;\cong\; \bd\bc^{(k+1)} \fk P\;\cong\; \bd^{(k)} \fk P \]
are isomorphic, we can form the following special pushouts
\begin{center} \xext=1000 \yext=650
\begin{picture}(\xext,\yext)(\xoff,\yoff)
 \setsqparms[1`-1`-1`1;1000`500]
 \putsquare(0,70)[\bc^{(k+1)}P^*`\bc^{(k+1)}P+_{\bc^{(k)}P}\fk P^*`\bc^{(k)}P^*`\fk P^*;
 \kappa_1`\bc^*_{\bc^{(k+1)} P}`\kappa_2`\bd^{(k),*}_{\sfk P}]
\end{picture}
\end{center}
and
\begin{center} \xext=1000 \yext=650
\begin{picture}(\xext,\yext)(\xoff,\yoff)
 \setsqparms[1`-1`-1`1;1000`500]
 \putsquare(0,70)[\bc^{(k+1)}\fk P^*`\bc^{(k+1)}\fk P+_{\bc^{(k)}P}\fk P^*`\bc^{(k)}P^*`\fk P^*;
 \kappa'_1`\bc^*_{\bc^{(k+1)}\fk P}`\kappa'_2`\bd^{(k),*}_{\sfk P}]
\end{picture}
\end{center}
We describe in detail the positive opetopic cardinals we have just defined. Their faces are described in the following table:

 \[ \begin{array}{|c|c|c|c|c|}\hline
    dim & \fkj P         & \fk P & P'=\bc^{(k+1)} P+_{\bc^{(k)}P}\fk P & P''=\bc^{(k+1)}\fk P+_{\bc^{(k)}P}\fk P \\ \hline

    n & \{ \bp_n \}& \{ \bp_n \}&\{ \bp_n \}&\{ \bp_n \} \\

    n\!\! -\!\! 1 &\{ \bq_{n-1}, \bp_{n-1} \} &\{\bq_{n-1}, \bp_{n-1} \} & \{\bq_{n-1}, \bp_{n-1} \} &\{ \bq_{n-1}, \bp_{n-1} \}\\
    \vdots &  \vdots &  \vdots &  \vdots &  \vdots \\

    k\!\! +\!\! 1   & \{ \bq_{k\! +\! 1}, \bp_{k\! +\! 1} \} &\{\bq_{k\! +\! 1}, \bp_{k\! +\! 1} \} &     \{ \br_{k\! +\! 1},\bq_{k\! +\! 1}, \bp_{k\! +\! 1} \}
    & \{\br_{k\! +\! 1}, \bq_{k\! +\! 1}, \bp_{k\! +\! 1} \} \\

    k  & \partial(\bp_{k+1}) &\{\bq_{k}, \bp_k \} & \delta(\bp_{k+1})\cup \{\bq_{k}, \bp_{k} \} & \{ \br_{k}, \bq_{k}, \bp_{k} \}\\

    k\!\! -\!\! 1 & P_{k-1} &  \partial(\bp_{k}) &  P_{k-1} &  \partial(\bp_{k})\\

    \vdots &  \vdots &  \vdots &  \vdots &  \vdots \\

    l & P_l & P_l & P_l & P_l\\\hline
  \end{array}
\]
$0\leq l <k$. The functions $\gamma$ and $\delta$ in $P'=\bc^{(k+1)} P+_{\bc^{(k)}P}\fk P$ are given by the following formulas

\[ \gamma^{P'}(x)  = \left\{ \begin{array}{ll}
        \bp_{l-1} & \mbox{if $x=\bq_l\;$ {\rm and} $\;k\leq l<n$}  \\
         \bq_k & \mbox{if $x=\br_{k+1}$}  \\
         \gamma^P(x) & \mbox{otherwise.}  \\
                                    \end{array}
                \right. \]

\[ \delta^{P'}(x)  = \left\{ \begin{array}{ll}
        \{ \bq_{l-1}\} & \mbox{if $x=\bq_{l}\;$ {\rm and} $\; k < l<n$}  \\

         & \mbox{or $x=\bp_l\;$ {\rm and} $\; k< l\leq n$}  \\
         \delta^P(\bp_{k+1}) & \mbox{if $x=\br_{k+1}$}  \\
         \delta^P(\bp_{k}) & \mbox{if $x=\bq_{k}$}  \\
         \delta^P(x) & \mbox{otherwise.}  \\
                                    \end{array}
                \right. \]
The functions $\gamma$ and $\delta$ in $P''=\bc^{(k+1)}\fk P+_{\bc^{(k)}P}\fk P$ are given by the following formulas

\[ \gamma^{P''}(x)  = \left\{ \begin{array}{ll}
        \bp_{l-1} & \mbox{if $x=\bq_l\;$ {\rm and} $\;k\leq l<n$}  \\
         \bq_k & \mbox{or $x=\br_{k+1}$}  \\
         \bp_{k-1} & \mbox{or $x=\br_k$}  \\
         \gamma^P(x) & \mbox{otherwise.}  \\
                                    \end{array}
                \right. \]

\[ \delta^{P''}(x)  = \left\{ \begin{array}{ll}
        \{ \bq_{l-1}\} & \mbox{if $x=\bq_{l}\;$ {\rm and} $\; k < l<n$}  \\

         & \mbox{or $x=\bp_l\;$ {\rm and} $\; k< l\leq n$}  \\
         \{ \br_{k}\} & \mbox{if $x=\br_{k+1}$}  \\
         \delta^P(\bp_{k}) & \mbox{if $x\in \{ \br_{k},\bq_{k},\bp_{k}  \}$  }  \\
         \delta^P(x) & \mbox{otherwise.}  \\
                                    \end{array}
                \right. \]

Now we shall define some $\o$-functors between some positive opetopic cardinals just defined. To describe their meaning, let us fix a special pullback preserving functors from $F:(\pOpeCardo)^{op}\lra Set$.

The $\o$-functors denoted by letter $\bmu$ are interpreted as operation that `globularize'  cells. We have two of them. The first one
$$\bmu^{S^*}:\alpha^{n,*}\lra S^*$$
was already introduced at the beginning of this section for any positive opetopic cardinal $S$. The second is the $\o$-functor
\[ \bmu_{\sfk} : \fk P^* \lra \fkj P^* \]
such that
\[ \bmu_{\sfk} (X)  = \left\{ \begin{array}{ll}
        (X-\{ \bq_k\})\cup\delta(\bp_{k+1}) & \mbox{if $\bq_k\in X$}  \\
        X & \mbox{otherwise.}
                                    \end{array}
                \right. \]
for $X\in \fk P^* $.

The fact that these operations are interpreted as globularization of cells can be explained as follows. The function
\[ F(\bmu_{\sfk}) : F(\fkj P^*) \lra F(\fk P^*), \]
takes a $(k+1)$-globular $n$-cell $a\in F(\fkj P^*)$ and returns a $k$-globular $n$-cell $a\cdot \bmu_{\sfk} = F(\bmu_{\sfk})(a)\in F(\fk P^*)$. Intuitively, $F(\bmu_{\sfk})$  is composing the $k$-domain of $a$ leaving the rest `untouched'. So it is a `one-step globularization'. On the other hand, the function
\[  F(\bmu^{S^*}):  F(S^*)\lra F(\alpha^{n,*}) \]
is taking an $n$-cell $b\in F(S^*)$ of an arbitrary shape $S^*$ of dimension $n$, and it is returning a globular $n$-cell $b\cdot \bmu^{S^*} \in F(\alpha^{n,*})$. This time $F(\bmu^{S^*})$ is composing all the domains and codomains in the cell $b$ as much as possible, so that there is nothing left to be composed. This is the `full globularization'.

We need a separate notation for the $\o$-functor $\bmu_k : \bc^{(k+1)}\fk P^* \lra \bc^{(k+1)} P^*$ such that

\[ \bmu_k (X)  = \left\{ \begin{array}{ll}
        \bc^{(k+1)} P & \mbox{if $X=\bc^{(k+1)}\fk P$}  \\
        \bd^{(k)} P & \mbox{if $X=\bd^{(k)}\fk P$}  \\
        X & \mbox{otherwise,}
                                    \end{array}
                \right. \]
for $X\in \bc^{(k+1)}\fk P^* $. It is a version of $\bmu_{\sfk}$. The $\o$-functor
\[ \bio_P :P^* \lra \bd P^* \]
is given by
\[ \bio_P(X)  = \left\{ \begin{array}{ll}
        \bd P & \mbox{if $\bc P\subseteq X$}  \\
        X & \mbox{otherwise,}
                                    \end{array}
                \right. \]
for $X\in P^*$. $\bio_P$ is a kind of degeneracy map and it is interpreted as `a kind of identity'.  For a cell $t\in F(\bd P^*)$, $t\cdot \bio_P\in  F(P^*)$ is `identity on $t$' but with the codomain composed. The $\o$-functor
\[ \bbeta_k: \bc^{(k)}P^*\lra  \bd^{(k)}P^*\]
such that
\[ \bbeta_k(X)  = \left\{ \begin{array}{ll}
        \bd^{(k)} P & \mbox{if $X=\bc^{(k)} P$}  \\
        X & \mbox{otherwise.}
                                    \end{array}
                \right. \]
for $X\in P^*$, is the operation of `composition of all the cells at the top' leaving the rest untouched. The map $\bbeta_{n-1}$ is equal to the composition
\begin{center} \xext=1200 \yext=100
\begin{picture}(\xext,\yext)(\xoff,\yoff)
 \putmorphism(0,50)(1,0)[\bc P^*`P^*`\bc^*_P]{600}{1}a
 \putmorphism(600,50)(1,0)[\phantom{P^* }`\bd P^*`\bio_P]{600}{1}a
\end{picture}
\end{center}
The following two $\o$-functors
\[ {[\bd^{(k),*}_P\circ \bio_{\bc^{(k+1)}P},\bmu_P ]} : \bc^{(k+1)} P+_{\bc^{(k)}P}\fk P^* \lra P^* \]
and
\[ {[\bd^{(k),*}_{\sfk P}\circ \bio_{\bc^{(k+1)}\fk P},1_{\sfk P} ]} : \bc^{(k+1)}\fk P+_{\bc^{(k)}P}\fk P^* \lra \fk P^* \]
are defined as the unique $\o$-functors making the following diagrams
\begin{center} \xext=2000 \yext=1200
\begin{picture}(\xext,\yext)(\xoff,\yoff)
 \setsqparms[1`-1`-1`1;1000`500]
 \putsquare(0,50)[\bc^{(k+1)}P^*`\bc^{(k+1)}P+_{\bc^{(k)}P}\fk P^*`\bc^{(k)}P^*`\fk P^*;
 \kappa_1`\bc^*_{\bc^{(k+1)} P}`\kappa_2`\bd^{(k),*}_{\sfk P}]
 \putmorphism(1000,550)(1,0)[\phantom{\bc^{(k+1)}P+_{\bc^{(k)}P}\fk P^*}`\fkj P^*`{[}\bd^{(k)}\!\circ \bio ,\bmu {]}]{1200}{1}a
 \put(1070,1070){\makebox(100,100){$\bd^{(k)}P^*$}}
 \putmorphism(200,700)(2,1)[``\bio_{\bc^{(k+1)}P^*}]{600}{1}l
 \putmorphism(1150,1150)(2,-1)[``\bd^{(k),*}_P]{1200}{1}r
 \putmorphism(1300,130)(2,1)[``\bmu_{\sfk}]{600}{1}r
\end{picture}
\end{center}
and
\begin{center} \xext=2000 \yext=1200
\begin{picture}(\xext,\yext)(\xoff,\yoff)
 \setsqparms[1`-1`-1`1;1000`500]
 \putsquare(0,50)[\bc^{(k+1)}\fk P^*`\bc^{(k+1)}\fk P+_{\bc^{(k)}P}\fk P^*`\bc^{(k)}P^*`\fk P^*;
 \kappa'_1`\bc^*_{\bc^{(k+1)}\fk P}`\kappa'_2`\bd^{(k),*}_{\sfk P}]
 \putmorphism(1000,550)(1,0)[\phantom{\bc^{(k+1)}P+_{\bc^{(k)}P}\fk P^*}`\fk P^*`{[} \bd^{(k)}\!\circ \bio ,1_{\sfk P^*} {]}]{1200}{1}a
 \put(1070,1070){\makebox(100,100){$\bd^{(k)}\fk P^*$}}
 \putmorphism(200,700)(2,1)[``\bio_{\bc^{(k+1)}\fk P^*}]{600}{1}l
 \putmorphism(1150,1150)(2,-1)[``\bd^{(k),*}_{\fk P}]{1200}{1}r
 \putmorphism(1300,130)(2,1)[``1_{\sfk P^*}]{600}{1}r
\end{picture}
\end{center}
commute in $\pOpeCardo$.

Finally, we introduce two maps that are a kind of binary composition combined with wiskering. The first
\[ \bcomp_{\sfkj} : \fkj P^* \lra \bc^{(k+1)} P+_{\bc^{(k)}P}\fk P \]
is given by
\[ \bcomp_{\sfkj}(X)  = \left\{ \begin{array}{ll}
        X\cup\{\br_{k+1}, \bq_{k} \} & \mbox{if $X\cap\{ \bq_{k+1}, \bp_{k+1}\} \neq\emptyset$}  \\
        X & \mbox{otherwise,}
                                    \end{array}
                \right. \]
for $X\in \fkj P^*$. The other composition map is
$$\bcomp_{k} : \fk P^* \lra \bc^{(k+1)} \fk P+_{\bc^{(k)}P}\fk P $$
given by the same defining formula as $\bcomp_{\sfkj}$, for $X\in \fk P^*$.

In the following diagram all the morphisms that we introduced above are displayed. Most of the subscripts of the morphisms are suppressed for clarity of the picture.
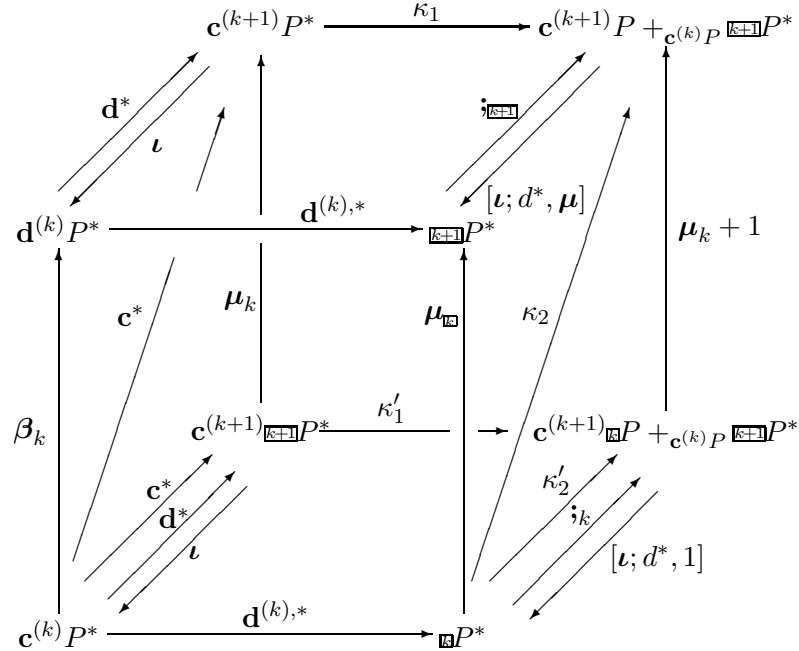
\begin{center} \xext=2500 \yext=2200
\begin{picture}(\xext,\yext)(\xoff,\yoff)
 \setsqparms[1`-1`-1`1;1400`1400]
 \putsquare(0,0)[\bd^{(k)}P^*`\fkj P^*`\bc^{(k)}P^*`\fk P^*;
 `\bbeta_k``]
 \put(1270,1050){\makebox(100,100){$\bmu_{\sfk}$}}
 \put(700,20){\makebox(100,100){$\bd^{(k),*}$}}
 \put(900,1420){\makebox(100,100){$\bd^{(k),*}$}}
 \setsqparms[1`0`-1`0;1400`1400]
 \putsquare(700,700)[\bc^{(k+1)}P^*`\bc^{(k+1)}P+_{\bc^{(k)}P}\fk P^*` \bc^{(k+1)}\fk P^*`\bc^{(k+1)}\fk P+_{\bc^{(k)}P}\fk P^*; \kappa_1``\bmu_k+1`]

 \putmorphism(160,250)(1,1)[``]{300}{1}l
 \put(300,450){\makebox(100,100){$\bc^*$}}
 \putmorphism(240,190)(1,1)[``]{300}{1}l
 \put(350,350){\makebox(100,100){$\bd^*$}}
 \putmorphism(280,140)(1,1)[``]{300}{-1}r
 \put(420,230){\makebox(100,100){$\bio$}}

 \putmorphism(1560,250)(1,1)[``]{300}{1}l
 \put(1670,480){\makebox(100,100){$\kappa'_2$}}
 \putmorphism(1640,170)(1,1)[``]{300}{1}l
 \put(1755,360){\makebox(100,100){$\bcomp_k$}}
 \putmorphism(1700,120)(1,1)[``]{300}{-1}r
 \put(2030,210){\makebox(100,100){$[d^*\!\circ\bio ,\!1]$}}

 \putmorphism(1410,1600)(1,1)[``]{340}{1}l
 \put(1470,1770){\makebox(100,100){$\bcomp_{\sfkj}$}}
 \putmorphism(1460,1550)(1,1)[``]{340}{-1}r
 \put(1600,1450){\makebox(100,100){$[ d^*\!\!\circ\bio,\!\bmu]$}}

 \putmorphism(70,1600)(1,1)[``]{340}{1}l
 \put(150,1770){\makebox(100,100){$\bd^*$}}
 \putmorphism(110,1550)(1,1)[``]{340}{-1}r
 \put(290,1640){\makebox(100,100){$\bio$}}

 \putmorphism(1450,250)(1,3)[``]{500}{1}l
 \putmorphism(500,1600)(1,3)[``]{50}{1}l
 \put(50,250){\line(1,3){350}}
 \put(200,1050){\makebox(100,100){$\bc^*$}}
 \put(1600,1050){\makebox(100,100){$\kappa_2$}}


 \put(580,1100){\makebox(100,100){$\bmu_k$}}
 \put(700,800){\line(0,1){550}}
 \put(700,1450){\vector(0,1){550}}

 \put(1100,730){\makebox(100,100){$\kappa'_1$}}
 \put(900,700){\line(1,0){450}}
 \put(1450,700){\vector(1,0){100}}
\end{picture}
\end{center}
The above cube contains two special pushouts mentioned above. The following Lemma describes some other commutations.

\begin{lemma}\label{principal pushout3}
With the notation as above we have, for $k\geq 1$,
\begin{enumerate}
  \item \label{eq1} $\bio_{\bc^{(k+1)}P^*}\circ \bc^*_{\bc^{(k+1)}P}=\bbeta_k$,
  \item \label{eq2} $\kappa_1 \circ \bd^*_{\bc^{(k+1)}P}= (\bcomp_{\sfkj})\circ \bd^{(k),*}_P$,
  \item \label{eq3} $\bio_{\bc^{(k+1)}P^*}\circ \bd^*_{\bc^{(k+1)}P}=1_{\bd^{(k)}P^*}$,
  \item \label{eq4} $(\bcomp_{\sfkj})\circ \bmu_{\sfk} =  (\bmu_k+1_{\sfk P})\circ (\bcomp_k)$,
  \item \label{eq5} $  \bbeta_k\circ \bio_{\bc^{(k+1)}\fk P^*} =  \bio_{\bc^{(k+1)} P^*}\circ \bmu_k$,
  \item \label{eq6} $ [\bd^{(k),*}_{\sfk P}\circ \bio_{\bc^{(k+1)}\sfk P^*}   ,1_{\sfk P^*} ]\circ (\bcomp_k) = 1_{\sfk P^*}$,
  \item \label{eq7} $[ \bd^{(k),*}_P\circ \bio_{\bc^{(k+1)}P^*},\bmu_{\sfk}  ]\circ (\bcomp_{\sfkj})  =  1_{\fkj P^*}$.
\end{enumerate}

\end{lemma}

{\it Proof.}~ Routine check. $~~\Box$

\begin{lemma}\label{principal pushout4}
Let $F:(\pOpeCardo)^{op}\lra Set$ be a special pullback preserving functor, $P$ a positive opetope of dimension $n$. Then, for any $0\leq k < n$, $F$ preserves the pullback in
$(\pOpeCardo)^{op}$
\begin{center} \xext=800 \yext=700
\begin{picture}(\xext,\yext)(\xoff,\yoff)
 \setsqparms[1`-1`-1`1;800`500]
 \putsquare(0,70)[\bd^{(k)}P^*`\fkj P^*`\bc^{(k)}P^*`\fk P^*;
 \bd^{(k),*}_{\sfkj P}`\bbeta_k`\bmu_{\fk}`\bd^{(k),*}_{\sfk P}]
\end{picture}
\end{center}
\end{lemma}

{\it Proof.}~ Let $F:(\pOpeCardo)^{op}\lra Set$ be a special pullback preserving functor, $P$ a positive opetope of dimension $n$. We need to show that the square
\begin{center} \xext=1000 \yext=720
\begin{picture}(\xext,\yext)(\xoff,\yoff)
 \setsqparms[-1`1`1`-1;1000`500]
 \putsquare(0,120)[F(\bd^{(k)}P^*)`F(\fkj P^*)`F(\bc^{(k)}P^*)`F(\fk P^*);  F(\bd^{(k),*}_{\sfkj P})`F(\bbeta_k)`F(\bmu_{\fk})`F(\bd^{(k),*}_{\sfk P})]
\end{picture}
\end{center}
is a pullback in $Set$. Let us fix $t\in F(\bd^{(k)}P^*)$ and $a\in F(\fk P^*)$ such that
$$t\cdot \bbeta_k=a\cdot \bd^{(k),*}_{\sfk P}.$$
We will check that it is a pullback, by showing existence and uniqueness of an element $b\in F(\fkj P^*)$ such that
$$b=t\cdot \bd^{(k),*}_{\sfkj P} \hskip 5mm{\rm and}\hskip 5mm b=a\cdot \bmu_{\fk}.$$

{\em Existence.}  Put $t' =t\cdot \bio_{\bc^{(k+1)}P^*}$. By Lemma \ref{principal pushout3}.\ref{eq1}, we have
\[ t'\cdot \bc^*_{\bc^{(k+1)}P}= t\cdot \bio_{\bc^{(k+1)}P^*}\cdot \bc^*_{\bc^{(k+1)}P} =t\cdot \bbeta_k=a\cdot \bd^{(k),*}_{\sfk P}.\]
Since $F$ preserves special pullbacks and $\bc^{(k+1)} P+_{\bc^{(k)}P}\fk P$ is a special pullback in $(\pOpeCardo)^{op}$, we have an element
\[ \lk t',a \rk\in F(\bc^{(k+1)} P)\times_{F(\bc^{(k)}P)}F(\fk P)  \cong F(\bc^{(k+1)} P+_{\bc^{(k)}P}\fk P) \]
such that
\[ \lk t',a \rk\cdot \kappa_1=t'\hskip 5mm {\rm and} \hskip 5mm \lk t,a \rk\cdot\kappa_2=a.\]
We put $b= \lk t',a \rk \cdot \diamond_{\sfkj}\in F(\fkj P^*)$. We have

\[ b\cdot \bd^{(k),*}_{\sfkj P} =  ({\rm def\; of}\; b)\]
\[ = (\lk t' ,a \rk \cdot \diamond_{\sfkj}) \cdot \bd^{(k),*}_{\sfkj P} = (F {\rm\; presheaf}) \]
\[ = \lk t' ,a \rk \cdot (\diamond_{\sfkj} \circ \bd^{(k),*}_{\sfkj P}) = ({\rm Lemma}\; 14.4.2) \]
\[ = \lk t' ,a \rk \cdot (\kappa_1\circ \bd^*_{\bc^{(k+1)}P}) =  (F {\rm\; presheaf})\]
\[ = (\lk t',a \rk \cdot \kappa_1)\cdot \bd^*_{\bc^{(k+1)}P} = (F {\rm\; pres.\; special\; pb's})  \]
\[ = t'\cdot \bd^*_{\bc^{(k+1)}P} = (F {\rm\; pres.\; special\; pb's, \; def\;} t' )  \]
\[ = (t\cdot\bio_{\bc^{(k+1)}P^*})  \cdot \bd^*_{\bc^{(k+1)}P} = (F {\rm\; presheaf})\]
\[ = t\cdot(\bio_{\bc^{(k+1)}P^*}  \circ \bd^*_{\bc^{(k+1)}P}) = ({\rm Lemma}\; 14.4.3)\]
\[ = t\cdot\ 1_{\bd^{(k)}P^*}  = (F {\rm\; presheaf})\]
\[ =t \]
and
\[ b\cdot \bmu_{\fk}= ({\rm def\; of}\; b)\]
\[ =(\lk t',a\rk\cdot  \diamond_{\sfkj})\cdot \bmu_{\fk} = (F {\rm\; presheaf})\]
\[ =\lk t',a\rk\cdot  (\diamond_{\sfkj}\circ \bmu_{\fk}) = ({\rm Lemma}\; 14.4.4) \]
\[ =\lk t',a\rk\cdot  ((\bmu_k+1)\circ \diamond_{k}) = (F {\rm\; presheaf}) \]
\[ =(\lk t',a\rk\cdot  (\bmu_k+1))\cdot \diamond_{k} = (F {\rm\; pres.\; special\; pb's})  \]
\[ =\lk t'\cdot \bmu_k,a\rk \cdot \diamond_{k} = ({\rm def\;} t' ) \]
\[ =\lk (t\cdot\bio_{\bc^{(k+1)}P^*}) \cdot \bmu_k,a\rk \cdot \diamond_{k} = (F {\rm\; presheaf})  \]
\[ =\lk (t\cdot(\bio_{\bc^{(k+1)}P^*} \circ \bmu_k),a\rk \cdot \diamond_{k} = ({\rm Lemma}\; 14.4.5) \]
\[ =\lk (t\cdot(\bbeta_k\circ \bio_{\bc^{(k+1)}\fk P^*}),a\rk \cdot \diamond_{k} = (F {\rm\; presheaf}) \]
\[ =\lk (t\cdot\bbeta_k)\cdot \bio_{\bc^{(k+1)}\fk P^*},a\rk \cdot \diamond_{k} = ({\rm assumption\; on}\; a\; {\rm and}\; t)  \]
\[ =\lk (a\cdot \bd^{(k),*}_{\sfk P})\cdot \bio_{\bc^{(k+1)}\fk P^*},a\rk \cdot \diamond_{k} = (F {\rm\; presheaf}) \]
\[ =\lk a\cdot (\bd^{(k),*}_{\sfk P}\circ \bio_{\bc^{(k+1)}\fk P^*}),a\rk \cdot \diamond_{k} =  (F {\rm\; pres.\; special\; pb's})   \]
\[ =(a \cdot [\bd^{(k),*}_{\sfk P}\circ \bio_{\bc^{(k+1)}\sfk P^*}   ,1_{\sfk P^*} ]) \cdot \diamond_{k} = (F {\rm\; presheaf}) \]
\[ =a \cdot ([\bd^{(k),*}_{\sfk P}\circ \bio_{\bc^{(k+1)}\sfk P^*}   ,1_{\sfk P^*} ] \circ \diamond_{k}) = ({\rm Lemma}\; 14.4.6)  \]
\[ =a \cdot 1_{\sfk P^*} = (F {\rm\; presheaf}) \]
\[ =a\]

{\em Uniqueness.}  Now suppose that we have two elements $b,b'\in F(\fkj P^*)$ such that $\bd^{(k)}(b)=t=\bd^{(k)}(b')$ and $\bmu(b)=a=\bmu(b')$. Then, using Lemma \ref{principal pushout3} \ref{eq7} and the assumption, we have (we won't mention that we use the fact that $F$ is a sheaf anymore)

\[ b = ({\rm Lemma}\; 14.4.7) \]
\[ = b\cdot ([ \bd^{(k),*}_P\circ \bio_{\bc^{(k+1)}P^*},\bmu_{\sfk}  ]\circ (\bcomp_{\sfkj})) = \]
\[ = \lk (b\cdot  \bd^{(k),*}_P)\cdot \bio_{\bc^{(k+1)}P^*}, b\cdot \bmu_{\sfk}  \rk \cdot (\bcomp_{\sfkj}) = ({\rm assumption\; on}\; b\; {\rm and}\; b')\]
\[ = \lk (b'\cdot  \bd^{(k),*}_P)\cdot \bio_{\bc^{(k+1)}P^*}, b'\cdot \bmu_{\sfk}  \rk \cdot (\bcomp_{\sfkj}) = \]
\[ = b'\cdot ([ \bd^{(k),*}_P\circ \bio_{\bc^{(k+1)}P^*},\bmu_{\sfk}  ]\circ (\bcomp_{\sfkj})) =({\rm Lemma}\; 14.4.7) \]
\[ = b'.  \]
So the element with these properties is unique. $~~\Box$

\begin{lemma}\label{principal pushout5}
Let $F:(\pOpeCardo)^{op}\lra Set$ be a special pullback preserving functor, $P$ a positive opetope of dimension $n$. Then, for any $0\leq k < n$, $F$ preserves the following pullback in $(\pOpeCardo)^{op}$ \samepage{\begin{equation}\label{principal
pushout6}
\end{equation}
\begin{center} \xext=800 \yext=450
\begin{picture}(\xext,\yext)(\xoff,\yoff)
 \setsqparms[1`-1`-1`1;800`500]
 \putsquare(0,70)[\bd^{(k)}P^*`\fkj P^*`\alpha^{k,*}`\alpha^{n,*};
 \bd^{(k),*}_{\sfkj P}`\bmu^{\bd^{(k)}P^*}`\bmu^{\fkj P^*}`\bd^{(k),*}_{\alpha^n}]
\end{picture}
\end{center}}
\end{lemma}

{\it Proof.}~The proof is by double induction, on dimension $n$ of the positive opetopic $P$, and $k<n$.

Note that if $k=0$, then, for any $n>0$,  the vertical arrows in (\ref{principal pushout6}) are isomorphisms, so any functor from $(\pOpeCardo)^{op}$ sends (\ref{principal pushout6}) to a pullback. This shows in particular that the Lemma holds for $n=1$. As we already mentioned, if $k=n-1$, the square (\ref{principal pushout6}) is an arbitrary special pushout.

Thus, we assume that $F:(\pOpeCardo)^{op}\lra Set$ is a special pullback preserving functor, and that $P$ is a positive opetope of dimension $n$, $0\leq k<n$, $F$ preserves the pullback (\ref{principal pushout6}). Moreover, for $m<n$ and the positive opetope $Q$ of dimension $m$, $F$ preserves the principal pullback in $(\pOpeCardo)^{op}$:
\begin{center} \xext=2600 \yext=700
\begin{picture}(\xext,\yext)(\xoff,\yoff)
 \setsqparms[1`-1`-1`1;700`500]
 \putsquare(0,50)[\bd Q^*`Q^*`\alpha^{m-1,*}`\alpha^{m,*};
 \bd_Q^*`\bmu^{\bd Q}`\bmu^Q`\bd_{\alpha^{m}}^*]

 \setsqparms[1`-1`-1`1;900`500]
 \putsquare(1700,70)[\bd^{(m-1)}Q^*`\mbox{\tiny{\frame{$m$}}} Q^*` \alpha^{m-1,*}`\alpha^{m,*}; \bd^{(m-1),*}_{\mbox{\tiny{\frame{$m$}}} Q}`\bmu^{\bd^{(m-1)}Q^*}`
 \bmu^{\mbox{\tiny{\frame{$_{m}$}}} Q^*}`\bd^{(m-1),*}_{\alpha^{m}}]
 \put(1000,300){\makebox(100,100){$=$}}
\end{picture}
\end{center}
We shall show that $F$ preserves the pullback
\begin{equation}\label{principal pushout7}
\end{equation}
\begin{center} \xext=800 \yext=450
\begin{picture}(\xext,\yext)(\xoff,\yoff)
 \setsqparms[1`-1`-1`1;800`500]
 \putsquare(0,70)[\bd^{(k+1)}P^*`\fkd P^*`\alpha^{k+1,*}`\alpha^{n,*}; \bd^{(k+1),*}_{\sfkd P}`\bmu^{\bd^{(k+1)}P^*}`\bmu^{\fkd P^*}`\bd^{(k+1),*}_{\alpha^n}]
\end{picture}
\end{center}
in $(\pOpeCardo)^{op}$, as well.

In the following diagram (most of the subscripts and some superscripts were suppressed for clarity):
\begin{center} \xext=2000 \yext=1900
\begin{picture}(\xext,\yext)(\xoff,\yoff)

\put(1500,1500){\makebox(100,100){$I$}}
\put(1500,600){\makebox(100,100){$II$}}
\put(500,600){\makebox(100,100){$III$}}

 \setsqparms[0`-1`0`1;1000`500]
 \putsquare(1000,300)[\bd^{(k)}P^*``\alpha^{k,*}`\alpha^{n,*}; `\bmu``]
 \put(1500,350){\makebox(100,100){$\bd^{(k),*}$}}

 \putmorphism(1230,950)(2,1)[``]{480}{1}a
 \put(1500,950){\makebox(100,100){$\bd^{(k),*}$}}
 \putmorphism(1000,1300)(0,-1)[\phantom{\bc^{(k+1)}P^*}`\phantom{\bd^{(k)}P^*}`\bd^*]{500}{-1}l
 \putmorphism(2000,1300)(0,-1)[\phantom{\fkj P^*}`\phantom{\alpha^{n,*}}`\bmu]{1000}{-1}r
 \putmorphism(0,300)(1,0)[\alpha^{k+1,*}`\phantom{\alpha^{k,*}}`\bd^*]{1000}{-1}a

 \put(80,380){\line(1,2){450}}
 \put(530,1280){\vector(1,0){250}}
 \put(400,950){\makebox(100,100){$\bmu$}}

 \put(1100,40){\makebox(100,100){$\bd^{(k+1),*}$}}
 \put(0,0){\line(0,1){250}}
 \put(0,0){\line(1,0){2000}}
 \put(2000,0){\vector(0,1){250}}
 \put(0,380){\line(0,1){1420}}
 \put(0,1800){\vector(1,0){700}}
 \put(40,1250){\makebox(100,100){$\bmu$}}

  \setsqparms[1`-1`-1`1;1000`500]
 \putsquare(1000,1300)[\bd^{(k+1)}P^*`\fkd P^*`\bc^{(k+1)}P^*`\fkj P^*;
 \bd^{(k+1),*}`\bbeta_{k+1}`\bmu`]
  \put(1500,1320){\makebox(100,100){$\bd^{(k+1),*}$}}
\end{picture}
\end{center}
all the squares and triangles commute. Moreover, $F$ sends the squares $I$, $II$, $III$ to pullbacks in $Set$: $I$ by Lemma \ref{principal pushout4}, $II$ by inductive hypothesis for $k$, $III$ by inductive hypothesis since $dim(\bc^{(k+1)}P)<n$.

Let $f: X\lra F(\bd^{(k+1)}P^*)$ and $g: X\lra F(\alpha^{n,*})$ be functions such that
\[ F(\bmu^{\bd^{(k+1)}P^*})\circ f = F(\bd^{(k+1),*}_{\alpha^n})\circ g. \]
Since $F$ applied to $II$ is a pullback in $Set$, and all squares and triangles in the above diagram commute, there is a unique function $h_1:X\lra F(\fkj P^*)$ such that
\begin{equation}\label{h1-cond}
F(d^{(k),*}_{\fkj P})\circ  h_1=F(\bd^*_{\bc^{(k+1)}P})\circ F(\bbeta_{k+1})\circ f \hskip 5mm \mathrm{ and }\hskip 5mm F(\bmu^{\sfkj P^*})\circ  h_1= g.
\end{equation}
To get a unique function $h_2: X\lra F(\fkd P^*)$ such that
\begin{equation}\label{h2-cond}
 F(\bd^{k+1,*}_{\sfkd P})\circ  h_2 = f \hskip 5mm \mathrm{ and }\hskip 5mm   F(\bmu_{\sfkj})\circ h_2= h_1,
\end{equation}
we use the fact that $F$ sends $III$ to a pullback in $Set$. The application of $F$ to the diagram above will give the following diagram in $Set$, where we added the additional functions $f$,
$g$, $h_1$, and $h_2$:
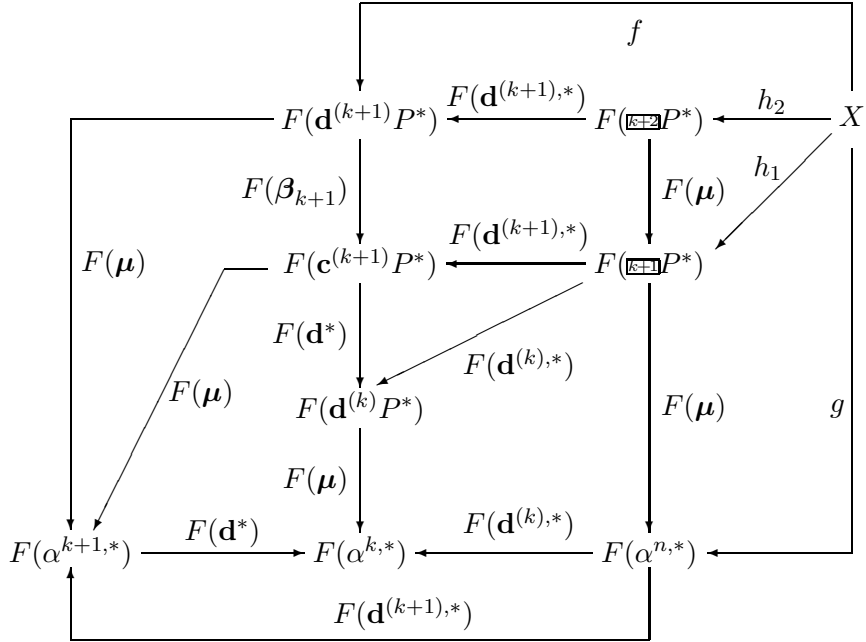
\begin{center} \xext=2500 \yext=2200
\begin{picture}(\xext,\yext)(\xoff,\yoff)
 \setsqparms[0`1`0`-1;1000`500]
 \putsquare(1000,300)[F(\bd^{(k)}P^*)``F(\alpha^{k,*})`F(\alpha^{n,*});
 `F(\bmu)``]
 \put(1500,350){\makebox(100,100){$F(\bd^{(k),*})$}}

 \putmorphism(1200,950)(2,1)[``]{430}{-1}a
 \put(1500,900){\makebox(100,100){$F(\bd^{(k),*})$}}
 \putmorphism(1000,1300)(0,-1)[\phantom{F(\bc^{(k+1)}P^*)}`\phantom{F(\bd^{(k)}P^*)}`F(\bd^*)]{500}{1}l
 \putmorphism(2000,1300)(0,-1)[\phantom{F(\fkj P^*)}`\phantom{F(\alpha^{n,*})}`F(\bmu)]{1000}{1}r
 \putmorphism(0,300)(1,0)[F(\alpha^{k+1,*})`\phantom{F(\alpha^{k,*})}`F(\bd^*)]{1000}{1}a

 \put(1100,40){\makebox(100,100){$F(\bd^{(k+1),*})$}}
 \put(2000,0){\line(0,1){250}}
 \put(0,0){\line(1,0){2000}}
 \put(0,0){\vector(0,1){250}}
 \put(0,1800){\vector(0,-1){1420}}
 \put(0,1800){\line(1,0){700}}
 \put(100,1250){\makebox(100,100){$F(\bmu)$}}

  \put(530,1280){\vector(-1,-2){450}}
 \put(530,1280){\line(1,0){150}}
 \put(400,800){\makebox(100,100){$F(\bmu)$}}

 \putmorphism(2000,1800)(1,0)[\phantom{F(\fkd P^*)}`X`h_2]{700}{-1}a
 \put(2630,1750){\vector(-1,-1){400}}
 \put(2360,1580){\makebox(100,100){$h_1$}}

 \put(2700,1700){\line(0,-1){1400}}
 \put(2600,750){\makebox(100,100){$g$}}
 \put(2700,300){\vector(-1,0){500}}

 \put(2700,1900){\line(0,1){300}}
 \put(2700,2200){\line(-1,0){1700}}
 \put(1900,2050){\makebox(100,100){$f$}}
 \put(1000,2200){\vector(0,-1){300}}

  \setsqparms[-1`1`1`-1;1000`500]
 \putsquare(1000,1300)[F(\bd^{(k+1)}P^*)`F(\fkd P^*)`F(\bc^{(k+1)}P^*)`F(\fkj P^*);
 F(\bd^{(k+1),*})`F(\bbeta_{k+1})`F(\bmu)`]
  \put(1500,1350){\makebox(100,100){$F(\bd^{(k+1),*})$}}
\end{picture}
\end{center}

Thus in order to verify that
 \[ F(\bbeta_{k+1})\circ f=F(\bd^{(k+1),*}_{\fkj P})\circ h_1\]
and to get $h_2$ satisfying (\ref{h2-cond}), it is enough to verify that
\begin{equation}\label{h2-III-1}
  F(\bmu^{\bc^{(k+1)}P^*})\circ F(\bbeta_{k+1})\circ f  =F(\bmu^{\bc^{(k+1)}P^*})\circ F(\bd^{(k+1),*}_{\fkj P}) \circ h_1
\end{equation}
and
\begin{equation}\label{h2-III-2}
  F(\bd^*_{\bc^{(k+1)}P}) \circ F(\bbeta_{k+1})\circ f= F(\bd^*_{\bc^{(k+1)}P}) \circ F(\bd^{(k+1),*}_{\sfkj P}) \circ h_1
\end{equation}

For (\ref{h2-III-1}), we have
\[ F(\bmu^{\bc^{(k+1)}P^*}) \circ F(\bbeta_{k+1})\circ f  = \]
\[ = F(\bmu^{\bd^{(k+1)}P^*})\circ f= \]
\[ = F(\bd^{(k+1),*}_{\alpha^n}) \circ g = \]
\[ = F(\bd^{(k+1),*}_{\alpha^n}) \circ F(\bmu^{\sfkj P^*}) \circ h_1   = \]
\[ =F(\bmu^{\bc^{(k+1)}P^*})\circ  F(\bd^{(k+1)^*}_{\sfkj P})\circ h_1, \]
and for (\ref{h2-III-2}), we have
\[  F(\bd^*_{\bc^{(k+1)}P}) \circ F(\bbeta_{k+1}) \circ f = \]
\[ =F(d^{(k),*}_{\fkj P})\circ h_1 = \]
\[ = F(\bd^*_{\bc^{(k+1)}P})\circ F(\bd^{(k+1),*}_{\sfkj P})\circ h_1. \]
By uniqueness of both $h_1$ and $h_2$, $h_2$ is the unique function such that
\[ \bd^{(k+1),*}_{\sfkd P}\circ h_2 = f \hskip 5mm \mathrm{ and }\hskip 5mm  \bmu^{\sfkd P^*} \circ h_2 = g,\]
i.e., $F$ sends (\ref{principal pushout7}) to a pullback in $Set$, as required. $~~\Box$

\section{A full nerve of $\o$-categories}\label{sec-nerve}

Let $\cS$ denote the category of simple $\o$-categories, described in \cite{MZ}. It was proved there that any simple $\o$-category is isomorphic to one of form $(\alpha^{\vec{u}})^*$  for some ud-vector $\vec{u}$. In fact what we need here is that any simple $\o$-category can be obtained from those of form $(\alpha^n)^*$, with $n\in\o$, via special pushouts. For more details the reader should consult \cite{MZ}.

As every simple $\o$-category is a positive opetopic cardinal (considered as an $\o$-category), we have a full inclusion functor
\[ \bk : \cS \lra \pOpeCardo\]
whose essential image is spanned by the opetopic cardinals whose all faces are globular.

In \cite{MZ} we have shown that $sPb(\cS^{op},Set)$, the category of special pullbacks preserving functors from the dual of $\cS$ to $Set$, is equivalent to the category $\o$-categories. We have in fact an adjoint equivalence

\begin{center} \xext=2000 \yext=300
\begin{picture}(\xext,\yext)(\xoff,\yoff)
\putmorphism(0,200)(1,0)[\phantom{\oC}`
\phantom{sPb(\cS^{op},Set)}`\widetilde{(-)}]{2000}{-1}a
\putmorphism(0,100)(1,0)[\phantom{\oC}`
\phantom{sPb(\cS^{op},Set)}`\widehat{(-)}=\oC(\simeq, -
)]{2000}{1}b
\putmorphism(0,150)(1,0)[\oC`sPb(\cS^{op},Set)`]{2000}{0}b
\end{picture}
\end{center}
where
\[ \widehat{C}: \cS^{op} \lra Set \]
is given by
\[  \widehat{C}(A) =\oC(A,C),\]
where $A$ is a simple $\o$-category. We have

\begin{proposition}\label{sp-nevre1}
The adjunction
\begin{center} \xext=1500 \yext=300
\begin{picture}(\xext,\yext)(\xoff,\yoff)
 \putmorphism(0,200)(1,0)[\phantom{Set^{\cS^{op}}}`\phantom{Set^{(\pOpeCardo)^{op}}}`Ran_\bk]{1500}{1}a

 \putmorphism(0,100)(1,0)[\phantom{Set^{\cS^{op}}}`\phantom{Set^{(\pOpeCardo)^{op}}}`\bk^*]{1500}{-1}b
\putmorphism(0,150)(1,0)[Set^{\cS^{op}}`Set^{(\pOpeCardo)^{op}}`]{1500}{0}b
\end{picture}
\end{center}
restricts to an equivalence of categories.
\begin{center} \xext=2000 \yext=300
\begin{picture}(\xext,\yext)(\xoff,\yoff)
 \putmorphism(0,200)(1,0)[\phantom{sPb(\cS^{op},Set)}`\phantom{sPb((\pOpeCardo)^{op},Set)}`Ran_\bk]{2000}{1}a
 \putmorphism(0,100)(1,0)[\phantom{sPb(\cS^{op},Set)}`\phantom{sPb((\pOpeCardo)^{op},Set)}`\bk^*]{2000}{-1}b
\putmorphism(0,150)(1,0)[sPb(\cS^{op},Set)`sPb((\pOpeCardo)^{op},Set)`]{2000}{0}b
\end{picture}
\end{center}
where $sPb((\pOpeCardo)^{op},Set)$ is the category of the special pullbacks preserving functors.
\end{proposition}

{\it Proof.}~ First we shall describe the adjunction in  details.

The counit. Let $G$ be a functor in $sPb(\cS^{op},Set)$ , $A$ a simple $\o$-category. We have a functor
\begin{center} \xext=1100 \yext=200
\begin{picture}(\xext,\yext)(\xoff,\yoff)
 \putmorphism(0,50)(1,0)[(\bk\da A)^{op}`\cS^{op}`\pi^A]{600}{1}a
 \putmorphism(600,50)(1,0)[\phantom{\cS^{op}}`Set`G]{500}{1}a
\end{picture}
\end{center}
with the limit, say $(Lim (G\circ \pi^A),\sigma^A)$. Then the counit $(\varepsilon_G)_A$ is
\begin{center} \xext=700 \yext=200
\begin{picture}(\xext,\yext)(\xoff,\yoff)
 \putmorphism(0,50)(1,0)[(\varepsilon_G)_A:(Ran_\bk(G)\circ \bk)(A)=Lim (G\circ \pi^A)`G(A)`\sigma^A_{1_A}. ]{1500}{1}a
\end{picture}
\end{center}

As $\bk$ is full and faithful\footnote{This conditions translates to the fact that $1_A$ is the initial object in $(\bk\da A)^{op}$ and therefore that we have an iso $\sigma^A_{1_A}: Lim G\circ \pi^A\cong G\circ \pi^A(1_A)=G(A)$.}, for any $G$, $\varepsilon_G$ is an iso. Thus $\varepsilon$ is an iso.

The unit. Let $F$ be a functor in $sPb((\pOpeCardo)^{op},Set)$, $T$ a positive opetopic cardinal. We have a functor
\begin{center} \xext=1100 \yext=200
\begin{picture}(\xext,\yext)(\xoff,\yoff)
 \putmorphism(0,50)(1,0)[(\bk\da T^*)^{op}`\cS^{op}`\pi^{T^*}]{600}{1}a
 \putmorphism(600,50)(1,0)[\phantom{\cS^{op}}`\pOpeCardo`\bk]{500}{1}a
 \putmorphism(1100,50)(1,0)[\phantom{\pOpeCardo}`Set`F]{500}{1}a
\end{picture}
\end{center}
with the limit, say $(Lim F\circ\bk\circ \pi^{T^*},\sigma^{T^*})$. Then the unit $(\eta_F)_{T^*}$ is the unique morphism into the limit:
\begin{center} \xext=2000 \yext=700
\begin{picture}(\xext,\yext)(\xoff,\yoff)
 \putmorphism(0,600)(1,0)[F(T^*)`Ran_\bk(F\circ \bk)(T^*)`(\eta_F)_{T^*}]{1800}{1}a
 \put(1800,460){\makebox(100,100){$\parallel$}}
 \settriparms[1`1`1-;400]
 \putAtriangle(1400,0)[Lim (F\circ\bk\circ \pi^{T^*})`F(A)`F(B);\sigma_a^T`\sigma_b^T`F(f)]
 \putmorphism(0,500)(3,-1)[``F(a)]{1500}{1}l
 \put(250,500){\line(4,-1){1100}}
 \putmorphism(1400,200)(4,-1)[``]{900}{1}l
 \put(800,380){\makebox(100,100){$F(b)$}}
\end{picture}
\end{center}
where the triangle in $\pOpeCardo$
\begin{center} \xext=600 \yext=350
\begin{picture}(\xext,\yext)(\xoff,\yoff)
 \settriparms[-1`-1`-1;300]
 \putAtriangle(0,0)[T^*`A`B;a`b`f]
\end{picture}
\end{center}
commutes.

Note that, as $F$ preserves special pullbacks, and any simple $\o$-category can be  obtained from those of form $\alpha^n$ with $n\in\o$, we can restrict the limiting cone $(Lim F\circ\bk\circ \pi^{T^*},\sigma^{T^*})$ to the objects of form $\alpha^n$, with $n\in\o$.

After this observation we shall prove, by induction on size of a positive opetopic cardinal $T$, that $(\eta_F)_{T^*}$ is an iso.

If $dim(T)\leq 1$, then $(\eta_F)_{T^*}$ is obviously an iso.

Suppose $T$ is not principal, i.e., we have $a\in Sd(T)$, for some $k\in\o$. By inductive hypothesis the morphisms
\[ (\eta_F)_{(T^{\da a})^*},\;\;\; (\eta_F)_{\bc^{(k)}(T^{\da a})^*},\;\;\;(\eta_F)_{(T^{\ua a})^*}\]
are isos, and the square
\begin{center}
\xext=650 \yext=500
\begin{picture}(\xext,\yext)(\xoff,\yoff)
 \setsqparms[1`-1`-1`1;650`450]
 \putsquare(0,50)[T^{\da a}`T`\bc^{(k)}(T^{\da a})`T^{\ua a};```]
\end{picture}
\end{center}
is a special pushout (see Proposition \ref{o-cat-comp}) which is sent by $F$ to a pullback. Hence the morphism
\[ (\eta_F)_{T^*}= (\eta_F)_{(T^{\da a})^*}\times (\eta_F)_{(T^{\ua a})^*} \]
 is indeed an iso in this case.

If $T$ is principal and $T=(\alpha^n)^*$, then the category $(\bk\da(\alpha^n)^*)^{op}$ has the initial object $1_{(\alpha^n)^*}$, so morphism
\[ (\eta_F)_{(\alpha^n)^*} :F((\alpha^n)^*) \lra Ran_\bk(F\circ \bk)((\alpha^n)^*)\]
 is an iso.

 Finally, let assume that $T(=P)$ is any  positive opetope of dimension $n$. Thus, by Corollary \ref{principal pushout1}, we have a principal pushout
\begin{center} \xext=700 \yext=650
\begin{picture}(\xext,\yext)(\xoff,\yoff)
 \setsqparms[1`-1`-1`1;700`500]
 \putsquare(0,50)[(\bd P)^*`P^*`(\alpha^{n-1})^*`(\alpha^{n})^*;
 \bd_P^*`\bmu^{\bd P}`\bmu^P`\bd_{\alpha^n}^*]
\end{picture}
\end{center}
which, by Theorem \ref{principal pushout2}, is preserved by $F$. By induction hypothesis the morphisms
\[ (\eta_F)_{(T^{\da a})^*},\;\;\; (\eta_F)_{\bc^{(k)}(T^{\da a})^*},\;\;\;(\eta_F)_{(T^{\ua a})^*}\]
are isos, so we have that the morphism
\[ (\eta_F)_{P^*}=(\eta_F)_{(\bd P)^*}\times
(\eta_F)_{(\alpha^n)^*}\] is an iso, as well. $~~\Box$

\begin{corollary}\label{sp-nevre2} We have a  commuting triangle of adjoint equivalences
\begin{center} \xext=2000 \yext=850
\begin{picture}(\xext,\yext)(\xoff,\yoff)

 \put(0,800){\makebox(500,100){$\oC$}}
 \put(0,100){\makebox(500,100){$sPb((\pOpeCardo)^{op},Set)$}}
 \put(1500,100){\makebox(500,100){$sPb(\cS^{op},Set) $}}

 \putmorphism(200,850)(0,1)[``\widehat{(-)}]{650}{1}l
 \putmorphism(300,850)(0,1)[``\widetilde{(-)}]{650}{-1}r

 \putmorphism(200,850)(2,-1)[``\widehat{(-)}]{1400}{1}l
 \putmorphism(400,850)(2,-1)[``\widetilde{(-)}]{1400}{-1}r

 \putmorphism(790,200)(1,0)[``\bk^*]{610}{1}a
 \putmorphism(790,120)(1,0)[``Ran_\bk ]{610}{-1}r
\end{picture}
\end{center}
In particular, the categories $\oC $ and  $sPb((\pOpeCardo)^{op},Set) $ are equivalent.
\end{corollary}

{\it Proof.}~ It is enough to show that in the above diagram $\bk^*\circ\widehat{(-)}=\widehat{(-)} $. But this is clear.  $~~\Box$

\section{A monadic adjunction}\label{sec_adjunction}
In this section we show that  the inclusion functor $\be :\pPoly \lra\oC$ has a right adjoint which is monadic.

First we will give an outline of the proof. Consider the following diagram of categories and functors
\begin{center} \xext=2400 \yext=850
\begin{picture}(\xext,\yext)(\xoff,\yoff)
 \put(1800,750){\makebox(500,100){$\oC $}}
 \put(1800,50){\makebox(500,100){$sPb((\pOpeCardo)^{op},Set) $}}
 \put(0,750){\makebox(500,100){$\pPoly $}}
 \put(0,50){\makebox(500,100){$sPb((\pOpeCard)^{op},Set) $}}
 \putmorphism(200,800)(0,1)[``\widehat{(-)}]{650}{1}l
 \putmorphism(300,800)(0,1)[``\widetilde{(-)}]{650}{-1}r
 \putmorphism(2000,800)(0,1)[``\widehat{(-)}]{650}{1}l
 \putmorphism(2100,800)(0,1)[``\widetilde{(-)}]{650}{-1}r
 \putmorphism(810,120)(1,0)[``Lan_\bj]{650}{1}a
 \putmorphism(810,50)(1,0)[``\bj^*]{650}{-1}b
 \putmorphism(420,780)(1,0)[``\be]{1450}{1}a
\end{picture}
\end{center}
where $\be$ is just an inclusion of positive-to-one polygraphs into $\o$-categories and $\bj=(-)^*:\pOpeCard \lra \pOpeCardo$ is an essentially surjective embedding. We have already shown (Proposition \ref{presheaf_eq2}, Corollary \ref{sp-nevre2}) that the vertical functors constitute two adjoint equivalences. The proof that $\be$ has a right adjoint takes few steps. We begin by presenting $Lan_\bj$ as a family representable functor (or local right adjoint). Then we check that $Lan_\bj$ is well defined, i.e.,  we will check that the functor $Lan_\bj(F)$, the left Kan extension of special pullbacks preserving functor $F$, preserves special pullbacks. Next we shall check that the above square commutes, i.e., $\widehat{(-)}\circ \be = Lan_\bj\circ\widehat{(-)}$. This will reduce the problem of monadicity of $\oC$ over $\widehat{\pOpe}$  to verification whether $\bj^*$, the left adjoint to $Lan_\bj$, is monadic. The last statement is verified directly checking assumptions of Beck's monadicity theorem.

We describe the left Kan extension along $\bj$ in a more convenient than the usual way, c.f. \cite{CWM}, as a family representable functor. We have

\begin{proposition}\label{Lan_j-description}
The functor of left Kan extension
\[ Lan_\bj : \widehat{\pOpeCard}  \lra \widehat{\pOpeCardo} \]
along the functor
\[ \bj: \pOpeCard\ra \pOpeCardo \]
is defined, for $F\in \widehat{\pOpeCard}$, as follows. For a
positive opetopic cardinal $S$, we have
\begin{center} \xext=1000 \yext=150
\begin{picture}(\xext,\yext)(\xoff,\yoff)
 \putmorphism(0,50)(1,0)[Lan_\bj(F)(S^*)=\coprod_{a:S^*\ra T^*\;
 {\rm inner}}\, F(T)`F(T)`\kappa^{S^*}_a]{1500}{-1}a
\end{picture}
\end{center}
where coproduct is taken over all up to iso inner maps in $\pOpeCardo$ with the domain $S^*$, with the coprojections as shown.

If $h:S_1^*\lra S_2^*$ is an $\o$-functor and $a_2:S^*_2\ra T^*_2$ is inner, by Lemma \ref{factorization}, we can form a diagram
\begin{center} \xext=500 \yext=480
\begin{picture}(\xext,\yext)(\xoff,\yoff)
 \setsqparms[1`1`1`1;500`400]
 \putsquare(0,30)[S^*_1`S^*_2`T^*_1`T^*_2;
 h`a_1`a_2`(h')^*]
\end{picture}
\end{center}
with $a_1$ inner and $h'$ a map of positive opetopic cardinals, i.e., the map $(h')^*$ is an outer map. $Lan_\bj(h)$ is so defined that, for any $h$, $h'$, $a_1$, $a_2$ as above, the diagram
\begin{center} \xext=1500 \yext=700
\begin{picture}(\xext,\yext)(\xoff,\yoff)
 \setsqparms[-1`1`1`-1;1500`500]
 \putsquare(0,50)[ Lan_\bj(F)(S_2^*)=\coprod_{a_2:S_2^*\ra T_2^*\,{\rm inner}}\, F(T_2)`F(T_2)` Lan_\bj(F)(S_1^*)=\coprod_{a_1:S_1^*\ra T_1^*\,{\rm inner}}\, F(T_1)`F(T_1); \kappa^{S_2^*}_{a_2}`Lan_\bj(F)(h)`F(h')`\kappa^{S_1^*}_{a_1}]
\end{picture}
\end{center}
commutes.
\end{proposition}
{\it Proof.}~ Fix $F$ in $sPb((\pOpeCardo)^{op},Set)$ for the whole proof.

Let $S$ be a positive opetopic cardinal. Then $Lan_\bj(F)(S)$ is the colimit of the following functor
\begin{center} \xext=1300 \yext=200
\begin{picture}(\xext,\yext)(\xoff,\yoff)
 \putmorphism(0,50)(1,0)[\bj^{op}\da S`(\pOpeCard)^{op}`\pi^S]{700}{1}a
 \putmorphism(700,50)(1,0)[\phantom{(\pOpeCard)^{op}}`Set`F]{600}{1}a
\end{picture}
\end{center}
i.e., $Lan_\bj(F)(S)=(F\circ \pi^S,\sigma^F)$.  A map $f:a\lra b$ in $\bj^{op}\da S$, is a commuting triangle
\begin{center} \xext=600 \yext=350
\begin{picture}(\xext,\yext)(\xoff,\yoff)
 \settriparms[1`1`1;300]
 \putAtriangle(0,0)[S^*`T_1^*`T_2^*;a`b`f^*]
\end{picture}
\end{center}
in $\pOpeCardo$, and hence by Lemma \ref{factorization} we can take the inner-outer factorizations, c.f. \ref{sec_inner-outer}, of both $a=(a'')^*\circ a'$ and $b=(b'')^*\circ b'$, with $a'$ and $b'$ inner maps. Then, again by Lemma\ref{factorization}, there is a morphism  $f':a'\lra b'$ in $\pOpeCard$ which must be an iso.  In this way we get a commuting diagram
\begin{center} \xext=1400 \yext=800
\begin{picture}(\xext,\yext)(\xoff,\yoff)
 \setsqparms[0`1`1`1;1400`400]
 \putsquare(0,50)[T_3^*`T_4^*`T_1^*`T_2^*;
 `(a'')^*`(b'')^*`f^*]
 \settriparms[1`1`0;700]
 \putAtriangle(0,50)[S^*`\phantom{T_1^*}`\phantom{T_2^*};a`b`]
 \putmorphism(250,600)(3,1)[``]{160}{-1}l
 \put(250,650){\makebox(100,100){$a'$}}
 \putmorphism(1150,600)(-3,1)[``]{160}{-1}r
 \put(1050,650){\makebox(100,100){$b'$}}
 \put(100,500){\line(1,0){250}}
 \put(500,500){\line(1,0){400}}
 \putmorphism(1000,500)(1,0)[``]{350}{1}a
 \put(700,350){\makebox(100,100){$(f')^*$}}
\end{picture}
\end{center}
in $\pOpeCardo$, which corresponds to the following part of the colimiting cocone
\begin{center} \xext=1400 \yext=900
\begin{picture}(\xext,\yext)(\xoff,\yoff)
 \setsqparms[0`-1`-1`-1;1400`400]
 \putsquare(0,50)[F(T_3)`F(T_4)`F(T_1)`F(T_2);
 `F(a'')`F(b'')`F(f)]
 \settriparms[-1`-1`0;700]
 \putAtriangle(0,50)[\phantom{Lan_\bj(F)(S^*)}`\phantom{F(T_1)}`\phantom{F(T_2};\sigma^F_a`\sigma^F_b`]
 \put(700,740){\makebox(100,100){$Lan_\bj(F)(S^*)$}}
 \putmorphism(250,600)(3,1)[``]{160}{1}l
 \put(250,640){\makebox(100,100){$\sigma^F_{a'}$}}
 \putmorphism(1150,600)(-3,1)[``]{160}{1}r
 \put(1050,640){\makebox(100,100){$\sigma^F_{b'}$}}
 \put(400,500){\vector(-1,0){250}}
 \put(500,500){\line(1,0){400}}
 \put(1000,500){\line(1,0){250}}
 \put(700,380){\makebox(100,100){$F(f')$}}
\end{picture}
\end{center}
Thus if there is a morphism $f:a\ra b$ between two objects in $\bj^{op}\da S$ , we have a commuting diagram
\begin{center} \xext=600 \yext=350
\begin{picture}(\xext,\yext)(\xoff,\yoff)
 \settriparms[1`1`1;300]
 \putAtriangle(0,50)[a'`a`b;a''`b''\circ f`f]
\end{picture}
\end{center}
in $\bj^{op}\da S$ with $a'$ being the inner part of both $a$ and $b$. There are no other comparison maps between these objects. But this say that in fact
\begin{center} \xext=1000 \yext=150
\begin{picture}(\xext,\yext)(\xoff,\yoff)
 \putmorphism(0,50)(1,0)[Lan_\bj(F)(S^*)=\coprod_{a:S^*\ra T^*\;{\rm inner}}\, F(T),`F(T)`\kappa^{S^*}_a]{1500}{-1}a
\end{picture}
\end{center}
where the coproduct is taken over all (up to iso) inner maps with the domain $S^*$, with the coprojections as shown.

To define $Lan_\bj(F)$ on morphisms, fix an $\o$-functor $h:S_1^*\lra S_2^*$ in $\pOpeCard$  and an inner map $a_2:S^*_2\ra T^*_2$. By Lemma \ref{factorization}, we can form a diagram
\begin{center} \xext=500 \yext=550
\begin{picture}(\xext,\yext)(\xoff,\yoff)
 \setsqparms[1`1`1`1;500`400]
 \putsquare(0,50)[S^*_1`S^*_2`T^*_1`T^*_2;
 h`a_1`a_2`(h')^*]
\end{picture}
\end{center}
with $a_1$ inner and $(h')^*$ outer. $Lan_\bj(h)$ is so defined that, for any $h'$, $a_1$, $a_2$ as above, the diagram
\begin{center} \xext=1500 \yext=650
\begin{picture}(\xext,\yext)(\xoff,\yoff)
 \setsqparms[-1`1`1`-1;1500`500]
 \putsquare(0,50)[
 Lan_\bj(F)(S_2^*)=\coprod_{a_2:S_2^*\ra T_2^*\,{\rm inner}}\, F(T_2)`F(T_2)`
 Lan_\bj(F)(S_1^*)=\coprod_{a_1:S_1^*\ra T_1^*\,{\rm inner}}\, F(T_1)`F(T_1);
 \kappa^{S_2^*}_{a_2}`Lan_\bj(F)(h)`F(h')`\kappa^{S_1^*}_{a_1}]
\end{picture}
\end{center}
commutes. This shows that the functor $Lan_\bj$ is a family representable functor.  $~~\Box$

\begin{lemma}\label{Lan j}
The functor of the left Kan extension along $\bj$ restricts to
\[ Lan_\bj : sPb((\pOpeCard)^{op},Set)  \lra sPb((\pOpeCardo)^{op},Set) \]
the above categories, i.e., whenever $F :(\pOpeCard)^{op}\lra Set$ preserves special pullbacks, so does $Lan_\bj(F): (\pOpeCardo)^{op}\lra Set$. Moreover, $Lan_\bj$ is the left adjoint to
 \[ \bj^*: sPb((\pOpeCardo)^{op},Set) \lra  sPb((\pOpeCard)^{op},Set).  \]
\end{lemma}

{\it Proof.}~Note that once we will prove the first part of the statement, the part following `moreover' will follow immediately.

Fix $F$ in $sPb((\pOpeCard)^{op},Set)$ for the whole proof. We shall use the description $Lan_\bj(F)$ from Proposition \ref{Lan_j-description} to show that $Lan_\bj(F)$ preserves special pushouts.

So assume that $S_1$ and $S_2$ are positive opetopic cardinals such that
$$\bc^{(k)}(S_1)=\bd^{(k)}(S_2),$$
i.e., we have a pushout
\begin{center} \xext=500 \yext=600
\begin{picture}(\xext,\yext)(\xoff,\yoff)
 \setsqparms[1`-1`-1`1;650`450]
 \putsquare(0,100)[S_1`S_1{\oplus}_k S_2`\bc^{(k)}(S_1)`S_2;
 \kappa_1`\bc^{(k)}_{S_1}`\kappa_2`\bd^{(k)}_{S_2}]
\end{picture}
\end{center}
in $\pOpeCard$. We need to show that the square
\begin{center} \xext=1500 \yext=850
\begin{picture}(\xext,\yext)(\xoff,\yoff)
 \setsqparms[-1`1`1`-1;1500`600]
 \putsquare(0,100)[Lan_\bj(F)(S_1)`Lan_\bj(F)(S_1{\oplus}_k S_2)`Lan_\bj(F)(\bc^{(k)}(S_1))`Lan_\bj(F)(S_2);
 Lan_\bj(F)(\kappa_1)`Lan_\bj(F)(\bc^{(k)}_{S_1})`Lan_\bj(F)(\kappa_2)`Lan_\bj(F)(\bd^{(k)}_{S_2})]
\end{picture}
\end{center}
is a pullback in $Set$, i.e., that the square
\begin{center} \xext=1800 \yext850
\begin{picture}(\xext,\yext)(\xoff,\yoff)
\setsqparms[-1`1`1`-1;1800`600]
\putsquare(0,100)[\coprod_{a:S_1^*\ra T^*\,{\rm inner}}\, F(T)`\coprod_{a:(S_1{;}_k S_2)^*\ra T^*\,{\rm inner}}\, F(T)` \coprod_{a:(\bc^{(k)}(S_1))^*\ra T^*\,{\rm inner}}\,F(T)`\coprod_{a:S_2^*\ra T^*\,{\rm inner}}\, F(T);Lan_\bj(F)(\kappa_1)`Lan_\bj(F)(\bc^{(k)}_{S_1})`Lan_\bj(F)(\kappa_2)`Lan_\bj(F)(\bd^{(k)}_{S_2})]
\end{picture}
\end{center}
is a pullback in $Set$. So suppose we have
\begin{center} \xext=1200 \yext=450
\begin{picture}(\xext,\yext)(\xoff,\yoff)
 \putmorphism(0,300)(1,0)[x_1\in F(T_1)`\coprod_{a:S_1^*\ra T^*\,{\rm inner}}\, F(T)`\kappa^{s^*}_{a_1}]{1200}{1}a
  \putmorphism(0,50)(1,0)[x_2\in F(T_2)`\coprod_{a:S_2^*\ra T^*\,{\rm inner}}\, F(T)`\kappa^{s^*}_{a_2}]{1200}{1}a
\end{picture}
\end{center}
such that
\[ Lan_\bj(F)(\bc^{(k)}_{S_1})(x_1) = Lan_\bj(F)(\bd^{(k)}_{S_2})(x_2),\]
i.e., we have a commuting diagram in $\pOpeCard$
\begin{center} \xext=2400 \yext700
\begin{picture}(\xext,\yext)(\xoff,\yoff)
 \setsqparms[-1`1`1`-1;1200`500]
 \putsquare(0,50)[S_1^*`(\bc^{(k)}(S_1))^*=(\bd^{(k)}(S_2))^*`T_1^*`T^*;(\bc^{(k)}_{S_1})^*`a_1``f_1^*]

 \setsqparms[1`0`1`1;1200`500]
 \putsquare(1200,50)[\phantom{(\bc^{(k)}(S_1))^*=(\bd^{(k)}(S_2))^*}`S_2^*` \phantom{T^*}`T_2^*;(\bd^{(k)}_{S_2})^*`a_0`a_2`f_2^*]
\end{picture}
\end{center}
such that
\[ F(f_1)(x_1)=F(f_2)(x_2). \]
By Proposition \ref{inner map},
\[ \bc^{(k)} a_1=a_0=\bd^{(k)} a_2,\;\;\;
 f_1^* =(\bc^{(k)}_{T_1})^*,\;\;\; f_2^* =(\bd^{(k)}_{T_2})^* \]
and the square
\begin{center} \xext=500 \yext=600
\begin{picture}(\xext,\yext)(\xoff,\yoff)
 \setsqparms[1`-1`-1`1;650`450]
 \putsquare(0,50)[T_1`T_1{\oplus}_k T_2`T`T_2;
 \kappa'_1`f_1`\kappa'_2`f_2]
\end{picture}
\end{center}
is a special pushout. We have a commuting diagram
\begin{center} \xext=2000 \yext=1350
\begin{picture}(\xext,\yext)(\xoff,\yoff)
 \setsqparms[1`1`1`1;850`850]
 \putsquare(0,50)[\bc^{(k)}(S_1)^*`S_2^*`T^*`T_2^*;`a_0``]
 \setsqparms[1`0`1`0;850`850]
 \putsquare(510,390)[S_1^*`(S_1{\oplus}_k S_2)^*`T_1^*`(T_1{\oplus}_k T_2)^*;``a_1\oplus_ka_2`]
 \putmorphism(180,180)(3,2)[``]{150}{1}l
 \putmorphism(180,1030)(3,2)[``]{150}{1}l
 \putmorphism(1030,180)(3,2)[``]{120}{1}l
 \putmorphism(1030,1030)(3,2)[``]{120}{1}l
 \put(600,390){\line(1,0){200}}
 \putmorphism(860,390)(1,0)[``]{300}{1}a
 \put(510,1170){\line(0,-1){200}}
 \putmorphism(510,920)(0,1)[``]{520}{1}a
  \put(420,700){\makebox(100,100){$a_1$}}
  \put(750,460){\makebox(100,100){$a_2$}}
\end{picture}
\end{center}
where the bottom square is the above square, the top square is the one we formed earlier. All the horizontal morphisms are outer. Since $a_1$ and $a_2$ are inner, $a_1(S_1)=T_1$ and  $a_2(S_2)=T_2$, we have
\[ (a_1\oplus_k a_2)(S_1\oplus_kS_2)=a_1(S_1)\oplus_ka_2(S_2)=T_1\oplus_kT_2,\]
i.e., $a_1\oplus a_2: (S_1{\oplus}_k S_2)^*\lra (T_1{\oplus}_k T_2)^*$ is inner, as well. So in fact all vertical morphisms in the above diagram are inner.

Suppose we have another inner map $u$ and outer maps $\kappa''_1$, $\kappa''_2$ so that the squares
\begin{center} \xext=2400 \yext800
\begin{picture}(\xext,\yext)(\xoff,\yoff)
 \setsqparms[1`1`1`1;700`500]
 \putsquare(0,50)[S_1^*`(S_1{\oplus}_k S_2)^*`T_1^*`U^*;
 \kappa_1^*`a_1`u`{\kappa''_1}^*]

 \setsqparms[-1`0`1`-1;700`500]
 \putsquare(700,50)[\phantom{(S_1{;}_k S_2)^*}`S_2^*`
 \phantom{U}`T_2^*;\kappa_2^*``a_2`{\kappa''_2}^*]
\end{picture}
\end{center}
commute.  A diagram chasing shows that
\[ {\kappa''_1}^*\circ f^*_1\circ a_1 = {\kappa''_2}^*\circ f^*_2\circ a_1 \]
As inner-outer factorization is essentially unique, it follows that
\[  {\kappa''_1}^*\circ f^*_1 = {\kappa''_2}^*\circ f^*_2 \]
By the universal property of the pushout $(T_1{\oplus}_k T_2)^*$, we have an $\o$-functor $$v: (T_1{\oplus}_k T_2)^*\lra U^*$$ such that
\[ \kappa''_1=u\circ \kappa'_1,\;\; \kappa''_2=u\circ \kappa'_2 \]
Then again, by a diagram chasing, we get
\[ u\circ \kappa_i = v\circ (a_1{\oplus_k} a_2)\circ \kappa_i, \]
for $i=1,2$. Hence, by universal property of the pushout $(S_1{\oplus}_k S_2)^*$, we have that $u=v\circ (a_1\oplus_ka_2)$. But both $u$ and $(a_1\oplus_ka_2)$ are inner so by uniqueness of factorization, see Lemma \ref{factorization}, $v$ must be an iso, as well. This means that we need to show that there is an
\begin{center} \xext=1200 \yext=200
\begin{picture}(\xext,\yext)(\xoff,\yoff)
 \putmorphism(0,50)(1,0)[x\in F(T_1\oplus_kT_2)`
 \coprod_{a:(S_1\oplus_kS_2)^*\ra T^*\,{\rm inner}}\, F(T)`
 \kappa^{(S_1\oplus_kS_2)^*}_{(a_1\oplus_ka_2)}]{1600}{1}a
\end{picture}
\end{center}
such that
\[ Lan_\bj(F)(\kappa_1)(x)=x_1,\;\;\; Lan_\bj(F)(\kappa_2)(x)=x_2.   \]
But $F$ sends special pushouts in $\pOpeCard$ to pullbacks in $Set$ so the square
\begin{center} \xext=750 \yext=600
\begin{picture}(\xext,\yext)(\xoff,\yoff)
 \setsqparms[-1`1`1`-1;750`450]
 \putsquare(0,50)[F(T_1)`F(T_1{\oplus}_k T_2)`F(T)`F(T_2);
 F(\kappa'_1)`F(f_1)`F(\kappa'_2)`F(f_2)]
\end{picture}
\end{center}
is a pullback in $Set$. Thus indeed there is a unique $x\in F(T_1{\oplus}_k T_2)$ such that $F(\kappa'_i)(x)=x_i$ for $i=1,2$. This shows that $Lan_\bj(F)$ preserves special pullbacks. $~~\Box$

\begin{lemma}\label{Lan_j-commutes}
The following square
\begin{center} \xext=1600 \yext=700
\begin{picture}(\xext,\yext)(\xoff,\yoff)
\setsqparms[1`1`1`1;1600`600]
 \putsquare(0,50)[\pPoly `\oC`sPb((\pOpeCard)^{op},Set) `sPb((\pOpeCardo)^{op},Set);\be`\widehat{(-)}`\widehat{(-)}`Lan_\bj]
\end{picture}
\end{center}
commutes, up to an isomorphism.
\end{lemma}

{\it Proof.}~ We shall define two natural transformations $\phi$ and $\psi$ which are mutually inverse, i.e., for a positive-to-one polygraph $Q$ we define
\begin{center} \xext=1500 \yext=320
\begin{picture}(\xext,\yext)(\xoff,\yoff)
\putmorphism(200,130)(1,0)[Lan_\bj\,(\pPoly((-)^*,Q))`\oC((-)^*,Q)`]{1500}{0}a
\putmorphism(200,190)(1,0)[\phantom{Lan_\bj\,(\pPoly((-)^*,Q))}`\phantom{\oC((-)^*,Q)}`\phi_Q]{1500}{1}a
\putmorphism(200,70)(1,0)[\phantom{Lan_\bj\,(\pPoly((-)^*,Q))}`\phantom{\oC((-)^*,Q)}`\psi_Q]{1500}{-1}b
\end{picture}
\end{center}

Let $a:S^*\lra T^*$ be an inner map and  $g:T^*\ra Q$ be a polygraph map, i.e., $g$ is in the following coproduct
\begin{center} \xext=1500 \yext=160
\begin{picture}(\xext,\yext)(\xoff,\yoff)
\putmorphism(200,30)(1,0)[g\in\pPoly(T^*,Q)`\coprod_{S^*\ra R^*\,{\rm inner}}\,\pPoly(R^*,Q)`\kappa^{S^*}_a]{1500}{1}a
\end{picture}
\end{center}
Then we put
\[ \phi_Q(g)=g\circ a .\]

On the other hand, for an $\o$-functor $f:S^*\ra Q\in \oC(S^*,Q)$, by Proposition \ref{type1}.4, we have a factorization
\begin{center}
\xext=600 \yext=300 \adjust[`I;I`;I`;`I]
\begin{picture}(\xext,\yext)(\xoff,\yoff)
 \settriparms[1`1`-1;300]
 \putVtriangle(0,0)[S^*`Q`T_{f(S)}^*;f`f^{in}`\tau_{f(S)}]
\end{picture}
\end{center}
Then we put
\begin{center} \xext=2200 \yext=160
\begin{picture}(\xext,\yext)(\xoff,\yoff)
\putmorphism(200,30)(1,0)[\psi_Q(f)=\tau_{f(S)}\in\pPoly(T_{f(S)}^*,Q)`\coprod_{S^*\ra
R^*\,{\rm inner}}\,\pPoly(R^*,Q)`\kappa^{S^*}_{f^{in}}]{1800}{1}a
\end{picture}
\end{center}
The fact that these transformations are mutually inverse follows from the fact that the above factorization is essentially unique.

The verifications that these transformations are natural is left for the reader. $~~\Box$

We have the following

\begin{theorem}\label{jStar monadic}
The functor
\[ \bj^* :  sPb((\pOpeCardo)^{op},Set)  \lra  sPb((\pOpeCard)^{op},Set) \]
is monadic.
\end{theorem}

{\it Proof.}~ We are going to verify Beck's conditions for monadicity. As $\bj$ is essentially surjective, $\bj^*$ is conservative. By Lemma \ref{Lan j}, the adjunction $Lan_\bj\dashv j^*$ restricts to the above categories.  So $\bj^*$ has a left adjoint. It remains to show that $sPb((\pOpeCardo)^{op},Set)$ has coequalizers of $\bj^*$-contractible coequalizer pairs and that $\bj^*$ preserves them.

To this aim, let assume that we have a parallel pair of morphisms in $sPb((\pOpeCardo)^{op},Set)$
\begin{center}
\xext=1000 \yext=300
\begin{picture}(\xext,\yext)(\xoff,\yoff)
 \putmorphism(0,120)(1,0)[A`B`]{600}{0}a
 \putmorphism(0,40)(1,0)[\phantom{A}`\phantom{B}`G]{600}{1}b
 \putmorphism(0,200)(1,0)[\phantom{A}`\phantom{B}`F]{600}{1}a
\end{picture}
\end{center}
such that
\begin{center}
\xext=1000 \yext=300
\begin{picture}(\xext,\yext)(\xoff,\yoff)
 \putmorphism(0,120)(1,0)[A((-)^*)`B((-)^*)`]{900}{0}a
 \putmorphism(0,40)(1,0)[\phantom{A((-)^*)}`\phantom{B((-)^*)}`G_{(-)^*}]{900}{1}b
 \putmorphism(0,200)(1,0)[\phantom{A((-)^*)}`\phantom{B((-)^*)}`F_{(-)^*}]{900}{1}a
 \putmorphism(0,100)(1,0)[\phantom{A((-)^*)}`\phantom{B((-)^*)}`t]{900}{-1}a
 \putmorphism(900,120)(1,0)[\phantom{B((-)^*)}`Q`q]{900}{1}a
 \putmorphism(900,60)(1,0)[\phantom{B((-)^*)}`\phantom{Q}`s]{900}{-1}b
\end{picture}
\end{center}
is a split coequalizer in $sPb((\pOpeCard)^{op},Set)$, i.e., the following equations
\[q\circ  s= 1_Q, \hskip 5mm q\circ  G_{(-)^*} = q\circ F_{(-)^*}, \hskip 5mm F_{(-)^*}\circ t=1_{B((-)^*)}, \hskip 5mm G_{(-)^*}\circ  t=s\circ q \]
hold.

We are going to construct a special pullbacks preserving  functor
\[ C:(\pOpeCardo)^{op} \lra  Set \]
and a natural transformation
\[ H:B\lra C \]
so that the diagram in $sPb((\pOpeCardo)^{op},Set)$
\begin{center}
\xext=1000 \yext=300
\begin{picture}(\xext,\yext)(\xoff,\yoff)
 \putmorphism(0,120)(1,0)[A`B`]{600}{0}a
 \putmorphism(0,40)(1,0)[\phantom{A}`\phantom{B}`G]{600}{1}b
 \putmorphism(0,200)(1,0)[\phantom{A}`\phantom{B}`F]{600}{1}a
 \putmorphism(600,120)(1,0)[\phantom{B}`C`H]{600}{1}a
\end{picture}
\end{center}
is a coequalizer, and $H_{(-)^*} = q$.

The functor $C$ on a morphism $f:T_1^*\lra T_2^*$ is defined as in the diagram
\begin{center} \xext=900 \yext=900
\begin{picture}(\xext,\yext)(\xoff,\yoff)
 \setsqparms[0`1`-1`1;800`400]
 \putsquare(0,100)[Q(T_1)`Q(T_2)`B(T_1^*)`B(T_2^*);`s_{T_1}`q_{T_2}`B(f)]
 \putmorphism(0,750)(1,0)[C(T_1^*)`C(T_2^*)`C(f)]{800}{1}a
 \put(0,580){\makebox(0,100){$\|$}}
 \put(750,580){\makebox(100,100){$\|$}}
\end{picture}
\end{center}
i.e., $C(T_i)=Q(T_i)$, for $i=1,2$ and $C(f)=q_{T_2}\circ B(f)\circ s_{T_1}$.

The natural transformation $H$ is given by
\[ H_{T^*}=q_T, \]
for $T\in\pOpeCard$.

It remains to verify that
\begin{enumerate}
  \item $C$ is a functor;
  \item $H$ is a natural transformation;
  \item  $C((-)^*)=Q$;
   \item  $H_{(-)^*}=q$;
  \item $C$ preserves the special pullbacks;
  \item $H$ is a coequalizer.
\end{enumerate}

Ad 1. Let
\begin{center} \xext=1000 \yext=150
\begin{picture}(\xext,\yext)(\xoff,\yoff)
 \putmorphism(0,50)(1,0)[T_1^*`T_2^*`f]{500}{-1}a
 \putmorphism(500,50)(1,0)[\phantom{T_2^*}`T_3^*`g]{500}{-1}a
\end{picture}
\end{center}
be a pair of morphisms in $\pOpeCardo$. We calculate
\[ C(g)\circ C(f) =q_{T_3}\circ B(g)  \circ s_{T_2} \circ q_{T_2} \circ B(f)\circ   s_{T_1} = \]
\[ =  q_{T_3} \circ B(g) \circ G_{T_2^*} \circ t_{T_2} \circ B(f) \circ s_{T_1} = \]
\[ =  q_{T_3} \circ G_{T_3^*} \circ A(g) \circ t_{T_2} \circ B(f) \circ  s_{T_1} = \]
\[ =  q_{T_3} \circ F_{T_3^*} \circ A(g) \circ t_{T_2} \circ B(f) \circ  s_{T_1} = \]
\[ = q_{T_3}  \circ B(g) \circ F_{T_2^*} \circ t_{T_2} \circ B(f) \circ s_{T_1}= \]
\[ = q_{T_3} \circ B(g) \circ 1_{B(T_2)^*} \circ B(f) \circ  s_{T_1} = \]
\[ = q_{T_3} \circ B(g) \circ B(f) \circ s_{T_1} = \]
\[ =q_{T_3} \circ B(f\circ g) \circ s_{T_1} =  C(f\circ g),  \]
i.e., $C$ preserves compositions. If $T$ is a positive opetopic cardinal,  we also have
\[ C(1_{T^*})= q_T \circ B(1_{T^*}) \circ s_T = q_T\circ s_T=1_{Q(T)}=1_{C(T^*)}, \]
i.e., $C$ preserves identities, as well.

Ad 2. Let $f:T^*_2\lra T^*_1$ be a morphism in $\pOpeCardo$. We have
\[ H_{T^*_2}\circ B(f)= q_{T_2} \circ B(f) = \]
\[ = q_{T_2} \circ B(f) \circ F(T^*_1) \circ t_{T_1}= \]
\[ = q_{T_2} \circ F(T^*_2) \circ A(f)\circ t_{T_1}= \]
\[ = q_{T_2} \circ G(T^*_2) \circ A(f) \circ  t_{T_1}= \]
\[ = q_{T_2} \circ B(f) \circ G(T^*_1) \circ t_{T_1} = \]
\[ = q_{T_2} \circ B(f) \circ s_{T_1} \circ q_{T_1}= \]
\[ = C(f) \circ q_{T_1} = C(f)\circ  H_{T^*_1}, \]
i.e., $H$ is a natural transformation.

Ad 3. Let $f:T_2\lra T_1$ be a morphism in $\pOpeCard$. Thus $q$ is natural with respect to $f$. So we have
\[ C(f^*)= q_{T_2} \circ B(f^*) \circ s_{T_1} =  Q(f) \circ q_{T_1} \circ s_{T_1} =
Q(f)\circ 1_{T_1}=Q(f), \]
i.e., $C_{(-)^*}=Q$. 

Ad 4. $H_{(-)^*}=q$ holds by definition.

Ad 5. Since special pullbacks involve only the outer morphisms (i.e., those that comes from $\pOpeCard$), and $Q$ preserves special pullbacks, so does $C$.

Ad 6. Finally, we shall show that $H$ is a coequalizer. Let $p:B\lra Z$ be a natural transformation in $sPb((\pOpeCardo)^{op},Set)$ such that $pF=pG$. We put $k=s;p:C\lra Z$ , so that we have a diagram
 \begin{center}
\xext=1000 \yext=700
\begin{picture}(\xext,\yext)(\xoff,\yoff)
 \putmorphism(0,520)(1,0)[A`B`]{600}{0}a
 \putmorphism(0,440)(1,0)[\phantom{A}`\phantom{B}`G]{600}{1}b
 \putmorphism(0,600)(1,0)[\phantom{A}`\phantom{B}`F]{600}{1}a
 \putmorphism(600,520)(1,0)[\phantom{B}`C`H]{600}{1}a
 \putmorphism(1200,500)(0,-1)[\phantom{Q}`Z`k=s;p]{500}{1}r
 \putmorphism(500,500)(3,-2)[``p]{750}{1}l
\end{picture}
\end{center}
We need to verify that $k$ is a natural transformation in $sPb((\pOpeCardo)^{op},Set)$, such that $p=H;k$. Then the uniqueness of $k$ will follow from the fact that $q$ is a split epi. Let $f:T^*_2\lra T^*_1$ be a morphism in $\pOpeCardo$. Then

\[ k_{T^*_2}\circ C(f) = k_{T^*_2} \circ q_{T_2} \circ B(f) \circ s_{T_1}= \]
\[ = p_{T^*_2} \circ s_{T_2} \circ q_{T_2} \circ B(f) \circ s_{T_1}= \]
\[ = p_{T^*_2} \circ G_{T_2^*} \circ t_{T_2} \circ B(f) \circ s_{T_1}= \]
\[ = p_{T^*_2} \circ F_{T_2^*} \circ t_{T_2} \circ B(f) \circ s_{T_1}= \]
\[ = p_{T^*_2} \circ B(f) \circ s_{T_1}= \]
\[ = D(f) \circ p_{T^*_1} \circ s_{T_1} = D(f)\circ k_{T^*_1}, \]
i.e., $k$ is a natural transformation and hence $H$ is indeed a coequalizer of $F$ and $G$ in $sPb((\pOpeCardo)^{op},Set)$, as required. $~~\Box$

Combining the above theorem with Corollaries \ref{xxx} and \ref{sp-nevre2} we get

\begin{theorem}\label{monadicity}
The nerve functor
\[ \widehat{(-)} : \oC  \lra  sPb((\pOpeCard)^{op},Set) \]
sending the $\o$-category $C$ to the presheaf
\[ \oC((-)^*,C):(\pOpeCard)^{op}\lra Set \]
is monadic.
\end{theorem}

We also have

\begin{proposition}\label{connected lims0}
The functor
\[ Lan_\bj : sPb((\pOpeCard)^{op},Set)  \lra sPb((\pOpeCardo)^{op},Set) \]
preserves connected limits.
\end{proposition}

{\it Proof.}~ This is a consequence of Lemma \ref{Lan_j-description}, where $Lan_\bj : \widehat{\pOpeCard} \lra \widehat{\pOpeCardo}$ is described as a family representable functor. In particular, it preserves connected limits. The above functor is a restriction of a family representable functor to the category of functors preserving special pullbacks. Since limits commute with limits, the functor
\[ Lan_\bj : sPb((\pOpeCard)^{op},Set)  \lra sPb((\pOpeCardo)^{op},Set) \]
preserves the connected limits, as well.  $~~\Box$

From Propositions \ref{presheaf_eq2}, \ref{connected lims0}, Lemma \ref{Lan_j-commutes} and Corollary \ref{sp-nevre2} we get immediately

\begin{theorem}\label{connected lims}
The embedding functor
\[ \be : \pPoly  \lra \oC \]
preserves connected limits. $~~\Box$
\end{theorem}

\section{More on monadic adjunctions and distributive laws}\label{sec-more-on-monadic-adj}

We have shown that $\oC$ is monadic over $\pPoly$ with the free functor being the embedding $\pPoly\ra \oC$. We also know that the category of positive-to-one polygraphs is equivalent to the category of presheaves on $\pOpe$ and to the subcategory of special pullback preserving functors $sPb((\pOpeCard)^{op},Set)$ of the presheaf category $\widehat{\pOpeCard}$. Because of the last equivalence we shall freely use the notation $X(Q)$ when $X$ is a presheaf on $\pPoly$ and $Q$ is a positive opetopic cardinal, not only a positive opetope. In this section we shall describe explicitly the whole strongly cartesian monad $(T_\o,\eta_\o,\mu_\o)$ on $\widehat{\pOpe}$ whose category of algebras is equivalent to $\oC$. We also show this monad decomposes into two other strongly cartesian monads of `pure composition' $(T_c,\eta_c,\mu_c)$ and of `adding identities $(T_\iota,\eta_\iota,\mu_\iota)$, in analogy with the decomposition of the strongly cartesian free monoid monad $T_{mon}$ into free semigroup monad and pointed set monad, c.f. \cite{TTT}, p. 258. In particular, the nerve theorem, c.f.\cite{W}, \cite{BMW}, applies.

\subsection*{The monad $T_\o$}

We write $P\raob(q)Q$ to indicate that the map $q$ is inner, i.e., it is an $\o$-functor between $\o$-categories $P^*$ and $Q^*$ so that $q(P)=Q$.
Let $u:X\ra Y$ be a map of presheaves on $\pOpe$, and $f: S'\ra S$ be a positive opetope.  Then $T_\o(X)(S)$ is given by the coproduct

\[   T_\o(X)(S) = \coprod_{S\rao(q) Q} X(Q) = \{ \lk x, q \rk | S\raob(q) Q \in \pOpeCard, \;\; x:Q\ra X \in X(Q)   \},  \]
with coprojections
\[ \kappa^1_q: X(Q) \lra  \coprod_{S\rao(q) Q} X(Q) =T_\o(X)(S),  \]
\[  X(Q) \ni x\mapsto \lk x,q \rk.  \]
Moreover
\[ T_\o(X)_f = T_\o(X)(S) \lra T_\o(X)(S'),  \]
\[ \lk x,q\rk \mapsto \lk x\circ f',q'\rk \]
where $(q',f')$ is the inner-outer factorization of $q\circ f$
\begin{center} \xext=500 \yext=500
\begin{picture}(\xext,\yext)(\xoff,\yoff)
\setsqparms[1`1`1`1;500`400]
\putsquare(0,50)[S'`S`Q'`Q;f`q'`q`f']
 \end{picture}
\end{center}
and
\[ T_\o(u)_S = T_\o(X)(S)\lra T_\o(Y)(S), \]
\[ \lk x,q\rk \mapsto \lk u\circ x,q\rk. \]
The iteration $T_\o^2(X)(S)$ is given by the coproduct
\[   T^2_\o(X)(S) = \coprod_{S\rao(q') R} \;\; \coprod_{R\rao(q) Q}   X(Q) = \]
\[ = \{ \lk x, q,q' \rk | S\raob(q') R, R\raob(q) Q \in \pOpeCard,\;\; x\in X(Q)   \},  \]
with coprojections
\[ \kappa^2_{q,q'}: X(R) \lra  \coprod_{S\rao(q') Q} \;\; \coprod_{Q\rao(q) R}   X(R) =T_\o^2(X)(S),  \]
\[  X(R) \ni x \mapsto \lk x, q, q' \rk.  \]

The unit is given by
\[ ((\eta_\o)_X)_S =\kappa^1_{id_S}: X(S) \lra T_\o(X)(S) = \coprod_{S\lrao(q) Q} X(Q) \]
\[ X(S) \ni x \mapsto \lk x, id_S \rk  \]
and the multiplication is the unique map commuting with the following coprojections
\begin{center} \xext=1900 \yext=600
\begin{picture}(\xext,\yext)(\xoff,\yoff)
 \putmorphism(0,500)(1,0)[T_\o^2(X)(S) = \coprod_{S\rao(q') Q} \;\; \coprod_{Q\rao(q) R}   X(R)`\coprod_{S\lrao(k) Q} X(Q)=T_\o(X)(S)`((\mu_\o)_X)_S )]{2000}{1}a
\put(900,150){\vector(-3,1){600}}
\put(1100,150){\vector(3,1){600}}
\put(350,180){$\kappa^2_{q,q'}$}
\put(1500,180){$\kappa^1_{q\circ q'}$}
\put(900,30){$X(R)$}
\end{picture}
\end{center}
i.e.,
\[ \lk x, q,q' \rk \mapsto \lk x, q\circ q' \rk.  \]

The factorization of a morphism
\[ f: P \lra T_\o(1) \]
through an inner map $q_f$ is
\begin{center} \xext=1200 \yext=100
\begin{picture}(\xext,\yext)(\xoff,\yoff)
 \putmorphism(0,50)(1,0)[P`T_\o(Q_f)`q_f]{600}{1}b
 \putmorphism(600,50)(1,0)[\phantom{T_\o(Q)}`T_\o(1),`]{600}{1}b
\end{picture}
\end{center}
where $q_f = f_P(id_P)$.

\subsection*{The comparison functor $K$}

Now we will describe the comparison functors between the category of $\o$-categories $\oC$ and the category of $T_\o$-algebras $\widehat{\pOpe}^{T_\o}$, i.e., we shall define the functor $K$ in the diagram

\begin{center} \xext=1900 \yext=600
\begin{picture}(\xext,\yext)(\xoff,\yoff)
 \putmorphism(0,500)(1,0)[\oC`\widehat{\pOpe}^{T_\o}`]{2000}{0}a
 \putmorphism(0,550)(1,0)[\phantom{\oC}`\phantom{\widehat{\pOpe}^{T_\o}}`L]{2000}{-1}a
 \putmorphism(0,450)(1,0)[\phantom{\oC}`\phantom{\widehat{\pOpe}^{T_\o}}`K]{2000}{1}b

\put(700,150){\vector(-3,1){600}}
\put(300,350){\vector(3,-1){600}}
\put(250,180){$F_\o$}
\put(750,220){$N$}

\put(1150,150){\vector(3,1){600}}
\put(1950,350){\vector(-3,-1){600}}
\put(1200,220){$F^{T_\o}$}
\put(1680,180){$U^{T_\o}$}
\put(900,30){$\widehat{\pOpe}$}
\end{picture}
\end{center}
and its left adjoint $L$. If $H:\cC \ra \cC'$ is an $\o$-functor, and $f:S\ra S'$ a morphism in $\pOpe$, then
\[ K(\cC)(S) = \oC(S^*,\cC), \]
and
\[  K(\cC)(f):  K(\cC)(S') \ra K(\cC)(S) \]
\[ h: S'^*\ra \cC \mapsto h\circ f.\]
Moreover, the $T_\o$-algebra map
\[  \xi_\cC:  T_\o(K(\cC)) \lra K(\cC), \]
for $S\in \pOpe$
\[ (\xi_\cC)_S:  T_\o(K(\cC))(S) \lra K(\cC)(S),   \]
is given by
\[ \lk S\raob(q)S',S'^*\rab(h)\cC\rk \mapsto h\circ q : S^*\ra \cC. \]
In particular, for $P, S\in \pOpe$,
\[ K(P^*)(S) = \coprod_{S\rao(h)S'} \pOpe(S',P) \cong \{\lk h,k\rk | h:S\raobp S',\; k:S'\ra P\in \pOpe \}\cong T_\o(P)(S). \]
For $\o$-functor $f:P^*\ra Q^*$ in $\pOpe_\o$, the map
\[ K(f) : T_\o(P)(S) \lra T_\o(Q)(S) \]
is given by
\[ \lk q:S\raobp S', x:S'\ra P\rk \mapsto \lk q':S\raobp S'', x:S''\ra P\rk, \]
where $(q',x')$ is the inner-outer factorization of the map in $\pOpe_\o$
\begin{center} \xext=1800 \yext=150
\begin{picture}(\xext,\yext)(\xoff,\yoff)
 \putmorphism(0,50)(1,0)[S^*`S'^*`q]{600}{1}a
 \putmorphism(600,50)(1,0)[\phantom{S'^*}`P^*`x^*]{600}{1}a
 \putmorphism(1200,50)(1,0)[\phantom{P^*}`Q^*`f]{600}{1}a
\end{picture}
\end{center}

\subsection*{The comparison functor $L$}
Below we describe explicitly the left adjoint, the essential inverse to the functor $K$.

For a $T_\o$-algebra $(X,\xi:T_\o(X)\ra X)$, the $\o$-category $L((X,\xi)$ is defined as follows. The set of $n$-cells is
\[ L(X,\xi)_n = X(\alpha^n), \]
and the map composing $n$-cells compatible over dimension $k$,
\[ m_{n,k,n} : X(\alpha^n)\times_{X(\alpha^k)}X(\alpha^n)\cong X(\alpha^{n,k,n})\lra X(\alpha^n), \]
is the composition of the maps

\begin{center} \xext=1800 \yext=150
\begin{picture}(\xext,\yext)(\xoff,\yoff)
 \putmorphism(0,50)(1,0)[X(\alpha^{n,k,n})`T_\o(X)(\alpha^n)`\kappa_{m_{n,k,n}}]{900}{1}a
 \putmorphism(900,50)(1,0)[\phantom{T_\o(X(\alpha^n))}`T_\o(\alpha^n)`\xi_{\alpha^n}]{900}{1}a
\end{picture}
\end{center}

Moreover, as there is a unique map from $\alpha^n$ to any opetopic cardinal $Q$ of dimension less or equal $n$ $\alpha^n\raobp Q$, we have
\[ L(T_\o(P))_n = T_\o(P)(\alpha^n) = \coprod_{din (S')\leq n} \pOpe(S',P). \]

\subsection*{Generic maps and the generic closure.}

Recall the notion of a $T_\o$-generic map from \cite{W}, \cite{BMW}.

Let $g: P\ra T_\o(D)$ be a map in $\widehat{\pOpe}$. Then we have a unique extension $\overline{g}$ as in the diagram
\begin{center}
\xext=800 \yext=500 \adjust[`I;I`;I`;`I]
\begin{picture}(\xext,\yext)(\xoff,\yoff)
\settriparms[1`1`1;400]
  \putVtriangle(0,50)[P`T_\o(P)`T_\o(D);\eta`g`\overline{g}]
\end{picture}
\end{center}
The morphism $g: P\ra T_\o(D)$ is $T_\o$-generic iff $z$ is an isomorphism, where $$\overline{g}_P(id_P)=\lk P\raob(q)P', P'\stackrel{z}{\ra} D \rk.$$

To see this, consider maps $w\ra T_\o(X)$, $v:Q\ra Y$ and $u:X\ra Y$ in $\widehat{\pOpe}$, so that the diagram
 \begin{center} \xext=800 \yext=500
\begin{picture}(\xext,\yext)(\xoff,\yoff)
  \setsqparms[1`1`1`1;800`400]
\putsquare(0,50)[P`T_\o(X)`T_\o(Q)`T_\o(Y);w`g`T_\o(u)`T_\o(v)]
\end{picture}
\end{center}
commutes. Then if $w(id_P) = \lk q',x'\rk$ with $x= x'\circ q':P\ra X$, and $g(id_P)=\lk q,id_Q \rk$, then the above commutation is equivalent to the commutation of the square
 \begin{center} \xext=600 \yext=500
\begin{picture}(\xext,\yext)(\xoff,\yoff)
  \setsqparms[1`1`1`1;600`400]
\putsquare(0,50)[P`X`Q`Y;x`q`u`v]
\end{picture}
\end{center}
Since $(q,v)$ is an inner-outer factorization and $u$ is outer, there is an outer map $d:Q\ra X$ making the square commute. Hence $T_\o(d)$ is the lift in the previous square showing that $g$ is indeed generic. Moreover, if $z$ as above is not an isomorphism then we do not have a lifting $d$ in general, and hence $g$ is not a generic morphism.

Thus, since any positive opetopic cardinal is a codomain of a $T_\o$-generic map whose domain is a positive opetope, in fact $\alpha^n$, for some $n\in \o$, see page \pageref{def_alpha_n}, the $T_\o$-generic closure of $\pOpe$ consists of all opetopic cardinals, i.e., $\pOpeCard$ is the $T_\o$-closure of $\pOpe$.

\subsection*{$T_\o$ as a local right adjoint}
The functor $T_\o$ is a local right adjoint as it is familialy representable. Below we describe explicitly the left adjoint $L_{\o,1}$ to the functor
\[ T_{\o,1} : \widehat{\pOpe}  \lra \widehat{\pOpe}_{T_\o(1)}.\]
Let $\pi: X\ra T_\o(1)$ be a morphism in $\widehat{\pOpe}_{T_\o(1)}$. We have a functor
\[ D_X: \int_{\pOpe} X \lra \widehat{\pOpe} \]
such that the image of a morphism
\begin{center}
\xext=800 \yext=500 \adjust[`I;I`;I`;`I]
\begin{picture}(\xext,\yext)(\xoff,\yoff)
\settriparms[1`1`1;400]
  \putVtriangle(0,50)[P`P'`X;f`h`h']
\end{picture}
\end{center}
in $\int_{\pOpe} X$ is the map $D_X(f): D_X(P,h)\ra D_X(P',h')$ defined from the diagram

 \begin{center} \xext=600 \yext=1500
\begin{picture}(\xext,\yext)(\xoff,\yoff)
  \setsqparms[1`1`1`1;800`800]
\putsquare(0,450)[P`P'`D_X(P,h)`D_X(P',h');f`q`q'`D_X(f)]

\settriparms[0`1`1;400]
  \putVtriangle(0,50)[\phantom{D_X(P,h)}`\phantom{D_X(P',h')}`T_\o(1);`\bar{h}`\bar{h}']
  \putVtriangle(0,850)[\phantom{P}`\phantom{P'}`X;`h`h']

\put(400,750){\line(0,-1){200}}
\put(400,300){\vector(0,-1){150}}
 \put(440,500){$\pi$}
\end{picture}
\end{center}
where $(q,\bar{h})$ is an inner-outer factorization of $\pi\circ h$, $(q',\bar{h}')$ is an inner-outer factorization of $\pi\circ h'$, and $D_X(f)$ is the unique map making the whole diagram commute.  Then $L_{\o,1}(X,\pi)$ is the colimit of the functor $D_X$.

\subsection*{The distributive law.}

If we replace in the above formulas the inner maps by inner epis (inner monos), we still get strongly cartesian monad $(T_{\o,\iota},\eta_{\o,\iota},\mu_{\o,\iota})$ ($(T_{\o,c},\eta_{\o,c},\mu_{\o,c})$) on $\widehat{\pOpe}$. These monads do compose to the monad $(T_\o,\eta_\o,\mu_\o)$. This is because there is a (cartesian) distributive law
\[  \lambda_\o :  T_{\o,c}\circ T_{\o,\iota} \lra  T_{\o,\iota}\circ T_{\o,c}. \]
For $X$ in $\widehat{\pOpe}$ and $S$ in $\pOpe$, both $((T_{\o,c}\circ T_{\o,\iota})_X)_S$ and $((T_{\o,\iota}\circ T_{\o,c})_X)_S$ are given by double coproducts with coprojections as displayed

\[ \sigma^{c,\iota}_{m,e} : X(R) \lra \coprod_{S\rao(m) Q\;mono } \;\; \coprod_{Q\rao(e) R\;epi}   X(R) = ((T_{\o,c}\circ T_{\o,\iota})_X)_S \]

\[ \sigma^{\iota,c}_{e',m'} : X(R) \lra \coprod_{S\rao(e') Q'\;epi} \;\; \coprod_{Q'\rao(m') R\;mono}   X(R) = ((T_{\o,\iota}\circ T_{\o,c})_X)_S \]
The component
\[  ((\lambda_\o)_X)_S :  ((T_{\o,c}\circ T_{\o,\iota})_X)_S \lra ((T_{\o,\iota}\circ T_{\o,c})_X)_S  \]
of distributive law $\lambda_\o$ is the unique map commuting, making all the following triangles
\begin{center} \xext=1900 \yext=600
\begin{picture}(\xext,\yext)(\xoff,\yoff)
 \putmorphism(0,500)(1,0)[\coprod_{S\rao(m) Q\;mono } \;\; \coprod_{Q\rao(e) R\;epi}   X(R)`\coprod_{S\rao(e') Q'\;epi} \;\; \coprod_{Q'\rao(m') R\;mono} X(R)` ((\lambda_\o)_X)_S]{2000}{1}a
\put(900,150){\vector(-3,1){600}}
\put(1100,150){\vector(3,1){600}}
\put(350,180){$\sigma^{c,\iota}_{m,e}$}
\put(1500,180){$\sigma^{\iota,c}_{e',m'}$}
\put(900,30){$X(R)$}
\end{picture}
\end{center}
commute, where $e'$, $m'$ is the epi-mono factorization of the inner map $e\circ m$, i.e., we have the following square of inner epi's and mono's in $\pOpe_\o$
\begin{center} \xext=400 \yext=500
\begin{picture}(\xext,\yext)(\xoff,\yoff)
  \setsqparms[2`3`3`2;400`400]
\putsquare(0,50)[S`Q`Q'`R;m`e'`e`m']
\end{picture}
\end{center}
that commutes.

All the above constructions and considerations can be to truncated to the level $n$, for $n\in \o$. Thus $\nC$ is monadic over $\widehat{pOpe_n}$ and the corresponding monad $(T_n,\eta_n,\mu_n)$, the $n$-truncation of the monad $(T_\o,\eta_\o,\mu_\o)$ is strongly cartesian decomposing into two strongly cartesian monads $(T_{n,\iota},\eta_{n,\iota},\mu_{n,\iota})$ and $(T_{n,c},\eta_{n,c},\mu_{n,c})$, related by a distributive law
\[\lambda_n: T_{n,c}\circ T_{n,\iota}\lra T_{n,\iota}\circ T_{n,c}\]  so that we have
\[T_n=T_{n,\iota}\circ T_{n,c},\]
as monads.

\section{Appendix: a definition of positive-to-one polygraphs}

The category of positive-to-one polygraphs $\pPoly$, the replete subcategory of the category of $\o$-categories $\oC$,  is a limit of a tower of positive-to-one polygraphs of dimension $n$ $\pPoly_n$,  replete subcategories of the category of $n$-categories $\nC$.  Thus we shall describe a diagram
\begin{center} \xext=1000 \yext=2300
\begin{picture}(\xext,\yext)(\xoff,\yoff)
\setsqparms[1`1`1`1;1000`600]
\putsquare(0,50)[\pPoly_1`\Cat`\pPoly_0`Set;\varphi_1`tr_0`tr_0`\varphi_0]
\putsquare(0,1000)[\pPoly_{n+1}`\njC`\pPoly_n`\nC;\varphi_{n+1}`tr_n`tr_n`\varphi_n]
\setsqparms[1`1`1`0;1000`400]
\putsquare(0,1800)[\pPoly`\oC``;\varphi_\o```]
\put(0,800){$\vdots$}
\put(1000,800){$\vdots$}
\put(0,1700){$\vdots$}
\put(1000,1700){$\vdots$}
\end{picture}
\end{center}
Having defined the subcategory $\pPoly_n$ of $\nC$, we define the category $\CP_{n+1}$ to be a non-full  subcategory of $(n+1)$-categories $\njC$ such that an $(n+1)$-category $C$ is an object of $\CP_{n+1}$ iff its truncation to $\nC$ is in  $\pPoly_n$, and an $(n+1)$-functor $f:C\ra D$ is $\CP_{n+1}$ iff its truncation to $\nC$ is in  $\pPoly_n$

The first two stages of the above tower are built as follows.
\begin{center} \xext=3200 \yext=1000
\begin{picture}(\xext,\yext)(\xoff,\yoff)
\setsqparms[1`1`1`1;1200`800]
\putsquare(1200,50)[\pPoly_{1}`\CP_1`\pPoly_0`\CP_0;\varphi_1``tr_0`\varphi_0]
\setsqparms[-1`0`0`-1;1200`800]
\putsquare(0,50)[Set`\phantom{\pPoly_{1}}`Set`\phantom{\pPoly_0};D_{1}```D_0]

\setsqparms[1`0`1`1;800`800]
\putsquare(2400,50)[\phantom{\CP_1}`\Cat`\phantom{\CP_0}`Set;\psi_1``tr_0`\psi_0]

\put(400,400){$Set\da D_0$}
\put(-100,250){$|-|_{1}$}
\put(930,250){$t_0$}
\put(480,330){\vector(-2,-1){400}}
\put(650,330){\vector(2,-1){400}}

\put(1500,710){$F_{1}$}
\put(1300,600){\vector(4,1){750}}
\put(940,500){\line(4,1){190}}

\put(1500,460){$U_{1}$}
\put(1300,500){\line(4,1){750}}
\put(1130,450){\vector(-4,-1){200}}

\put(500,650){$\overline{(-)}^{1}$}
\put(670,520){\vector(1,1){270}}
\end{picture}
\end{center}
The categories $\pPoly_0$ and $\CP_0$  are the category of sets, the functors $\varphi_0$ and $\psi_0$ are the identity on $Set$, and the functor $D_0$ sends set $X$ to its product $X\times X$. $\CP_1$ is $\Cat$ and $\psi_1$ is an identity on $\Cat$. Having this data, we form a comma category $Set\da D_0$ with projections
\[ |-|_1: Set\da D_0\ra Set,\;\; {\rm and}  tr_0: Set\da D_0\ra \pPoly_0.\]
There is a forgetful functor from
$U_1: \CP_1=\Cat \ra Set\da U_0$ sending category $(C,d,c,\circ,i)$ to  $(C,d,c)$. The functor $U_1$ has a right adjoint $F_1$ of a free category on a graph.
Then $\pPoly_1$ is the replete subcategory of $\Cat$ such that we have
\begin{center} \xext=2000 \yext=300
\begin{picture}(\xext,\yext)(\xoff,\yoff)
 \putmorphism(0,50)(1,0)[Set\da U_0`\pPoly_1`\overline{(-)}^1]{1000}{1}a
 \putmorphism(1000,50)(1,0)[\phantom{\pPoly_1}`\Cat`\varphi_1]{1000}{1}a
\end{picture}
\end{center}
is a full and faithful/essentially surjective factorization. Thus $\pPoly_1$ is a category of free categories over graphs with functors that send generating morphisms to generating morphisms. Finally, we define the functor
\[ D_1: \pPoly_1 \lra Set \]
\[ X\mapsto  \{ \lk x,y\rk\in X_1 |\;  x\| y\;  y\;{\rm is\;a\; ideterminate \; and}\; x\;{\rm is\;a\; non-identity\;morphism} \},\]
i.e., sending a free category $X$ to the set of parallel pairs of morphisms $\lk x,y\rk$ so that they are parallel ($\bd(x)=\bd(y)$ and $\bc(x)=\bc(y)$) and $y$ is a generating morphism and $x$ is a non-identity morphism.

The inductive step in the construction of the above tower is similar to the second one.
\begin{center} \xext=3300 \yext=1000
\begin{picture}(\xext,\yext)(\xoff,\yoff)
\setsqparms[1`1`1`1;1200`800]
\putsquare(1200,50)[\pPoly_{n+1}`\CP_{n+1}`\pPoly_n`\CP_n;\varphi_{n+1}``tr_n`\varphi_n]
\setsqparms[-1`0`0`-1;1200`800]
\putsquare(0,50)[Set`\phantom{\pPoly_{n+1}}`Set`\phantom{\pPoly_n};D_{n+1}```D_n]

\setsqparms[1`0`1`1;900`800]
\putsquare(2400,50)[\phantom{\CP_{n+1}}`\njC`\phantom{\CP_n}`\nC;\psi_{n+1}``tr_n`\psi_n]

\put(400,400){$Set\da D_n$}
\put(-100,250){$|-|_{n+1}$}
\put(930,250){$t_n$}
\put(480,330){\vector(-2,-1){400}}
\put(650,330){\vector(2,-1){400}}

\put(1500,730){$F_{n+1}$}
\put(1300,600){\vector(4,1){750}}
\put(940,500){\line(4,1){190}}

\put(1500,460){$U_{n+1}$}
\put(1300,500){\line(4,1){750}}
\put(1130,450){\vector(-4,-1){200}}

\put(500,650){$\overline{(-)}^{n+1}$}
\put(670,520){\vector(1,1){270}}
\end{picture}
\end{center}
Having defined the subcategory $\pPoly_n$ of $\nC$ and the functor $D_n:\pPoly_n\ra Set$, we built the other parts of the above diagram.
We form a comma category $Set\da D_n$ with projections
\[ |-|_{n+1}: Set\da D_n\ra Set,\;\;\; {\rm and}\;\;\;  tr_n: Set\da D_n\ra \pPoly_n.\] There is a forgetful functor from $U_{n+1}: \CP_{n+1}\ra Set\da U_n$ sending $(n+1)$-category $C$ in $\CP_{n+1}$ to its $n$-truncation together the set of $(n+1)$-cells and functions assigning their domains and codomains. The functor $U_{n+1}$ has a left adjoint $F_{n+1}$ of a free $(n+1)$-category. Then $\pPoly_{n+1}$ is the replete subcategory of $\CP_{n+1}$ such that
\begin{center} \xext=2000 \yext=300
\begin{picture}(\xext,\yext)(\xoff,\yoff)
 \putmorphism(0,50)(1,0)[Set\da D_n`\pPoly_{n+1}`\overline{(-)}^{n+1}]{1000}{1}a
 \putmorphism(1000,50)(1,0)[\phantom{\pPoly_{n+1}}`\CP_{n+1}`\varphi_{n+1}]{1000}{1}a
\end{picture}
\end{center}
is a full and faithful/essentially surjective factorization. Thus $\pPoly_{n+1}$ is a category of free $(n+1)$-categories whose generators/indeterminates are cells of positive-to-one shapes $(n+1)$-functors that send indeterminates to indeterminates. Finally, we define the functor
\[ D_{n+1}: \pPoly_{n+1} \lra Set \]
\[ X\mapsto  \{ \lk x,y\rk\in X_1 |\;  x\| y,\;\;  y\;{\rm is\;a\;ideterminate\; and}\; x\;{\rm is\;a\; non-identity\;cell} \},\]
i.e., sending a free category $X$ to the set of parallel pairs of cells $\lk x,y\rk$, ($\bd(x)=\bd(y)$ and $\bc(x)=\bc(y)$),   so that $y$ is an indeterminate
and $x$ is a non-identity map.


\begin{thebibliography}{CWMW} \frenchspacing
\bibitem[BD]{BaezDolan}
J. Baez, J. Dolan,  {\em Higher-dimensional algebra III: n-Categories and the algebra of opetopes}. Advances in Math. 135 (1998), pp. 145-206.

\bibitem[TTT]{TTT}
M. Barr, C. Wells, {\em Toposes, Triples and Theories}, Springer-Verlag, New York, 1985. Republished in:
Reprints in Theory and Applications of Categories, No. 12 (2005) pp. 1-287.

\bibitem[B]{Batanin} M. Batanin, {\em Monoidal globular categories as a natural environment for the theory of weak n-categories}. Advances in Math. 136 (1998), 39-103.

\bibitem[Be]{Berger} C. Berger, {\em A Cellular Nerve for Higher Categories}. Adv. in Mathematics 169,  (2002), pp. 118-175.

\bibitem[BMW]{BMW} C. Berger, P-A. Melli\'{e}s, M. Weber , {\em Monads with arities and their associsted theories}  J. Pure and Applied Alg.  216 (2012), pp. 2029-2048..

\bibitem[Bu]{Bu} A. Burroni, {\em Higher-dimensional word problems with applications to equational logic}, Theoretical Computer Science 115 (1993), no. 1, pp. 43–62.

\bibitem[CJ]{CJ} A. Carboni, P.T. Johnstone {\em Connected limits, familial representability and Artin glueing}, Math. Struct. Computer Science 5 (1995), pp. 441–459.

\bibitem[H]{H} V. Harnik, {\em Private communication.}

\bibitem[HMP]{HMP} C. Hermida, M. Makkai, J. Power, {\em On weak higher dimensional categories, I} Parts 1,2,3, J. Pure and Applied Alg.  153 (2000), pp. 221-246, 157 (2001), pp. 247-277, 166 (2002), pp. 83-104.

\bibitem[J]{Joyal} A. Joyal, {\em Disks, Duality and $\Theta$-categories}. Preprint, (1997).

\bibitem[CWM]{CWM} S. MacLane, {\em Categories for Working Mathematician}, Springer-Verlag, New York, (1971).

\bibitem[SGL]{SGL} S. MacLane, I.Moerdijk, {\em Sheaves in Geometry and Logic: a first introduction to topos theory}, Springer-Verlag, New York, (1992).

\bibitem[MM]{MM} M.Makkai, {\em The multitopic omega-category of all multitopic omega-categories.} Preprint (1999).

\bibitem[MZ]{MZ} M. Makkai and M. Zawadowski, {\em Disks and duality}. Theory and Applications of Categories, 8(7), 2001, pp. 114-243.

\bibitem[M]{M} F. Metayer, {\em Private communication.}

\bibitem[W]{W} M. Weber, {\em Famillial 2-functors and parametric right adjoints} Theory and Applications of Categories, 18(22), 2007, pp. 665-732.

\bibitem[Z]{Z} M. Zawadowski, {\em Positive face structures and many-to-positive computads}, arXiv:0708.2658v1, pp. 1-77.
\end{thebibliography}
\end{document}